\documentclass[12pt,a4paper,reqno]{article}

\pdfoutput=1 

\usepackage{amsmath}
\usepackage{amsfonts}
\usepackage{amssymb}
\usepackage{amsthm}

\usepackage{enumerate}
\usepackage{geometry}
\usepackage{graphicx}
\usepackage{hyperref}
\usepackage[utf8]{inputenc}
\usepackage{mathbbol}
\usepackage[numbers,square]{natbib}
\usepackage{sidecap}
\usepackage{stmaryrd}
\usepackage{subfig}
\usepackage{tikz}

\usetikzlibrary{arrows}
\usetikzlibrary{matrix}

\newcommand*{\doi}[1]{\href{http://dx.doi.org/\detokenize{#1}}{doi: \detokenize{#1}}}

\makeatletter
\@ifundefined{chapter}
   {\newtheorem{theorem}{Theorem}[section]}
   {}
\makeatother

\newtheorem*{exercise*}{Exercise}

\AtBeginDocument{
    
    \let\obar\bar
    
    \let\otilde\tilde

    \let\bar\undefined
    \let\hat\undefined
    \let\tilde\undefined

    \newcommand{\bar}[1]{\ensuremath{\overline{#1}}}
    \newcommand{\hat}[1]{\ensuremath{\widehat{#1}}}
    \newcommand{\tilde}[1]{\ensuremath{\widetilde{#1}}}
    
}

\def\ba{\begin{align}}
\def\ea{\end{align}}

\newcommand{\basis}[1]{\ensuremath{\mathrm{#1}}}
\newcommand{\coord}[1]{\ensuremath{#1}}

\newcommand{\bx}{\ensuremath{\basis{x}}}
\newcommand{\by}{\ensuremath{\basis{y}}}
\newcommand{\bz}{\ensuremath{\basis{z}}}
\newcommand{\bu}{\ensuremath{\basis{u}}}
\newcommand{\bv}{\ensuremath{\basis{v}}}
\newcommand{\bw}{\ensuremath{\basis{w}}}
\newcommand{\bq}{\ensuremath{\basis{x}}}
\newcommand{\bp}{\ensuremath{\basis{p}}}
\newcommand{\bl}{\ensuremath{\basis{\lambda}}}

\newcommand{\cz}{\ensuremath{\coord{z}}}

\newcommand{\cv}{\ensuremath{\coord{v}}}

\newcommand{\cc}{\ensuremath{\coord{c}}}
\newcommand{\cq}{\ensuremath{\coord{q}}}
\newcommand{\cp}{\ensuremath{\coord{p}}}
\newcommand{\cl}{\ensuremath{\coord{\lambda}}}

\newcommand{\dz}{\ensuremath{\dot{\cz}}}

\newcommand{\cd}{\ensuremath{\dot{\cc}}}
\newcommand{\qd}{\ensuremath{\dot{\cq}}}
\newcommand{\pd}{\ensuremath{\dot{\cp}}}

\newcommand{\oq}{\ensuremath{\obar{\cq}}}
\newcommand{\op}{\ensuremath{\obar{\cp}}}
\newcommand{\oz}{\ensuremath{\obar{\cz}}}

\newcommand{\rkq}{\ensuremath{\coord{Q}}}
\newcommand{\rkv}{\ensuremath{\dot{\rkq}}}
\newcommand{\rkp}{\ensuremath{\coord{P}}}
\newcommand{\rkf}{\ensuremath{\dot{\rkp}}}

\newcommand{\mbb}[1]{\ensuremath{\mathbb{#1}}}
\newcommand{\mcal}[1]{\ensuremath{\mathcal{#1}}}

\newcommand{\rsp}{\ensuremath{\mathbb{R}}}

\newcommand{\mf}[1]{\ensuremath{\mathcal{#1}}}
\newcommand{\tb}[2][\empty]{\ensuremath{\mathsf{T}_{#1}     #2}}
\newcommand{\cb}[2][\empty]{\ensuremath{\mathsf{T}_{#1}^{*} #2}}

\newcommand{\ext}{\ensuremath{\mathsf{d}}}

\newcommand{\eps}{\ensuremath{\epsilon}}
\newcommand{\veps}{\ensuremath{\varepsilon}}
\newcommand{\phy}{\ensuremath{\varphi}}

\newcommand{\abs}[1]{\ensuremath{\left \vert #1 \right \vert}}

\newcommand{\pair}[2]{\ensuremath{\left < #1 \, , \, #2 \right >}}

\DeclareMathOperator{\id}{id}

\newcommand{\unity}{\ensuremath{\mathbb{1}}}
\newcommand{\identity}{\ensuremath{\mathbb{I}}}

\DeclareMathAlphabet{\mathpzc}{OT1}{pzc}{m}{it}

\renewcommand{\basis}[1]{\ensuremath{#1}}
\renewcommand{\bq}{\ensuremath{\basis{q}}}
\renewcommand{\bv}{\ensuremath{\basis{\qd}}}

\newcommand{\bxd}{\ensuremath{\dot{\bx}}}
\newcommand{\byd}{\ensuremath{\dot{\by}}}
\newcommand{\bud}{\ensuremath{\dot{\bu}}}

\begin{document}

\title{Projected Variational Integrators for\\Degenerate Lagrangian Systems}

\author{
\large{Michael Kraus}\\
\small{(michael.kraus@ipp.mpg.de)}
\vspace{.5em}\\
\normalsize{Max-Planck-Institut f\"ur Plasmaphysik}\\
\normalsize{Boltzmannstra\ss{}e 2, 85748 Garching, Deutschland}%
\vspace{.1em}\\
\small{and}
\vspace{.1em}\\
\normalsize{Technische Universit\"at M\"unchen, Zentrum Mathematik}\\
\normalsize{Boltzmannstra\ss{}e 3, 85748 Garching, Deutschland}%
\vspace{1em}\\
}

\date{August 24, 2017}

\maketitle

\begin{abstract}

We propose and compare several projection methods applied to variational integrators for degenerate Lagrangian systems, whose Lagrangian is of the form $L = \vartheta(\bq) \cdot \bv - H(\bq)$ and thus linear in velocities.
While previous methods for such systems only work reliably in the case of $\vartheta$ being a linear function of $\bq$, our methods are long-time stable also for systems where $\vartheta$ is a nonlinear function of $\bq$.
We analyse the properties of the resulting algorithms, in particular with respect to the conservation of energy, momentum maps and symplecticity.
In numerical experiments, we verify the favourable properties of the projected integrators and demonstrate their excellent long-time fidelity.
In particular, we consider a two-dimensional Lotka--Volterra system, planar point vortices with position-dependent circulation and guiding centre dynamics.

\end{abstract}

\newpage

\tableofcontents

\newpage

%%%%%%%%%%%%%%%%%%%%%%%%%%%%%%

\section{Introduction}

In various areas of physics, we are confronted with degenerate Lagrangian systems, whose Lagrangian is of the form $L = \vartheta(\bq) \cdot \bv - H(\bq)$, that is $L$ is linear in the velocities $\bv$.
Examples for such systems are planar point vortices, certain Lotka--Volterra models and guiding centre dynamics.

In order to derive structure-preserving integrators~\cite{HairerLubichWanner:2006} for such systems, it seems natural to apply the variational integrator method~\cite{MarsdenWest:2001, LewMata:2016, Rowley:2002, Qin:2008, Qin:2009, Tyranowski:2014, Tyranowski:2014thesis}.
This, however, does not immediately lead to stable integrators, as the resulting numerical methods will in general be multi-step methods which are subject to parasitic modes~\cite{Ellison:2016, Ellison:2015}. Moreover, we face the initialization problem, that is how to determine initial data for all previous time steps used by the method without introducing a large error into the solution.

A potential solution to the first problem is provided by the discrete fibre derivative while a solution to the second problem is provided by the continuous fibre derivative. Using the discrete fibre derivative, we can rewrite the discrete Euler--Lagrange equations in position-momentum form, which constitutes a one-step method for numerically computing the phasespace trajectory in terms of the generalized coordinates $\cq$ together with their conjugate momenta $\cp$. The resulting system can be solved, as in general the discrete Lagrangian will not be degenerate even though the continuous Lagrangian is.
The continuous fibre derivative $(\cq (t), \cp (t) = \partial L / \partial \bv (\cq (t), \qd (t)))$ can then be used in order to obtain an initial value for the conjugate momenta as functions of the coordinates as $\cp (t) = \vartheta (\cq (t))$.
Tyranowski et al.~\cite{Tyranowski:2014, Tyranowski:2014thesis} show that this is a viable strategy when $\vartheta$ is a linear function. Unfortunately, for the general case of $\vartheta$ being a nonlinear function, this idea does in general not lead to stable integrators as the numerical solution will drift away from the constraint submanifold defined by the continuous fibre derivative, $\phi(\cq(t), \cp(t)) = \cp (t) - \vartheta (\cq (t)) = 0$. 

A standard solution for such problems is to project the solution back to the constraint submanifold after each time step~\cite{HairerWanner:1996, Hairer:2001, HairerLubichWanner:2006}. This, however, renders the integrator non-symmetric (assuming the variational integrator itself is symmetric), which leads to growing errors in the solution and consequently a drift in the total energy of the system.
Improved long-time stability is achieved by employing a symmetric projection~\cite{Hairer:2000, Hairer:2001, Chan:2004}, where the initial data is perturbed away from the constraint submanifold before the variational integrator is applied and then projected back to the manifold.
While these projection methods are standard procedures for holonomic and non-holonomic constraints, there are only few traces in the literature on their application to Dirac constraints $\phi(\cq,\cp)=0$.
Some authors consider general differential algebraic systems of index two~\cite{HairerLubichRoche:1989, Ascher:1991, Ascher:1992, Chan:2002, Chan:2004, Jay:2003, Jay:2006}, but do not go into the details of Lagrangian or Hamiltonian systems. 
Or they consider symplectic integrators for Hamiltonian systems subject to holonomic constraints~\cite{Leimkuhler:1994, Jay:1996, Jay:2007eg, HairerLubichWanner:2006}, which simplifies the situation dramatically compared to the case of Dirac constraints.
Most importantly, we are not aware of any discussion of the influence of such a projection on the symplecticity of the algorithm assuming that the underlying numerical integrator is symplectic. As symplecticity is a crucial property of Lagrangian and Hamiltonian systems, which is often important to preserve in a numerical simulation, we will analyse all of the proposed methods regarding its preservation.
We will find that the well-known projection methods, both standard projection and symmetric projection, are not symplectic. However, we can introduce small modifications to the symmetric projection method which make it symplectic.

\pagebreak

The outline of the paper is as follows.
In Section~\ref{sec:degenerate_lagrangian_systems} we provide an overview of degenerate Lagrangian systems and Dirac constraints and their various formulations and discuss symplecticity and momentum maps.
In Section~\ref{sec:vi} we review the discrete action principle leading to variational integrators and problems that arise when this method is applied to degenerate Lagrangians.
This is followed by a discussion of the proposed projection methods in Section~\ref{sec:projection} and numerical experiments in Section~\ref{sec:numerical_experiments}.

\section{Degenerate Lagrangian Systems}
\label{sec:degenerate_lagrangian_systems}

Degenerate Lagrangian systems have attracted quite some interest in the geometric mechanics literature~\cite{Kunzle:1969uc, Gotay:1978tz, Gotay:1979wp, Gotay:1980to, Cantrijn:1985ue, Carinena:1990gs, Carinena:2003hm} due to their interesting properties.
They are also relevant for practical applications like the study of population models, point vortex dynamics or reduced charged particle models like the guiding centre system.
In the following, we will consider degenerate Lagrangian systems characterized by a Lagrangian that is linear or singular in the velocities. In particular, we consider the class of systems whose Lagrangian is of the form
\begin{align}\label{eq:degenerate_lagrangian}
L (\bq, \bv) = \vartheta(\bq) \cdot \bv - H(\bq) .
\end{align}
The Lagrangian $L$ is a function on the tangent bundle $\tb{\mf{M}}$,
\begin{align}\label{eq:lagrangian_tangent_bundle}
L : \tb{\mf{M}} \rightarrow \rsp ,
\end{align}
where $\mf{M}$ denotes the configuration manifold of the system which is assumed to be of dimension $d$. The cotangent bundle of the configuration manifold $\mf{M}$ is denoted by $\cb{\mf{M}}$.
Further, we denote the coordinates of a point $m \in \mf{M}$ by $\bq(m) = (\bq^{1} (m), \hdots, \bq^{d} (m))$ and similarly coordinates of points in $\tb{\mf{M}}$ by $(\bq^{i}, \bv^{i})$ and coordinates of points in $\cb{\mf{M}}$ by $(\bq^{i}, \bp^{i})$.
In the following, we will always assume the existence of a global coordinate chart, so that $\mf{M}$ can be identified with the Euclidean space $\rsp^{d}$. For simplicity, we often use short-hand notation where we write $\bq$ to refer to both a point in $\mf{M}$ as well as its coordinates. Similarly, we often denote points in the tangent bundle $\tb{\mf{M}}$ by $(\bq, \bv)$.
In local coordinates, the Lagrangian~\eqref{eq:lagrangian_tangent_bundle} is thus written as a map $(\bq, \bv) \mapsto L(\bq, \bv)$.

In Equation~\eqref{eq:degenerate_lagrangian}, $\vartheta = \vartheta_{i} (\bq) \, \ext \bq^{i}$ is a differential one-form $\vartheta : \mf{M} \rightarrow \cb{\mf{M}}$, whose components $\vartheta_{i} : \mf{M} \rightarrow \rsp$ are general, possibly nonlinear functions of $\bq$, some of which (but not all) could be identically zero.
For details on differential forms, tangent and cotangent bundles the interested reader may consult any modern book in mathematical physics or differential geometry. We recommend~\cite{Dray:2014, BaezMuniain:1994, Darling:1994, Frankel:2011} for more physics oriented accounts and~\cite{Lee:2012, Lee:2009, Tu:2011, Morita:2001} for more mathematics oriented accounts. In the following we assume a basic understanding of these concepts. To see their usefulness for classical mechanics we refer to~\cite{AbrahamMarsden:1978, MarsdenRatiu:2002, Holm:2009}.

\subsection{Hamilton's Action Principle}

The evolution of Lagrangian systems is described by curves $\cq$ on $\mf{M}$.
To make this precise, let us fix two points $\cq_{1}, \cq_{2} \in \mf{M}$ and an interval $[t_{1}, t_{2}] \subset \rsp$ and define the path space connecting $\cq_{1}$ and $\cq_{2}$ as
\begin{align}
\label{eq:tangent_bundle_path_space}
\mf{Q} ( t_{1}, t_{2}, \cq_{1}, \cq_{2} ) = \big\{ \cq : [t_{1}, t_{2}] \rightarrow \mf{M} \; \big\vert \; \text{$\cq$ is a $C^{2}$ curve} \; \cq(t_{1}) = \cq_{1}, \, \cq(t_{2}) = \cq_{2} \big\} .
\end{align}
Elements $\cq$ of $\mf{Q} ( t_{1}, t_{2}, \cq_{1}, \cq_{2} )$ map the time interval $[t_{1}, t_{2}]$ to curves on $\mf{M}$, whereby the first and last points, $\cq(t_{1})$ and $\cq(t_{2})$, take fixed values, $\cq_{1}$ and $\cq_{2}$, respectively.
Such a curve $\cq : [t_{1}, t_{2}] \rightarrow \mf{M}$ with $\cq(t) = (q^{1} (t), \hdots, \cq^{d} (t))$ can be lifted to a curve $\hat{\cq} : [t_{1}, t_{2}] \rightarrow \tb{\mf{M}}$, which in coordinates is given by
\begin{align}
\hat{\cq} (t) = \bigg( \cq^{1} (t), \hdots, \cq^{d} (t), \dfrac{d\cq^{1}}{dt} (t), \hdots, \dfrac{d\cq^{d}}{dt} (t) \bigg) .
\end{align}
In the following the derivative of the curve with respect to the parameter $t$ is denoted by $\qd = d\cq/dt$.
This constitutes slight abuse of notation as we also denote the tangent bundle coordinates that way (note that not all curves in the tangent bundle are lifts of curves in the configuration manifold), but it should be clear at any time if we refer to the derivative of the curve or to the coordinates.

\begin{figure}[tb]
	\centering
	\includegraphics[width=.5\textwidth]{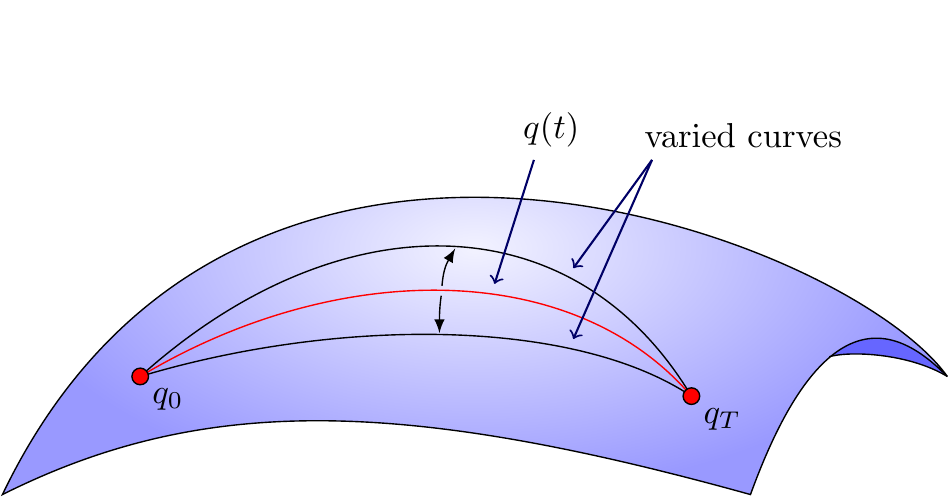}
	\caption{Variations of the trajectory $\cq$.}
	\label{fig:variation_continuous}
\end{figure}

In order to determine the equations of motion of a Lagrangian system, we consider infinitesimal variations of the action integral
\begin{align}\label{eq:action_integral}
\mcal{A} [\cq] = \int_{0}^{T} L (\cq (t), \qd (t)) \, dt ,
\end{align}
where, without loss of generality and in order to simplify the discrete treatment, we choose an interval $[0,T]$.
Infinitesimal variations (c.f., Figure~\ref{fig:variation_continuous}) are defined in terms of $C^{2}$ maps $\cc : [0,T] \rightarrow \rsp^{d}$ which vanish at the boundaries of the interval $[0,T]$, that is $\cc(0) = 0$ and $\cc(T) = 0$, and are such that $c(t) \in \tb[q(t)]{\mf{M}}$ for $0 \leq t \leq T$ with $\tb[q(t)]{\mf{M}}$ the tangent space to $\mf{M}$ at the point $q(t)$.
In coordinates, we can write
\begin{align}
\cc (t) = \big( \cc^{1} (t), \hdots, \cc^{d} (t) \big) .
\end{align}
We can now define a one-parameter family of trajectories $\cq^{\eps} \in \mf{Q} ( 0, T, \cq_{1}, \cq_{2} )$, for which
\begin{align}
\cc = \dfrac{d}{d\eps} \cq^{\eps} \bigg\vert_{\eps = 0} ,
\end{align}
with $\eps \in (-r,+r)$ and $r \in \rsp_{+}$. That is $\cq^{\eps}$ is a differentiable mapping $(-r,+r) \times [0,T] \rightarrow \mf{M}$, such that $\cq^{\eps} (0) = q (0)$ and $\cq^{\eps} (T) = q (T)$ for all values of $\eps$ and $q^{0} (t) = q(t)$ for all values of $t$.
If $\mf{Q}$ is a vector space, the simplest example of such a family of trajectories is given by
\begin{align}\label{eq:one_parameter_family_of_curves}
\cq^{\eps} (t) = \big( \cq^{1} (t) + \eps \cc^{1} (t), \hdots, \cq^{d} (t) + \eps \cc^{d} (t) \big) .
\end{align}
On a general manifold, the corresponding expressions are usually more complicated.
However, if the manifold is embedded in $\rsp^{d}$, each member of the family of trajectories can be expanded in a power series, whose leading-order terms are those given in~\eqref{eq:one_parameter_family_of_curves}.

Hamilton's principle states that in order to determine the equations of motion, we need to find a curve $\cq = \cq^{0}$, such that the action~\eqref{eq:action_integral} takes a stationary point with respect to all curves $\cq^{\eps}$.
A necessary condition for $\cq$ making the action stationary is that 
\begin{align}\label{eq:hamiltons_action_principle_1}
\dfrac{d}{d\eps} \mcal{A} [\cq^{\eps}] \bigg\vert_{\eps = 0}
\nonumber
&= \dfrac{d}{d\eps} \int_{0}^{T} L (\cq^{\eps} (t), \qd^{\eps} (t)) \, dt \bigg\vert_{\eps = 0}
\\
&= \int_{0}^{T} \bigg[ \dfrac{\partial L}{\partial \bq} \big( \cq (t), \qd (t) \big) \cdot \cc (t) + \dfrac{\partial L}{\partial \bv} \big( \cq (t), \qd (t) \big) \cdot \cd (t) \bigg] \, dt = 0 .
\end{align}
Let us note that the derivatives of $L$ with respect to $\bq$ and $\bv$ are sometimes also written as $D_{1}$ and $D_{2}$, respectively, where $D_{i}$ denotes the slot-derivative with respect to the $i$th argument of $L$.
Assuming that the operations of computing variations of $\cq$ and computing the time derivative of $\cq$ commute (which is a fair assumption, see e.g. \citet{SaletanCromer:1971}), we can integrate the second term of~\eqref{eq:hamiltons_action_principle_1} by parts to obtain
\begin{multline}\label{eq:hamiltons_action_principle_2}
\dfrac{d}{d\eps} \mcal{A} [q^{\eps}] \bigg\vert_{\eps = 0}
= \int_{0}^{T} \bigg[ \dfrac{\partial L}{\partial \bq} \big( \cq (t), \qd (t) \big) - \dfrac{d}{dt} \bigg( \dfrac{\partial L}{\partial \bv} \big( \cq (t), \qd (t) \big) \bigg) \bigg] \cdot \cc (t) \, dt \\
 + \bigg[ \dfrac{\partial L}{\partial \bv} \big( \cq (t), \qd (t) \big) \cdot \cc (t) \bigg]_{0}^{T} = 0 .
\end{multline}
As the infinitesimal variations $\cc$ are required to vanish at the boundaries, the boundary term vanishes, and as the variations $\cc$ are otherwise arbitrary, vanishing of the variations of $\mcal{A}$ implies vanishing of the term in square brackets in the integrand.
This leads us to the Euler--Lagrange equations, that is the equations of motion,
\begin{align}\label{eq:euler_lagrange_equations_general}
\dfrac{\partial L}{\partial \bq} \big( \cq (t), \qd (t) \big) - \dfrac{d}{dt} \bigg( \dfrac{\partial L}{\partial \bv} \big( \cq (t), \qd (t) \big) \bigg) = 0 .
\end{align}
For the Lagrangian~\eqref{eq:degenerate_lagrangian} the Euler--Lagrange equations yield
\begin{align}\label{eq:euler_lagrange_equations}
\nabla \vartheta (\cq (t)) \cdot \qd (t) - \nabla H(\cq (t)) - \dot{\vartheta} (\cq (t)) = 0 ,
\end{align}
which after computing the time-derivative of $\vartheta$ can be written as
\begin{align}\label{eq:euler_lagrange_equations_matrix}
\obar{\Omega}^{T} (\cq (t)) \, \qd (t) &= \nabla H(\cq (t)) &
& \text{with} &
\obar{\Omega}_{ij} &= \dfrac{\partial \vartheta_{j}}{\partial \bq^{i}} - \dfrac{\partial \vartheta_{i}}{\partial \bq^{j}} .
\end{align}
The skew-symmetric matrix $\obar{\Omega}$ plays an important role as it holds the components of the symplectic form $\obar{\omega}$ on $\mf{M}$. 
Let us note that in principle $\obar{\Omega}$ can be of odd dimension, in which case the corresponding two-form $\obar{\omega}$ is degenerate and therefore does not classify as a symplectic form. However, most degenerate Lagrangians of the form~\eqref{eq:degenerate_lagrangian}, including all of the examples discussed in Section~\ref{sec:numerical_experiments}, originate from some kind of coordinate transformation of a canonical system, possibly followed by some reduction procedure, which always results in a system of even degree. We therefore assume in the following that the system under consideration is of even dimension and hence symplectic.
Details of the symplectic structure will be discussed in Section~\ref{sec:symplecticity}.

Equation~\eqref{eq:euler_lagrange_equations_matrix} has the structure of a noncanonical Hamiltonian system on $\mf{M}$, characterized by the skew-symmetric matrix $\obar{\Omega}$ and the Hamiltonian $H$.
For such noncanonical Hamiltonian systems no general geometric integrators are known.
In principal it is possible to use the Darboux theorem to find canonical coordinates and then apply some canonical symplectic integrator.
In practice, however, the construction of such Darboux coordinates tends to be a non-trivial task and is often possible only locally but not globally.
Our strategy will thus be to reformulate the system as a canonical Hamiltonian system by adding canonical conjugate momenta, thus doubling the size of the solution space, and restricting the numerical solution to the physical subspace of this extended solution space.
The geometrical foundation of this procedure is the theory of Dirac constraints.

\subsection{Dirac Constraints}

Degenerate systems of the form~\eqref{eq:degenerate_lagrangian} can also be formulated in terms of the phasespace trajectory $(\cq, \cp)$ in the cotangent bundle $\cb{\mf{M}}$, subject to a primary constraint in the sense of Dirac, determined by the function $\phi : \cb{\mf{M}} \rightarrow \rsp^{d}$, given by
\begin{align}\label{eq:dirac_constraint}
\phi (\bq, \bp) = \bp - \vartheta(\bq) = 0 ,
\end{align}
and originating from the fibre derivative $\mbb{F} L : \tb{\mf{M}} \rightarrow \cb{\mf{M}}$,
\begin{align}\label{eq:fibre_derivative_general}
\mbb{F} L (\bv_{\bq}) \cdot \bw_{\bq} = \dfrac{d}{d\eps} \bigg\vert_{\eps=0} L(\bv_{\bq} + \eps \bw_{\bq}) ,
\end{align}
where $\bv_{\bq} = (\bq, \bv)$ and $\bw_{\bq} = (\bq, \bw)$ denote two points in $\tb{\mf{M}}$ which share the same base point $\bq$ and are thus elements of the same fibre of $\tb{\mf{M}}$.
By acting point-wise for each $t$, the fibre derivative maps the curve $(\cq, \qd)$ in the tangent bundle $\tb{\mf{M}}$ into the curve $(\cq, \cp)$ in the cotangent bundle $\cb{\mf{M}}$,
\begin{align}\label{eq:fibre_derivative}
(\cq (t), \cp (t)) = \left( \cq (t), \dfrac{\partial L}{\partial \bv} (\cq (t), \qd (t)) \right) = (\cq (t), \vartheta (\cq (t))) ,
\end{align}
where the last equality follows for Lagrangians of the form~\eqref{eq:degenerate_lagrangian}.
The Dirac constraint arising from the degenerate Lagrangian restricts the dynamics to the submanifold
\begin{align}\label{eq:constraint_submanifold}
\Delta = \big\{ (\bq, \bp) \in \cb{\mf{M}} \, \big\vert \, \phi (\bq, \bp) = 0  \big\} \subset \cb{\mf{M}} .
\end{align}
In the preceding and the following, we assume that the Lagrangian is degenerate in all velocity components, that is, the Lagrangian is either linear or singular in each component of $\bv$, so that
\begin{align}
\dfrac{\partial^{2} L}{\partial \bv^{i} \, \partial \bv^{j}} = 0
\hspace{3em}
\text{for all $1 \leq i,j \leq d$.}
\end{align}
For instructive reasons, however, assume for a moment that the Lagrangian is degenerate in only $m<d$ components of $\bv$ and, e.g., quadratic in the other $d-m$ components. That is to say we can write
\begin{align}
\cp (t)
= \big( \beta_{1} (\cq (t), \qd (t)), \hdots, \beta_{d-m} (\cq (t), \qd (t)), \, \vartheta_{d-m+1} (\cq (t)), \hdots, \vartheta_{d} (\cq (t)) \big)^{T} ,
\end{align}
where
\begin{align}
\dfrac{\partial L}{\partial \bv^{i}} (\cq (t), \qd (t)) = \begin{cases}
\beta_{i} (\cq (t), \qd (t)) & 1 \leq i \leq d - m , \\
\vartheta_{i} (\cq (t)) & d-m+1 \leq i \leq d . \\
\end{cases}
\end{align}
We can then denote coordinates in $\Delta$ by $(\bq^{i}, \basis{\pi}^{j})$ with $1 \leq i \leq d$ and $1 \leq j \leq d-m$, where the $\pi^{i}$ denote those momenta which are ``free'', i.e., not determined by the Dirac constraint. The inclusion map $i : \Delta \rightarrow \cb{\mf{M}}$ can then be written as
\begin{align}\label{eq:dirac_inclusion_genereal}
i : (\bq, \basis{\pi}) \mapsto (\bq, \basis{\pi}, \vartheta (\bq)) .
\end{align}
In the fully degenerate case, however, we have $m=d$, so that the configuration manifold $\mf{M}$ and the constraint submanifold $\Delta$ are isomorphic and we can label points in $\Delta$ by the same $\bq$ we use to label points in $\mf{M}$.
The inclusion map $i : \Delta \rightarrow \cb{\mf{M}}$ simplifies accordingly and reads
\begin{align}\label{eq:dirac_inclusion}
i : \bq \mapsto (\bq, \vartheta (\bq)) ,
\end{align}
where it is important to keep in mind that $\bq$ denotes a point in $\Delta$.
The inverse operation is given by the projection $\pi_{\Delta} : \cb{\mf{M}} \rightarrow \Delta$, defined such that $\pi_{\Delta} \circ i = \id$.

As we are lacking a general framework for constructing structure-preserving numerical algorithms for noncanonical Hamiltonian systems on $\mf{M}$, we will construct such algorithms on $i(\Delta)$. This can be achieved by using canonically symplectic integrators on $\cb{\mf{M}}$ and assuring that their solution stays on $i(\Delta)$. To this end we will employ various projection methods, as discussed in Section~\ref{sec:projection}.

\subsection{Augmented Hamiltonian Approach}

Hamilton's form of the equations of motion for a degenerate Lagrangian system~\eqref{eq:degenerate_lagrangian} can be derived from the phasespace action
\begin{align}\label{eq:phasespace_action}
\obar{\mcal{A}} [q,p,\lambda] = \int_{0}^{T} \left[ \cp (t) \cdot \qd (t) - \obar{H} (\cq (t), \cp (t), \cl (t)) \right] dt ,
\end{align}
with the augmented Hamiltonian $\obar{H} : \cb{\mf{M}} \times \rsp^{d} \rightarrow \rsp$ defined as
\begin{align}\label{eq:augmented_hamiltonian}
\obar{H} (\bq, \bp, \lambda) = H(\bq) + \phi (\bq,\bp) \cdot \bl .
\end{align}
Applying Hamilton's principle of stationary action to~\eqref{eq:phasespace_action} results in the following index two differential-algebraic system of equations (see e.g.~\citet{HairerLubichRoche:1989} for a definition of the notion of index),
\begin{subequations}\label{eq:euler_lagrange_dae}
\begin{align}
\label{eq:euler_lagrange_dae_1}
\qd (t) &= \phi_{\bp}^{T} (\cq (t), \cp (t)) \, \cl (t) , \\
\label{eq:euler_lagrange_dae_2}
\pd (t) &= - H_{\bq} (\cq (t)) - \phi_{\bq}^{T} (\cq (t), \cp (t)) \, \cl (t) , \\
\label{eq:euler_lagrange_dae_3}
0 &= \phi (\cq (t), \cp (t)) .
\end{align}
\end{subequations}
Here, subscripts $\bq$ and $\bp$ denote partial derivatives with respect to the coordinates of $\cb{\mf{M}}$. For the constraint function $\phi$ these derivatives are explicitly written as
\begin{align}
\phi_{\bq} &= \begin{pmatrix}
\dfrac{\partial \phi_{1}}{\partial \bq^{1}} & \dfrac{\partial \phi_{1}}{\partial \bq^{2}} & \hdots & \dfrac{\partial \phi_{1}}{\partial \bq^{d}} \\
\vdots & \vdots & \ddots & \vdots \\
\dfrac{\partial \phi_{d}}{\partial \bq^{1}} & \dfrac{\partial \phi_{d}}{\partial \bq^{2}} & \hdots & \dfrac{\partial \phi_{d}}{\partial \bq^{d}} \\
\end{pmatrix} &
& \text{and} &
\phi_{\bp} &= \begin{pmatrix}
\dfrac{\partial \phi_{1}}{\partial \bp^{1}} & \dfrac{\partial \phi_{1}}{\partial \bp^{2}} & \hdots & \dfrac{\partial \phi_{1}}{\partial \bp^{d}} \\
\vdots & \vdots & \ddots & \vdots \\
\dfrac{\partial \phi_{d}}{\partial \bp^{1}} & \dfrac{\partial \phi_{d}}{\partial \bp^{2}} & \hdots & \dfrac{\partial \phi_{d}}{\partial \bp^{d}} \\
\end{pmatrix} ,
\end{align}
so that in $\phi_{\bq}^{T} \lambda$ and $\phi_{\bp}^{T} \lambda$, the components of $\lambda$ are contracted with the components of $\phi$, not with the derivatives.
As the Hamiltonian $H$ does not depend on $\bp$ and $\phi = \bp - \vartheta(\bq)$, we find that the first equation reduces to $\qd (t) = \cl (t)$, that is the Lagrange multiplier takes the role of the velocity.
Denoting the trajectory in the cotangent bundle $\cb{\mf{M}}$ by $\cz=(\cp, \cq)$, the equations of motion~\eqref{eq:euler_lagrange_dae} can be rewritten more compactly as
\begin{subequations}\label{eq:euler_lagrange_dae_z}
\begin{align}
\label{eq:euler_lagrange_dae_z1}
\dz (t) &= \Omega^{-1} \, \nabla H (\cz (t)) + \Omega^{-1} \, \nabla \phi^{T} (\cz (t)) \, \cl (t) , \\
\label{eq:euler_lagrange_dae_z2}
0 &= \phi (\cz (t)) ,
\end{align}
\end{subequations}
where $\nabla$ denotes the derivatives with respect to $\bz = (\bq, \bp)$.

\subsection{Hamilton--Pontryagin Principle}
\label{sec:hamilton_pontryagin_principle}

The phasespace action principle of the previous section is equivalent to the Hamilton--Pontryagin principle~\cite{Yoshimura:2006a, Yoshimura:2006b} on $\tb{\mf{M}} \oplus \cb{\mf{M}}$, given by
\begin{align}\label{eq:hamilton_pontryagin_action}
\delta \int_{0}^{T} \big[ L (\cq (t), \cv (t)) + \cp (t) \cdot ( \qd (t) - \cv (t) ) \big] \, dt = 0 .
\end{align}
Here, the dynamics of the system are described by the evolution of $(\cq, \cv, \cp)$, which constitute a trajectory in the Pontryagin bundle $\tb{\mf{M}} \oplus \cb{\mf{M}}$.
If the Lagrange multiplier $\cl$ is replaced by the velocity $\cv$ it is easy to verify that the Lagrangian is related to the augmented Hamiltonian~\eqref{eq:augmented_hamiltonian} by $L (\cq, \cv) = \cp \cdot \cv - \obar{H} (\cq, \cp, \cv)$, and that the Hamilton--Pontryagin principle~\eqref{eq:hamilton_pontryagin_action} is equivalent to the phasespace action principle $\delta \obar{\mcal{A}} (\cq, \cp, \cv) = 0$ with the augmented action $\obar{\mcal{A}}$ given in~\eqref{eq:phasespace_action}.
Computing variations of~\eqref{eq:hamilton_pontryagin_action}, where $\cq$, $\cv$ and $\cp$ are all varied independently and only restricted in that the variations of $\cq$ have to vanish at the endpoints, we obtain the implicit Euler--Lagrange equations,
\begin{align}
\qd (t) &= \cv (t) , &
\cp (t) &= \dfrac{\partial L}{\partial \basis{v}} (\cq (t), \cv (t)) , &
\pd (t) &= \dfrac{\partial L}{\partial \bq} (\cq (t), \cv (t)) ,
\end{align}
which are easily seen to be equivalent to~\eqref{eq:euler_lagrange_equations_general} and~\eqref{eq:euler_lagrange_dae}.
Here, the Dirac constraint $\phi(\cq (t), \cp (t)) = 0$ appears quite naturally as one of the equations of motion, which suggests that the Hamilton--Pontryagin principle might be the natural starting point for the discretization of degenerate Lagrangian systems.
That this is not necessarily the case will be discussed in Section~\ref{sec:discrete_hamilton_pontryagin_principle}.

\subsection{Symplecticity}
\label{sec:symplecticity}

Our aim is to construct methods which retain the symplecticity of the integrator as well as its momentum maps.
Care has to be taken, when stating that the variational integrator and the projection are symplectic.
The continuous system preserves two symplectic forms, the canonical two-form $\omega$ on $\Delta$, but also the noncanonical two-form $\obar{\omega}$ on $\mf{M}$ defined by
\begin{align}\label{eq:noncanonical_symplectic_matrix}
\obar{\omega} &= \ext \vartheta = \tfrac{1}{2} \obar{\Omega}_{ij} \, \ext \bq^{i} \wedge \ext \bq^{j} &
& \text{with} &
\obar{\Omega}_{ij} &= \dfrac{\partial \vartheta_{j}}{\partial \bq^{i}} - \dfrac{\partial \vartheta_{i}}{\partial \bq^{j}} .
\end{align}
The matrix $\obar{\Omega}$ is the noncanonical symplectic matrix which we already encountered in the equations of motion~\eqref{eq:euler_lagrange_equations_matrix}.
The function $\vartheta$ is interpreted as a one-form on $\mf{M}$, in coordinates given by $\vartheta = \vartheta_{i} (\bq) \, \ext \bq^{i}$.
In principle it is possible that $\obar{\omega}$ is degenerate, namely in the case of a system of odd dimension $d$. Then $\obar{\omega}$ is not a symplectic form but a presymplectic form. Most of the following discussion also holds in this case. However, in almost all examples of practical relevance the configuration space is even-dimensional. For this reason, we will always refer to $\obar{\omega}$ as symplectic form.

Note that $\obar{\omega}$ is not the symplectic form $\omega_{L}$ on $\tb{\mf{M}}$ originating from the boundary terms in the action principle~\eqref{eq:hamiltons_action_principle_2}.
Besides leading to the equations of motion, the variational principle provides a direct and natural way to derive the fundamental geometric structures of classical mechanics.
For this derivation, the boundary conditions $\cc (0) = \cc (T) = 0$ are relaxed, while the time interval $[0,T]$ is kept fixed.
Thus the variational principle reads
\begin{align}\label{eq:action_principle_cartan_form}
\dfrac{d}{d\eps} \mcal{A} [\cq^{\eps}] \bigg\vert_{\eps = 0}
&= \int_{0}^{T} \bigg[ \dfrac{\partial L}{\partial \bq} \big( \cq (t), \qd (t) \big) \cdot \cc (t) + \dfrac{\partial L}{\partial \bv} \big( \cq (t), \qd (t) \big) \cdot \cd (t) \bigg] \, dt  + \bigg[ \dfrac{\partial L}{\partial \bv} \cdot \cc \bigg]_{t_{1}}^{t_{2}} .
\end{align}
where the variations $\cc (t)$ do not vanish at the boundary point, so that the last term on the right hand side does not vanish.
This last term corresponds to a linear pairing of the function $\partial L / \partial \bv$, which in general is a function of $(\cq, \qd)$, with the tangent vector $\cc$.
The boundary term in~\eqref{eq:hamiltons_action_principle_2} can be written as $\pair{\Theta_{L}}{\delta q}\vert_{0}^{T}$, where $\Theta_{L}$ is the so-called Lagrangian one-form or Cartan one-form, in coordinates given by
\begin{align}\label{eq:cartan_one_form}
\Theta_{L} = \dfrac{\partial L}{\partial \bv^{i}} \, \ext \bq^{i} ,
\end{align}
One could be tempted to regard $\partial L / \partial \bv$ as a one-form on $\mf{M}$ as it only has a component in $\ext \bq$. The same way $\cc$ could be regarded as a tangent vector on $\mf{M}$.
However, in general $\partial L / \partial \bv$ is a function of $(\bq, \qd)$ and therefore clearly a function on $\tb{\mf{M}}$.
The exterior derivative of the Lagrangian one-form gives the Lagrangian two-form, also referred to as the symplectic two-form,
\begin{align}
\label{eq:mechanics_lagrangian_symplectic_two_form}
\omega_{L} = \ext \Theta_{L} ,
\end{align}
given in coordinates by
\begin{align}
\omega_{L}
&= \dfrac{\partial^{2} L}{\partial \bq^{i} \, \partial \bv^{j}} \, d \bq^{i} \wedge d \bq^{j}
 + \dfrac{\partial^{2} L}{\partial \bv^{i} \, \partial \bv^{j}} \, d \bv^{i} \wedge d \bq^{j} .
\end{align}
For more details on this derivation see e.g. \citet{MarsdenRatiu:2002}.
As the Lagrangian~\eqref{eq:degenerate_lagrangian} is degenerate, so is the corresponding symplectic matrix $\Omega_{L}$, which can be written in block form as
\begin{align}
\Omega_{L} =
\begin{pmatrix}
\obar{\Omega}  & \mbb{0} \\
\mbb{0} & \mbb{0} \\
\end{pmatrix} .
\end{align}
We recognize the upper left block, which corresponds to the noncanonical symplectic matrix $\obar{\Omega}$ on $\mf{M}$. When we discuss symplecticity in the following, we are always referring to the noncanonical symplectic form $\obar{\omega}$ or its matrix representation $\obar{\Omega}$.

Preserving the noncanonical symplectic form $\obar{\omega}$ on $\mf{M}$ is equivalent to preserving the canonical symplectic form $\omega = \ext \Theta$ on the embedding of $\Delta$ in $\cb{\mf{M}}$.
Denoting coordinates on $\cb{\mf{M}}$ by $\bz = (\bq, \bp)$, the canonical one-form $\Theta$ and the symplectic two-form $\omega$ can be written in coordinates as
\begin{align}
\Theta = p_{i} \, \ext q^{i} , \qquad
\omega = \tfrac{1}{2} \Omega_{ij} \, \ext \bz^{i} \wedge \ext \bz^{j} = \ext \bp_{i} \wedge \ext \bq^{i} ,
\end{align}
with $\Omega$ the canonical symplectic matrix, given by
\begin{align}\label{eq:canonical_symplectic_matrix}
\Omega = \begin{pmatrix}
\mathbb{0} &           -  \mathbb{1} \\
\mathbb{1} & \hphantom{-} \mathbb{0} 
\end{pmatrix} .
\end{align}
On the constraint submanifold we have that $\bp_{i} = \vartheta_{i} (\bq)$, and therefore $\ext \bp_{i} = \ext \vartheta_{i} (\bq)$, so that $\omega$ restricted to $\Delta$ reads
\begin{align}
\omega \vert_{\Delta}
\nonumber
 = \dfrac{\partial \vartheta_{i}}{\partial \bq^{j}} \, \ext \bq^{j} \wedge \ext \bq^{i}
&= \dfrac{1}{2} \dfrac{\partial \vartheta_{i}}{\partial \bq^{j}} \, \ext \bq^{j} \wedge \ext \bq^{i}
 + \dfrac{1}{2} \dfrac{\partial \vartheta_{j}}{\partial \bq^{i}} \, \ext \bq^{i} \wedge \ext \bq^{j} \\
&= \dfrac{1}{2} \left( \dfrac{\partial \vartheta_{j}}{\partial \bq^{i}} - \dfrac{\partial \vartheta_{i}}{\partial \bq^{j}} \right) \ext \bq^{i} \wedge \ext \bq^{j}
 = \obar{\omega} .
\end{align}
Using the inclusion~\eqref{eq:dirac_inclusion}, we can write $\obar{\omega} = i^{*} \omega$.
The preceding arguments thus suggest that in order to construct a numerical algorithm that preserves the noncanonical symplectic form $\obar{\omega}$ on $\mf{M}$, a viable strategy could be to construct a canonically symplectic algorithm on $\cb{\mf{M}}$ whose solution stays on the constraint submanifold $\Delta$.

\subsection{Noether Theorem and Conservation Laws}

On of the most influential results of classical mechanics in the 20th century is the correspondence of point-symmetries of the Lagrangian and conservation laws of the Euler--Lagrange equations established by Emmy Noether~\cite{Noether:1918, KosmannSchwarzbach:2010}.
In the following we will summarize her famous theorem.

Consider a Lagrangian system $L : \tb{\mf{M}} \rightarrow \rsp$ and a one-parameter group of transformations $\{ \sigma^{\eps} : \eps \in B_{0}^{r} , \, \sigma_{0} = \id \}$, where $B_{0}^{r}$ denotes the open ball with radius $r > 0$ centred at $0$.
We denote the transformed trajectory by $\cq^{\eps} = \sigma^{\eps} \circ \cq$ and its time derivative by $\qd^{\eps} = d (\sigma^{\eps} \circ \cq) / dt$ such that $\cq^{0} = q$ and $\qd^{0} = \qd$.
We have a symmetry if the transformation $\sigma^{\eps}$ leaves the Lagrangian $L$ invariant, that is
\begin{align}\label{eq:noether_invariance_condition}
L \big( \cq^{\eps} (t), \qd^{\eps} (t) \big) = L \big( \cq (t), \qd (t) \big)
\hspace{3em}
\text{for all $\eps$ and all $\cq$.}
\end{align}
Taking the $\eps$ derivative of~\eqref{eq:noether_invariance_condition}, we obtain the infinitesimal invariance condition,
\begin{align}\label{eq:noether_infinitesimal_invariance_condition}
\dfrac{d}{d\eps} L \big( \cq^{\eps} (t), \qd^{\eps} (t) \big) \bigg\vert_{\eps=0} = 0 ,
\end{align}
which is equivalent to (\ref{eq:noether_invariance_condition}).
Explicitly computing this $\eps$ derivative, we obtain
\begin{align}\label{eq:noether_infinitesimal_invariance_condition_2}
\dfrac{d}{d\eps} L \big( \cq^{\eps} (t), \qd^{\eps} (t) \big) \bigg\vert_{\eps=0}
&= \dfrac{\partial L}{\partial \bq} \big( \cq (t), \qd (t) \big) \cdot \dfrac{d \sigma^{\eps}}{d\eps} \bigg\vert_{\eps=0} + \dfrac{\partial L}{\partial \bv} \big( \cq (t), \qd (t) \big) \cdot \dfrac{d \dot{\sigma}^{\eps}}{d\eps} \bigg\vert_{\eps=0} = 0 .
\end{align}
Denoting by $V$ the vector field with flow $\sigma^{\eps}$, defined as follows,
\begin{align}
V &= \dfrac{d \sigma^{\eps}}{d\eps} \bigg\vert_{\eps=0} ,
\end{align}
and if $\cq$ solves the Euler--Lagrange equations~\eqref{eq:euler_lagrange_equations_general}, we can rewrite~\eqref{eq:noether_infinitesimal_invariance_condition_2} as
\begin{align}\label{eq:noether_infinitesimal_invariance_condition_3}
\bigg[ \dfrac{d}{dt} \dfrac{\partial L}{\partial \bv} \big( \cq (t), \qd (t) \big) \bigg] \cdot V \big( \cq (t) \big) + \dfrac{\partial L}{\partial \bv} \big( \cq (t), \qd (t) \big) \cdot \dfrac{dV}{dt} \big( \cq (t), \qd (t) \big) &= 0 .
\end{align}
The time derivative of the vector field $V$ is simply computed by the chain rule, so that assuming that the transformation $\sigma^{\eps}$ does not explicitly depend on time it is given by
\begin{align}
\dfrac{dV}{dt} &= \bv^{j} \, \dfrac{\partial V^{i}}{\partial \bq^{j}} \dfrac{\partial}{\partial \bv^{i}} .
\end{align}
The expression in~\eqref{eq:noether_infinitesimal_invariance_condition_3} amounts to a total time derivative of the so-called Noether current, which constitutes the preserved quantity,
\begin{align}\label{eq:noether_theorem}
\dfrac{d}{dt} \big[ P \big( \cq (t), \qd (t) \big) \big] &= 0 , &
P  \big( \cq (t), \qd (t) \big) &= \dfrac{\partial L}{\partial \bv} \big( \cq (t), \qd (t) \big) \cdot V (\cq (t)) .
\end{align}
Thus, the momentum $\partial L / \partial \bv$ in the direction $V$ is conserved along solutions $\cq$ of the Euler--Lagrange equations~\eqref{eq:euler_lagrange_equations_general} obtained from $L$ for all times $t$.
In Section~\ref{sec:numerical_experiments}, we will apply the Noether theorem several times in order to determine the conservation laws for the various examples we will consider.

\section{Variational Integrators}\label{sec:vi}

Variational integrators can be seen as the Lagrangian equivalent of symplectic integrators for Hamiltonian systems. Instead of discretizing the equations of motion, the action integral is discretized, followed by the application of a discrete version of Hamilton's principle of stationary action. This leads to discrete Euler--Lagrange equations (the discrete equations of motion) at once. 
The evolution map that corresponds to the discrete Euler--Lagrange equations is what is called a variational integrator. Such a numerical scheme preserves a discrete symplectic form which originates from the boundary terms in the variation of the discrete action.

The seminal work in the development of a discrete equivalent of classical mechanics was presented by \citet{Veselov:1988, Veselov:1991}.
His method, based on a discrete variational principle, leads to symplectic integration schemes that automatically preserve constants of motion~\cite{WendlandtMarsden:1997, MarsdenWendlandt:1997}.
Comprehensive reviews of variational integrators and discrete mechanics can be found in \citet{MarsdenWest:2001} and \citet{LewMata:2016}, including thorough accounts on the historical development preceding and following the work of \citeauthor*{Veselov:1988}.

In the following we collect some material on variational integrators, specifically on the discrete action principle, the position-momentum form, and on variational Runge--Kutta methods, before discussing the problems that arise when trying to apply the method to degenerate Lagrangians.

\subsection{Discrete Action Principle}

Time will be discretized uniformly, i.e., the time step $h \in \rsp_{+}$ is constant. We thus split the interval $[0,T]$ into a finite sequence of times $\{ t_{n} = nh \; \vert \; n = 0, \hdots, N \}$, where $h = T/N$, so that $t_{N} = T$.
Let us denote the configuration of the discrete system at time $t_{n}$ by $\cq_{n}$, so that $\cq_{n} \approx \cq (t_{n})$, where $\cq (t_{n})$ is the configuration of the continuous system at time $t_{n}$.
Then a discrete trajectory can be written as $\cq_{d} = \{ \cq_{n} \}_{n=0}^{N}$.

The discrete Lagrangian is defined as an approximation of the time integral of the continuous Lagrangian over the interval $\mcal{I}_{n} = (t_{n}, t_{n+1})$, i.e.,
\begin{align}\label{eq:vi_finite_quadrature}
L_{d} (\cq_{n}, \cq_{n+1}) \approx \int_{\mcal{I}_{n}} L \big( \cq_{n,n+1} (t) , \, \qd_{n,n+1} (t) \big) \, dt ,
\end{align}
where $\cq_{n,n+1}$ denotes the solution of the Euler--Lagrange equations~\eqref{eq:euler_lagrange_equations_general} in $\mcal{I}_{n}$.
The specific expression of the discrete Lagrangian is determined by the polynomial approximation of the trajectory and the quadrature rule used to approximate the integral.
The discrete action then becomes merely a sum over the time index of discrete Lagrangians,
\begin{align}\label{eq:vi_finite_action}
\mcal{A}_{d} [\cq_{d}] = \sum \limits_{n=0}^{N-1} L_{d} (\cq_{n}, \cq_{n+1}) .
\end{align}
Using linear interpolation between $\cq_{n}$ and $\cq_{n+1}$ to describe the discrete trajectory, thus approximating $\cq_{n,n+1}$ by
\begin{align}
\cq_{n,n+1} (t) &\approx \cq_{n} \, \dfrac{t_{n+1} - t}{t_{n+1} - t_{n}} + \cq_{n+1} \, \dfrac{t - t_{n}}{t_{n+1} - t_{n}} ,
\end{align}
the velocity $\qd_{n,n+1}$ will be approximated by a simple finite-difference expression, namely
\begin{align}
\qd_{n,n+1} (t) &\approx \dfrac{\cq_{n+1} - \cq_{n}}{t_{n+1} - t_{n}} .
\end{align}
As we assume the time step to be constant, in the following we will just write $h$ instead of $t_{n+1} - t_{n}$.
The quadrature approximating the integral in~\eqref{eq:vi_finite_quadrature} is most often realized by either the trapezoidal rule, leading to the discrete Lagrangian
\begin{align}\label{eq:vi_finite_trapezoidal}
L_{d}^{\text{tr}} (\cq_{n}, \cq_{n+1}) = \dfrac{h}{2} \, L \bigg( \cq_{n}, \dfrac{\cq_{n+1} - \cq_{n}}{h} \bigg) + \dfrac{h}{2} \, L \bigg( \cq_{n+1}, \dfrac{\cq_{n+1} - \cq_{n}}{h} \bigg) ,
\end{align}
or the midpoint rule, leading to the discrete Lagrangian
\begin{align}\label{eq:vi_finite_midpoint}
L_{d}^{\text{mp}} (\cq_{n}, \cq_{n+1}) = h \, L \bigg( \dfrac{\cq_{n} + \cq_{n+1}}{2}, \dfrac{\cq_{n+1} - \cq_{n}}{h} \bigg) .
\end{align}
The configuration manifold of the discrete system is still $\mf{M}$, but the discrete state space is $\mf{M} \times \mf{M}$ instead of $\tb{\mf{M}}$,
such that the discrete Lagrangian $L_{d}$ is a function
\begin{align}\label{eq:vi_finite_discrete_lagrangian}
L_{d} : \mf{M} \times \mf{M} \rightarrow \rsp .
\end{align}

\begin{figure}[tb]
	\centering
	\includegraphics[width=.5\textwidth]{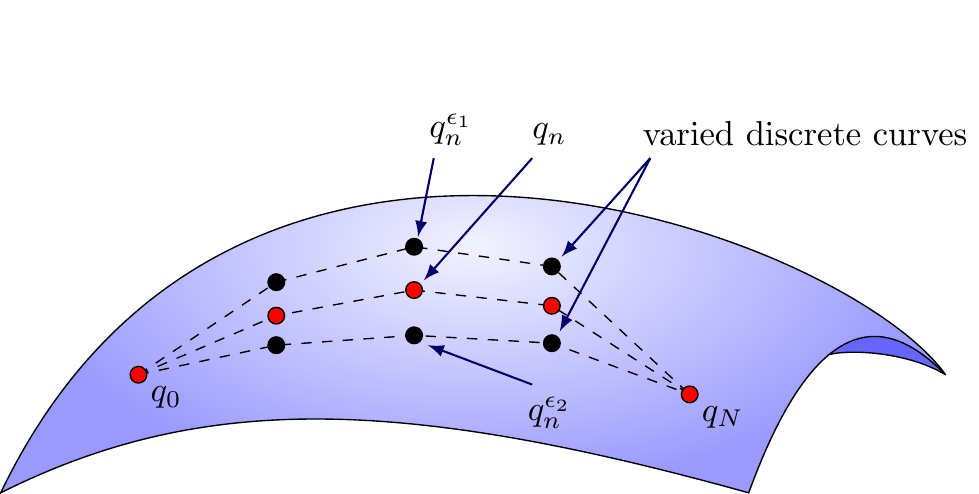}
	\caption{Variations of the discrete trajectory $\cq_{d} = \{ \cq_{n} \}_{n=0}^{N}$.}
	\label{fig:variation_discrete}
\end{figure}

The discrete equations of motion are determined in the same way as the continuous equations of motion~\eqref{eq:euler_lagrange_equations_general}, that is by applying Hamilton's principle of stationary action.
Infinitesimal variations of the discrete trajectories $\cq_{d}$ are given in terms of maps $\cc_{d} : \{ t_{n} \}_{n=0}^{N} \rightarrow \rsp^{d}$, which vanish at $t_{0}$ and $t_{N}$, that is $\cc_{d}(t_{0}) = 0$ and $\cc_{d}(t_{N}) = 0$, and are such that $\cc_{d}(t_{n}) \in \tb[q_{d} (t_{n})]{\mf{M}}$, where $\cc_{d}(t_{n}) = \cc_{n}$ and $q_{d} (t_{n}) = q_{n}$ and we denote such maps by $\cc_{d} = \{ \cc_{n} \}_{n=0}^{N}$.
Discrete one-parameter families of trajectories $\cq_{d}^{\eps} = \{ \cq_{n}^{\eps} \}_{n=0}^{N}$, are defined by
\begin{align}
\cc_{d} = \dfrac{d}{d\eps} \cq_{d}^{\eps} \bigg\vert_{\eps = 0} ,
\end{align}
with the simplest example (c.f., Figure~\ref{fig:variation_discrete}) given by
\begin{align}
\cq_{d}^{\eps} = \{ \cq_n + \eps \cc_n \}_{n=0}^{N} .
\end{align}
Such trajectories are elements of the discrete path space, defined as
\begin{align}
\mf{Q}_{d} (t_{0}, t_{N}, \cq_{0}, \cq_{N}) = \big\{ \cq_{d} : \{ t_{n} \}_{n=0}^{N} \rightarrow \mf{M} \; \big\vert \; \cq_{d} (t_{0}) = \cq_{0}, \, \cq_{d} (t_{N}) = \cq_{N} \big\} .
\end{align}
A necessary condition for the discrete trajectory $\cq_{d} = \cq_{d}^{0}$ making the discrete action~\eqref{eq:vi_finite_action} stationary with respect to all curves $\cq^{\eps}$, is that 
\begin{align}
\dfrac{d}{d\eps} \mcal{A}_{d} [\cq_{d}^{\eps}] \bigg\vert_{\eps = 0}
= \dfrac{d}{d\eps} \sum \limits_{n=0}^{N-1} L_{d} (\cq_{n}^{\eps}, \cq_{n+1}^{\eps}) \bigg\vert_{\eps = 0}
= 0 .
\end{align}
Computing the $\eps$-derivative of the discrete action explicitly, we obtain
\begin{align}\label{eq:vi_finite_variation_01}
\dfrac{d}{d\eps} \mcal{A}_{d} [\cq_{d}^{\eps}] \bigg\vert_{\eps = 0}
&= \sum \limits_{n=0}^{N-1} \big[ D_{1} L_{d} (\cq_{n}, \cq_{n+1}) \cdot \cc_{n} + D_{2} L_{d} (\cq_{n}, \cq_{n+1}) \cdot \cc_{n+1} \big] ,
\end{align}
where $D_{i}$ denotes the slot-derivative with respect to the $i$th argument of $L_{d}$.
What follows corresponds to a discrete integration by parts, i.e., a reordering of the summation.
The $n=0$ term is removed from the first part of the sum and the $n=N-1$ term is removed from the second part,
\begin{align}\label{eq:vi_finite_variation_02}
\dfrac{d}{d\eps} \mcal{A}_{d} [\cq_{d}^{\eps}] \bigg\vert_{\eps = 0}
\nonumber
&= D_{1} L_{d} (\cq_{0}, \cq_{1}) \cdot c_{0} + \sum \limits_{n=1}^{N-1} D_{1} L_{d} (\cq_{n}, \cq_{n+1}) \cdot \cc_{n} \\
&+ \sum \limits_{n=0}^{N-2} D_{2} L_{d} (\cq_{n}, \cq_{n+1}) \cdot \cc_{n+1} + D_{2} L_{d} (\cq_{N-1}, \cq_{N}) \cdot \cc_{N} .
\end{align}
As the variations at the endpoints, $\cc_{0}$ and $\cc_{N}$, vanish, the corresponding terms in the above sum also vanish.
At last, the summation range of the second sum is shifted upwards by one with the arguments of the discrete Lagrangian adapted correspondingly, so that
\begin{align}\label{eq:vi_finite_variation_03}
\dfrac{d}{d\eps} \mcal{A}_{d} [\cq_{d}^{\eps}] \bigg\vert_{\eps = 0}
&=  \sum \limits_{n=1}^{N-1} \big[ D_{1} L_{d} (\cq_{n}, \cq_{n+1}) + D_{2} L_{d} (\cq_{n-1}, \cq_{n}) \big] \cdot \cc_{n} .
\end{align}
Hamilton's principle of least action requires the variation of the discrete action $\mcal{A}_{d}$ to vanish for any choice of $\cc_{n}$. Consequently, the expression in the square brackets of~\eqref{eq:vi_finite_variation_03} has to vanish. This defines the discrete Euler--Lagrange equations
\begin{align}\label{eq:discrete_euler_lagrange_equations}
D_{1} L_{d} (\cq_{n}, \cq_{n+1}) + D_{2} L_{d} (\cq_{n-1}, \cq_{n}) = 0 .
\end{align}
Denoting coordinates on $\mf{M} \times \mf{M}$ by $(\bx_{1}^{1}, \hdots, \bx_{1}^{d}, \bx_{2}^{1}, \hdots, \bx_{2}^{d})$, this can also be written as
\begin{align}\label{eq:discrete_euler_lagrange_equations2}
\dfrac{\partial L_{d}}{\partial \bx_{1}} (\cq_{n}, \cq_{n+1}) + \dfrac{\partial L_{d}}{\partial \bx_{2}} (\cq_{n-1}, \cq_{n}) = 0 .
\end{align}
The discrete Euler--Lagrange equations~\eqref{eq:discrete_euler_lagrange_equations} define an evolution map
\begin{align}\label{eq:vi_finite_evolution_map}
F_{L_{d}}
\; : \; \mf{M} \times \mf{M} \rightarrow \mf{M} \times \mf{M}
\; : \; ( \cq_{n-1}, \cq_{n} ) \mapsto ( \cq_{n}, \cq_{n+1} ) .
\end{align}
Starting from two configurations, $\cq_{0} \approx \cq (t_{0})$ and $\cq_{1} \approx \cq (t_{1} = t_{0} + h)$, the successive solution of the discrete Euler--Lagrange equations (\ref{eq:discrete_euler_lagrange_equations}) for $\cq_{2}$, $\cq_{3}$, etc., up to $\cq_{N}$, determines the discrete trajectory $\cq_{d}$.

For the class of degenerate Lagrangians~\eqref{eq:degenerate_lagrangian} under consideration, the prescription of two sets of initial conditions, $\cq_{0}$ and $\cq_{1}$, is not natural. The continuous Euler--Lagrange equations are ordinary differential equations of first order and therefore need only one set of initial conditions in order to solve the equations. In practice, we face the problem that there is no unique way of determining a second set of initial conditions. All methods will introduce some error that will propagate to the solution and eventually most often lead to a break down of the solution.

\subsection{Position-Momentum Form}

A viable way around this problem appears to be to use the discrete fibre derivative to reformulate the discrete Euler--Lagrange equations~\eqref{eq:discrete_euler_lagrange_equations} in position-momentum form. For regular Lagrangians $L$ this is equivalent to rewriting the continuous Euler--Lagrange equations in the form of Hamilton's equations by using the continuous fibre derivative and Legendre transform.

Given a point $(\cq_{n}, \cq_{n+1})$ in $\mf{M} \times \mf{M}$ and a discrete Lagrangian $L_{d} : \mf{M} \times \mf{M} \rightarrow \rsp$, we define two discrete fibre derivatives, $\mbb{F}^{-} L_{d}$ and $\mbb{F}^{+} L_{d}$, in analogy to the continuous case~\eqref{eq:fibre_derivative} by
\begin{subequations}\label{eq:discrete_fibre_derivative_2}
\begin{align}
\mbb{F}^{-} L_{d} : (\cq_{n}, \cq_{n+1}) &\mapsto (\cq_{n},   \cp_{n})   = \big( \cq_{n}  , - D_{1} L_{d} (\cq_{n}, \cq_{n+1}) \big) , \\
\mbb{F}^{+} L_{d} : (\cq_{n}, \cq_{n+1}) &\mapsto (\cq_{n+1}, \cp_{n+1}) = \big( \cq_{n+1},   D_{2} L_{d} (\cq_{n}, \cq_{n+1}) \big) .
\end{align}
\end{subequations}
With this, the discrete Euler--Lagrange equations~\eqref{eq:discrete_euler_lagrange_equations} can be written as
\begin{align}
\mbb{F}^{+} L_{d} (\cq_{n-1}, \cq_{n}) - \mbb{F}^{-} L_{d} (\cq_{n}, \cq_{n+1}) = 0 ,
\end{align}
which motivates the introduction of the position-momentum form $\obar{F}_{L_{d}} : \cb{\mf{M}} \rightarrow \cb{\mf{M}}$ of the variational integrator~\eqref{eq:vi_finite_evolution_map} by
\begin{subequations}\label{eq:position_momentum_variational_integrator}
\begin{align}
\label{eq:position_momentum_variational_integrator_1}
\cp_{n  } &=           -  D_{1} L_{d} (\cq_{n}, \cq_{n+1}) , \\
\label{eq:position_momentum_variational_integrator_2}
\cp_{n+1} &= \hphantom{-} D_{2} L_{d} (\cq_{n}, \cq_{n+1}) .
\end{align}
\end{subequations}
Given $(\cq_{n}, \cp_{n})$, Equation~\eqref{eq:position_momentum_variational_integrator_1} can be solved for $\cq_{n+1}$. This is generally a nonlinearly implicit equation that has to be solved by some iterative technique like Newton's method.
Equation~\eqref{eq:position_momentum_variational_integrator_2} is an explicit function, so to obtain $\cp_{n+1}$ we merely have to plug in $\cq_{n}$ and $\cq_{n+1}$.
The corresponding Hamiltonian evolution map is
\begin{align}\label{eq:vi_finite_evolution_map_position_momentum}
\obar{F}_{L_{d}}
\; : \; \cb{\mf{M}} \rightarrow \cb{\mf{M}}
\; : \; ( \cq_{n}, \cp_{n} ) \mapsto ( \cq_{n+1}, \cp_{n+1} ) .
\end{align}
In terms of the discrete fibre derivatives it can be equivalently expressed as
\begin{subequations}\label{eq:position_momentum_fibre_derivative}
\begin{align}
\label{eq:position_momentum_fibre_derivative_1}
\obar{F}_{L_{d}} &= \big( \mbb{F}^{-} L_{d} \big) \circ F_{L_{d}} \circ \big( \mbb{F}^{-} L_{d} \big)^{-1} , \\
\label{eq:position_momentum_fibre_derivative_2}
\obar{F}_{L_{d}} &= \big( \mbb{F}^{+} L_{d} \big) \circ F_{L_{d}} \circ \big( \mbb{F}^{+} L_{d} \big)^{-1} , \\
\label{eq:position_momentum_fibre_derivative_3}
\obar{F}_{L_{d}} &= \big( \mbb{F}^{+} L_{d} \big) \circ \big( \mbb{F}^{-} L_{d} \big)^{-1} .
\end{align}
\end{subequations}
In position-momentum form, the variational integrator can be initialized by prescribing an initial position $\cq_{0}$ in conjunction with the corresponding momentum $\cp_{0} = \vartheta (\cq_{0})$. We thus have an exact initialization mechanism, as $\cp_{0}$ constitutes a well-defined second set of initial conditions which can be exactly determined.
Starting with an initial position $\cq_{0}$ and an initial momentum $\cp_{0}$, the repeated solution of~\eqref{eq:position_momentum_variational_integrator} gives the same discrete trajectory $\cq_{d} = \{ \cq_{n} \}_{n=0}^{N}$ as~\eqref{eq:discrete_euler_lagrange_equations}.

\subsection{Discrete Symplectic Structure}

In the following we will show why variational integrators can be considered symplectic integrators and shed some light on the relation with symplectic integrators.
As in the continuous case, we can obtain a discrete Lagrangian one-form by computing the variation of the action for varying endpoints,
\begin{align}
\dfrac{d}{d\eps} \mcal{A}_{d} [\cq_{d}^{\eps}] \bigg\vert_{\eps = 0}
\nonumber
&= \sum \limits_{n=0}^{N-1} \big[ D_{1} \, L_{d} (\cq_{n}, \cq_{n+1}) \cdot \cc_{n} + D_{2} \, L_{d} (\cq_{n}, \cq_{n+1}) \cdot \cc_{n+1} \big] \\
\nonumber
&= \sum \limits_{n=1}^{N-1} \big[ D_{1} \, L_{d} (\cq_{n}, \cq_{n+1}) + D_{2} \, L_{d} (\cq_{n-1}, \cq_{n}) \big] \cdot \cc_{n} \\
&\hspace{3em}
+ D_{1} \, L_{d} (\cq_{0}, \cq_{1}) \cdot \cc_{0} + D_{2} \, L_{d} (\cq_{N-1}, \cq_{N}) \cdot \cc_{N} .
\end{align}
The two latter terms originate from the variation at the boundaries. They form the discrete counterpart of the Lagrangian one-form.
However, there are two boundary terms that define two distinct one-forms on $\mf{M} \times \mf{M}$,
{%
\setlength{\arraycolsep}{2pt}%
\begin{align}
\begin{split}
\begin{array}{ll}
\Theta_{L_{d}}^{-} ( q_{0}  , q_{1} ) \cdot ( \cc_{0}   , \cc_{1} ) &\equiv           -  D_{1} L_{d} (q_{0},   q_{1}) \cdot \cc_{0} , \\
\Theta_{L_{d}}^{+} ( q_{N-1}, q_{N} ) \cdot ( \cc_{N-1} , \cc_{N} ) &\equiv \hphantom{-} D_{2} L_{d} (q_{N-1}, q_{N}) \cdot \cc_{N} .
\end{array}
\end{split}
\end{align}
}%
In general, these one-forms are defined as
\begin{align}\label{eq:vi_finite_discrete_one_form}
\begin{split}
\Theta_{L_{d}}^{-} ( \cq_{n} , \cq_{n+1} ) &\equiv           -  D_{1} L (\cq_{n}, \cq_{n+1}) \, \ext \cq_{n} , \\
\Theta_{L_{d}}^{+} ( \cq_{n} , \cq_{n+1} ) &\equiv \hphantom{-} D_{2} L (\cq_{n}, \cq_{n+1}) \, \ext \cq_{n+1} .
\end{split}
\end{align}
As $\ext L_{d} = \Theta_{L_{d}}^{+} - \Theta_{L_{d}}^{-}$ and $\ext^{2} L_{d} = 0$ one observes that
\begin{align}
\ext \Theta_{L_{d}}^{+} = \ext \Theta_{L_{d}}^{-}
\end{align}
such that the exterior derivative of both discrete one-forms defines the same discrete Lagrangian two-form or discrete symplectic form
\begin{align}\label{eq:vi_finite_discrete_two_form}
\omega_{L_{d}}
&= \ext \Theta_{L_{d}}^{+}
 = \ext \Theta_{L_{d}}^{-}
 = \dfrac{\partial^{2} L_{d} (\cq_{n}, \cq_{n+1})}{\partial \cq_{n } \, \partial \cq_{n+1}} \, \ext \cq_{n} \wedge \ext \cq_{n+1} &
& \text{(no summation over $n$)} .
\end{align}
Now consider the exterior derivative of the discrete action~\eqref{eq:vi_finite_action}.
Upon insertion of the discrete Euler--Lagrange equations (\ref{eq:vi_finite_symplectic_preservation}), it becomes
\begin{align}\label{eq:vi_finite_symplectic_ext_action}
\ext \mcal{A}_{d}
&= D_{1} L_{d} (\cq_{0}, \cq_{1}) \cdot \ext \cq_{0} + D_{2} L_{d} (\cq_{N-1}, \cq_{N}) \cdot \ext \cq_{N}
= \Theta_{L_{d}}^{+} (\cq_{N-1}, \cq_{N}) - \Theta_{L_{d}}^{-} (\cq_{0}, \cq_{1}) .
\end{align}
On the right hand side we find the just defined Lagrangian one-forms (\ref{eq:vi_finite_discrete_one_form}).
Taking the exterior derivative of (\ref{eq:vi_finite_symplectic_ext_action}) gives
\begin{align}\label{eq:vi_finite_symplectic_preservation}
\omega_{L_{d}} (\cq_{0}, \cq_{1}) = \omega_{L_{d}} (\cq_{N-1}, \cq_{N}) ,
\end{align}
where $\cq_{N-1}$ and $\cq_{N}$ are connected with $\cq_{0}$ and $\cq_{1}$ through the discrete Euler--Lagrange equations (\ref{eq:discrete_euler_lagrange_equations}).
Therefore, (\ref{eq:vi_finite_symplectic_preservation}) implies that the discrete symplectic structure $\omega_{L_{d}}$ is preserved while the system advances from $t=0$ to $t=Nh$ according to the discrete equations of motion (\ref{eq:discrete_euler_lagrange_equations}).
As the number of time steps $N$ is arbitrary, the discrete symplectic form $\omega_{L_{d}}$ is preserved at all times of the simulation.
Note that this does not automatically imply that the continuous symplectic structure $\omega_{L}$ is preserved by the variational integrator.
However, as can be seen by comparing~\eqref{eq:vi_finite_discrete_one_form} and~\eqref{eq:position_momentum_variational_integrator}, the discrete one-forms~\eqref{eq:vi_finite_discrete_one_form} correspond to the canonical one-form $\cp_{i} \, \ext \cq^{i}$ under pullback by the discrete fibre derivatives~\eqref{eq:discrete_fibre_derivative_2}. Thus conservation of the discrete symplectic form $\omega_{L_{d}}$ by the discrete Euler--Lagrange equations~\eqref{eq:discrete_euler_lagrange_equations} on $\mf{M} \times \mf{M}$ is equivalent to conservation of the canonical symplectic form $\Omega$ by the position momentum form~\eqref{eq:position_momentum_variational_integrator} on $\cb{\mf{M}}$.

\subsection{Variational Runge--Kutta Methods}

The derivation of higher-order variational integrators in either standard or position-momentum form is rather cumbersome.
A convenient framework for the derivation of integrators of arbitrary order is provided by variational Runge--Kutta methods, which can be seen as a generalization of the position-momentum form.
These methods constitute a special family of symplectic-partitioned Runge--Kutta methods for Lagrangian systems, which are of the form
\begin{subequations}\label{eq:vprk}
\begin{align}
\label{eq:vprk_1}
\rkp_{n,i} &= \dfrac{\partial L}{\partial \bv} (\rkq_{n,i}, \rkv_{n,i}) , &
\rkf_{n,i} &= \dfrac{\partial L}{\partial \bq} (\rkq_{n,i}, \rkv_{n,i}) , \\
\label{eq:vprk_2}
\rkq_{n,i} &= \cq_{n} + h \sum \limits_{j=1}^{s} a_{ij} \, \rkv_{n,j} , &
\rkp_{n,i} &= \cp_{n} + h \sum \limits_{j=1}^{s} \bar{a}_{ij} \, \rkf_{n,j} , \\
\label{eq:vprk_3}
\cq_{n+1} &= \cq_{n} + h \sum \limits_{i=1}^{s} b_{i}  \, \rkv_{n,i} , &
\cp_{n+1} &= \cp_{n} + h \sum \limits_{i=1}^{s} \bar{b}_{i}  \, \dot{P}_{n,i} ,
\end{align}
\end{subequations}
with coefficients satisfying the symplecticity conditions,
\begin{align}\label{eq:vprk_coefficients}
b_{i} \bar{a}_{ij} + \bar{b}_{j} a_{ji} &= b_{i} \bar{b}_{j} &
& \text{and} &
\bar{b}_{i} &= b_{i} .
\end{align}
Here, $s$ denotes the number of internal stages, $a_{ij}$ and $\bar{a}_{ij}$ are the coefficients of the Runge--Kutta method and $b_{i}$ and $\bar{b}_{i}$ the corresponding weights.
Note that while $\rkv_{n,i}$ and $\rkf_{n,i}$ represent velocities and forces at the internal stages, they strictly speaking do not correspond to the time derivatives of $\rkq_{n,i}$ and $\rkp_{n,i}$, respectively.
As $(\rkq_{n,i}, \rkp_{n,i})$ is nothing else than a point in $\cb{\mf{M}}$, the concept of a time derivative of these quantities does not make any sense.
Instead, $\rkq_{n,i}$ and $\rkv_{n,i}$ as well as $\rkp_{n,i}$ and $\rkf_{n,i}$, respectively, denote independent degrees of freedom, which are however related by~\eqref{eq:vprk_2}.

\citet{MarsdenWest:2001} show that variational Runge--Kutta methods~\eqref{eq:vprk} correspond to the position-momentum form~\eqref{eq:position_momentum_variational_integrator} of the discrete Lagrangian
\begin{align}
L_{d} (\cq_{n}, \cq_{n+1}) = h \sum \limits_{i=1}^{s} b_{i} \, L \big( \rkq_{n,i}, \rkv_{n,i} \big) .
\end{align}
They can also be obtained from a discrete action principle similar to the Hamilton--Pontryagin principle presented in Section~\eqref{sec:hamilton_pontryagin_principle}.
For discretizations of Gauss--Legendre type, like the midpoint Lagrangian~\eqref{eq:vi_finite_midpoint}, this is achieved by extremizing a discrete action of the following form~\cite[Section VI.6.3]{HairerLubichWanner:2006} (see also~\cite{BouRabee:2008gx, OberBloebaum:2016}),
\begin{multline}\label{eq:vprk_action_gauss}
\mathcal{A}_{d} =
\sum \limits_{n=0}^{N-1} \Bigg\lgroup
h \sum \limits_{i=1}^{s} b_{i} \, \bigg[ L \big( \rkq_{n,i}, \rkv_{n,i} \big)
+ \rkf_{n,i} \cdot \bigg( \rkq_{n,i} - \cq_{n} - h \sum \limits_{j=1}^{s} a_{ij} \, \rkv_{n,j} \bigg) \bigg] \\
- \cp_{n+1} \cdot \bigg( \cq_{n+1} - \cq_{n} - h \sum \limits_{i=1}^{s} b_{i} \, \rkv_{n,i} \bigg)
\Bigg\rgroup .
\end{multline}
Here, the definition of the generalized coordinates at the internal stages $\rkq_{n,i}$ and the update rule determining $\cq_{n+1}$ are added as constraints with the corresponding momenta $\rkp_{n,i}$ and $\cp_{n+1}$ taking the role of Lagrange multipliers.
Requiring stationarity of the discrete action~\eqref{eq:vprk_action_gauss} for arbitrary variations of $\cq_{n}$, $\cp_{n+1}$, $\cq_{n,i}$, $\rkv_{n,i}$ and $\rkf_{n,i}$, we recover~\eqref{eq:vprk} with the conditions~\eqref{eq:vprk_coefficients} automatically satisfied.
For discretizations of Lobatto--IIIA type like the trapezoidal Lagrangian~\eqref{eq:vi_finite_trapezoidal}, where the first internal stage coincides with the solution at the previous time step, the velocities $\rkv_{n,i}$ are not linearly independent and the discrete action~\eqref{eq:vprk_action_gauss} needs to be augmented by an additional constraint to take this dependence into account (for details see~\citet{OberBloebaum:2016}),
\begin{multline}\label{eq:vprk_action_lobatto}
\mathcal{A}_{d} =
\sum \limits_{n=0}^{N-1} \Bigg\lgroup
h \sum \limits_{i=1}^{s} b_{i} \, \bigg[ L \big( \rkq_{n,i}, \rkv_{n,i} \big)
+ \rkf_{n,i} \cdot \bigg( \rkq_{n,i} - \cq_{n} - h \sum \limits_{j=1}^{s} a_{ij} \, \rkv_{n,j} \bigg) \bigg] \\
- \cp_{n+1} \cdot \bigg( \cq_{n+1} - \cq_{n} - h \sum \limits_{i=1}^{s} b_{i} \, \rkv_{n,i} \bigg)
+ \mu_{n} \cdot \bigg( \sum \limits_{i=1}^{s} d_{i} \rkv_{n,i} \bigg)
\Bigg\rgroup .
\end{multline}
Requiring stationarity of~\eqref{eq:vprk_action_lobatto}, we obtain a modified system of equations,
\begin{subequations}\label{eq:vprk_lobatto}
\begin{align}
\rkp_{n,i} &= \dfrac{\partial L}{\partial \bv} (\rkq_{n,i}, \rkv_{n,i}) , &
\rkf_{n,i} &= \dfrac{\partial L}{\partial \bq} (\rkq_{n,i}, \rkv_{n,i}) , \\
\rkq_{n,i} &= q_{n} + h \sum \limits_{j=1}^{s} a_{ij} \, \rkv_{n,j} , &
\rkp_{n,i} &= p_{n} + h \sum \limits_{j=1}^{s} \bar{a}_{ij} \, \rkf_{n,j} - \mu_{n} \dfrac{d_{i}}{b_{i}} , \\
\cq_{n+1} &= q_{n} + h \sum \limits_{i=1}^{s} b_{i}  \, \rkv_{n,i} , &
\cp_{n+1} &= p_{n} + h \sum \limits_{i=1}^{s} \bar{b}_{i}  \, \rkf_{n,i} , \\
0 &= \sum \limits_{i=1}^{s} d_{i} \rkv_{n,i} ,
\end{align}
\end{subequations}
accounting for the linear dependence of the $\dot{Q}_{n,i}$ and consequently also of the $P_{n,i}$.
The particular values of $d_{i}$ depend on the number of stages $s$ and the definition of the $\cq_{n,i}$~\cite{OberBloebaum:2016}. For two stages, we have $d_{1} = - d_{2}$, so that we can choose, for example, $d_{1} = 1$ and $d_{2} = -1$, and~\eqref{eq:vprk_lobatto} becomes equivalent to the variational integrator of the trapezoidal Lagrangian~\eqref{eq:vi_finite_trapezoidal}. For three stages, we can choose $d = (\tfrac{1}{2}, -1, \tfrac{1}{2})$, and for four stages we can use $d = (+1, -\sqrt{5}, +\sqrt{5}, -1)$.

\subsection{Variational Integrators and Degenerate Lagrangians}

In this section, we want to give an overview of some approaches for discretizing degenerate Lagrangian systems which fail and discuss why they fail.

\begin{figure}[tb]
	\begin{center}
		\includegraphics[width=.32\textwidth]{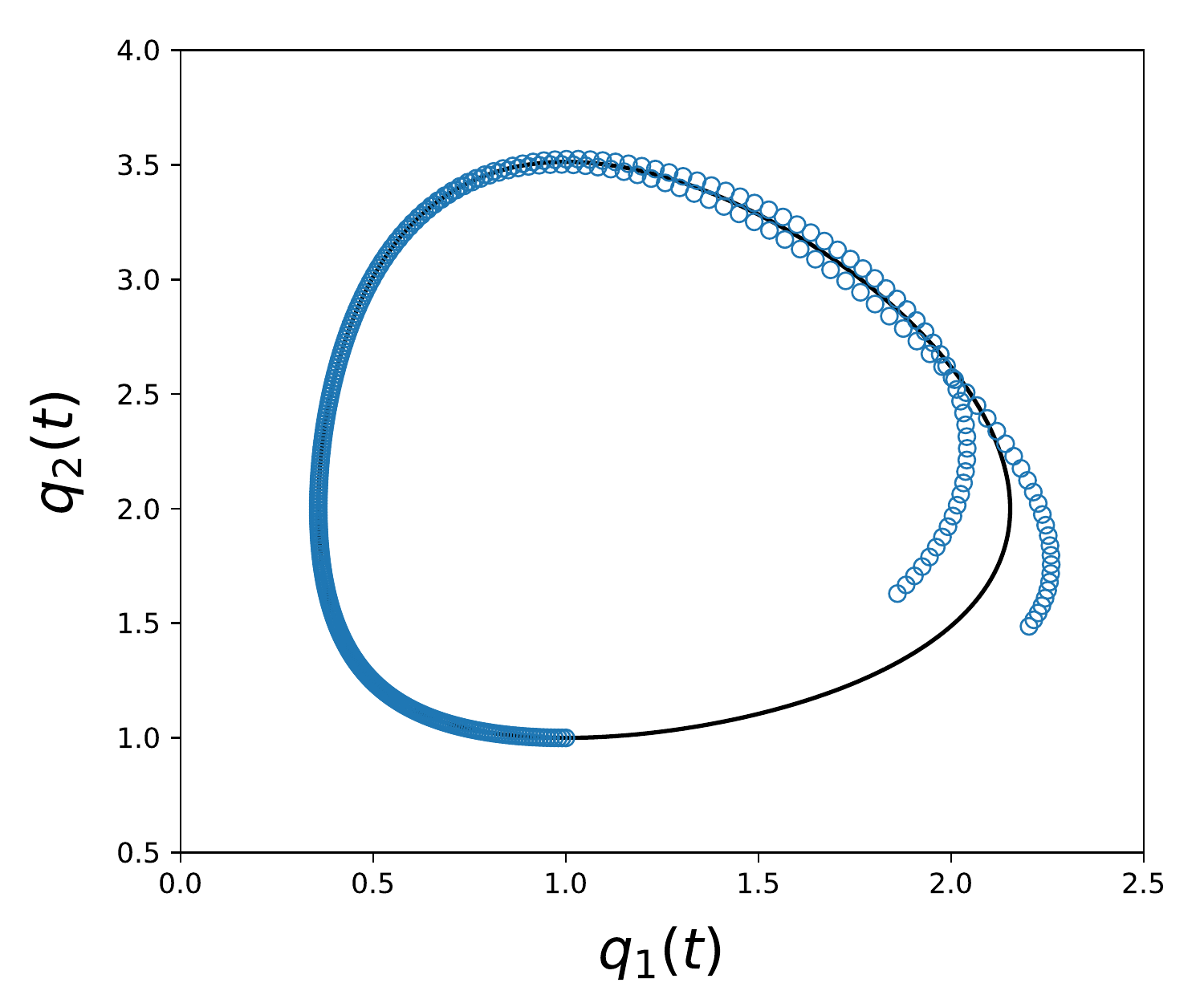}
		\includegraphics[width=.32\textwidth]{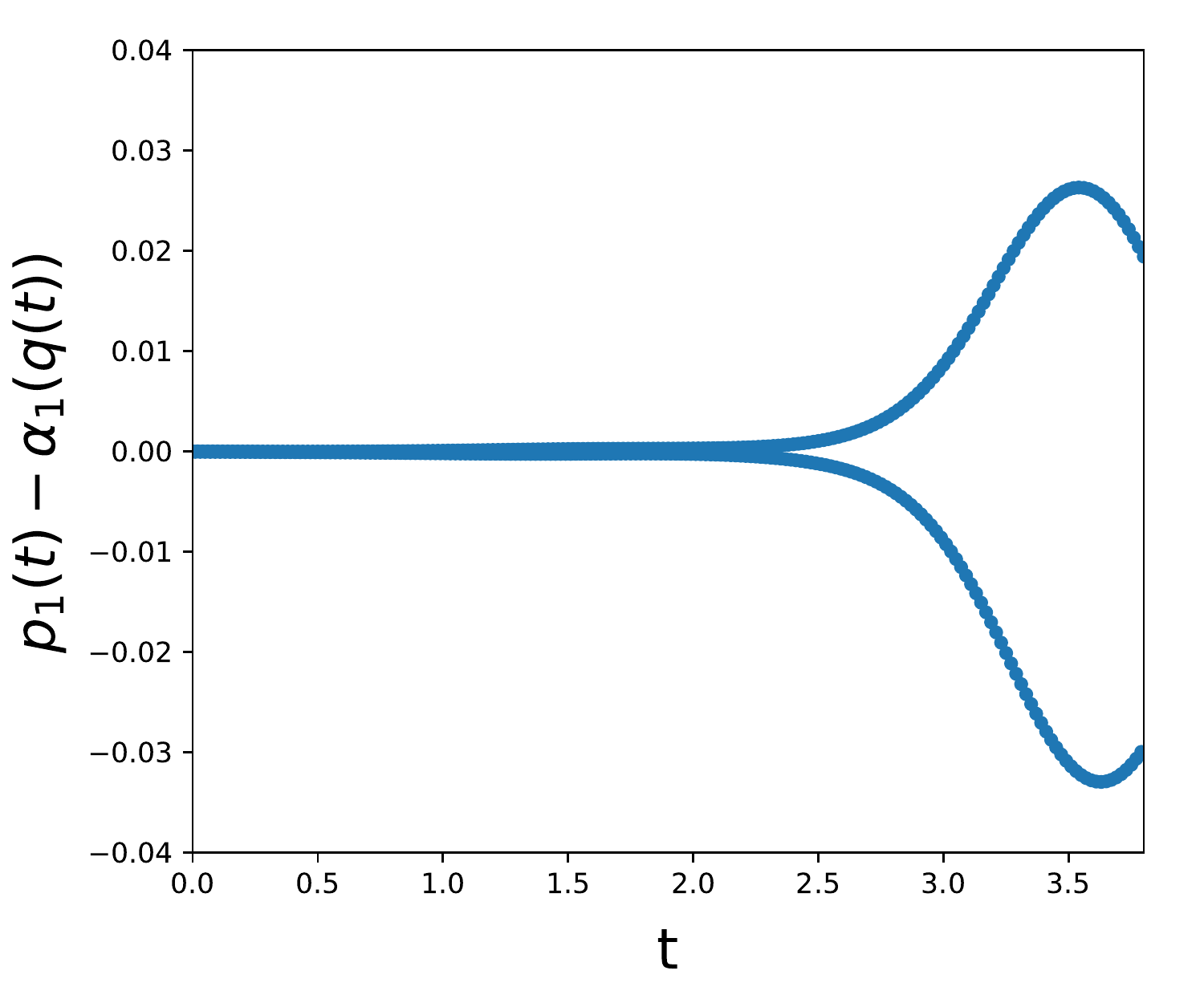}
		\includegraphics[width=.32\textwidth]{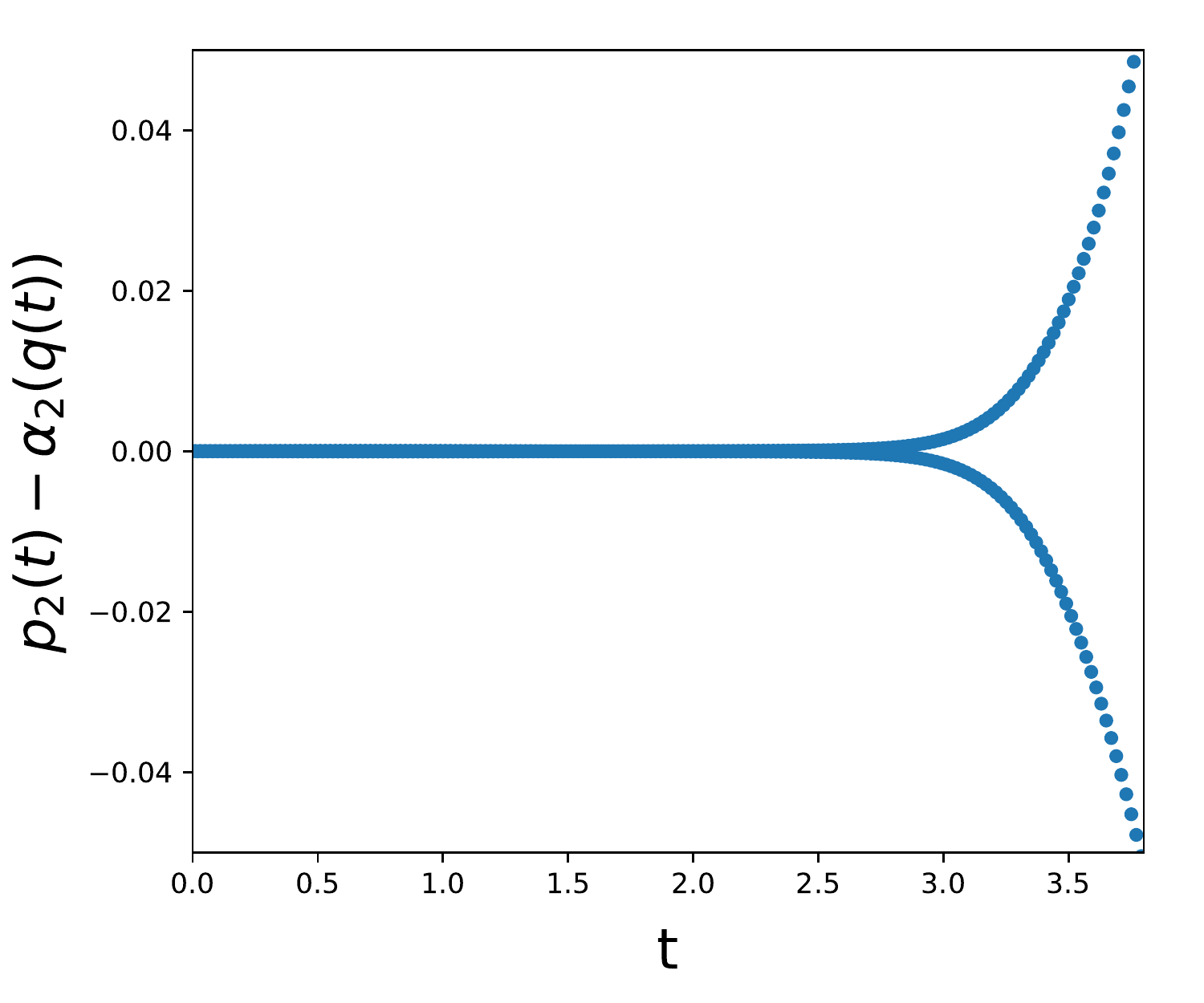}
		\caption{Lotka-Volterra problem in 2D from Section~\ref{sec:app_lotka_volterra_2d}, computed with the trapezoidal variational integrator for a time step of $h = 0.01$. The solution (blue dots) quickly deviates from the reference solution (black line) as the momenta $p_{1}$ and $p_{2}$ diverge from the functions $\vartheta_{1} (q)$ and $\vartheta_{2} (q)$.}
		\label{fig:lotka_volterra_2d_tr}
	\end{center}
\end{figure}

The most obvious option for obtaining a geometric integrator for any Lagrangian systems is to directly discretize the Lagrangian~\eqref{eq:degenerate_lagrangian} and compute the corresponding discrete Euler--Lagrange equations~\eqref{eq:discrete_euler_lagrange_equations}, followed by a discrete fibre derivative~\eqref{eq:discrete_fibre_derivative_2} in order to obtain the position-momentum form~\eqref{eq:position_momentum_variational_integrator} of a variational integrator.
Indeed, it has been shown by \citeauthor{Tyranowski:2014}~\cite{Tyranowski:2014, Tyranowski:2014thesis} that this is a viable strategy for those cases where $\vartheta$ is a linear function.
For cases where $\vartheta$ is a nonlinear function, however, we observe in simulations with such integrators that the numerical solution will in general not satisfy the constraint~\eqref{eq:dirac_constraint}.
Thus the discrete trajectory $\{ (q_{n}, p_{n}) \}_{n=0}^{N}$ will drift away from the constraint submanifold~\eqref{eq:constraint_submanifold}, i.e., even though $(q_{0}, p_{0}) \in i(\Delta)$ we usually find that $(q_{n}, p_{n}) \notin i(\Delta)$ for $n \ge 1$. Whence the solution becomes unphysical, c.f., Figure~\ref{fig:lotka_volterra_2d_tr}.
In the standard form of the variational integrators, which are multi-step methods, this behaviour can be explained in terms of parasitic modes~\cite{Ellison:2015, Ellison:2016}.

In some cases, the solution stays close to the constraint submanifold, which is to say that $\abs{\phi (q_{n}, p_{n})}$ although not zero at least stays bounded for very long times. In such cases variational integrators might still be a viable solution method. However, in general it is not clear to which extend this behaviour depends on the initial conditions. It is easily perceivable that the deviation from the constraint submanifold is bounded for some initial conditions but not for others. And indeed, we observe such behaviour for the example of guiding centre dynamics, described in Section~\ref{app:guiding_centre_dynamics}, where certain particles are found to stay close to the constraint submanifold for very long times while other particles diverge further and further from the constraint submanifold as the simulation proceeds until it eventually crashes.

Let us note that for variational Runge--Kutta methods~\eqref{eq:vprk} and~\eqref{eq:vprk_lobatto} the constraint is automatically satisfied for the internal stages of the method, so that $\phi(Q_{n,i}, P_{n,i}) = 0$ for all $1 \leq i \leq s$, but not for the solution at the next time step $(q_{n+1}, p_{n+1})$. This is due to the fact that the internal stages obey discrete versions of the equations of motion, whereas the final step merely amounts to a numerical quadrature.
It appears, though, that the drift-off problem could be avoided by choosing particular coefficient matrices $a_{ij}$ and $\bar{a}_{ij}$ and weights $b_{i}$ and $\bar{b}_{i}$ in~\eqref{eq:vprk} or~\eqref{eq:vprk_lobatto} such that the last internal stage corresponds to the solution at the next time step, that is $Q_{n,s} = q_{n+1}$ and $P_{n,s} = p_{n+1}$.
It turns out, however, that such a choice is incompatible with the symplecticity conditions~\eqref{eq:vprk_coefficients}, i.e., variational Runge--Kutta methods for which the coefficients and weights are such that the Dirac constraint is automatically satisfied do not exist~\cite{Jay:1996}.

It seems then natural to augment the discrete action~\eqref{eq:vprk_action_gauss} or~\eqref{eq:vprk_action_lobatto} with the constraint evaluated at the solution at the next time step via a Lagrange multiplier similar to~\eqref{eq:phasespace_action}.
This approach and why it fails will be discussed in some detail in the next section.

Yet another seemingly natural approach is to apply a discrete version of the Hamilton--Pontryagin principle~\eqref{eq:hamilton_pontryagin_action} as proposed by \citet{LeokOhsawa:2011}.
We will see in Section~\ref{sec:discrete_hamilton_pontryagin_principle} that the resulting integrator is exactly equivalent to the position-momentum form~\eqref{eq:position_momentum_variational_integrator} and therefore shares the same problems.

After discussing in more detail why these approaches fail, we will present several projection methods for enforcing the constraint $\phi (q_{n+1}, p_{n+1}) = 0$ in Section~\ref{sec:projection}.
In these, we take the solution of a variational integrator and project it to the constraint submanifold~\eqref{eq:constraint_submanifold}.
Although these methods are not strictly-speaking variational or geometric, they lead to useful long-time stable integration algorithms.

\subsection{Augmented Variational Runge--Kutta Methods}

For degenerate Lagrangian systems, we see that the constraint $p = \vartheta(q)$ is automatically enforced at the internal stages. But as the constraint is not enforced at the time steps $n$, the solution tends to drift away from the constraint submanifold~\eqref{eq:constraint_submanifold}.
It seems to be natural to add the constraint to the discrete action~\eqref{eq:vprk_action_gauss} or~\eqref{eq:vprk_action_lobatto}, e.g.,
\begin{multline}\label{eq:vprk_action_gauss_augmented}
\mathcal{A}_{d} =
\sum \limits_{n=0}^{N-1} \bigg\lgroup
h \sum \limits_{i=1}^{s} b_{i} \, \bigg[ L \big( Q_{n,i}, \dot{Q}_{n,i} \big)
- \dot{P}_{n,i} \cdot \bigg( Q_{n,i} - q_{n} - h \sum \limits_{j=1}^{s} a_{ij} \, \dot{Q}_{n,j} \bigg) \bigg]
\\
+ \obar{p}_{n+1} \cdot \bigg( q_{n+1} - q_{n} - h \sum \limits_{i=1}^{s} b_{i} \, \dot{Q}_{n,i} \bigg) + \phi (q_{n+1}, \obar{p}_{n+1}) \cdot \lambda_{n+1}
\bigg\rgroup .
\end{multline}
The momenta are denoted by $\obar{p}_{n+1}$ instead of $p_{n+1}$ in order to highlight the problem with this approach as discussed below.
Variation of~\eqref{eq:vprk_action_gauss_augmented} leads to a modified integrator, which upon defining
\begin{align}
p_{n} &= \obar{p}_{n} - \phi_{\bq}^{T} (q_{n}, \obar{p}_{n}) \, \lambda_{n} &
& \text{and} &
\obar{q}_{n} &= q_{n} - \phi_{\bp}^{T} (q_{n}, \obar{p}_{n}) \, \lambda_{n} ,
\end{align}
can be written as
\begin{subequations}\label{eq:vprk_augmented}
\begin{align}
P_{n,i} &= \dfrac{\partial L}{\partial \bv} (Q_{n,i}, \dot{Q}_{n,i}) , &
\dot{P}_{n,i} &= \dfrac{\partial L}{\partial \bq} (Q_{n,i}, \dot{Q}_{n,i}) , \\
Q_{n,i} &= q_{n} + h \sum \limits_{j=1}^{s} a_{ij} \, \dot{Q}_{n,j} , &
P_{n,i} &= p_{n} + h \sum \limits_{j=1}^{s} \bar{a}_{ij} \, \dot{P}_{n,j} , \\
\obar{q}_{n+1} &= q_{n} + h \sum \limits_{i=1}^{s} b_{i}  \, \dot{Q}_{n,i} , &
\obar{p}_{n+1} &= p_{n} + h \sum \limits_{i=1}^{s} \bar{b}_{i}  \, \dot{P}_{n,i} ,
\end{align}
with projection
\begin{align}
q_{n+1} &= \obar{q}_{n+1} + \phi_{\bp}^{T} (q_{n+1}, \obar{p}_{n+1}) \, \lambda_{n+1} , \\
p_{n+1} &= \obar{p}_{n+1} - \phi_{\bq}^{T} (q_{n+1}, \obar{p}_{n+1}) \, \lambda_{n+1} , \\
0 &= \phi (q_{n+1}, \obar{p}_{n+1}) ,
\end{align}
\end{subequations}
where we assume that $p_{0} = \vartheta(q_{0})$ so that $\obar{p}_{0} = p_{0}$ and $\lambda_{0} = 0$.
We observe that the constraint is enforced at $(q_{n+1}, \obar{p}_{n+1})$, that is the projected coordinate but the unprojected momentum. Practically, we fix the momentum $\obar{p}_{n+1}$ and change the coordinate $\cq_{n+1}$ until it matches the constraint $\phi (q_{n+1}, \obar{p}_{n+1})$.  Then, the momentum is shifted using $\lambda_{n+1}$ as determined from the projection of $\cq_{n+1}$. The result is that while $(q_{n+1}, \obar{p}_{n+1})$ is guaranteed to lie on the constraint submanifold, this is not the case for $(q_{n+1}, p_{n+1})$.

\subsection{Discrete Hamilton--Pontryagin Principle}
\label{sec:discrete_hamilton_pontryagin_principle}

As the Hamilton--Pontryagin principle of Section~\ref{sec:hamilton_pontryagin_principle} provides a very natural setting for degenerate Lagrangian systems, where the Dirac constraint appears as one of the equations of motion, it appears as an appropriate starting point for discretization.
To that end, consider the $(+)$ and $(-)$ discrete Lagrange--Pontryagin principles proposed by~\citet{LeokOhsawa:2011} and given by
\begin{align}
\delta \sum \limits_{n=0}^{N-1} \left[ L_{d} (q_{n}, q_{n}^{+}) + p_{n+1} \cdot (q_{n+1} - q_{n}^{+}) \right] = 0 ,
\end{align}
and
\begin{align}
\delta \sum \limits_{n=0}^{N-1} \left[ L_{d} (q_{n+1}^{-}, q_{n+1}) - p_{n} \cdot (q_{n} - q_{n+1}^{-}) \right] = 0 ,
\end{align}
respectively.
Computing the variations and assuming that the variations of $\cq$ are fixed at the endpoints, $\delta q_{0} = \delta q_{N} = 0$, we obtain
\begin{align}
q_{n+1} &= q_{n}^{+} , &
p_{n} &= - D_{1} L_{d} (q_{n}, q_{n}^{+}) , &
p_{n+1} &= D_{2} L_{d} (q_{n}, q_{n}^{+}) ,
\end{align}
as well as
\begin{align}
q_{n} &= q_{n+1}^{-} , &
p_{n} &= - D_{1} L_{d} (q_{n+1}^{-}, q_{n+1}) , &
p_{n+1} &= D_{2} L_{d} (q_{n+1}^{-}, q_{n+1}) ,
\end{align}
which in both cases are immediately recognized as being equivalent to the position-momentum form~\eqref{eq:position_momentum_variational_integrator} and therefore subject to the same instabilities.

\section{Projection Methods}
\label{sec:projection}

Projection methods are a standard technique for the integration of ordinary differential equations on manifolds~\cite{Hairer:2001, HairerLubichWanner:2006}.
The problem of constructing numerical integrators on manifolds with complicated structure is often difficult and thus avoided by embedding the manifold into a larger space with simple, usually Euclidean structure, where standard integrators can be applied.
Projection methods are then used to ensure that the solution stays on the correct subspace of the extended solution space, as that is usually not guaranteed by the numerical integrator itself.

In the standard projection method, a projection is applied after each step of the numerical algorithm.
Assuming that the initial condition lies in the manifold, the solution of the projected integrator will stay in the manifold.
The problem with this approach is that even though assuming that the numerical integrator is symmetric, the whole algorithm comprised of the integrator and the projection will not be symmetric.
This often leads to growing errors in the solution and consequently a drift in the total energy of the system.
This can be remedied by symmetrizing the projection~\cite{Hairer:2000, Hairer:2001, Chan:2004, HairerLubichWanner:2006}, where the initial data is first perturbed out of the constraint submanifold, before the numerical integrator is applied, and then projected back to the manifold.
This leads to very good long-time stability and improved energy behaviour.

While such projection methods, both standard and symmetric ones, are standard procedures for conserving energy, as well as holonomic and non-holonomic constraints, not much is known about their application to Dirac constraints. 
Some authors consider general differential algebraic systems of index two~\cite{HairerLubichRoche:1989, Ascher:1991, Chan:2002, Chan:2004, Jay:2003, Jay:2006}, the class to which the systems considered here belong, but a discussion of symplecticity seems to be mostly lacking from the literature, aside from some remarks on the conservation of quadratic invariants by the post-projection method of~\citet{Chan:2002}.

In the following, we apply several projection methods (standard, symmetric, symplectic, midpoint) to variational integrators in position-momentum form.
As it turns out, both the standard projection and the symmetric projection are not symplectic. The symmetric projection nevertheless shows very good long-time stability, as it can be shown to be pseudo-symplectic. The symplectic projection method, as the name suggests, is indeed symplectic, although in a generalized sense. The midpoint projection method is symplectic in the usual sense but only for particular integrators.

The general procedure is as follows.
We start with initial conditions $\cq_{n}$ on $\Delta$ (recall that for the particular Lagrangian~\eqref{eq:degenerate_lagrangian} considered here, the configuration manifold $\mf{M}$ and the constraint submanifold $\Delta$ are isomorphic, so that we can use the same coordinates on $\Delta$ as we use on $\mf{M}$).
We compute the corresponding momentum $p_{n}$ by the continuous fibre derivative~\eqref{eq:fibre_derivative}, which yields initial conditions $(q_{n}, p_{n} = \vartheta(q_{n}))$ on $\cb{\mf{M}}$ satisfying the constraint $\phi(q_{n}, p_{n}) = 0$. This corresponds to the inclusion map~\eqref{eq:dirac_inclusion}.
Then, we may or may not perturb these initial conditions off the constraint submanifold by applying a map $(q_{n}, p_{n}) \mapsto (\obar{q}_{n}, \obar{p}_{n})$ which is either the inverse $\mbb{P}^{-1}$ of a projection $\mbb{P} : \cb{\mf{M}} \rightarrow i(\Delta)$ or, in the case of the standard projection of Section~\ref{sec:standard_projection}, just the identity.
The perturbation is followed by the application of some canonically symplectic algorithm $\Psi_{h}$ on $\cb{\mf{M}}$, namely a variational integrator in position-momentum form~\eqref{eq:position_momentum_variational_integrator} or a variational Runge--Kutta method~\eqref{eq:vprk} or \eqref{eq:vprk_lobatto}, in which cases we have that $\Psi_{h} = \big( \mbb{F}^{+} L_{d} \big) \circ \big( \mbb{F}^{-} L_{d} \big)^{-1}$.
In general, the result of this algorithm, $(\obar{q}_{n+1}, \obar{p}_{n+1}) = \Psi_{h} (\obar{q}_{n}, \obar{p}_{n})$, will not lie on the constraint submanifold~\eqref{eq:constraint_submanifold}. Therefore we apply a projection $(\obar{q}_{n+1}, \obar{p}_{n+1}) \mapsto (q_{n+1}, p_{n+1})$ which enforces $\phi (q_{n+1}, p_{n+1}) = p_{n+1} - \vartheta(q_{n+1}) = 0$. As this final result is a point in $i(\Delta)$ it is completely characterized by the value $\cq_{n+1}$.

\begin{SCfigure}
	\includegraphics[width=.4\textwidth]{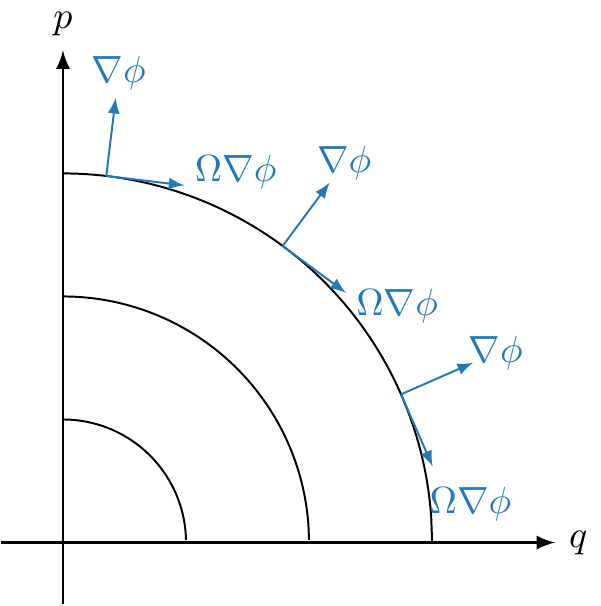}
	\hspace{1em}
	\caption{Gradient of the constraint function $\phi$ orthogonal and $\Omega$-orthogonal to constant surfaces of $\phi(\cq, \cp) = \cp - \sqrt{\cp_{0}^{2} - \cq^{2}}$ for $\cp_{0} \in \{ 1, 2, 3 \}$.}
	\label{fig:omega-orthogonal-projection}
\end{SCfigure}

Let us emphasize that in contrast to standard projection methods, where the solution is projected orthogonal to the constrained submanifold, along the gradient of $\phi$, here the projection has to be $\Omega$-orthogonal (c.f., Figure~\ref{fig:omega-orthogonal-projection}), where $\Omega$ is the canonical symplectic matrix~\eqref{eq:canonical_symplectic_matrix}. That is, denoting by $\lambda$ the Lagrange multiplier, the projection step is given by $\Omega^{-1} \nabla \phi^{T} \cl$ instead of an orthogonal projection $\nabla \phi^{T} \cl$. This appears quite natural when comparing with~\eqref{eq:euler_lagrange_dae}.

Let us also note that, practically speaking, the momenta $\cp_{n}$ and $\cp_{n+1}$ are merely treated as intermediate variables much like the internal stages of a Runge--Kutta method.
The Lagrange multiplier $\cl$, on the other hand, is determined in different ways for the different methods and can be the same or different in the perturbation and the projection. It thus takes the role of an internal variable only for the standard, symmetric projection and midpoint projection, but not for the symplectic projection.

\subsection{Projected Fibre Derivatives}

In the following, we will try to underpin the construction of the various projection methods with some geometric ideas.
We already mentioned several times that the position-momentum form of the variational integrator~\eqref{eq:position_momentum_variational_integrator} suffers from the problem that it does not preserve the constraint submanifold $\Delta$ defined in~\eqref{eq:constraint_submanifold}.
That is, even though it is applied to a point in $i(\Delta)$, it usually returns a point in $\cb{\mf{M}}$, but outside of $i(\Delta)$.
In order to understand the reason for this, let us define $\Delta_{\mf{M}}^{-}$ and $\Delta_{\mf{M}}^{+}$ as the subsets of $\mf{M} \times \mf{M}$ which are mapped into the constraint submanifold $i(\Delta)$ by the discrete fibre derivatives $\mbb{F}^{-} L_{d}$ and $\mbb{F}^{+} L_{d}$, respectively, i.e.,
\begin{subequations}\label{eq:constraint_submanifold_QxQ}
\begin{align}
\Delta_{\mf{M}}^{-} &= \{ (\cq_{n}, \cq_{n+1}) \in \mf{M} \times \mf{M} \, \big\vert \, \mbb{F}^{-} L_{d} (\cq_{n}, \cq_{n+1}) = (\cq_{n}, \cp_{n}) \in i(\Delta) \big\} , \\
\Delta_{\mf{M}}^{+} &= \{ (\cq_{n}, \cq_{n+1}) \in \mf{M} \times \mf{M} \, \big\vert \, \mbb{F}^{+} L_{d} (\cq_{n}, \cq_{n+1}) = (\cq_{n+1}, \cp_{n+1}) \in i(\Delta) \big\} ,
\end{align}
\end{subequations}
or more explicitly,
\begin{subequations}\label{eq:constraint_submanifold_QxQ_coordinates}
\begin{align}
\Delta_{\mf{M}}^{-} &= \{ (\cq_{n}, \cq_{n+1}) \in \mf{M} \times \mf{M} \, \big\vert \, - D_{1} L_{d} (\cq_{n}, \cq_{n+1}) = \vartheta(\cq_{n}) \big\} , \\
\Delta_{\mf{M}}^{+} &= \{ (\cq_{n}, \cq_{n+1}) \in \mf{M} \times \mf{M} \, \big\vert \, D_{2} L_{d} (\cq_{n}, \cq_{n+1}) = \vartheta(\cq_{n+1}) \big\} .
\end{align}
\end{subequations}
A sufficient condition for the position-momentum form of the variational integrator~\eqref{eq:position_momentum_variational_integrator} to preserve the constraint submanifold~\eqref{eq:constraint_submanifold} would be that $\Delta_{\mf{M}}^{-}$ and $\Delta_{\mf{M}}^{+}$ are identical. 
Slightly weaker necessary conditions can be formulated depending on the formulation of the position-momentum form in terms of the discrete Euler--Lagrange equations~\eqref{eq:discrete_euler_lagrange_equations} and the discrete fibre derivative~\eqref{eq:discrete_fibre_derivative_2}, c.f., Equation~\eqref{eq:position_momentum_fibre_derivative}.
For example, considering~\eqref{eq:position_momentum_fibre_derivative_3}, a necessary condition for the position-momentum form to preserve $\Delta$ is that the image of the inverse of $\mbb{F}^{-} L_{d}$, namely $\Delta_{\mf{M}}^{-}$, is in $\Delta_{\mf{M}}^{+}$,
\begin{align}
\big( \mbb{F}^{-} L_{d} \big)^{-1} (i(\Delta)) = \Delta_{\mf{M}}^{-} \subset \Delta_{\mf{M}}^{+} .
\end{align}
Further, from~\eqref{eq:position_momentum_fibre_derivative_1} and~\eqref{eq:position_momentum_fibre_derivative_2} it follows that the image of the variational integrator $F_{L_{d}}$ applied to $\Delta_{\mf{M}}^{-}$ must be in $\Delta_{\mf{M}}^{-}$ and the image of $F_{L_{d}}$ applied to $\Delta_{\mf{M}}^{+}$ must be in $\Delta_{\mf{M}}^{+}$,
\begin{align}
F_{L_{d}} \big( \Delta_{\mf{M}}^{-} \big) &\subset \Delta_{\mf{M}}^{-} , &
F_{L_{d}} \big( \Delta_{\mf{M}}^{+} \big) &\subset \Delta_{\mf{M}}^{+} .
\end{align}
None of these conditions can be guaranteed and they are in general not satisfied. Although $\Delta_{\mf{M}}^{-}$ and $\Delta_{\mf{M}}^{+}$ might have some overlap, they are usually not identical, and the variational integrator, applied to a point in $\Delta_{\mf{M}}^{-}$ or $\Delta_{\mf{M}}^{+}$, does not necessarily result in a point in $\Delta_{\mf{M}}^{-}$ or $\Delta_{\mf{M}}^{+}$, respectively.

In order to construct a modified algorithm which does preserve the constraint submanifold, we compose the discrete fibre derivatives $\mbb{F}^{\pm}$ with appropriate projections $\mbb{P}^{\pm}$,
\begin{align}
\label{eq:position_momentum_projection_1}
(\cq_{n}, \cp_{n  })
&= \big( \mbb{P}^{-} \circ \mbb{F}^{-} L_{d} \big) (\oq_{n}, \oq_{n+1})
 = \mbb{P}_{\cl_{n}^{-}}^{-} \big( \oq_{n}, -  D_{1} L_{d} (\oq_{n}, \oq_{n+1}) \big) , \\
\label{eq:position_momentum_projection_2}
(\cq_{n+1}, \cp_{n+1})
&= \big( \mbb{P}^{+} \circ \mbb{F}^{+} L_{d} \big) (\oq_{n}, \oq_{n+1})
 = \mbb{P}_{\cl_{n+1}^{+}}^{+} \big( \oq_{n+1}, D_{2} L_{d} (\oq_{n}, \oq_{n+1}) \big) , 
\end{align}
so that they take any point in $\mf{M} \times \mf{M}$ to the constraint submanifold $\Delta$.
The Lagrange multiplier $\lambda$ is indicated as subscript and implicitly determined by requiring that the constraint $\phi$ is satisfied by the projected values of $\cq$ and $\cp$.
These projected fibre derivatives will not be a fibre-preserving map anymore, but they will change both $\cq$ and $\cp$, mimicking the continuous equations~\eqref{eq:euler_lagrange_dae}.
Noting that the nullspace of $\mbb{P}_{\cl}$ is the span of $\Omega^{-1} \nabla \phi$, a natural candidate for the projection $\mbb{P}_{\cl}$ is given by
\begin{align}\label{eq:projector}
\mbb{P}_{\lambda}^{\pm} (\oq, \op) : (\cq, \cp) &= (\oq, \op) \pm h \, \Omega^{-1} \nabla \phi^{T} (\cq, \cp) \lambda , &
0 &= \phi(\cq, \cp) ,
\end{align}
so that $( \mbb{P}^{-} \circ \mbb{F}^{-} L_{d} ) (\oq_{n}, \oq_{n+1})$ explicitly reads
\begin{subequations}
\begin{align}
\cq_{n} &= \oq_{n} - h \, \phi_{\bp}^{T} (\cq_{n}, \cp_{n}) \cl_{n}^{-} , \\
\cp_{n} &= - D_{1} L_{d} (\oq_{n}, \oq_{n+1}) + h \, \phi_{\bq}^{T} (\cq_{n}, \cp_{n}) \cl_{n}^{-} , \\
0 &= \phi(\cq_{n}, \cp_{n}) ,
\end{align}
\end{subequations}
and $( \mbb{P}^{+} \circ \mbb{F}^{+} L_{d} ) (\oq_{n}, \oq_{n+1})$ explicitly reads
\begin{subequations}
\begin{align}
\cq_{n+1} &= \oq_{n+1} + h \, \phi_{\bp}^{T} (\cq_{n+1}, \cp_{n+1}) \cl_{n+1}^{+} , \\
\cp_{n+1} &= D_{2} L_{d} (\oq_{n}, \oq_{n+1}) - h \, \phi_{\bq}^{T} (\cq_{n+1}, \cp_{n+1}) \cl_{n+1}^{+} , \\
0 &= \phi(\cq_{n+1}, \cp_{n+1}) .
\end{align}
\end{subequations}
The signs in front of the projections have been chosen in correspondence with the signs of the discrete forces in~\citet[Chapter~3]{MarsdenWest:2001}.
With these projections we obtain all of the algorithms introduced in the following sections, except for the midpoint projection, in a similar fashion to the definition of the position-momentum form of the variational integrator~\eqref{eq:position_momentum_fibre_derivative_3}, as a map $\Delta \rightarrow \Delta$ which can formally be written as
\begin{align}\label{eq:projection_composition_map}
\Phi_{h} = \big( \pi_{\Delta} \circ \mbb{P}^{+} \circ \mbb{F}^{+} L_{d} \big) \circ \big( \pi_{\Delta} \circ \mbb{P}^{-} \circ \mbb{F}^{-} L_{d} \big)^{-1} .
\end{align}
In total, we obtain algorithms which map $\cq_{n}$ into $\cq_{n+1}$ via the steps
\begin{align}
\Delta
\xrightarrow{\pi_{\Delta}^{-1}}
i(\Delta)
\xrightarrow{(\mbb{P}^{-})^{-1}}
\cb{\mf{M}}
\xrightarrow{(\mbb{F}^{-} L_{d})^{-1}}
\mf{M} \times \mf{M}
\xrightarrow{\mbb{F}^{+} L_{d}}
\cb{\mf{M}}
\xrightarrow{\mbb{P}^{+}}
i(\Delta)
\xrightarrow{\pi_{\Delta}}
\Delta ,
\end{align}
where $\pi_{\Delta}^{-1}$ is identical to the inclusion~\eqref{eq:dirac_inclusion}.
The difference of the various algorithms lies in the choice of $\cl_{n}^{-}$ and $\cl_{n+1}^{+}$ as follows
\begin{center}
\begin{tabular}{|l|c|c|}
\hline
Projection & $\cl_{n}^{-}$ & $\cl_{n+1}^{+}$ \\
\hline
Standard   & $0$               & $\cl_{n+1}$ \\
Symplectic & $\cl_{n}$     & $R (\infty) \, \cl_{n+1\hphantom{/2}}$ \\
Symmetric  & $\cl_{n+1/2}$ & $R (\infty) \, \cl_{n+1/2}$ \\
Midpoint   & $\cl_{n+1/2}$ & $R (\infty) \, \cl_{n+1/2}$ \\
\hline
\end{tabular}
\end{center}
For the symmetric, symplectic and midpoint projections, it is important to adapt the sign in the projection according to the stability function $R(\infty)$ of the basic integrator (for details see e.g.~\citet{Chan:2004}).
For the methods we are interested in, namely Runge--Kutta methods, the stability function is given by $R(z) = 1 + z b^{T} (\identity - zA)^{-1} e$ with $e = (1, 1, ..., 1)^{T} \in \rsp^{s}$, and we have $\abs{R(\infty)}=1$ or, more specifically, for Gauss--Legendre methods $R(\infty) = (-1)^{s}$ and for partitioned Gauss--Lobatto IIIA--IIIB and IIIB--IIIA methods we have $R(\infty) = (-1)^{s-1}$.

Let us remark that for the standard projection, the basic integrator and the projection step can be applied independently.
Similarly, for the symplectic projection, the three steps, namely perturbation, numerical integrator, and projection, decouple and can be solved consecutively, as we use different Lagrange multipliers $\cl_{n}$ in the perturbation and $\cl_{n+1}$ in the projection.
For the symmetric projection and the midpoint projection, however, this is not the case. 
There, we used the same Lagrange multiplier $\cl_{n+1/2}$ in both the perturbation and the projection, so that the whole system has to be solved at once, which is more costly.
This also implies that for the projection methods where $\cl_{n}^{-}$ and $\cl_{n+1}^{+}$ are the same (possibly up to a sign due to $R(\infty)$), strictly speaking we cannot write the projected algorithm in terms of a composition of two steps as we did in~\eqref{eq:projection_composition_map}. Instead the whole algorithm has to be treated as one nonlinear map.
The idea of the construction of the methods is still the same, though. Only the midpoint projection of Section~\ref{sec:symplectic_midpoint_projection} needs special treatment. There, the operator $\mbb{P}_{\cl}$ is defined in a slightly more complicated way than in~\eqref{eq:projector}, using different arguments in the projection step, which does not quite fit the general framework outlined here.

\begin{SCfigure}[\sidecaptionrelwidth][tbh]
	\includegraphics[width=.5\textwidth]{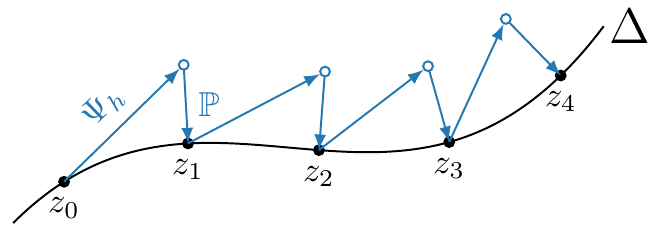}
	\caption{Illustration of the standard projection method. The solution is projected to the constraint submanifold $\Delta$ after each step of the numerical integrator $\Psi_{h}$.}
	\label{fig:standard-projection}
\end{SCfigure}

\subsection{Standard Projection}
\label{sec:standard_projection}

The standard projection method~\citep[section IV.4]{HairerLubichWanner:2006} is the simplest projection method.
Starting from $\cq_{n}$, we use the continuous fibre derivative~\eqref{eq:fibre_derivative} to compute $\cp_{n} = \vartheta (\cq_{n})$. Then we apply some symplectic one-step method $\Psi_{h}$ to $\cz_{n} = (\cq_{n}, \cp_{n})$ to obtain an intermediate solution $\oz_{n+1}$,
\begin{align}
\oz_{n+1} &= \Psi_{h} (\cz_{n}) ,
\end{align}
which is projected onto the constraint submanifold~\eqref{eq:constraint_submanifold} by
\begin{align}\label{eq:orthogonal_projection}
\cz_{n+1} &= \oz_{n+1} + h \, \Omega^{-1} \nabla \phi^{T} (\cz_{n+1}) \cl_{n+1} ,
\end{align}
enforcing the constraint
\begin{align}
0 &= \phi (\cz_{n+1}) .
\end{align}
This projection method, combined with the variational integrators~\eqref{eq:position_momentum_variational_integrator}, is not symmetric, and therefore not reversible. Moreover, it exhibits a drift of the energy, as has been observed before, e.g., for holonomic constraints~\cite{Hairer:2000, Hairer:2001, HairerLubichWanner:2006}.

\subsubsection*{Symplecticity}

In order to verify the symplecticity condition, we write the projection~\eqref{eq:orthogonal_projection} in terms of $(\cp, \cq)$, that is
\begin{subequations}
\begin{align}
\cp_{n+1}^{i} &= \op_{n+1}^{i} - h \, \phi_{\bq^{i}}^{k} (\cq_{n+1}, \cp_{n+1}) \, \cl_{n+1}^{k} , \\
\cq_{n+1}^{i} &= \oq_{n+1}^{i} + h \, \phi_{\bp^{i}}^{k} (\cq_{n+1}, \cp_{n+1}) \, \cl_{n+1}^{k} ,
\end{align}	
\end{subequations}
and assume that $\Psi_{h}$ is a symplectic integrator so that
\begin{align}\label{eq:orthogonal_symplectic_integrator}
\ext \cp_{n}^{i} \wedge \ext \cq_{n}^{i} = \ext \op_{n+1}^{i} \wedge \ext \oq_{n+1}^{i} .
\end{align}
We start by taking the exterior derivative of $\cp_{n+1}$ and $\cq_{n+1}$,
\begin{subequations}\label{eq:orthogonal_projection_differentials}
\begin{align}
\label{eq:orthogonal_projection_differential_p}
\ext \op_{n+1}^{i} &= \ext \cp_{n+1}^{i} + \ext \big( h \, \phi_{\bq^{i}}^{k} (\cq_{n+1}, \cp_{n+1}) \, \cl_{n+1}^{k} \big) , \\
\label{eq:orthogonal_projection_differential_q}
\ext \oq_{n+1}^{i} &= \ext \cq_{n+1}^{i} - \ext \big( h \, \phi_{\bp^{i}}^{k} (\cq_{n+1}, \cp_{n+1}) \, \cl_{n+1}^{k} \big) .
\end{align}	
\end{subequations}
Take the wedge product of the two equations,
\begin{align}
\ext \op_{n+1}^{i} \wedge \ext \oq_{n+1}^{i}
 = \ext \cp_{n+1}^{i} & \wedge \ext \cq_{n+1}^{i}
\nonumber
 - \ext \cp_{n+1}^{i} \wedge \ext \big( h \, \phi_{\bp^{i}}^{k} (\cq_{n+1}, \cp_{n+1}) \, \cl_{n+1}^{k} \big) \\
\nonumber
&+ \ext \big( h \, \phi_{\bq^{i}}^{k} (\cq_{n+1}, \cp_{n+1}) \, \cl_{n+1}^{k} \big) \wedge \ext \cq_{n+1}^{i} \\
&- \ext \big( h \, \phi_{\bq^{i}}^{k} (\cq_{n+1}, \cp_{n+1}) \, \cl_{n+1}^{k} \big) \wedge \ext \big( h \, \phi_{\bp^{i}}^{l} (\cq_{n+1}, \cp_{n+1}) \, \cl_{n+1}^{l} \big) .
\end{align}
The second and third term on the right-hand side become
\begin{align}
\nonumber
   \ext \big( h \, \cl_{n+1}^{k} \phi_{\bp^{i}}^{k} (\cq_{n+1} &, \cp_{n+1}) \big) \wedge \ext \cp_{n+1}^{i}
 + \ext \big( h \, \cl_{n+1}^{k} \phi_{\bq^{i}}^{k} (\cq_{n+1}  , \cp_{n+1}) \big) \wedge \ext \cq_{n+1}^{i} \\
\nonumber
&= \ext \big( h \, \cl_{n+1}^{k} \big[ \phi_{\bp^{i}}^{k} (\cq_{n+1}, \cp_{n+1}) \, \ext \cp_{n+1}^{i} + \phi_{\bq^{i}}^{k} (\cq_{n+1}, \cp_{n+1}) \, \ext \cq_{n+1}^{i} \big] \big) \\
&= \ext \big( h \, \cl^{k} \, \big[ \ext \phi^{k} (\cq_{n+1}, \cp_{n+1}) \big] \big)
 = 0 .
\end{align}
The term in square brackets vanishes as $\phi (\cq_{n+1}, \cp_{n+1}) = 0$ and therefore $\ext \phi (\cq_{n+1}, \cp_{n+1}) = 0$.
Further we have
\begin{align}
\ext \big( h \, \cl_{n+1}^{k} \phi_{\bq^{i}}^{k} (\cq_{n+1}, \cp_{n+1}) \big) 
\nonumber
&= h \, \cl_{n+1}^{k} \phi_{\bq^{i} \bq^{j}}^{k} (\cq_{n+1}, \cp_{n+1}) \, \ext \cq_{n+1}^{j} \\
\nonumber
&+ h \, \cl_{n+1}^{k} \phi_{\bq^{i} \bp^{j}}^{k} (\cq_{n+1}, \cp_{n+1}) \, \ext \cp_{n+1}^{j} \\
\nonumber
&+ h \, \phi_{\bq^{i}}^{k} (\cq_{n+1}, \cp_{n+1}) \, \ext \cl_{n+1}^{k} \\
&= - h \, \cl_{n+1}^{k} \vartheta_{k,ij} (\cq_{n+1}) \, \ext \cq_{n+1}^{j}
   - h \, \vartheta_{k,i} (\cq_{n+1}) \, \ext \cl_{n+1}^{k} ,
\end{align}
and
\begin{align}
\ext \big( h \, \cl_{n+1}^{l} \phi_{\bp^{i}}^{l} (\cq_{n+1}, \cp_{n+1}) \big)
\nonumber
&= h \, \cl_{n+1}^{l} \phi_{\bp^{i} \bq^{j}}^{l} (\cq_{n+1}, \cp_{n+1}) \, \ext \cq_{n+1}^{j} \\
\nonumber
&+ h \, \cl_{n+1}^{l} \phi_{\bp^{i} \bp^{j}}^{l} (\cq_{n+1}, \cp_{n+1}) \, \ext \cp_{n+1}^{j} \\
\nonumber
&+ h \, \phi_{\bp^{i}}^{l} (\cq_{n+1}, \cp_{n+1}) \, \ext \cl_{n+1}^{l} \\
&= h \, \ext \cl_{n+1}^{i} ,
\end{align}
where the terms involving $\phi_{\bp^{i} \bq^{j}}^{k}$ or $\phi_{\bq^{i} \bp^{j}}^{k}$ vanish as $\phi(p,q)$ is separable and the terms involving $\phi_{\bp^{i} \bp^{j}}^{k}$ vanish as $\phi$ is linear in $\bp$.
The wedge product of the two expressions becomes
\begin{multline}
- \ext \big( h \, \phi_{\bq^{i}}^{k} (\cq_{n+1}, \cp_{n+1}) \, \cl_{n+1}^{k} \big) \wedge \ext \big( h \, \phi_{\bp^{i}}^{l} (\cq_{n+1}, \cp_{n+1}) \, \cl_{n+1}^{l} \big) = \\
= \big( h \, \cl_{n+1}^{k} \vartheta_{k,ij} (\cq_{n+1}) \, \ext \cq_{n+1}^{j}
   + h \, \vartheta_{k,i} (\cq_{n+1}) \, \ext \cl_{n+1}^{k} \big) \wedge \big( h \, \ext \cl_{n+1}^{i} \big) = \\
= h^{2} \, \vartheta_{j,i} (\cq_{n+1}) \, \ext \cl_{n+1}^{j} \wedge \ext \cl_{n+1}^{i} .
\end{multline}
The result can be anti-symmetrized so that by using~\eqref{eq:noncanonical_symplectic_matrix} as well as~\eqref{eq:orthogonal_symplectic_integrator}, we obtain
\begin{multline}
\ext \cp_{n}^{i} \wedge \ext \cq_{n}^{i}
 = \ext \cp_{n+1}^{i} \wedge \ext \cq_{n+1}^{i}
 + \dfrac{h^{2}}{2} \obar{\Omega}_{ij} (\cq_{n+1}) \, \ext \cl_{n+1}^{i} \wedge \ext \cl_{n+1}^{j} \\
\vphantom{\dfrac{1}{2}}
 + h^{2} \, \cl_{n+1}^{k} \vartheta_{k,ij} (\cq_{n+1}) \, \ext \cq_{n+1}^{i} \wedge \ext \cl_{n+1}^{j} .
\end{multline}
Using that the constraint $\phi(\cq, \cp) = p - \vartheta(q) = 0$ holds for both $(\cq_{n}, \cp_{n})$ and $(\cq_{n+1}, \cp_{n+1})$, this can be rewritten as
\begin{multline}
\dfrac{1}{2} \obar{\Omega}_{ij} (\cq_{n}) \, \ext \cq_{n}^{i} \wedge \ext \cq_{n}^{j}
 = \dfrac{1}{2} \obar{\Omega}_{ij} (\cq_{n+1}) \, \ext \cq_{n+1}^{i} \wedge \ext \cq_{n+1}^{j}
 - \dfrac{h^{2}}{2} \obar{\Omega}_{ij} (\cq_{n+1}) \, \ext \cl_{n+1}^{i} \wedge \ext \cl_{n+1}^{j} \\
\vphantom{\dfrac{1}{2}}
 - h^{2} \, \cl_{n+1}^{k} \vartheta_{k,ij} (\cq_{n+1}) \, \ext \cq_{n+1}^{i} \wedge \ext \cl_{n+1}^{j} ,
\end{multline}
and we see that the noncanonical symplectic form~\eqref{eq:noncanonical_symplectic_matrix} is not preserved, but in each step accumulates an error $h^{2} \, \cl_{n+1}^{k} \vartheta_{k,ij} (\cq_{n+1}) \, \ext \cq_{n+1}^{i} \wedge \ext \cl_{n+1}^{j} + \tfrac{1}{2} h^{2} \, \obar{\Omega}_{ij} (\cq_{n+1}) \, \ext \cl_{n+1}^{i} \wedge \ext \cl_{n+1}^{j}$.
In numerical simulations, this error accumulation usually manifests itself in form of a drift of the solution and the energy.

\begin{figure}[bht]
\begin{center}
	\subfloat[$R(\infty)=+1$]{
		\includegraphics[width=.48\textwidth]{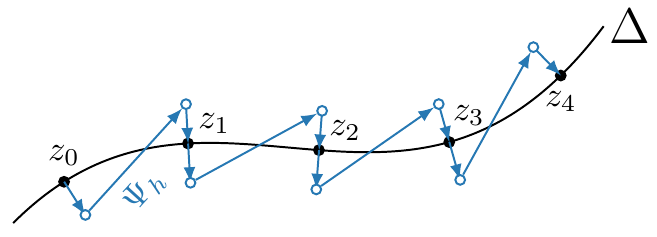}
	}
	\subfloat[$R(\infty)=-1$]{
		\includegraphics[width=.48\textwidth]{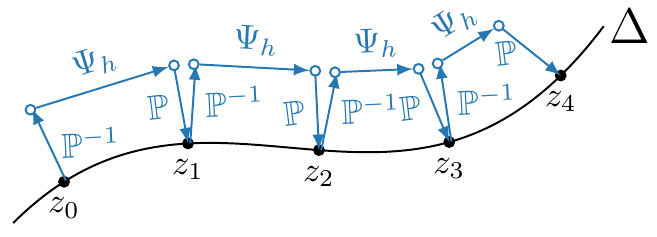}
	}

	\caption{Illustration of the symmetric projection method. The solution is first perturbed off the constraint submanifold $\Delta$, then one step of the numerical integrator $\Psi_{h}$ is performed, and the result is projected back onto $\Delta$.}
	\label{fig:symmetric-projection}
\end{center}
\end{figure}

\subsection{Symmetric Projection}
\label{sec:symmetric_projection}

To overcome the shortcomings of the standard projection, we consider a symmetric projection of the variational Runge--Kutta integrators following~\citet{Hairer:2000, Hairer:2001, Chan:2004}, c.f., Figure~\ref{fig:symmetric-projection} (see also~\cite[section V.4.1]{HairerLubichWanner:2006}). Here, one starts again by computing the momentum $\cp_{n}$ as a function of the coordinates $\cq_{n}$ according to the continuous fibre derivative, which can be expressed with the constraint function as
\begin{subequations}\label{eq:symmetric_symplectic_projection}
\begin{align}
0 &= \phi (\cz_{n}) .
\end{align}
Then the initial value $\cz_{n}$ is first perturbed,
\begin{align}\label{eq:symmetric_projection_pre}
\oz_{n} &= \cz_{n} + h \, \Omega^{-1} \nabla \phi^{T} (\cz_{n}) \, \cl_{n+1/2} , 
\end{align}
followed by the application of some one-step method $\Psi_{h}$,
\begin{align}
\oz_{n+1} &= \Psi_{h} (\oz_{n}) ,
\end{align}
and a projection of the result onto the constraint submanifold,
\begin{align}\label{eq:symmetric_projection_post}
\cz_{n+1} &= \oz_{n+1} + h \, R(\infty) \, \Omega^{-1} \nabla \phi^{T} (\cz_{n+1}) \cl_{n+1/2} ,
\end{align}
which enforces the constraint
\begin{align}
0 &= \phi (\cz_{n+1}) .
\end{align}
\end{subequations}
Here, it is important to note that Lagrange multiplier $\cl_{n+1/2}$ is the same in both the perturbation and the projection step, and to account for the stability function $R(\infty)$ of the basic integrator, as mentioned before.
The algorithm composed of the symmetric projection and some symmetric variational integrator in position-momentum form, constitutes a symmetric map
\begin{align}
\Phi_{h} : \cq_{n} \mapsto \cq_{n+1} ,
\end{align}
where, from a practical point of view, $\cp_{n}$, $\cp_{n+1}$ and $\cl_{n+1/2}$ are treated as intermediate variables.

\subsubsection*{Symplecticity}

In the following, we assume that $\abs{R(\infty)} = 1$. Then, the considerations of symplecticity for the symmetric projection follow along the very same lines as for the standard projection. 
In addition to the projection, we also have to consider the perturbation. Assuming the integrator $\Psi_{h}$ is such that
\begin{align}\label{eq:symmetric_symplectic_integrator}
\ext \op_{n}^{i} \wedge \ext \oq_{n}^{i} = \ext \op_{n+1}^{i} \wedge \ext \oq_{n+1}^{i} ,
\end{align}
we obtain
\begin{multline}\label{eq:symmetric_projection_symplecticity_condition_0}
   \ext \cp_{n}^{i} \wedge \ext \cq_{n}^{i}
 - \dfrac{h^{2}}{2} \obar{\Omega}_{ij} (\cq_{n}) \, \ext \cl^{i}_{n+1/2} \wedge \ext \cl^{j}_{n+1/2}
 - h^{2} \cl_{n+1/2}^{k} \vartheta_{k,ij} (\cq_{n}) \, \ext \cq_{n}^{i} \wedge \ext \cl_{n+1/2}^{j} = \\
= \ext \cp_{n+1}^{i} \wedge \ext \cq_{n+1}^{i}
 - \dfrac{h^{2}}{2} \obar{\Omega}_{ij} (\cq_{n+1}) \, \ext \cl^{i}_{n+1/2} \wedge \ext \cl^{j}_{n+1/2}
 - h^{2} \cl_{n+1/2}^{k} \vartheta_{k,ij} (\cq_{n+1}) \, \ext \cq_{n+1}^{i} \wedge \ext \cl_{n+1/2}^{j} .
\end{multline}
The symmetrically projected integrator admits a certain symmetry in the error terms and can be shown to be pseudo-symplectic~\cite{AubryChartier:1998}.
It is worth to go one step back, and reconsider the derivation that leads to~\eqref{eq:symmetric_projection_symplecticity_condition_0}. By the same considerations as for the standard projection, we obtain
\begin{multline}\label{eq:symmetric_projection_symplecticity_condition}
   \ext \cp_{n}^{i} \wedge \ext \cq_{n}^{i}
 - \ext \big( h \, \phi_{\bq^{i}}^{k} (\cp_{n}, \cq_{n}) \, \cl_{n+1/2}^{k} \big) \wedge \ext \big( h \, \phi_{\bp^{i}}^{l} (\cp_{n}, \cq_{n}) \, \cl_{n+1/2}^{l} \big) = \\
 = \ext \cp_{n+1}^{i} \wedge \ext \cq_{n+1}^{i}
 - \ext \big( h \, \phi_{\bq^{i}}^{k} (\cq_{n+1}, \cp_{n+1}) \, \cl_{n+1/2}^{k} \big) \wedge \ext \big( h \, \phi_{\bp^{i}}^{l} (\cq_{n+1}, \cp_{n+1}) \, \cl_{n+1/2}^{l} \big) .
\end{multline}
We see, that in general the symmetric projection is not symplectic unless
\begin{align}\label{eq:symmetric_projection_symplecticity_condition_1}
\ext \big( h \, \phi_{\bq^{i}}^{k} (\cp_{n}, \cq_{n}) \, \cl_{n+1/2}^{k} \big) &= \ext \big( h \, \phi_{\bq^{i}}^{k} (\cq_{n+1}, \cp_{n+1}) \, \cl_{n+1/2}^{k} \big)
\hspace{3em}
\text{for all $i$,}
\end{align}
as well as
\begin{align}\label{eq:symmetric_projection_symplecticity_condition_2}
\ext \big( h \, \phi_{\bp^{i}}^{k} (\cp_{n}, \cq_{n}) \, \cl_{n+1/2}^{k} \big) &= \ext \big( h \, \phi_{\bp^{i}}^{k} (\cq_{n+1}, \cp_{n+1}) \, \cl_{n+1/2}^{k} \big)
\hspace{3em}
\text{for all $i$,}
\end{align}
that is the initial perturbation is exactly the same as the final projection.
While the first condition~\eqref{eq:symmetric_projection_symplecticity_condition_1} is not obvious, the second condition is immediately seen to be satisfied for $\phi (\cq, \cp) = p - \vartheta (q)$, as $\phi_{\bp} = \unity$, so that~\eqref{eq:symmetric_projection_symplecticity_condition_2} reduces to $\ext \cl_{n+1/2} = \ext \cl_{n+1/2}$. 
If the first condition is not satisfied, though, the method is not symplectic. However, as the error terms to the symplecticity condition appear on both sides of~\eqref{eq:symmetric_projection_symplecticity_condition}, the accumulated error is much smaller than with the standard projection.

Again, using that the constraint $\phi(\cq, \cp) = p - \vartheta(q) = 0$ holds for both $(\cq_{n}, \cp_{n})$ and $(\cq_{n+1}, \cp_{n+1})$, the symplecticity condition~\eqref{eq:symmetric_projection_symplecticity_condition} can be rewritten as
\begin{multline}\label{eq:symmetric_projection_symplecticity_condition_3}
   \dfrac{1}{2} \obar{\Omega}_{ij} (\cq_{n}) \, \big( \ext \cq_{n}^{i} \wedge \ext \cq_{n}^{j}
 - h^{2} \, \ext \cl_{n+1/2}^{i} \wedge \ext \cl_{n+1/2}^{j} \big)
 - h^{2} \cl_{n+1/2}^{k} \vartheta_{k,ij} (\cq_{n}) \, \ext \cq_{n}^{i} \wedge \ext \cl_{n+1/2}^{j} = \\
 = \dfrac{1}{2} \obar{\Omega}_{ij} (\cq_{n+1}) \, \big( \ext \cq_{n+1}^{i} \wedge \ext \cq_{n+1}^{j}
 - h^{2} \, \ext \cl_{n+1/2}^{i} \wedge \ext \cl_{n+1/2}^{j} \big)
 - h^{2} \cl_{n+1/2}^{k} \vartheta_{k,ij} (\cq_{n+1}) \, \ext \cq_{n+1}^{i} \wedge \ext \cl_{n+1/2}^{j} .
\end{multline}
This formulation suggests the following construction.

\begin{SCfigure}[\sidecaptionrelwidth][bth]
	\includegraphics[width=.5\textwidth]{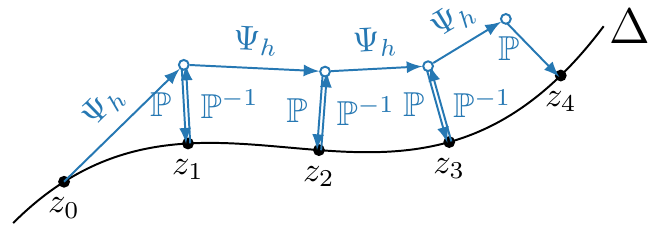}
	\hspace{.5em}
	\caption{Illustration of the post projection method. Starting on the constraint submanifold $\Delta$, the numerical integrator $\Psi_{h}$ moves the solution away from $\Delta$ in the first step. After each step, the solution is projected back onto $\Delta$, but the perturbation at the beginning of each consecutive step is exactly the inverse of the previous projection, so that, practically speaking, the solution is projected back onto $\Delta$ only for output purposes.}
	\label{fig:post-projection}
\end{SCfigure}

\subsection{Symplectic Projection}

If we modify the perturbation~\eqref{eq:symmetric_projection_pre} to use the Lagrange multiplier at the previous time step, $\cl_{n}$, instead of $\cl_{n+1}$, that is we replace~\eqref{eq:symmetric_symplectic_projection} by
\begin{subequations}\label{eq:symplectic_projection}
\begin{align}
\label{eq:symplectic_projection_pre_constraint}
0 &= \phi (\cz_{n}) , \\
\label{eq:symplectic_projection_pre}
\oz_{n\hphantom{+1}} &= \cz_{n} + h \, \Omega^{-1} \nabla \phi^{T} (\cz_{n}) \, \cl_{n} , \\
\oz_{n+1} &= \Psi_{h} (\oz_{n}) , \\
\label{eq:symplectic_projection_post}
\cz_{n+1} &= \oz_{n+1} + h \, R(\infty) \, \Omega^{-1} \nabla \phi^{T} (\cz_{n+1}) \cl_{n+1} , \\
\label{eq:symplectic_projection_post_constraint}
0 &= \phi (\cz_{n+1}) ,
\end{align}
\end{subequations}
the symplecticity condition~\eqref{eq:symmetric_projection_symplecticity_condition_3} is modified as follows,
\begin{multline}\label{eq:symplectic_projection_symplecticity_condition}
   \dfrac{1}{2} \obar{\Omega}_{ij} (\cq_{n}) \, \big( \ext \cq_{n}^{i} \wedge \ext \cq_{n}^{j}
 - h^{2} \, \ext \cl_{n}^{i} \wedge \ext \cl_{n}^{j} \big)
 - h^{2} \cl_{n}^{k} \vartheta_{k,ij} (\cq_{n}) \, \ext \cq_{n}^{i} \wedge \ext \cl_{n}^{j} = \\
 = \dfrac{1}{2} \obar{\Omega}_{ij} (\cq_{n+1}) \, \big( \ext \cq_{n+1}^{i} \wedge \ext \cq_{n+1}^{j}
 - h^{2} \, \ext \cl_{n+1}^{i} \wedge \ext \cl_{n+1}^{j} \big)
 - h^{2} \cl_{n+1/2}^{k} \vartheta_{k,ij} (\cq_{n+1}) \, \ext \cq_{n+1}^{i} \wedge \ext \cl_{n+1}^{j} ,
\end{multline}
implying the conservation of a modified symplectic form $\omega_{\bl}$ defined on an extended phasespace $\mf{M} \times \rsp^{d}$ with coordinates $(\bq, \bl)$ by
\begin{align}\label{eq:symplectic_projection_two_form}
\omega_{\cl}
 = \dfrac{1}{2} \obar{\Omega}_{ij} (\bq) \, \ext \bq^{i} \wedge \ext \bq^{j}
 - \dfrac{h^{2}}{2} \obar{\Omega}_{ij} (\bq) \, \ext \cl^{i} \wedge \ext \cl^{j}
 - h^{2} \cl^{k} \vartheta_{k,ij} (\cq) \, \ext \cq^{i} \wedge \ext \cl^{j} ,
\end{align}
with matrix representation
\begin{align}
\Omega_{\cl} =
\begin{pmatrix}
\obar{\Omega} & - h^{2} \cl \cdot \vartheta_{\bq\bq} \\
h^{2} \cl \cdot \vartheta_{\bq\bq} & - h^{2} \obar{\Omega} \\
\end{pmatrix} .
\end{align}
To this corresponds a modified one-form $\vartheta_{\cl}$, such that $\omega_{\cl} = \ext \vartheta_{\cl}$, given by
\begin{align}\label{eq:symplectic_projection_one_form}
\vartheta_{\cl} = ( \vartheta_{i} (q) - h \, \lambda^{k} \vartheta_{k,i} (q) ) \, ( \ext q^{i} - h \, \ext \lambda^{i} )
\end{align}
As noted by~\citet{Chan:2004}, the modified perturbation~\eqref{eq:symplectic_projection_pre_constraint}-\eqref{eq:symplectic_projection_pre} can be viewed as a change of variables from $(\bq, \bl)$ on $\mf{M} \times \rsp^{d}$ to $(\bq, \bp)$ on $\cb{\mf{M}}$, and the projection~\eqref{eq:symplectic_projection_post}-\eqref{eq:symplectic_projection_post_constraint} as a change of variables back from $(\bq, \bp)$ to $(\bq, \bl)$.
The symplectic form $\omega_{\bl}$ on $\mf{M} \times \rsp^{d}$ thus corresponds to the pullback of the canonical symplectic form $\omega$ on $\cb{\mf{M}}$ by this variable transformation.

Let us note that the sign in in front of the projection in~\eqref{eq:symplectic_projection_post}, given by the stability function of the basic integrator, has very important implications on the nature of the algorithm.
If it is the same as in~\eqref{eq:symplectic_projection_pre}, the character of the method is very similar to the symmetric projection method described before.
If the sign is the opposite of the one in~\eqref{eq:symplectic_projection_pre}, like for Gauss--Legendre Runge--Kutta methods with an odd number of stages, the perturbation reverses the projection of the previous step, so that we effectively apply the post-projection method of~\citet{Chan:2002}. That is, the projected integrator $\Phi_{h}$ is conjugate to the unprojected integrator $\Psi_{h}$ by
\begin{align}
\Phi_{h} = \mbb{P}^{-1} \circ \Psi_{h} \circ \mbb{P} ,
\end{align}
so that the following diagram commutes
\begin{center}
\begin{tikzpicture}
\matrix (m) [matrix of math nodes,row sep=3em,column sep=4em,minimum width=2em] {
\oq_{n} & \oq_{n+1} \\
(\cq_{n}, \cl_{n}) & (\cq_{n+1}, \cl_{n+1}) \\
};

\path[-stealth, line width=.4mm]
(m-1-1) edge node [above] {$\Psi_{h}$}       (m-1-2)
(m-2-1) edge node [left]  {$\mbb{P}^{-1}$}   (m-1-1)
(m-1-2) edge node [right] {$\mbb{P}$}        (m-2-2)
(m-2-1) edge node [below] {$\Phi_{h}$}       (m-2-2);
\end{tikzpicture}
\end{center}
and the projection is effectively only applied for the output of the solution, but the actual advancement of the solution in time happens outside of the constraint submanifold (c.f., Figure~\ref{fig:post-projection}).
In other words, applying $n$ times the algorithm $\Phi_{h}$ to a point $(\cq_{0}, 0)$ is equivalent to applying the perturbation $\mbb{P}^{-1}$, then applying $n$ times the algorithm $\Psi_{h}$ and projecting the result with $\mbb{P}$.

Potentially, this might degrade the performance of the algorithm. If the accumulated global error drives the solution too far away from the constraint submanifold, the projection step might not have a solution anymore. Interestingly, however, post-projected Gauss--Legendre Runge--Kutta methods retain their optimal order of $2s$~\cite{Chan:2002}. Moreover, for methods with an odd number of stages, the global error of the unprojected solution is $\mcal{O}(h^{s+1})$, compared to $\mcal{O}(h^{s})$ for methods with an even number of stages. In practice this seems to be at least part of the reason of the good long-time stability of these methods.

\subsubsection*{Symplecticity}

While in the continuous case, the symplectic form on $\tb{\mf{M}}$ is always degenerate, thus not symplectic but presymplectic, in the discretization this is changed. The discrete Lagrangian on $\mf{M} \times \mf{M}$ is in general not degenerate, thus the symplectic form on $\mf{M} \times \mf{M}$ is non-degenerate as well. Composing the usual position-momentum form with the projection to $\Delta \subset \cb{\mf{M}}$, thus enforcing $\phi(\cq, \cp)=0$ in the way outlined before, we effectively obtain an algorithm mapping $\mf{M} \times \rsp^{d}$ into $\mf{M} \times \rsp^{d}$ instead of the original variational integrator, which mapped $\mf{M} \times \mf{M}$ into $\mf{M} \times \mf{M}$. However, the new algorithm preserves a true symplectic form on $\mf{M} \times \rsp^{d}$, which is not the same as the presymplectic form of the continuous dynamics, and also not the same as the discrete symplectic form on $\mf{M} \times \mf{M}$.
This change of the presymplectic form to a symplectic form appears to be due to the initial ``non-conservation'' of degeneracy when discretizing the Lagrangian in conjunction with the projection.

\subsection{Midpoint Projection}
\label{sec:symplectic_midpoint_projection}

For certain variational Runge--Kutta methods, we can also modify the symmetric projection in a different way in order to obtain a symplectic projection, namely by evaluating the projection at the midpoint
\begin{align}
\oz_{n+1/2} &= (\oq_{n+1/2}, \op_{n+1/2}) , &
\oq_{n+1/2} &= \tfrac{1}{2} \big( \oq_{n} + \oq_{n+1} \big) , &
\op_{n+1/2} &= \tfrac{1}{2} \big( \op_{n} + \op_{n+1} \big) , 
\end{align}
so that the projection algorithm becomes
\begin{subequations}\label{eq:symplectic_midpoint_projection}
\begin{align}
\label{eq:symplectic_midpoint_projection_pre_constraint}
0 &= \phi (\cz_{n}) , \\
\label{eq:symplectic_midpoint_projection_pre}
\oz_{n\hphantom{+1}} &= \cz_{n} + h \, \Omega^{-1} \nabla \phi^{T} (\oz_{n+1/2}) \, \cl_{n+1/2} , \\
\oz_{n+1} &= \Psi_{h} (\oz_{n}) , \\
\label{eq:symplectic_midpoint_projection_post}
\cz_{n+1} &= \oz_{n+1} + h \, \Omega^{-1} \nabla \phi^{T} (\oz_{n+1/2}) \cl_{n+1/2} , \\
\label{eq:symplectic_projection_midpoint_post_constraint}
0 &= \phi (\cz_{n+1}) .
\end{align}
\end{subequations}
This method can be shown to be symplectic with respect to the original noncanonical symplectic form on $\mf{M}$ if the integrator $\Psi_{h}$ is a symmetric, symplectic Runge--Kutta method with an odd number of stages $s$, for which the central stage with index $(s+1)/2$ corresponds to $\oz_{n+1/2}$.
This is obviously the case for the implicit midpoint rule, that is the Gauss--Legendre Runge--Kutta method with $s=1$, but unfortunately not for higher-order Gauss--Legendre or for Gauss-Lobatto methods.
However, following~\citet{OevelSofroniou:1997} and~\citet{Zhao:2014}, higher-order methods similar to Gauss--Legendre methods but satisfying the requested property can be obtained. See for example the method with three stages given in Table~\ref{tab:srk3}.

\begin{table}[bth]
\begin{center}
	\begin{tabular}{c|ccc}
		$\tfrac{1}{2} - \tfrac{\sqrt{15}}{10}$ & $\tfrac{5}{36}$ & $\tfrac{2}{9}$ & $\tfrac{25}{180} - \tfrac{\sqrt{15}}{10}$ \\
		$\tfrac{1}{2}$ & $\tfrac{5}{36}$ & $\tfrac{2}{9}$ & $\tfrac{5}{36}$ \\
		$\tfrac{1}{2} + \tfrac{\sqrt{15}}{10}$ & $\tfrac{25}{180} + \tfrac{\sqrt{15}}{10}$ & $\tfrac{2}{9}$ & $\tfrac{5}{36}$ \\
		\hline
		& $\tfrac{5}{18}$ & $\tfrac{4}{9}$ & $\tfrac{5}{18}$ \\
	\end{tabular}
	\caption{Butcher Tableau of a symplectic Runge--Kutta method with three stages~\cite{Zhao:2014}, satisfying $Z_{2} = \oz_{n+1/2} = \tfrac{1}{2} ( \oz_{n} + \oz_{n+1} )$.}
	\label{tab:srk3}
\end{center}
\end{table}

\subsubsection*{Symplecticity}

In order to show symplecticity, we follow a similar path as before for the standard projection method.
We start by computing the exterior derivative of the perturbation and projection steps,
\begin{subequations}\label{eq:symplectic_midpoint_projection_differentials}
\begin{align}
\label{eq:symplectic_midpoin_projection_differential_p0}
\ext \cp_{n}^{i} &= \ext \op_{n}^{i} + \ext \big( h \, \phi_{\bq^{i}}^{k} (\oq_{n+1/2}, \op_{n+1/2}) \, \cl_{n+1/2}^{k} \big) , \\
\label{eq:symplectic_midpoin_projection_differential_q0}
\ext \cq_{n}^{i} &= \ext \oq_{n}^{i} - \ext \big( h \, \phi_{\bp^{i}}^{k} (\oq_{n+1/2}, \op_{n+1/2}) \, \cl_{n+1/2}^{k} \big) ,
\end{align}
and
\begin{align}
\label{eq:symplectic_midpoin_projection_differential_p1}
\ext \cp_{n+1}^{i} &= \ext \op_{n+1}^{i} - \ext \big( h \, \phi_{\bq^{i}}^{k} (\oq_{n+1/2}, \op_{n+1/2}) \, \cl_{n+1/2}^{k} \big) , \\
\label{eq:symplectic_midpoin_projection_differential_q1}
\ext \cq_{n+1}^{i} &= \ext \oq_{n+1}^{i} + \ext \big( h \, \phi_{\bp^{i}}^{k} (\oq_{n+1/2}, \op_{n+1/2}) \, \cl_{n+1/2}^{k} \big) .
\end{align}	
\end{subequations}
Then we compute the wedge products $\ext \cp_{n}^{i} \wedge \ext \cq_{n}^{i}$,
\begin{align}
\ext \cp_{n}^{i} \wedge \ext \cq_{n}^{i}
 = \ext \op_{n}^{i} \wedge \ext \oq_{n}^{i}
\nonumber
&+ \ext \big( h \, \phi_{\bq^{i}}^{k} (\oq_{n+1/2}, \op_{n+1/2}) \, \cl_{n+1/2}^{k} \big) \wedge \ext \oq_{n}^{i} \\
\nonumber
&- \ext \op_{n}^{i} \wedge \ext \big( h \, \phi_{\bp^{i}}^{k} (\oq_{n+1/2}, \op_{n+1/2}) \, \cl_{n+1/2}^{k} \big) \\
&- \ext \big( h \, \phi_{\bq^{i}}^{k} (\oq_{n+1/2}, \op_{n+1/2}) \, \cl_{n+1/2}^{k} \big) \wedge \ext \big( h \, \phi_{\bp^{i}}^{k} (\oq_{n+1/2}, \op_{n+1/2}) \, \cl_{n+1/2}^{k} \big) ,
\end{align}
and $\ext \cp_{n+1}^{i} \wedge \ext \cq_{n+1}^{i}$,
\begin{align}
\ext \cp_{n+1}^{i} \wedge \ext \cq_{n+1}^{i}
 = \ext \op_{n+1}^{i} \wedge \ext \oq_{n+1}^{i}
\nonumber
&- \ext \big( h \, \phi_{\bq^{i}}^{k} (\oq_{n+1/2}, \op_{n+1/2}) \, \cl_{n+1/2}^{k} \big) \wedge \ext \oq_{n+1}^{i} \\
\nonumber
&+ \ext \op_{n+1}^{i} \wedge \ext \big( h \, \phi_{\bp^{i}}^{k} (\oq_{n+1/2}, \op_{n+1/2}) \, \cl_{n+1/2}^{k} \big) \\
&- \ext \big( h \, \phi_{\bq^{i}}^{k} (\oq_{n+1/2}, \op_{n+1/2}) \, \cl_{n+1/2}^{k} \big) \wedge \ext \big( h \, \phi_{\bp^{i}}^{k} (\oq_{n+1/2}, \op_{n+1/2}) \, \cl_{n+1/2}^{k} \big) .
\end{align}
Now assume that the integrator $\Psi_{h}$ is symplectic and thus satisfies $\ext \op_{n+1}^{i} \wedge \ext \oq_{n+1}^{i} = \ext \op_{n}^{i} \wedge \ext \oq_{n}^{i}$, which allows us to insert the second equation into the first to obtain
\begin{align}\label{eq:symplectic_midpoint_projection_symplecticity_condition_1}
\ext \cp_{n+1}^{i} \wedge \ext \cq_{n+1}^{i}
 = \ext \cp_{n}^{i} \wedge \ext \cq_{n}^{i}
\nonumber
&- \ext \big( h \, \phi_{\bq^{i}}^{k} (\oq_{n+1/2}, \op_{n+1/2}) \, \cl_{n+1/2}^{k} \big) \wedge \ext ( \oq_{n}^{i} + \oq_{n+1}^{i} ) \\
&- \ext \big( h \, \phi_{\bp^{i}}^{k} (\oq_{n+1/2}, \op_{n+1/2}) \, \cl_{n+1/2}^{k} \big) \wedge \ext ( \op_{n}^{i} + \op_{n+1}^{i} ) .
\end{align}
Noting that $\oq_{n}^{i} + \oq_{n+1}^{i} = 2 \oq_{n+1/2}$ and $\op_{n}^{i} + \op_{n+1}^{i} = 2 \op_{n+1/2}$, we can rewrite the previous expression as
\begin{align}\label{eq:symplectic_midpoint_projection_symplecticity_condition_2}
\ext \cp_{n+1}^{i} \wedge \ext \cq_{n+1}^{i}
 = \ext \cp_{n}^{i} \wedge \ext \cq_{n}^{i}
\nonumber
&- \ext \big( h \, \cl_{n+1/2}^{k} \, \phi_{\bq^{i}}^{k} (\oq_{n+1/2}, \op_{n+1/2}) \, \ext \oq_{n+1/2}^{i} \big) \\
&- \ext \big( h \, \cl_{n+1/2}^{k} \, \phi_{\bp^{i}}^{k} (\oq_{n+1/2}, \op_{n+1/2}) \, \ext \op_{n+1/2}^{i} \big) .
\end{align}
The last two terms can be combined, so that the symplecticity condition reads
\begin{align}\label{eq:symplectic_midpoint_projection_symplecticity_condition_3}
\ext \cp_{n+1}^{i} \wedge \ext \cq_{n+1}^{i} = \ext \cp_{n}^{i} \wedge \ext \cq_{n}^{i}
- \ext \big( h \, \cl_{n+1/2}^{k} \, \ext \phi^{k} (\oq_{n+1/2}, \op_{n+1/2}) \big) .
\end{align}
The additional terms vanish under the assumption that $\oz_{n+1/2} = (\oq_{n+1/2}, \op_{n+1/2})$ is equivalent to one of the internal stages of the variational Runge--Kutta method.
We pointed out before that for the internal stages the Dirac constraint $\phi(\cq, \cp)=0$ is automatically satisfied by the first equation in~\eqref{eq:vprk}.
Therefore, if $\oz_{n+1/2}$ corresponds to one of the internal stages, we have that $\phi(\oq_{n+1/2}, \op_{n+1/2}) = 0$ and thus also $\ext \phi (\oq_{n+1/2}, \op_{n+1/2}) = 0$ so that
\begin{align}\label{eq:symplectic_midpoint_projection_symplecticity_condition_4}
\ext \cp_{n+1}^{i} \wedge \ext \cq_{n+1}^{i} = \ext \cp_{n}^{i} \wedge \ext \cq_{n}^{i} .
\end{align}
It is worth pointing out that this holds for arbitrary constraints $\phi(\bq, \bp) = 0$ and that we did not use the particular structure of~\eqref{eq:dirac_constraint} like separability or $\phi_{\bp} = \unity$.
Therefore, the midpoint projection method is applicable to arbitrary Hamiltonian systems with Dirac constraints, not just the degenerate Lagrangian systems discussed in this paper.

\section{Numerical Experiments}
\label{sec:numerical_experiments}

The projection methods described in the previous section have all been implemented in the \texttt{GeometricIntegrators.jl} package, which is a library of geometric integrators for ordinary differential equations and differential algebraic equations in the Julia programming language~\cite{Bezanson:2012, Bezanson:2017} freely available on GitHub~\cite{GeometricIntegrators}.
We use Newton's method with quadratic line search for solving the nonlinear systems and LU decomposition for solving the linear systems. The Jacobian is computed via automatic differentiation via the \texttt{ForwardDiff.jl} package~\cite{RevelsLubinPapamarkou:2016} and updated in every time step but only every five nonlinear iterations.
If possible, the numerical integration step and the projection step are solved separately (that is for the standard and symplectic projection, but not for the symmetric and midpoint projection).
The updates of the solution are computed using compensated summation (Kahan's algorithm) in order to reduce the propagation of round-off errors.

The examples we will consider are a two-dimensional Lotka--Volterra model, planar point vortices with varying circulation and guiding centre dynamics.
The first two examples are implemented in the \texttt{GeometricProblems.jl} package. The latter is implemented in the \texttt{ChargedParticleDynamics.jl} package. Both packages are also available on GitHub~\cite{GeometricProblems,ChargedParticleDynamics}.
Except for the first example, all systems possess Noether symmetries and some related conservation law, whose preservation will be monitored in the simulations.

We perform simulation with Gauss--Legendre Runge--Kutta methods with one to six stages as well as Gauss--Lobatto--IIIA, IIIB, IIIC, IIID and IIIE methods~\cite{Jay:2015hc} with two, three and four stages. 
Here, the referenced method always provides the coefficients $a$ and the coefficients $\bar{a}$ are chosen, such that the symplecticity condition~\eqref{eq:vprk_coefficients} is satisfied. That is Gauss--Lobatto--IIIA denotes the IIIA--IIIB pair, Gauss--Lobatto--IIIB denotes the IIIB--IIIA pair, and Gauss--Lobatto--IIIC denotes the IIIC--IIIC* pair. For the Gauss--Legendre as well as Gauss--Lobatto--IIID and IIIE methods, we have $\bar{a} = a$.
The Gauss--Lobatto--IIIC* method is sometimes also referred to as Gauss--Lobatto--III. Similar inconsistent naming is found for the IIID and IIIE methods. Here, we denote by the Gauss--Lobatto--IIID method the special case of the Gauss--Lobatto--IIIS method with $\sigma=1.0$ and by the IIIE method the special case of the IIIS method with $\sigma=0.5$. The Gauss--Lobatto--IIIS methods are interpolations of the IIIA, IIIB, IIIC and IIIC* methods with coefficients given by
\begin{align*}
a_{ij}^{S} (\sigma) = (1-\sigma) \big( a_{ij}^{A} + a_{ij}^{B} \big) + (\sigma - \tfrac{1}{2})  \big( a_{ij}^{C} + a_{ij}^{C*} \big) ,
\end{align*}
so that
\begin{align*}
a_{ij}^{D} = a_{ij}^{S} (1) = \tfrac{1}{2} (a_{ij}^{A} + a_{ij}^{B} ) ,
\end{align*}
and
\begin{align*}
a_{ij}^{E} = a_{ij}^{S} (\tfrac{1}{2}) = \tfrac{1}{2} (a_{ij}^{C} + a_{ij}^{C*} ) .
\end{align*}
We compare the results of the variational Runge--Kutta methods with simulations of Radau--IIA methods, which have the advantage that they automatically preserve the Dirac constraint but also have the disadvantage of dissipating energy.

For all methods, we perform simulations both without projection and with standard, symmetric, symplectic and midpoint projection. Due to the limited space we will only show some selected examples. The collection of all simulation results can be found in the documentation of the \texttt{GeometricExamples.jl} package~\cite{GeometricExamples}.

For most examples the simulations with the Gauss--Lobatto--IIIA, IIIB and IIIC methods break down after very few time steps. Even when reducing the time step by an order of magnitude, the IIIA, IIIB and IIIC methods perform rather poorly in almost all of the experiments. For the unprojected integrator this was already shown as a motivating example in Figure~\ref{fig:lotka_volterra_2d_tr} for the Lotka--Volterra model.
The origin of this behaviour is most likely related to the fact that for the IIIA, IIIB and IIIC methods, different Runge--Kutta coefficient $a_{ij}$ and $\obar{a}_{ij}$ are used for the integration of the trajectory $\cq$ and the conjugate momenta $p$. Even though the nodes of the stages $c_{i}$ are the same for both $\cq$ and $p$, the definition of the values at the nodes $\cq_{n,i}$ and $\rkp_{n,i}$ in terms of the corresponding vector fields $\rkv_{n,i}$ and $\rkf_{n,i}$ is different. While this is usually fine for regular problems, especially with separable Hamiltonians, it does not seem appropriate for degenerate problems where there is a functional relationship between the momenta and the position along the trajectory given at the internal stages by $\rkp_{n,i} = \vartheta (\rkq_{n,i})$ for $1 \leq i \leq s$. This particular property of degenerate systems suggests that the same coefficient matrices should be used for the definition of the internal stages of both $\cq$ and $p$.

For the two- and four-stage Gauss--Lobatto--IIIA, IIIB and IIIC methods, the symplectic projection amounts to a post-projection. Therefore, if the simulation without projection breaks down, so does the simulation with symplectic projection.
The midpoint projection is only symplectic for the Gauss--Legendre method with one stage and the SRK3 method whose tableau was given in Table~\ref{tab:srk3}. Nevertheless, we run experiments with this projection and all integrators to study the long time behaviour.

\subsection{Lotka--Volterra Model}
\label{sec:app_lotka_volterra_2d}

Lotka--Volterra models~\cite{Lotka:1925, Volterra:1927} are used in mathematical biology for modelling population dynamics of animal species, as well as many other fields where predator-prey and similar models appear.
The dynamics of the growth of two interacting species can be modelled by the following Lagrangian system~\cite{FernandezNunez:1998},
\begin{align}
L (\bq, \bv) = \left( \dfrac{\log \bq_{2}}{\bq_{1}} + \bq_{2} \right) \bv_{1} + \bq_{1} \bv_{2} - H (\bq) ,
\end{align}
with the Hamiltonian $H$ given by
\begin{align}
H (\bq) = a_{1} \bq_{1} + a_{2} \bq_{2} - b_{1} \log \bq_{1} - b_{2} \log \bq_{2} .
\end{align}
The noncanonical symplectic form~\eqref{eq:noncanonical_symplectic_matrix} is computed as
\begin{align}
\obar{\omega} = - \dfrac{1}{\bq_{1} \bq_{2}} \, \ext \bq_{1} \wedge \ext \bq_{2} .
\end{align}
In the position-momentum form, which is the basis for the variational Runge-Kutta methods we employ in the numerical experiments, we obtain the following functions for the momenta and forces,
\begin{subequations}
\begin{align}
\vartheta_{1} (\bq) &= \dfrac{\log \bq_{2}}{\bq_{1}} + \bq_{2} , &
f_{1} (\bq, \bv) &= \bv_{2} - \dfrac{\log \bq_{2}}{\bq_{1}^{2}} \, \bv_{1} - a_{1} + \dfrac{b_{1}}{\bq_{1}} , \\
\vartheta_{2} (\bq) &= \bq_{1} , &
f_{2} (\bq, \bv) &= \left( 1 + \dfrac{1}{\bq_{1} \bq_{2}} \right) \bv_{1} - a_{2} + \dfrac{b_{2}}{\bq_{2}} .
\end{align}
\end{subequations}
In the simulations, we use a time step of $h = 0.1$ and consider initial conditions $(q_{1,0}, q_{2,0}) = (1, 1)$ with parameters $(a_{1}, a_{2}, b_{1}, b_{2}) = (1, 1, 1, 2)$, which give a periodic solution.

\begin{figure}[p]
	\begin{center}
		\subfloat[GLRK1]{\label{fig:lotka_volterra_2d_vprk_pnone_glrk1}
			\includegraphics[height=2.5cm]{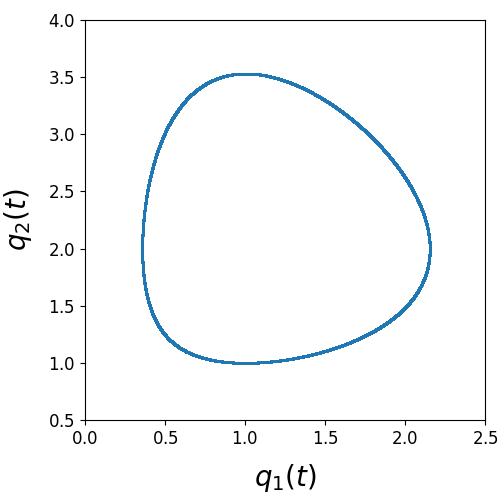}
			\includegraphics[height=2.5cm]{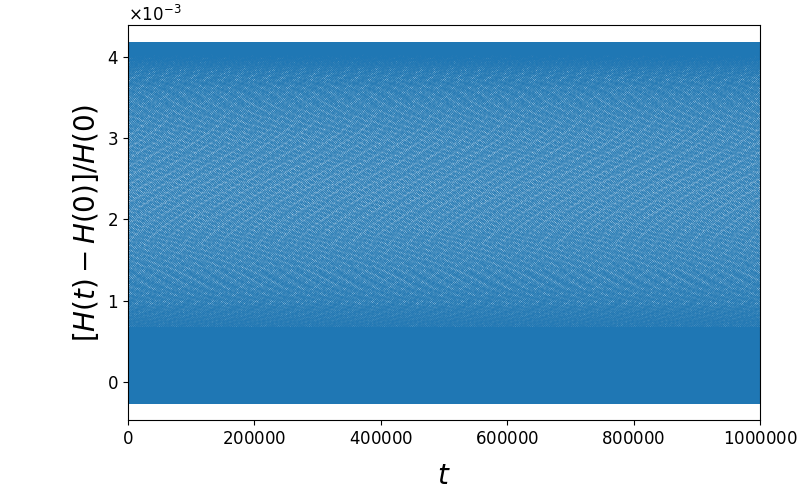}
			\includegraphics[height=2.5cm]{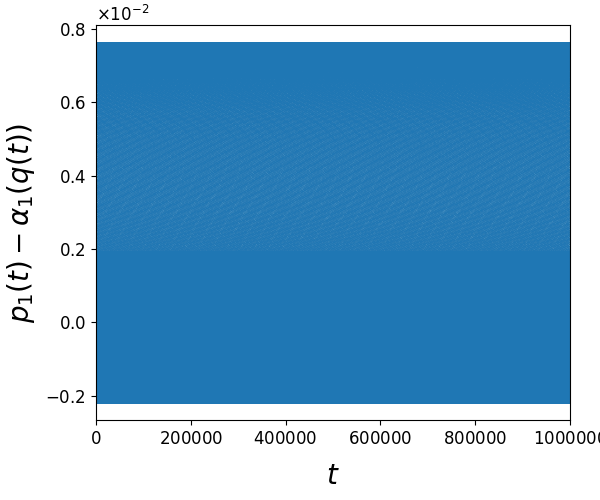}
			\includegraphics[height=2.5cm]{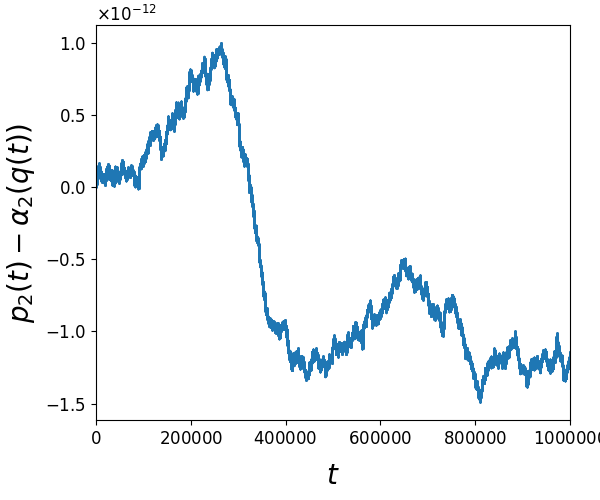}
		}

		\subfloat[GLRK2]{\label{fig:lotka_volterra_2d_vprk_pnone_glrk2}
			\includegraphics[height=2.5cm]{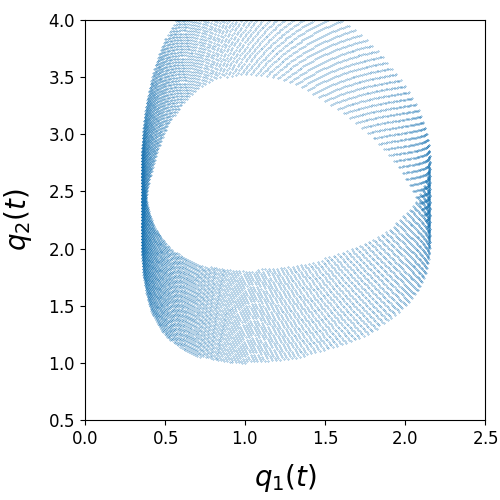}
			\includegraphics[height=2.5cm]{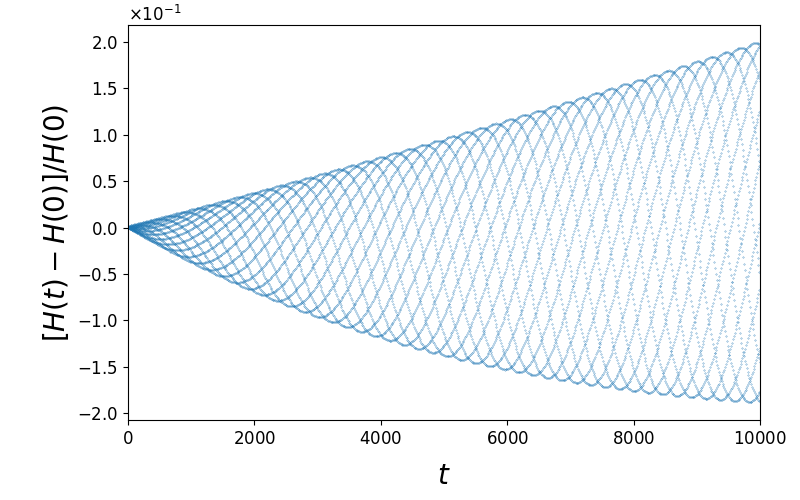}
			\includegraphics[height=2.5cm]{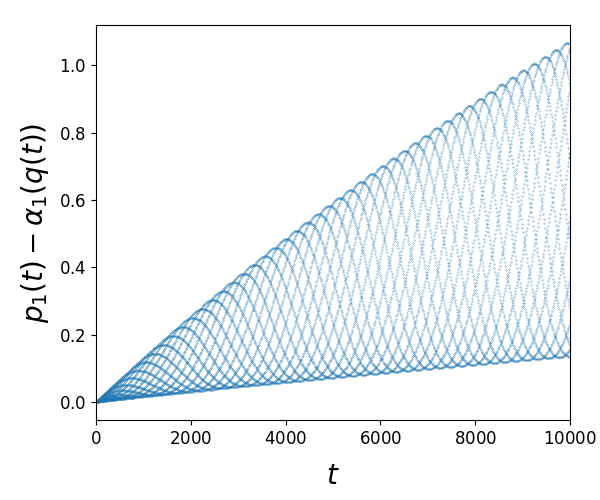}
			\includegraphics[height=2.5cm]{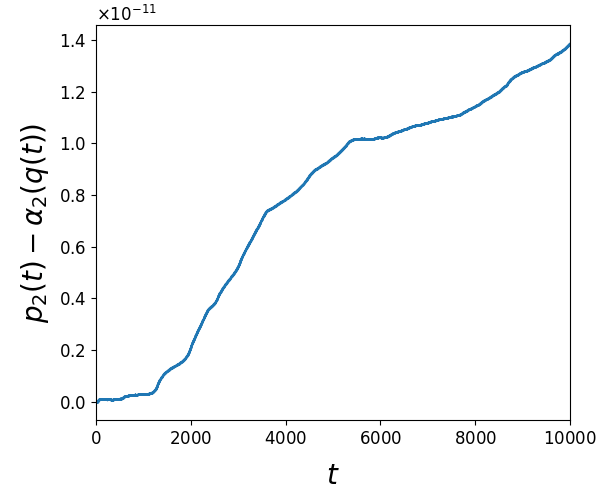}
		}
		
		\subfloat[GLRK3]{\label{fig:lotka_volterra_2d_vprk_pnone_glrk3}
			\includegraphics[height=2.5cm]{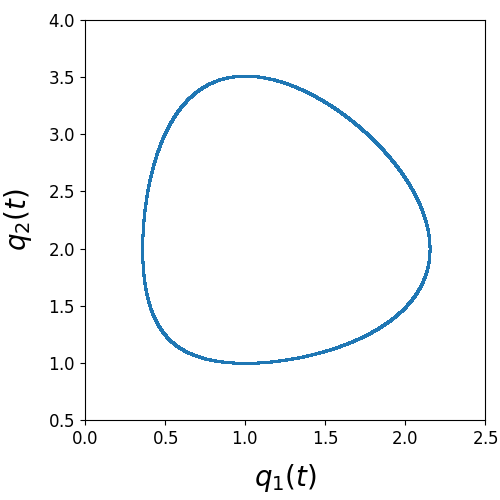}
			\includegraphics[height=2.5cm]{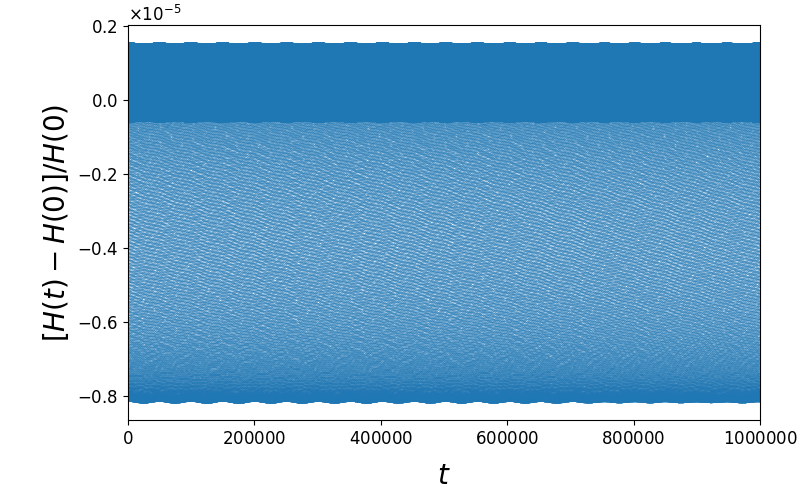}
			\includegraphics[height=2.5cm]{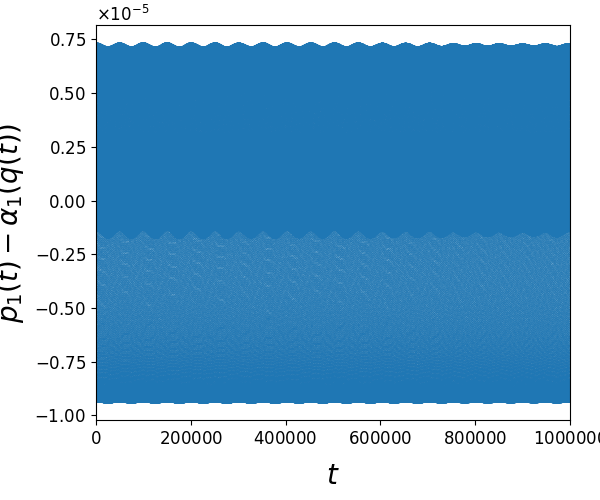}
			\includegraphics[height=2.5cm]{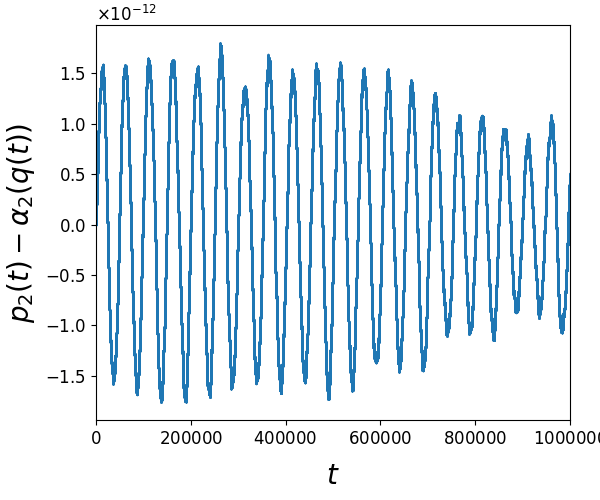}
		}
		
		\subfloat[GLRK4]{\label{fig:lotka_volterra_2d_vprk_pnone_glrk4}
			\includegraphics[height=2.5cm]{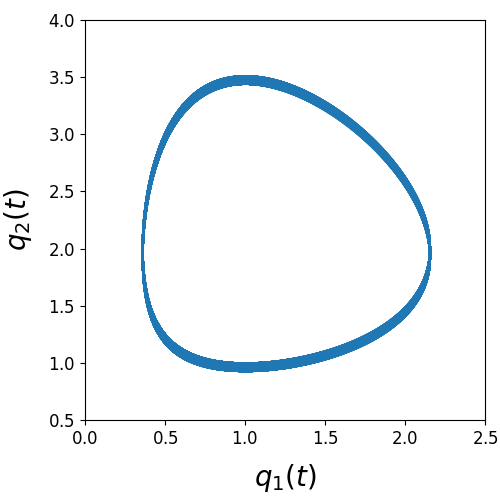}
			\includegraphics[height=2.5cm]{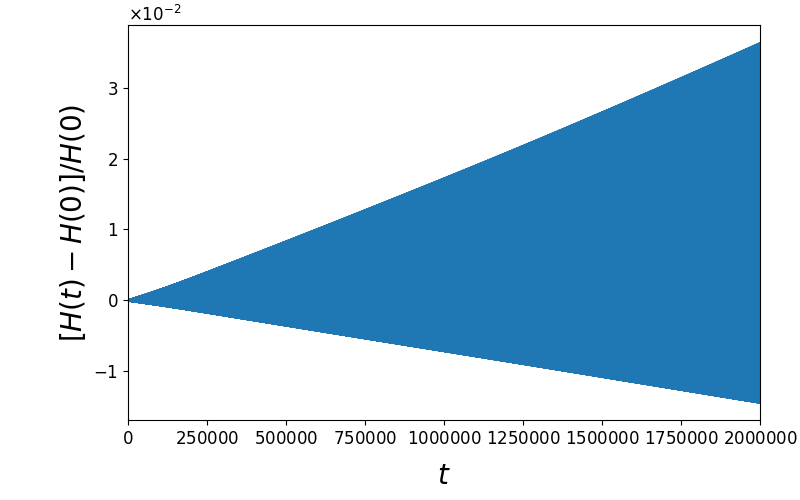}
			\includegraphics[height=2.5cm]{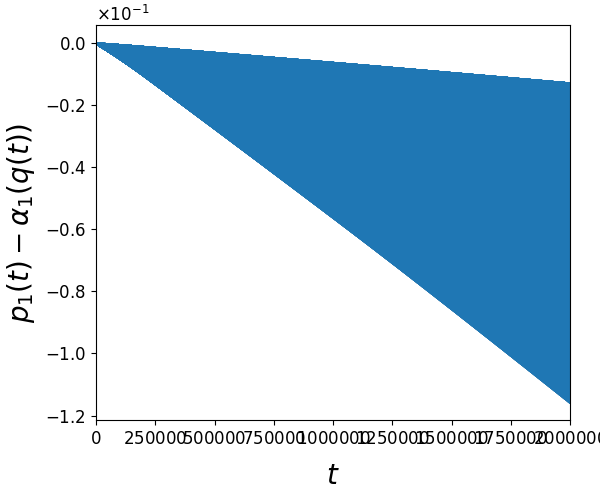}
			\includegraphics[height=2.5cm]{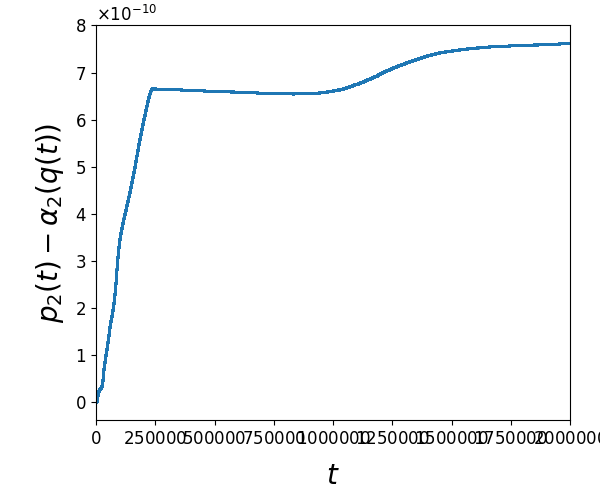}
		}
		
		\subfloat[SRK3]{\label{fig:lotka_volterra_2d_vprk_pnone_srk3}
			\includegraphics[height=2.5cm]{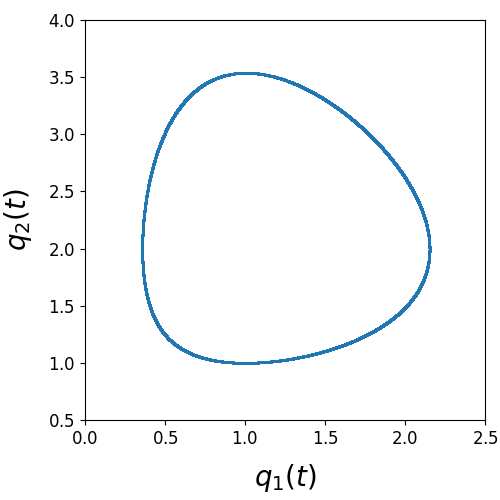}
			\includegraphics[height=2.5cm]{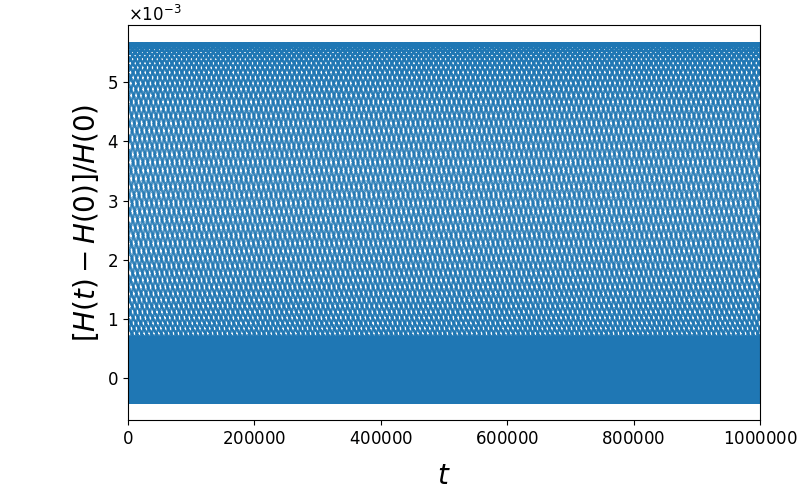}
			\includegraphics[height=2.5cm]{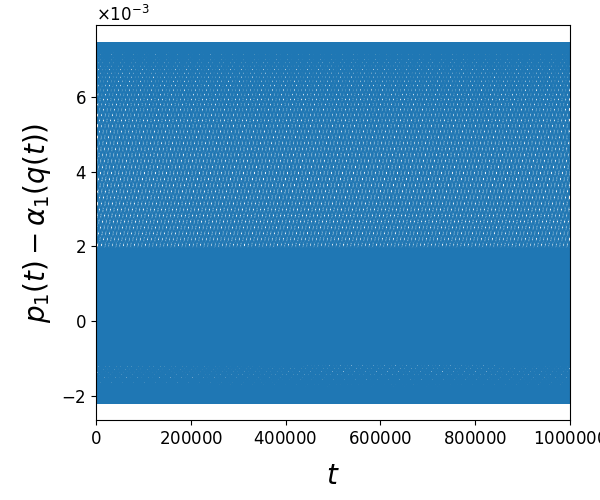}
			\includegraphics[height=2.5cm]{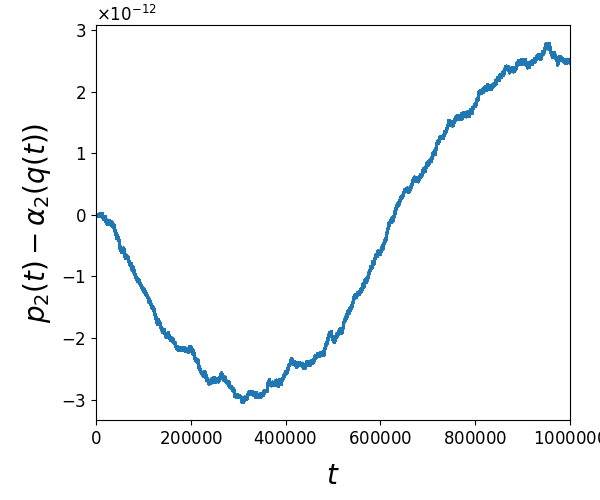}
		}
		
		\subfloat[RadIIA2]{\label{fig:lotka_volterra_2d_vprk_RadIIA2}
			\includegraphics[height=2.5cm]{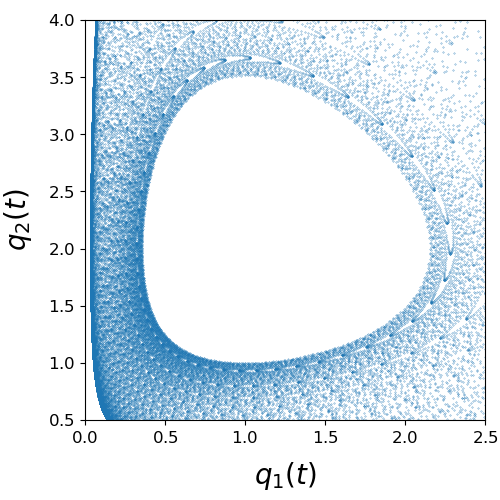}
			\includegraphics[height=2.5cm]{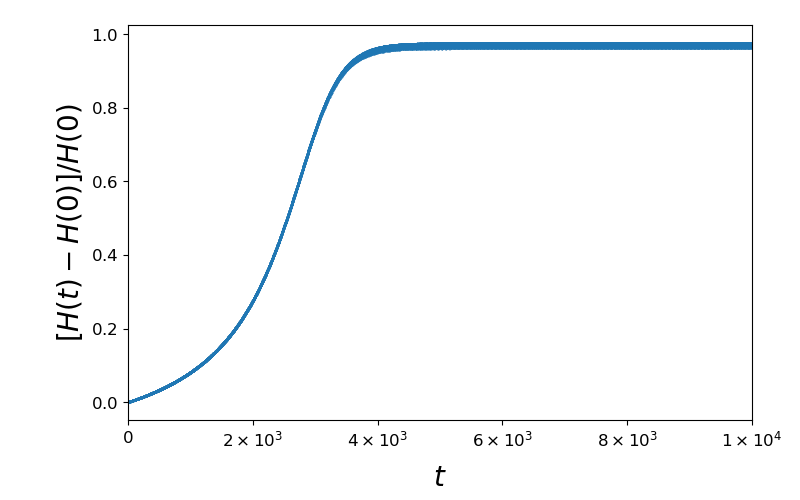}
			\includegraphics[height=2.5cm]{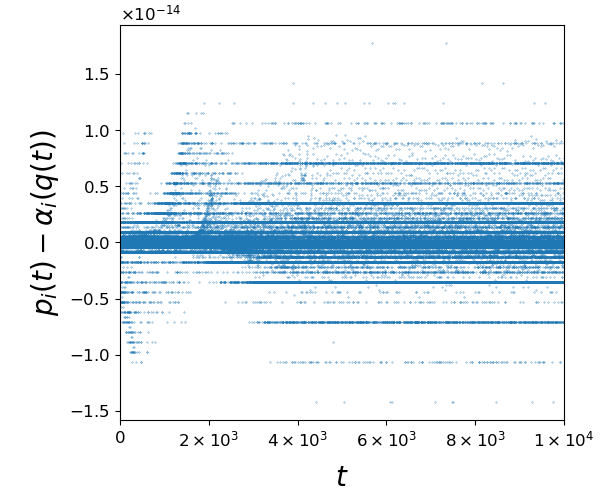}
			\includegraphics[height=2.5cm]{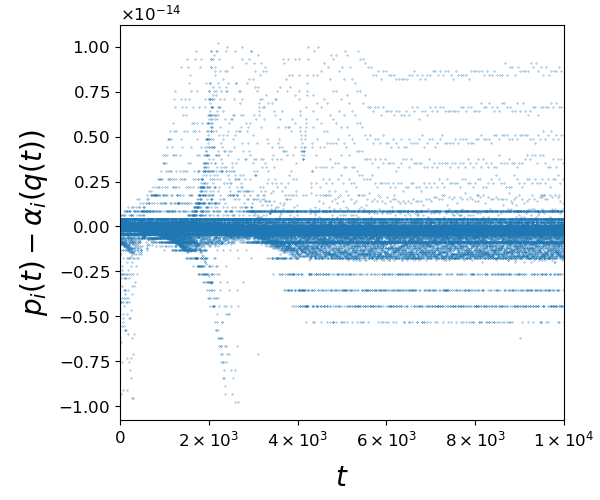}
		}
		
		\subfloat[RadIIA3]{\label{fig:lotka_volterra_2d_vprk_RadIIA3}
			\includegraphics[height=2.5cm]{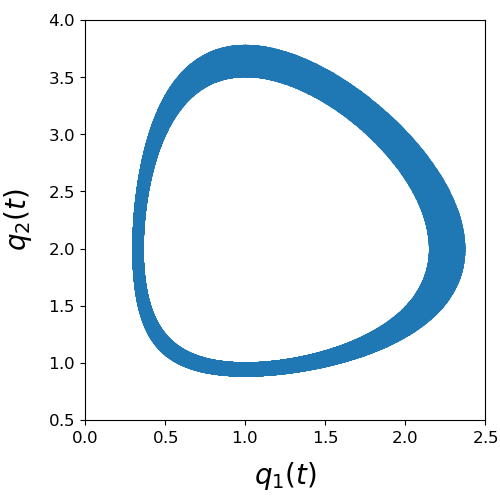}
			\includegraphics[height=2.5cm]{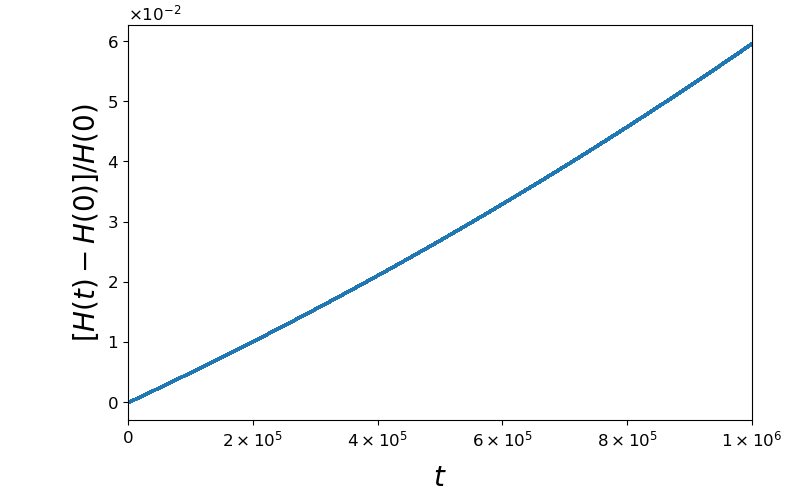}
			\includegraphics[height=2.5cm]{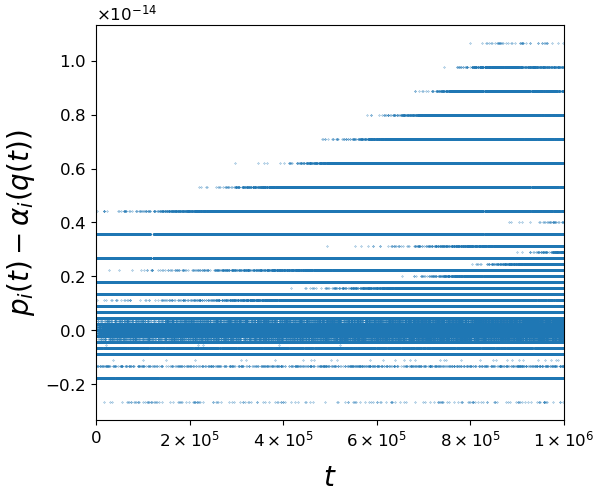}
			\includegraphics[height=2.5cm]{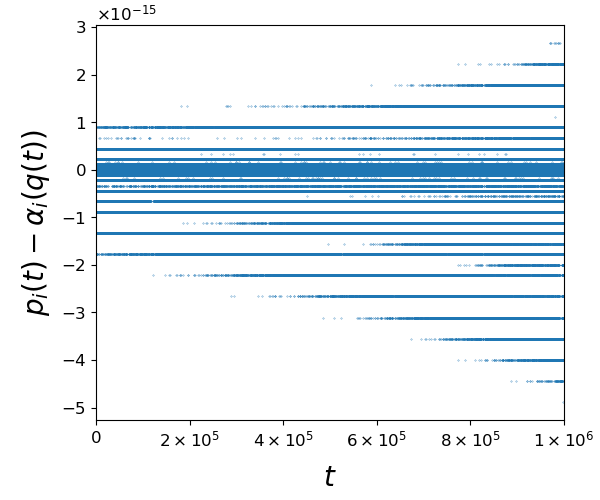}
		}
		
		\caption{Lotka--Volterra model with variational-partitioned Runge--Kutta and non-variational Radau methods.}
		\label{fig:lotka_volterra_2d_vprk_pnone1}
	\end{center}
\end{figure}

\begin{figure}[htp]
	\begin{center}
		\subfloat[GLRK3pStandard]{\label{fig:lotka_volterra_2d_vprk_pstandard_glrk3}
			\includegraphics[height=2.5cm]{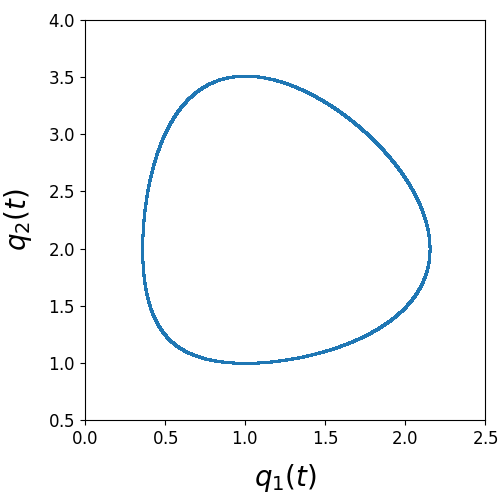}
			\includegraphics[height=2.5cm]{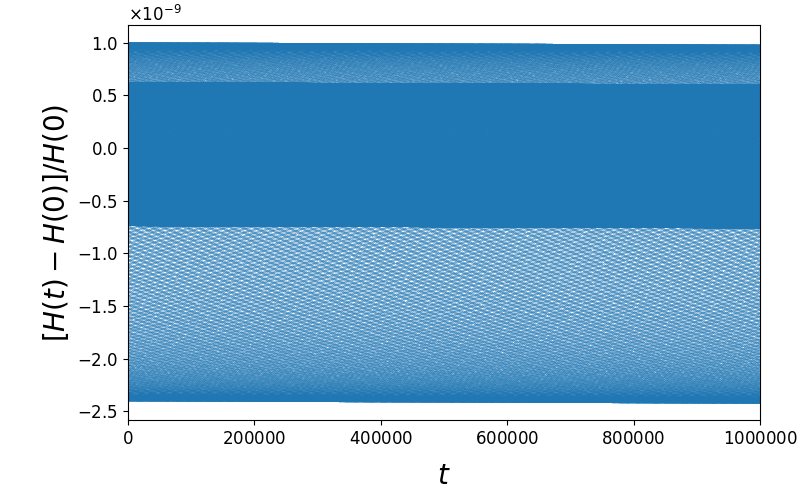}
		}
		\subfloat[GLRK3pSymmetric]{\label{fig:lotka_volterra_2d_vprk_psymmetric_glrk3}
			\includegraphics[height=2.5cm]{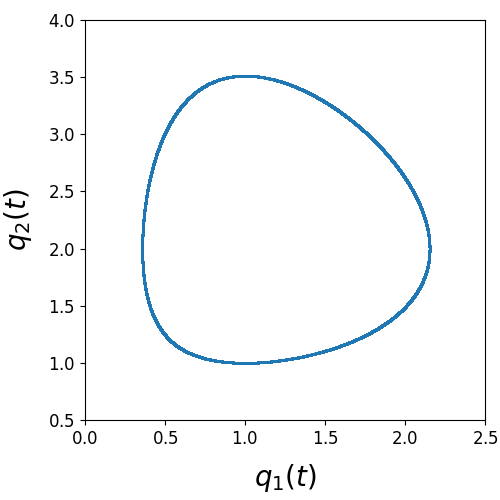}
			\includegraphics[height=2.5cm]{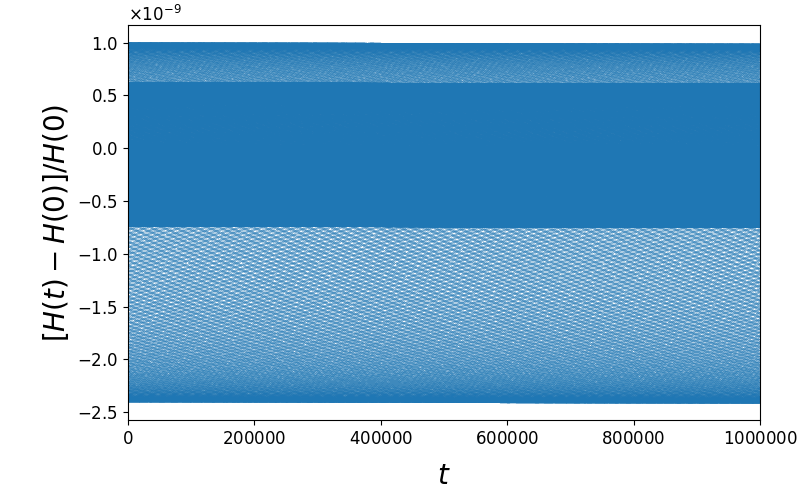}
		}
		
		\subfloat[GLRK4pStandard]{\label{fig:lotka_volterra_2d_vprk_pstandard_glrk4}
			\includegraphics[height=2.5cm]{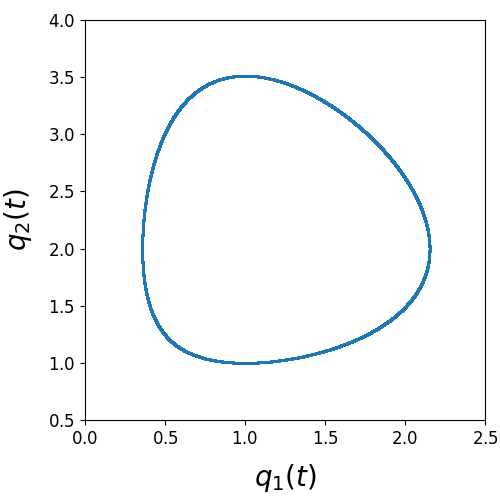}
			\includegraphics[height=2.5cm]{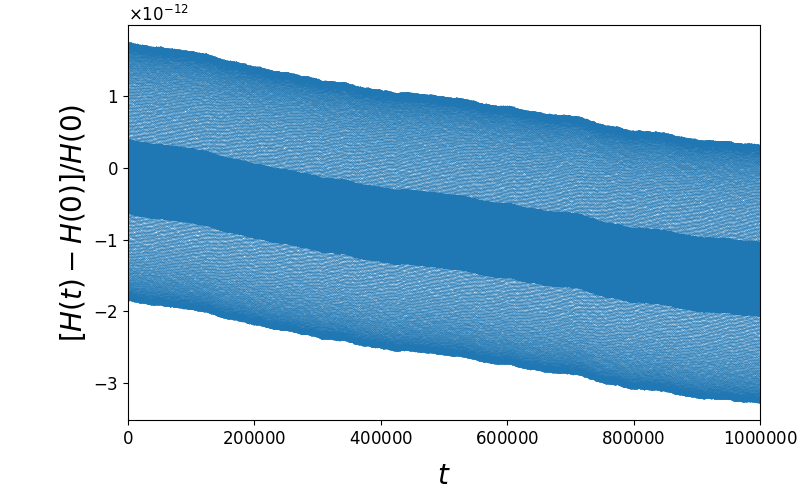}
		}
		\subfloat[GLRK4pSymmetric]{\label{fig:lotka_volterra_2d_vprk_psymmetric_glrk4}
			\includegraphics[height=2.5cm]{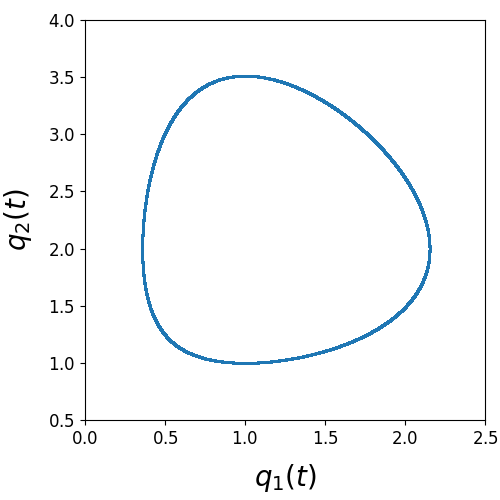}
			\includegraphics[height=2.5cm]{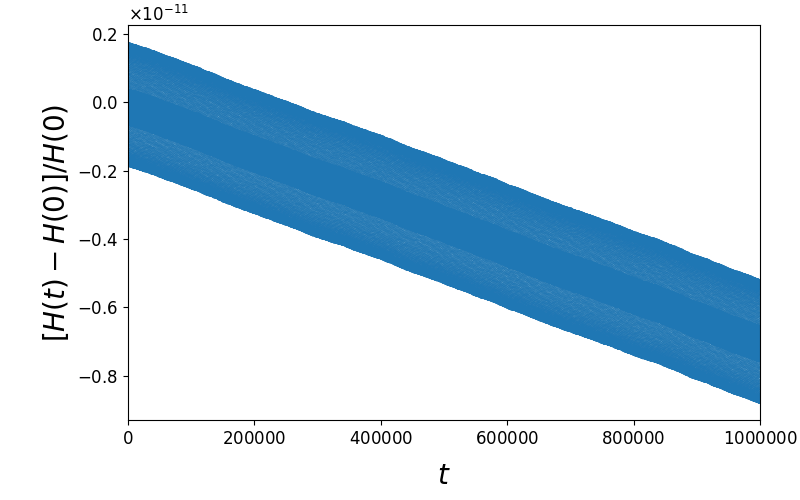}
		}

		\subfloat[SRK3pStandard]{\label{fig:lotka_volterra_2d_vprk_pstandard_srk3}
			\includegraphics[height=2.5cm]{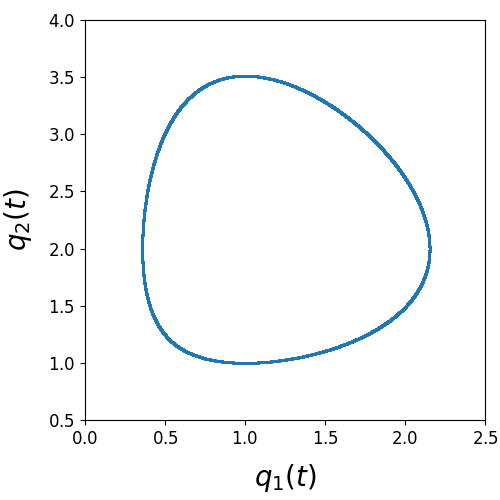}
			\includegraphics[height=2.5cm]{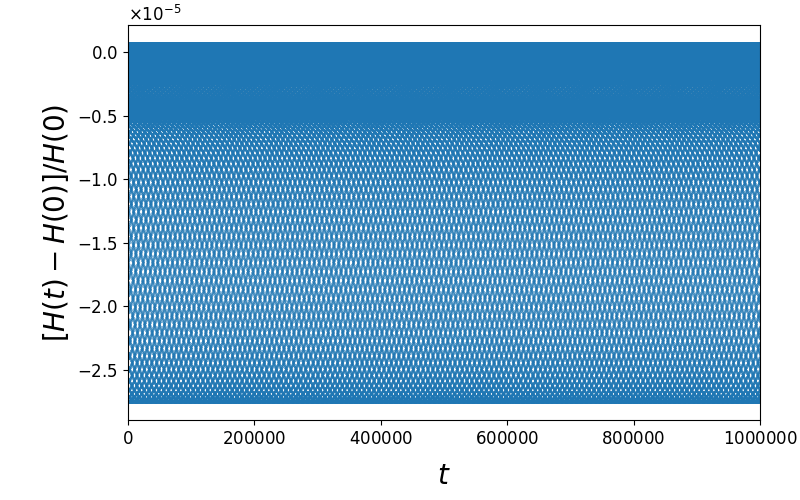}
		}
		\subfloat[SRK3pSymmetric]{\label{fig:lotka_volterra_2d_vprk_psymmetric_srk3}
			\includegraphics[height=2.5cm]{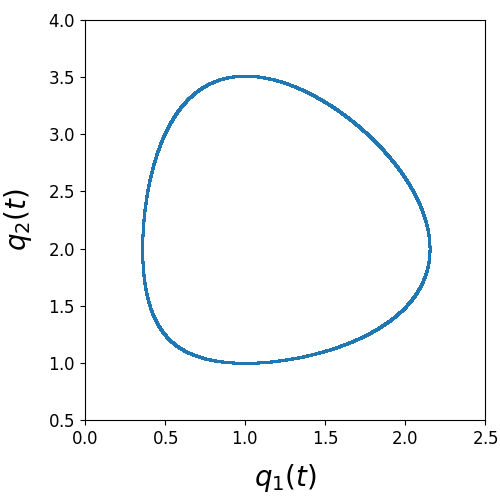}
			\includegraphics[height=2.5cm]{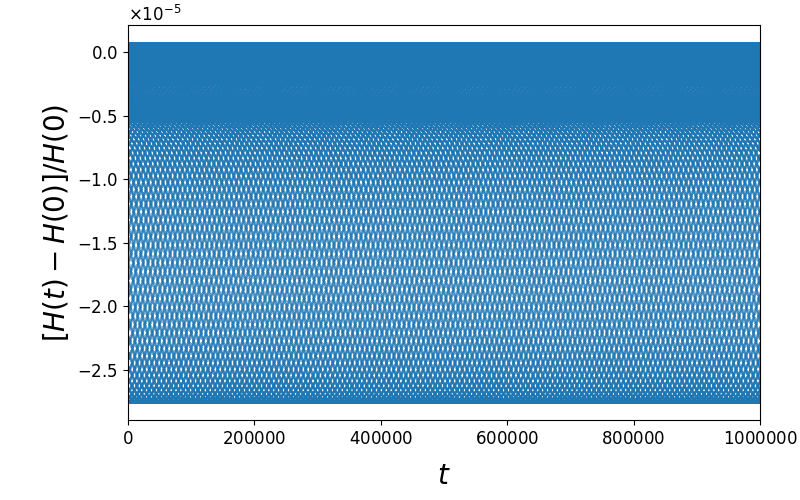}
		}
		
		\caption{Lotka--Volterra model with variational Runge--Kutta methods with standard and symmetric projection.}
		\label{fig:lotka_volterra_2d_vprk_pstandard}
	\end{center}
\end{figure}

\begin{figure}[htp]
	\begin{center}
		\subfloat[GLRK3pSymplectic]{\label{fig:lotka_volterra_2d_vprk_psymplectic_glrk3}
			\includegraphics[height=2.5cm]{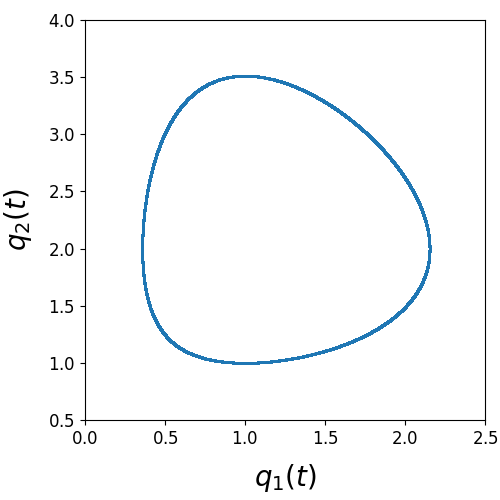}
			\includegraphics[height=2.5cm]{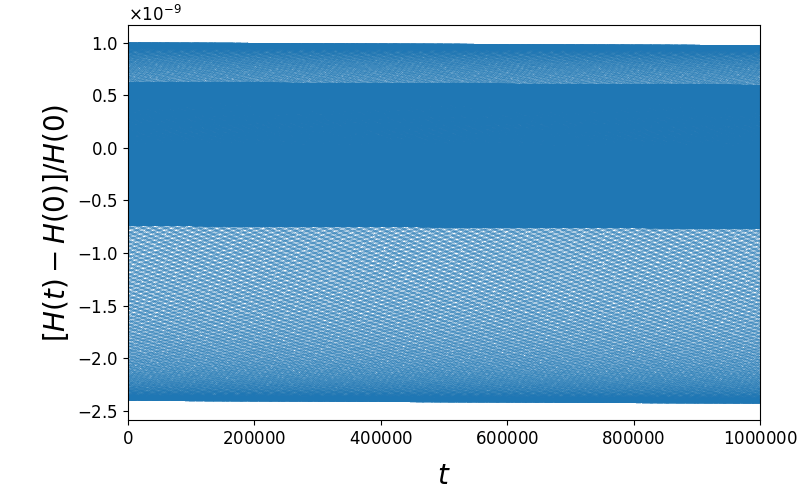}
		}
		\subfloat[GLRK3pMidpoint]{\label{fig:lotka_volterra_2d_vprk_pmidpoint_glrk3}
			\includegraphics[height=2.5cm]{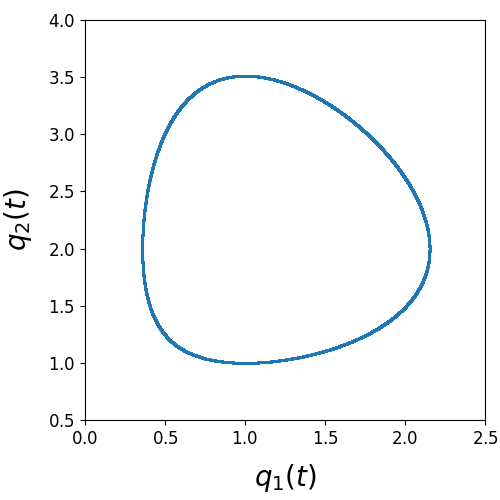}
			\includegraphics[height=2.5cm]{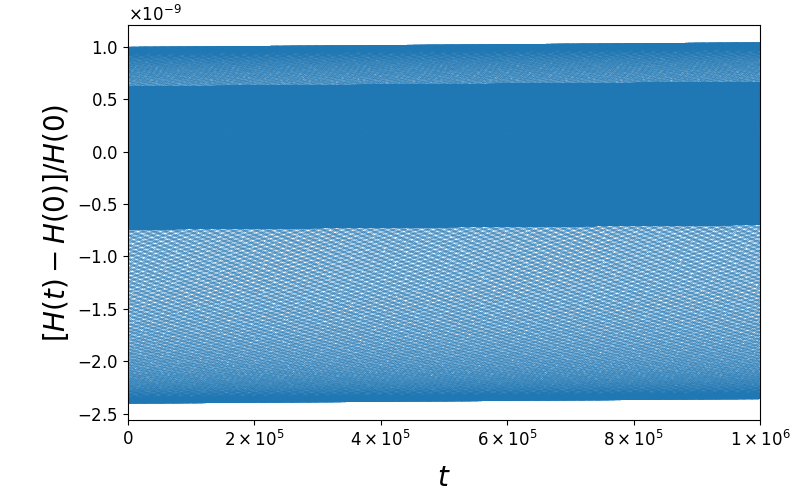}
		}
		
		\subfloat[GLRK4pSymplectic]{\label{fig:lotka_volterra_2d_vprk_psymplectic_glrk4}
			\includegraphics[height=2.5cm]{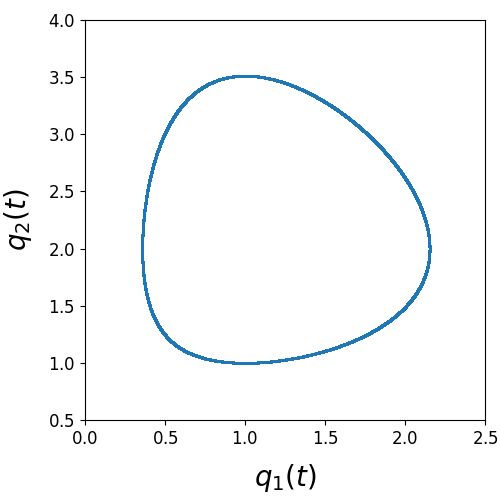}
			\includegraphics[height=2.5cm]{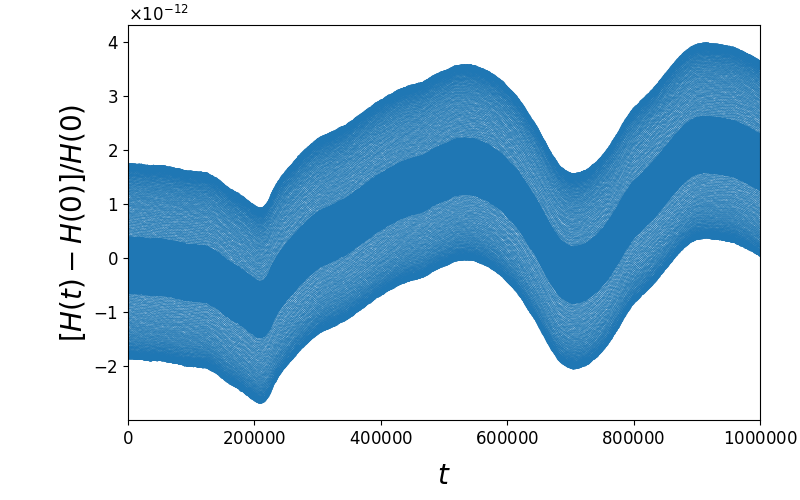}
		}
		\subfloat[GLRK4pMidpoint]{\label{fig:lotka_volterra_2d_vprk_pmidpoint_glrk4}
			\includegraphics[height=2.5cm]{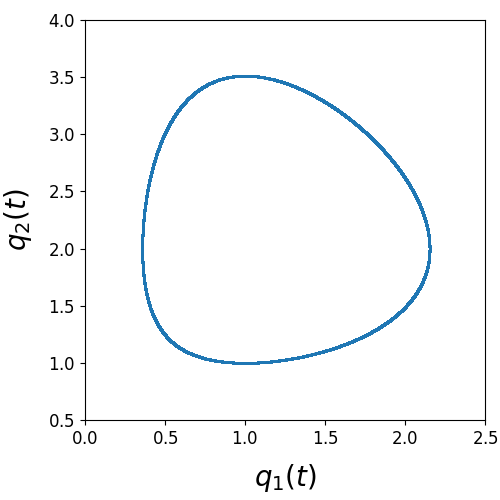}
			\includegraphics[height=2.5cm]{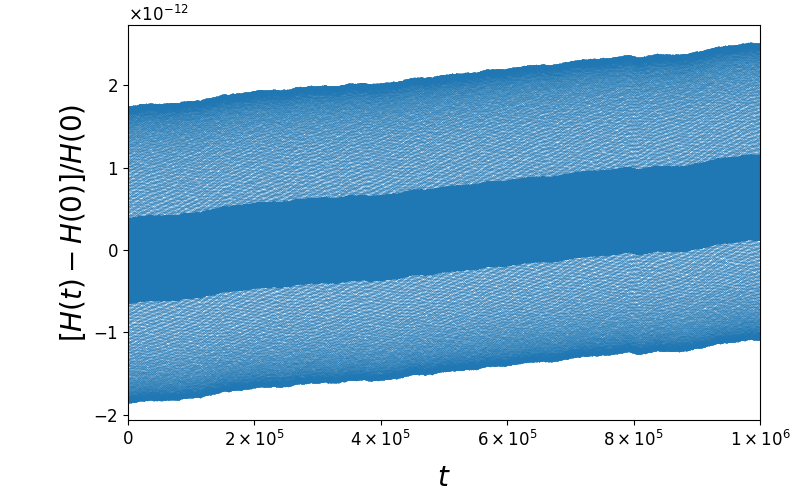}
		}
		
		\subfloat[SRK3pSymplectic]{\label{fig:lotka_volterra_2d_vprk_psymplectic_srk3}
			\includegraphics[height=2.5cm]{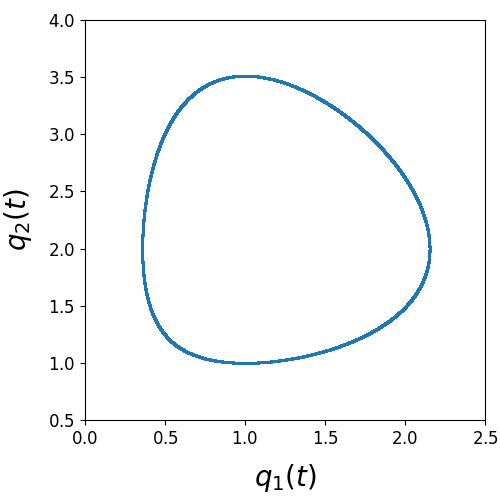}
			\includegraphics[height=2.5cm]{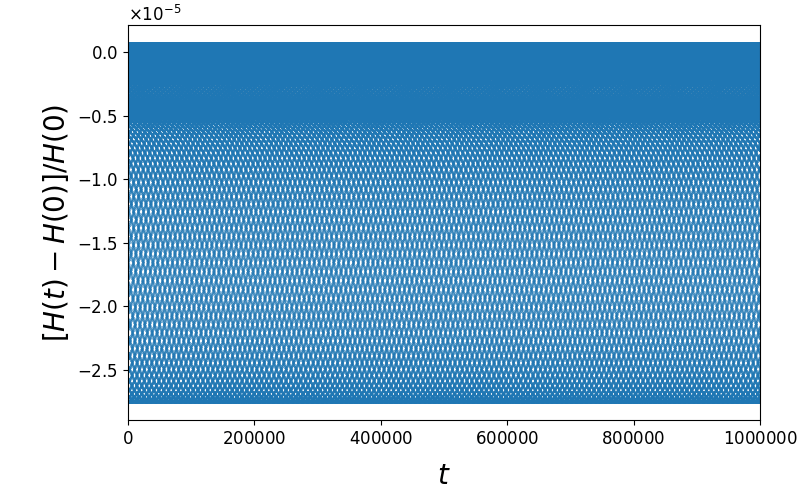}
		}
		\subfloat[SRK3pMidpoint]{\label{fig:lotka_volterra_2d_vprk_pmidpoint_srk3}
			\includegraphics[height=2.5cm]{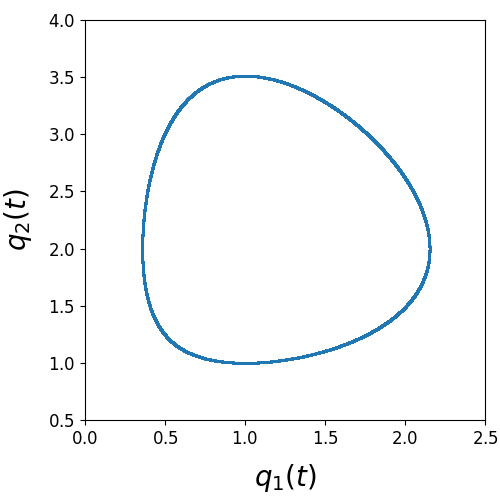}
			\includegraphics[height=2.5cm]{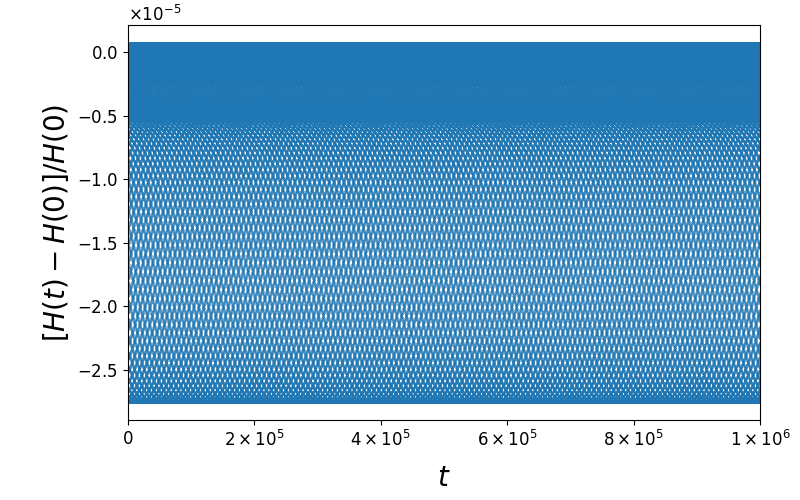}
		}
		
		\caption{Lotka--Volterra model with variational Runge--Kutta methods with symplectic and midpoint projection.}
		\label{fig:lotka_volterra_2d_vprk_psymplectic}
	\end{center}
\end{figure}

We make the following observations:
\begin{itemize}
\item The Gauss--Legendre Runge--Kutta methods with an odd number of stages  (Figure~\ref{fig:lotka_volterra_2d_vprk_pnone_glrk3}) as well as the SRK3 method (Figure~\ref{fig:lotka_volterra_2d_vprk_pnone_srk3}) are stable even without projection. Even though they do not preserve the Dirac constraint exactly, the error in the constraint oscillates about zero and the amplitude of that oscillation appears to be bounded or at least grows only slowly. A similar behaviour is observed for the Gauss--Lobatto--IIID and IIIE methods (not shown).

\item The Gauss--Legendre Runge--Kutta methods with an even number of stages (Figure~\ref{fig:lotka_volterra_2d_vprk_pnone_glrk4}) show an increasing error in the Dirac constraint and also in the energy, which eventually renders the simulation unstable (after about 250~000 time steps for the two-stage method and after about 1~000~000 time steps for the four-stage method).

\item The Gauss--Lobatto--IIIA, IIIB and IIIC methods (not shown) are unstable without projection. For the integrator with two stages, the simulation crashes after about 25 time steps. For the integrator with three stages, it crashes immediately on the first time step. Decreasing the time step to $h = 0.01$ both integrator run for a short period. The integrator with two stages crashes after about 350 time steps and the integrator with three stages after about 1.000 time steps.

\item The standard projection leads to very good results with all Gauss--Legendre methods (Figures~\ref{fig:lotka_volterra_2d_vprk_pstandard_glrk3}, \ref{fig:lotka_volterra_2d_vprk_pstandard_glrk4}), the SRK3 method (Figure~\ref{fig:lotka_volterra_2d_vprk_pstandard_srk3}), as well as the Gauss--Lobatto--IIID and IIIE methods (not shown), but not with the Gauss--Lobatto--IIIA, IIIB and IIIC methods (not shown), whose solution deteriorates quickly. We observe small drifts in the energy error, but over 10~000~000 time steps this drift is of the order of $10^{-12}$.

\item For the Gauss--Legendre methods (Figures~\ref{fig:lotka_volterra_2d_vprk_psymmetric_glrk3}, \ref{fig:lotka_volterra_2d_vprk_psymmetric_glrk4}), the SRK3 method (Figure~\ref{fig:lotka_volterra_2d_vprk_psymmetric_srk3}), and the Gauss--Lobatto--IIID and IIIE methods (not shown), the symmetric projection leads to similar results as the standard projection. In some cases the drift in the energy seems to be slightly larger than with the standard projection. This, however, is due to round-off errors (c.f., Section~\ref{sec:app_energy_drift}). The errors of the Gauss--Lobatto--IIIA, IIIB and IIIC methods (not shown) are smaller than with the standard projection, but there still is a substantial drift in the energy.

\item The symplectic projection (Figures~\ref{fig:lotka_volterra_2d_vprk_psymplectic_glrk3}, \ref{fig:lotka_volterra_2d_vprk_psymplectic_glrk4}, \ref{fig:lotka_volterra_2d_vprk_psymplectic_srk3}) leads to very good results with all Gauss--Legendre methods, the SRK3 method, as well as the Gauss--Lobatto--IIID and IIIE methods, comparable to the results obtained with the symmetric projection. For the Gauss--Lobatto--IIIA, IIIB and IIIC methods with an even number of stages, the symplectic projection corresponds to a post-projection method. For this reason and as these methods are unstable without projection, the symplectic projection is also unstable. Although for an odd number of stages, the projection does not correspond to a post-projection, simulations still tend to crash as quickly as without projection.

\item The midpoint projection (Figures~\ref{fig:lotka_volterra_2d_vprk_pmidpoint_glrk3}, \ref{fig:lotka_volterra_2d_vprk_pmidpoint_glrk4}, \ref{fig:lotka_volterra_2d_vprk_pmidpoint_srk3}) leads to good results with all Gauss--Legendre methods, the SRK3 method as well as the Gauss--Lobatto--IIID and IIIE methods, again comparable to the results obtained with the symmetric projection. However, it is only symplectic for the Gauss--Legendre method with one stage and the SRK3 method.

\item With the Radau methods, we observe exact conservation of the Dirac constraint, as expected, but dissipation of energy, which is related to large errors in the solution.
\end{itemize}

In summary, the numerical experiments for the Lotka--Volterra problem suggest that the Gauss--Legendre, the SRK and the Gauss--Lobatto--IIID and IIIE methods lead to good results with all projection methods, whereas (for the time step used) the the results with the Gauss--Lobatto--IIIA, IIIB and IIIC methods are never satisfactory, even with projection. For these methods, the time step needs to be reduced by at least a factor $10$ in order to obtain stable simulations. With such small time steps, however, the Gauss--Lobatto--IIIA, IIIB and IIIC methods are not competitive anymore and one should rather use the Gauss--Legendre or Gauss--Lobatto--IIID and IIIE methods.
For the Gauss--Legendre methods with an odd number of stages, the simulation appears to be stable even without projection, at least for very long times (10 million time steps), although the order of the integrators is decreased in this case (c.f., Section~\ref{sec:app_convergence}).
It was already reported by~\citet{Chan:2004} that for index two differential-algebraic equations Gauss--Legendre Runge--Kutta methods with an odd number of stages behave much better than those with an even number of stages, which is related to the stability function $R(\infty)$ being $+1$ for the former and $-1$ for the latter.

\subsection{Planar Point Vortices with Varying Circulation}
\label{sec:app_planar_point_vortices}

\begin{figure}[bt]
	\begin{center}
		\includegraphics[width=.5\textwidth]{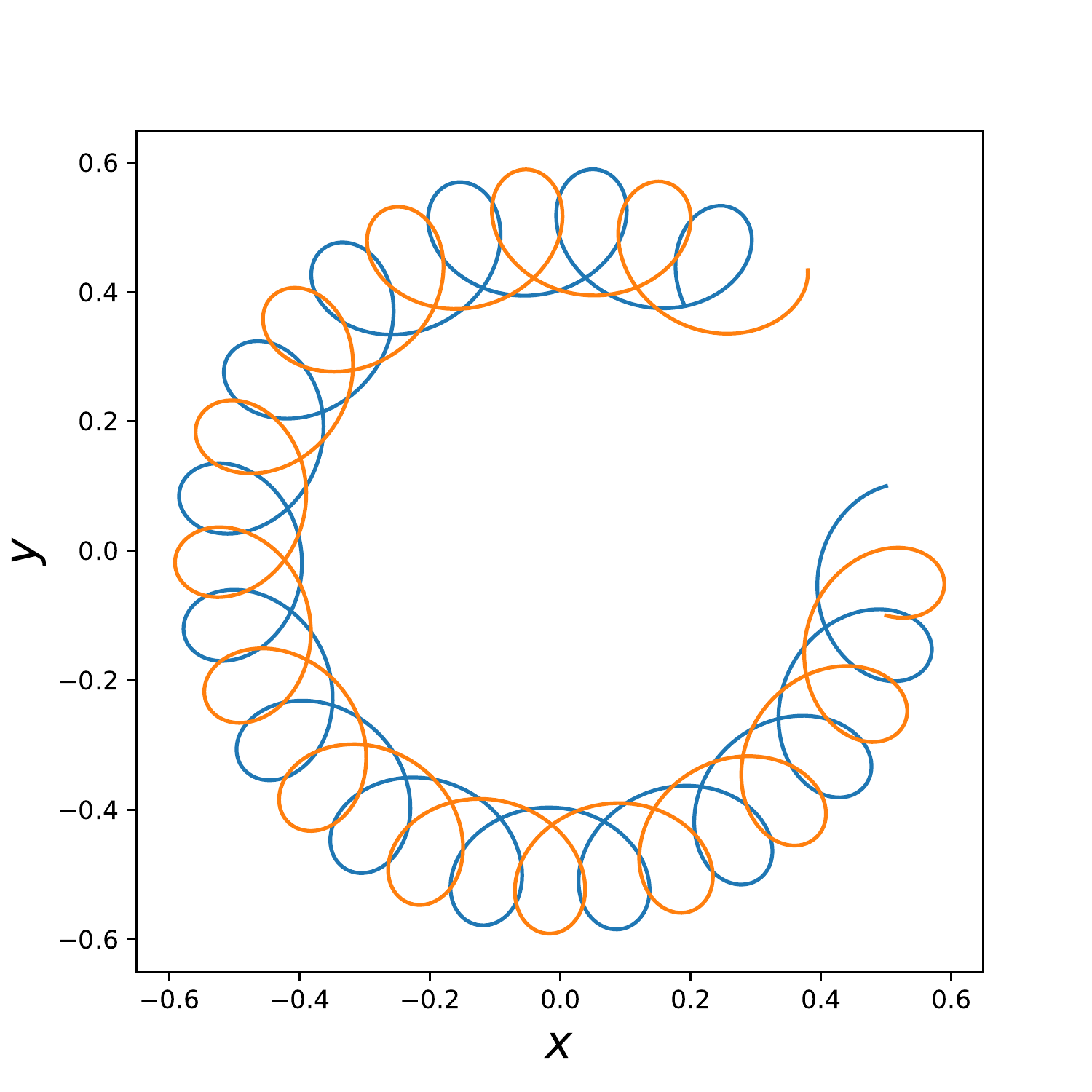}
		\caption{Two leapfrogging point vortices with position dependent circulation.}
		\label{fig:point_vortices}
	\end{center}
\end{figure}

Systems of planar point vortices~\cite{PullinSaffman:1991, Rowley:2001} provide a challenging problem for numerical integrators. Such systems are integrable for up to three vortices but produce chaotic behaviour for a minimum number of four vortices.
An interesting phenomenon is that of leapfrogging, which is usually observed only for two pairs of point vortices.
However, also one pair of point vortices can leapfrog by itself (see Figure~\ref{fig:point_vortices}) if the circulation is position dependent~\cite{Montaldi:2011}. In this case, the function $\vartheta$ in the Lagrange is nonlinear, hence this provides an interesting test case for our integrators.

We denote coordinates on $\mf{M}$ by $(\bx,\by) = (\bx^{1}, \hdots, \bx^{d}, \by^{1}, \hdots, \by^{d})$ and correspondingly coordinates on $\tb{\mf{M}}$ by $(\bx,\by,\bxd,\byd) = (\bx^{1}, \hdots, \bx^{d}, \by^{1}, \hdots, \by^{d}, \bxd^{1}, \hdots, \bxd^{d}, \byd^{1}, \hdots, \byd^{d})$. The coordinates on $\mf{M}$ are sometimes also collectively referred to by $q$.
The general form of the Lagrangian for point vortices is 
\begin{align}
L
= \dfrac{1}{2} \sum \limits_{i,j=1}^{d} \Gamma_{ij} (\bx^{i} \byd^{j} - \by^{i} \bxd^{j} )
- \dfrac{1}{4\pi} \sum \limits_{i \neq k}^{d} \sum \limits_{j \neq l}^{d} \Gamma_{ij} \Gamma_{kl} \, \log \big( (\bx^{i} - \bx^{k})^{2} + (\by^{j} - \by^{l})^{2} \big) ,
\end{align}
with $d$ the number of vortices and $\Gamma$ the matrix of vortex strengths, which is assumed to be of the form $\Gamma_{ij} = \gamma_{i} \delta_{ij}$, where $\gamma_{i}$ is the circulation of the $i$th vortex.
Here, we consider the special case of $\Gamma$ being position-dependent, specifically $\Gamma_{ij} (\bx^{i}, \by^{i}, \bx^{j}, \by^{j}) = \gamma_{i} \, S(\bx^{i}, \by^{i}) \, \delta_{ij}$, where we assume later on that $S$ is such that it has rotational symmetry.
For $d=2$ the Lagrangian thus becomes
\begin{align}
L &= \dfrac{1}{2} \big( \gamma^{1} \, S(\bx^{1}, \by^{1}) \, (\bx^{1} \byd^{1} - \by^{1} \bxd^{1} ) + \gamma^{2} \, S(\bx^{2}, \by^{2}) \, (\bx^{2} \byd^{2} - \by^{2} \bxd^{2} ) \big) 
- H (\bx, \by) , \\
H &= \dfrac{1}{2\pi} \gamma^{1} \gamma^{2} \, S(\bx^{1}, \by^{1}) \, S(\bx^{2}, \by^{2}) \, \log \big( (\bx^{1} - \bx^{2})^{2} + (\by^{1} - \by^{2})^{2} \big) .
\end{align}
The noncanonical symplectic form~\eqref{eq:noncanonical_symplectic_matrix} of this system reads
\begin{align}
\obar{\omega}
 = \sum \limits_{i=1}^{2} \gamma_{i} \, S(\bx^{i}, \by^{i}) \, \ext \bx^{i} \wedge \ext \by^{i}
 + \dfrac{1}{2} \sum \limits_{i=1}^{2} \gamma_{i} \left( \bx^{i} \, \dfrac{\partial S}{\partial \bx} (\bx^{i}, \by^{i}) + \by^{i} \, \dfrac{\partial S}{\partial \by} (\bx^{i}, \by^{i}) \right) \ext \bx^{i} \wedge \ext \by^{i} .
\end{align}
Assuming that the function $S$ is of the form $S(a,b) = s(a^{2} + b^{2})$ with some function $s : \rsp \rightarrow \rsp$, the Lagrangian is invariant under rotations of all coordinates by a constant angle $\vartheta$, that is, the following transformation of the coordinates,
\begin{align}
\sigma^{\eps} : 
\begin{pmatrix}
\bx^{i} \\
\by^{i} \\
\end{pmatrix}
\mapsto
\begin{pmatrix}
\bx^{i} \cos (\eps \vartheta) - \by^{i} \sin (\eps \vartheta) \\
\by^{i} \cos (\eps \vartheta) + \bx^{i} \sin (\eps \vartheta) \\
\end{pmatrix} ,
\end{align}
together with the corresponding transformation of the velocities, leaves the Lagrangian $L$ invariant.
The generating vector field is computed as
\begin{align}
V = \dfrac{d \sigma^{\eps}}{d\eps} \bigg\vert_{\eps=0}
  = - \by^{i} \dfrac{\partial}{\partial \bx^{i}} + \bx^{i} \dfrac{\partial}{\partial \by^{i}} ,
\end{align}
and the corresponding conserved quantity~\eqref{eq:noether_theorem} is obtained as
\begin{align}
P = \dfrac{\partial L}{\partial \bv} \cdot V
\nonumber
 &= \dfrac{1}{2} \sum \limits_{i=1}^{2} \gamma_{i} \, \big( \abs{\bx^{i}}^{2} + \abs{\by^{i}}^{2} \big) \, S(\bx^{i}, \by^{i}) \\
 &= q^{1} \vartheta_{2} (\bx, \by) - q^{2} \vartheta_{1} (\bx, \by)
  + q^{3} \vartheta_{4} (\bx, \by) - q^{4} \vartheta_{3} (\bx, \by)  .
\end{align}
We are particularly interested in the behaviour of this angular momentum under the various projection methods.

We consider the simple case of $s(r) = 1 + r$, so that $S(a,b) = 1 + a^{2} + b^{2}$ and the functions for the momenta are computed as
\begin{subequations}
\begin{align}
\vartheta_{1} (\bx,\by) &= - \tfrac{1}{2} \gamma^{1} \, \by^{1} \, S(\bx^{1}, \by^{1}) , &
\vartheta_{2} (\bx,\by) &= \tfrac{1}{2} \gamma^{1} \, \bx^{1} \, S(\bx^{1}, \by^{1}) , \\
\vartheta_{3} (\bx,\by) &= - \tfrac{1}{2} \gamma^{2} \, \by^{2} \, S(\bx^{2}, \by^{2}) , &
\vartheta_{4} (\bx,\by) &= \tfrac{1}{2} \gamma^{2} \, \bx^{2} \, S(\bx^{2}, \by^{2}) , 
\end{align}
\end{subequations}
and those for the forces as
\begin{subequations}
\begin{align}
f^{1} (\bx,\by,\bxd,\byd)
&= \tfrac{1}{2} \gamma^{1} \, (\bx^{1} \byd^{1} - \by^{1} \bxd^{1} ) \, S^{(\bx)} (\bx^{1}, \by^{1})
 + \tfrac{1}{2} \gamma^{1} \byd^{1} \, S(\bx^{1}, \by^{1})
 - \nabla_{1} H (\bx,\by) , \\
f^{2} (\bx,\by,\bxd,\byd)
&= \tfrac{1}{2} \gamma^{2} \, (\bx^{2} \byd^{2} - \by^{2} \bxd^{2} ) \, S^{(\bx)} (\bx^{2}, \by^{2})
 + \tfrac{1}{2} \gamma^{2} \byd^{2} \, S(\bx^{2}, \by^{2})
 - \nabla_{2} H (\bx,\by) , \\
f_{3} (\bx,\by,\bxd,\byd)
&= \tfrac{1}{2} \gamma^{1} \, (\bx^{1} \byd^{1} - \by^{1} \bxd^{1} ) \, S^{(\by)} (\bx^{1}, \by^{1})
 - \tfrac{1}{2} \gamma^{1} \bxd^{1} \, S(\bx^{1}, \by^{1})
 - \nabla_{3} H (\bx,\by) , \\
f_{4} (\bx,\by,\bxd,\byd)
&= \tfrac{1}{2} \gamma^{2} \, (\bx^{2} \byd^{2} - \by^{2} \bxd^{2} ) \, S^{(\by)} (\bx^{2}, \by^{2})
 - \tfrac{1}{2} \gamma^{2} \bxd^{2} \, S(\bx^{2}, \by^{2})
 - \nabla_{4} H (\bx,\by) ,
\end{align}
\end{subequations}
with the gradient of the Hamiltonian being
\begingroup
\allowdisplaybreaks
\begin{subequations}
\begin{align}
\nabla_{1} H (\bx,\by)
\nonumber
&= \dfrac{\gamma^{1} \gamma^{2}}{2\pi} \, S^{(\bx)} (\bx^{1}, \by^{1}) \, S(\bx^{2}, \by^{2}) \, \log \big( (\bx^{1} - \bx^{2})^{2} + (\by^{1} - \by^{2})^{2} \big) \\
&+ \dfrac{\gamma^{1} \gamma^{2}}{\pi} \, S(\bx^{1}, \by^{1}) \, S(\bx^{2}, \by^{2}) \, \dfrac{\bx^{1} - \bx^{2}}{(\bx^{1} - \bx^{2})^{2} + (\by^{1} - \by^{2})^{2}} , \\
\nabla_{2} H (\bx,\by)
\nonumber
&= \dfrac{\gamma^{1} \gamma^{2}}{2\pi} \, S^{(\bx)} (\bx^{2}, \by^{2}) \, S(\bx^{1}, \by^{1}) \, \log \big( (\bx^{1} - \bx^{2})^{2} + (\by^{1} - \by^{2})^{2} \big) \\
&- \dfrac{\gamma^{1} \gamma^{2}}{\pi} \, S(\bx^{1}, \by^{1}) \, S(\bx^{2}, \by^{2}) \, \dfrac{\bx^{1} - \bx^{2}}{(\bx^{1} - \bx^{2})^{2} + (\by^{1} - \by^{2})^{2}} , \\
\nabla_{3} H (\bx,\by)
\nonumber
&= \dfrac{\gamma^{1} \gamma^{2}}{2\pi} \, S^{(\by)} (\bx^{1}, \by^{1}) \, S(\bx^{2}, \by^{2}) \, \log \big( (\bx^{1} - \bx^{2})^{2} + (\by^{1} - \by^{2})^{2} \big) \\
&+ \dfrac{\gamma^{1} \gamma^{2}}{\pi} \, S(\bx^{1}, \by^{1}) \, S(\bx^{2}, \by^{2}) \, \dfrac{\by^{1} - \by^{2}}{(\bx^{1} - \bx^{2})^{2} + (\by^{1} - \by^{2})^{2}} , \\
\nabla_{4} H (\bx,\by)
\nonumber
&= \dfrac{\gamma^{1} \gamma^{2}}{2\pi} \, S^{(\by)} (\bx^{2}, \by^{2}) \, S(\bx^{1}, \by^{1}) \, \log \big( (\bx^{1} - \bx^{2})^{2} + (\by^{1} - \by^{2})^{2} \big) \\
&- \dfrac{\gamma^{1} \gamma^{2}}{\pi} \, S(\bx^{1}, \by^{1}) \, S(\bx^{2}, \by^{2}) \, \dfrac{\by^{1} - \by^{2}}{(\bx^{1} - \bx^{2})^{2} + (\by^{1} - \by^{2})^{2}} ,
\end{align}
\end{subequations}
\endgroup
where $S^{(\bx)}$ and $S^{(\by)}$ denote the $\bx$ and $\by$ derivative of $S$, respectively.

We use the time step $h = 0.1$, circulations $\gamma_{1} = \gamma_{2} = 0.1$ and initial conditions $\cq_{0} = (1.0, 0.1, 1.0, -0.1)$. This setup leads to a circular leapfrogging of the two point vortices. \\

We make the following observations:
\begin{itemize}
\item All methods except the two-stage Gauss--Lobatto--IIIA and IIIB method and all of the Gauss--Lobatto--IIIC methods are stable even without projection, although with reduced order (c.f., Section~\ref{sec:app_convergence}).

\item For all methods except the two-stage Gauss--Lobatto--IIIA and IIIB and the Gauss--Lobatto--IIIC methods method we observe that the angular momentum oscillates about its initial value where the amplitude of the oscillation seems bounded, that is the angular momentum seems to be preserved in a nearby sense, similar to the energy with symplectic integrators.

\item The standard projection worsens the result for all methods except the two-stage Gauss--Lobatto--IIIA and IIIB and the Gauss--Lobatto--IIIC methods method, which do not crash when applying the projection, but the projected methods still show large errors and do not provide satisfactory results.

\item The symplectic, symmetric and midpoint projections lead to very good results with almost all methods, restoring the original order of the methods and showing good long-time behaviour of both the energy and the angular momentum. There are some exceptions, however:

\begin{itemize}
\item The symplectic projection applied to the two-stage Gauss--Lobatto--IIIA and IIIB and the Gauss--Lobatto--IIIC methods is just as unstable as the corresponding unprojected methods.
\item Both, the symmetric and midpoint projection applied to all of the Gauss--Lobatto--IIIC methods lead to an improved behaviour compared to the unprojected case, but exhibit a strong drift in the energy.
\item The Gauss--Lobatto--IIID and IIIE methods with an even number of stages together with the midpoint projection exhibit a rather erratic behaviour in the energy error.
\end{itemize}

\item For the symmetric projection and the higher-order methods (e.g. Gauss--Legendre with four or more stages or Gauss--Lobatto-IIID with four stages), we observe a small drift in the angular momentum, but over 1~000~000 time steps this drift is of the order of $10^{-12}$. This drift is most likely caused by round-off errors (see Section~\ref{sec:app_energy_drift} for more details).
\end{itemize}

In summary, the numerical experiments suggest that the combination of almost all integration methods and all projection methods excluding the standard projection provide suitable integration algorithms for the point vortex example. Exceptions are the Gauss--Lobatto--IIIC methods with any of the projection methods and the combination of the midpoint projection with the two-stage Gauss--Lobatto--IIIA and IIIB method and Gauss--Lobatto--IIID and IIIE methods with an even number of stages.

\subsection{Guiding Centre Dynamics}
\label{app:guiding_centre_dynamics}

In plasma physics, the search for geometric integrators for guiding centre dynamics and gyrokinetics is currently of great interest.
As the Hamiltonian structure of the guiding centre system is noncanonical, there are practically no standard methods which can be easily applied.
As the guiding centre equations can also be obtained from a Lagrangian, the application of variational integrators seems natural and has recently been tried by various researchers~\cite{Qin:2008, Qin:2009, Ellison:2015, Ellison:2016}.
However, the guiding centre Lagrangian is degenerate, leading to all the problems discussed so far. We will see in the following if our projection methods can overcome these deficits.

Guiding centre dynamics~\cite{Northrop:1961cd} is a reduced version of charged particle dynamics, where the motion of the particle in a strong magnetic field $B$ is reduced to the motion of the guiding centre, that is the centre of the gyro motion of the particle about a magnetic field line.
The dynamics of the guiding centre can be described in terms of only four coordinates (as compared to six for the full motion of the charged particle), the position of the guiding centre $\bx$ and the parallel velocity $\bu$, where parallel refers to the direction of the magnetic field.
Denoting coordinates on $\mf{M}$ by $(\bx, \bu) = (\bx^{1}, \bx^{2}, \bx^{3}, \bu)$ and correspondingly coordinates on $\tb{\mf{M}}$ by $(\bx, \bu, \bxd, \bud) = (\bx^{1}, \bx^{2}, \bx^{3}, \bu, \bxd^{1}, \bxd^{2}, \bxd^{3}, \bud)$, the guiding centre Lagrangian~\cite{Littlejohn:1983, Cary:2009} can be written as
\begin{align}
L &= (A (\bx) + \bu \, b (\bx)) \cdot \bv - H(\bx,\bu) , &
H &= \tfrac{1}{2} \bu^{2} + \mu \abs{B (\bx)} ,
\end{align}
where $b = B / \abs{B}$ is the unit vector of the magnetic field $B = \nabla \times A$ with $A$ the magnetic vector potential and $\mu$ is the magnetic moment. The first term in $H$ denotes the parallel part of the kinetic energy and the second term the perpendicular part (parallel and perpendicular to the direction of the magnetic field). Here, we consider the case of only a magnetic field with vanishing electrostatic potential.

Denoting a curve in $\tb{\mf{M}}$ by $t \mapsto (x(t), u(t))$, the Euler--Lagrange equations~\eqref{eq:euler_lagrange_equations} are computed as follows,
\begin{align}
\nabla \vartheta^{T} \big( x(t), u(t) \big) \cdot \dot{x} (t) - \dot{\vartheta} \big( x(t), u(t) \big) &= \nabla H \big( x(t), u(t) \big) , \\
b \big( x(t) \big) \cdot \dot{x} (t) &= u (t) ,
\end{align}
with $\vartheta(\bx,\bu) = A (\bx) + \bu \, b (\bx)$ and the gradient denoting the derivative with respect to $\bx$.
This can be rewritten in an explicit form as
\begin{subequations}
\begin{align}
\dot{x} (t) &= \dfrac{u (t) \, \beta (x (t))}{b (x (t)) \cdot \beta (x (t))} + \dfrac{B (x (t))}{B (x (t)) \cdot \beta (x (t))} \times \nabla H (x (t),u (t)) , \\
\dot{u} (t) &= - \dfrac{\beta (x (t))}{b (x (t)) \cdot \beta (x (t))} \cdot \nabla H (x (t),u (t)) ,
\end{align}
\end{subequations}
where $\beta = \nabla \times \vartheta$.
The noncanonical symplectic form~\eqref{eq:noncanonical_symplectic_matrix} is given by
\begin{align}
\obar{\omega} = \tfrac{1}{2} ( \vartheta_{j,i} (\bx, \bu) - \vartheta_{i,j} (\bx, \bu) ) \, \ext \bx^{i} \wedge \ext \bx^{j} - b_{i} (\bx) \, \ext \bx^{i} \wedge \ext \bu .
\end{align}
Let us assume that the magnetic field $B$ is not uniform, but that both $A$ and $B$ do not depend on one of the coordinates, say $\bx_{3}$.
Than we have a symmetry for the transformation
\begin{align}
\sigma^{\eps} : \bx^{3} \mapsto \bx^{3} + \eps ,
\end{align}
with generating vector field
\begin{align}
V = \dfrac{d \sigma^{\eps}}{d\eps} \bigg\vert_{\eps=0}
  = \dfrac{\partial}{\partial \bx_{3}} .
\end{align}
The corresponding conserved momentum map~\eqref{eq:noether_theorem},
\begin{align}
P = \dfrac{\partial L}{\partial \bv} \cdot V
  = \vartheta_{3} (\bx,\bu) ,
\end{align}
which, depending on the actual form of $\vartheta$, can be quite complicated and is therefore a good test for our algorithms.
Although the basic integrator will preserve this toroidal momentum if the discrete Lagrangian preserves the corresponding symmetry, the projection could potentially modify its value.
The projection guarantees preservation of the constraint $p_{n+1} = \vartheta(q_{n+1})$ but it does not guarantee that $p_{n+1} = p_{n}$.

In the numerical experiments, we use toroidal coordinates $\bx = (R, Z, \phy)$, where $R$, $Z$ and $\phy$ denote the radial, vertical and toroidal direction, respectively.
For the magnetic field $B$ and the vector potential $A$ we will use analytic expressions following~\citet{Qin:2009}.
The vector potential is given as
\begin{align}
A_{R} &= \dfrac{B_{0} R_{0} Z}{2 R} , &
A_{Z} &= - \ln \bigg( \dfrac{R}{R_{0}} \bigg) \, \dfrac{B_{0} R_{0}}{2} , &
A_{\phy} &= - \dfrac{B_{0} r^{2}}{2 q_{0} R} .
\end{align}
The magnetic field $B = \nabla \times A$ is computed as
\begin{align}
B_{R   } &= - \dfrac{B_{0} Z}{q_{0} R} , &
B_{Z   } &=   \dfrac{B_{0} \, (R - R_{0})}{q_{0} R} , &
B_{\phy} &= - \dfrac{B_{0} R_{0}}{R} , &
\abs{B} &= \dfrac{B_{0} S}{q_{0} R} , 
\end{align}
and the normalized magnetic field as
\begin{align}
b_{R   } &= - \dfrac{Z}{S} , &
b_{Z   } &= \dfrac{R - R_{0}}{S} , &
b_{\phy} &= - \dfrac{q_{0} R_{0}}{S} .
\end{align}
Here, $R_{0}$ is the radial position of the magnetic axis, $B_{0}$ is the magnetic field at $R_{0}$, and $\cq_{0}$ is the safety factor, regarded as constant. In all of the examples, these constants are set to $R_{0} = 2$, $B_{0} = 5$ and $q = 2$, respectively.
The functions $r$ and $S$ are given by
\begin{align}
r &= \sqrt{ (R - R_{0})^{2} + Z^{2} } , &
S &= \sqrt{ r^{2} + q_{0}^{2} R_{0}^{2} } .
\end{align}
In toroidal coordinates, the functions for the momenta are
\begin{subequations}
\begin{align}
\vartheta_{1} (\bx,\bu) &= A_{R} (\bx) + \bu \, b_{R} (\bx) , &
\vartheta_{2} (\bx,\bu) &= A_{Z} (\bx) + \bu \, b_{Z} (\bx) , \\
\vartheta_{3} (\bx,\bu) &= R ( A_{\phy} (\bx) + \bu \, b_{\phy} (\bx) ) , &
\vartheta_{4} (\bx,\bu) &= 0 , 
\end{align}
\end{subequations}
and those for the forces are computed as
\begin{subequations}
\begin{align}
f_{1} (\bx,\bu,\bxd,\bud)
&= \vartheta_{1,1} (\bx,\bu) \, \bxd^{1} + \vartheta_{2,1} (\bx,\bu) \, \bxd^{2} + \vartheta_{3,1} (\bx,\bu) \, \bxd^{3}
 - \nabla_{1} H (\bx,\bu) , \\
f_{2} (\bx,\bu,\bxd,\bud)
&= \vartheta_{1,2} (\bx,\bu) \, \bxd^{1} + \vartheta_{2,2} (\bx,\bu) \, \bxd^{2} + \vartheta_{3,2} (\bx,\bu) \, \bxd^{3}
 - \nabla_{2} H (\bx,\bu) , \\
f_{3} (\bx,\bu,\bxd,\bud)
&= \vartheta_{1,3} (\bx,\bu) \, \bxd^{1} + \vartheta_{2,3} (\bx,\bu) \, \bxd^{2} + \vartheta_{3,3} (\bx,\bu) \, \bxd^{3}
 - \nabla_{3} H (\bx,\bu) , \\
f_{4} (\bx,\bu,\bxd,\bud)
&= b(x) \cdot \bxd
 - \nabla_{4} H (\bx,\bu) ,
\end{align}
\end{subequations}
with the gradient of the Hamiltonian being
\begin{subequations}
\begin{align}
\nabla_{1} H (\bx,\bu) &= \nabla_{1} \abs{B(\bx)} , &
\nabla_{2} H (\bx,\bu) &= \nabla_{2} \abs{B(\bx)} , \\
\nabla_{3} H (\bx,\bu) &= \nabla_{3} \abs{B(\bx)} , &
\nabla_{4} H (\bx,\bu) &= \bu .
\end{align}
\end{subequations}

\begin{figure}[tb]
	\begin{center}
		\subfloat[Deeply Passing]{
			\includegraphics[width=.24\textwidth]{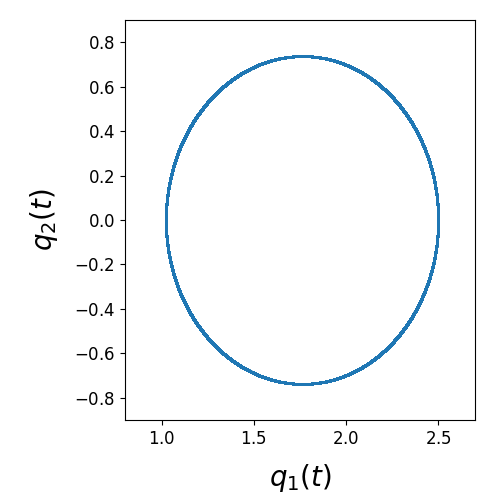}
		}
		\subfloat[Barely Passing]{
			\includegraphics[width=.24\textwidth]{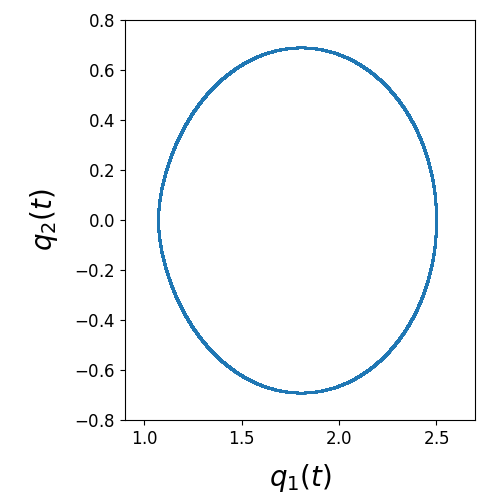}
		}
		\subfloat[Barely Trapped]{
			\includegraphics[width=.24\textwidth]{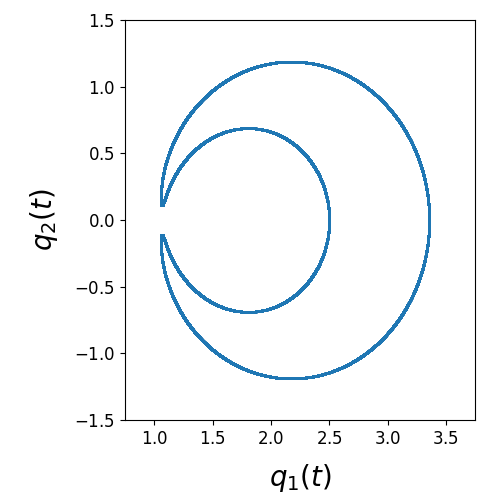}
		}
		\subfloat[Deeply Trapped]{
			\includegraphics[width=.24\textwidth]{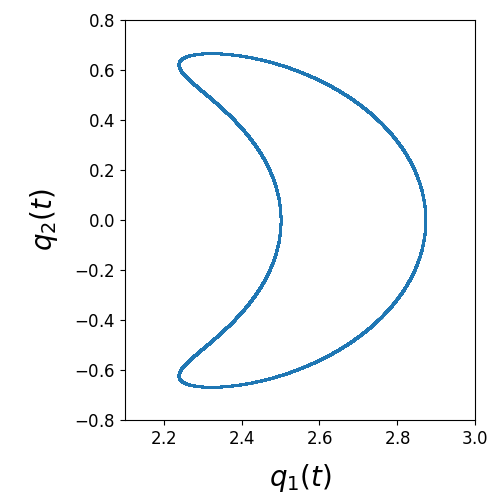}
		}
		
		\caption{Guiding centre test particles. Projection of the guiding centre trajectory onto the $(R,Z)$ plane.}
		\label{fig:guiding_center_4d}
	\end{center}
\end{figure}

We consider four different initial conditions, that of a deeply trapped particle, a barely trapped particle, a barely passing particle and a deeply passing particle (c.f., Figure~\ref{fig:guiding_center_4d}). The second and third case are particularly challenging. In all three cases we choose $(R,Z,\phy) = (2.5, 0, 0)$ and set the magnetic moment $\mu = 0.01$. The parallel velocity $u$ and the time step $h$ are chosen as follows:
\begin{center}
\begin{tabular}{|l|c|c|c|c|}
\hline
 & deeply trapped & barely trapped & barely passing & deeply passing \\
\hline
$u$ & $0.1$ & $0.3375$ & $0.3425$ & $0.5$ \\
\hline
$h$ & $5.0$ & $3.0$    & $2.5$    & $2.5$ \\
\hline
\end{tabular}
\end{center}

\begin{figure}[p]
	\begin{center}
		\subfloat[GLRK1]{\label{fig:guiding_centre_4d_barely_passing_vprk_pnone_glrk1}
			\includegraphics[height=2.5cm]{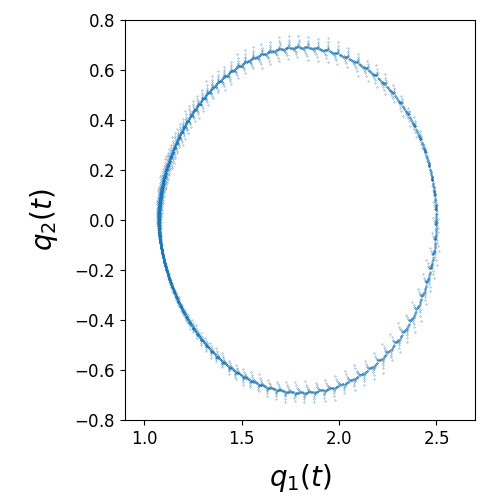}
			\includegraphics[height=2.5cm]{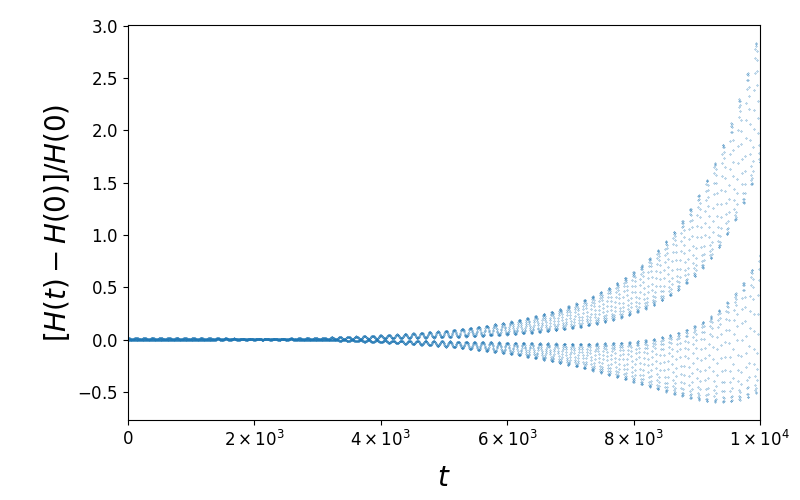}
			\includegraphics[height=2.5cm]{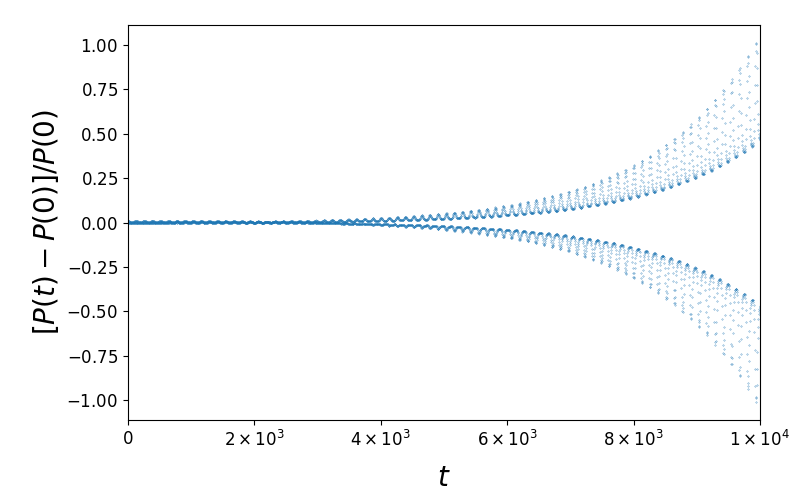}
		}
		
		\subfloat[GLRK2]{\label{fig:guiding_centre_4d_barely_passing_vprk_pnone_glrk2}
			\includegraphics[height=2.5cm]{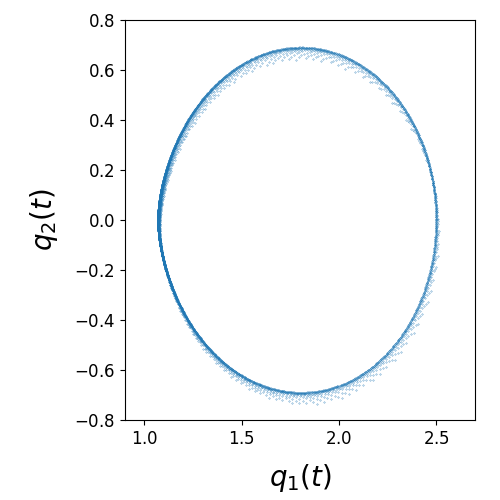}
			\includegraphics[height=2.5cm]{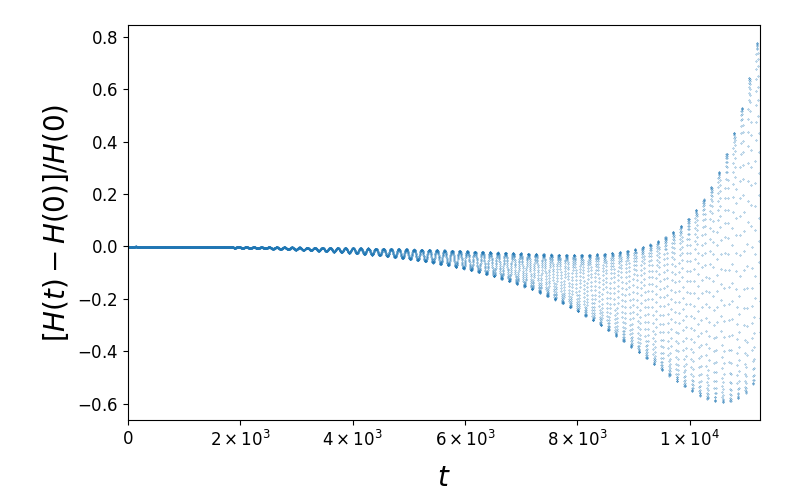}
			\includegraphics[height=2.5cm]{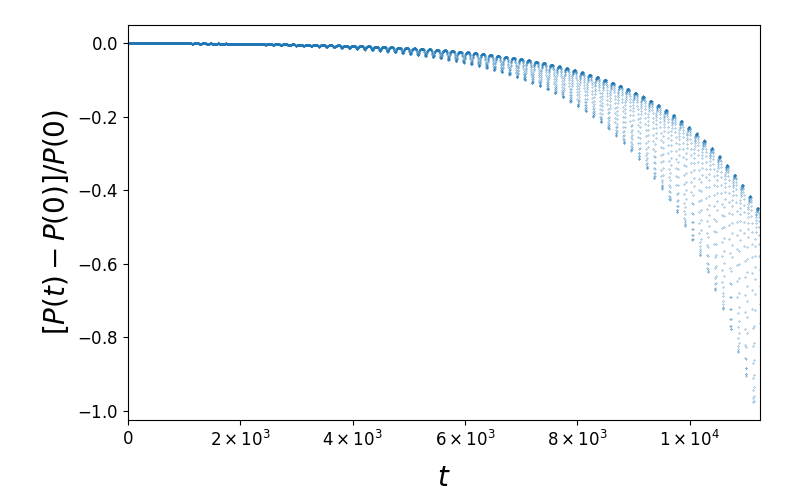}
		}

		\subfloat[GLRK3]{\label{fig:guiding_centre_4d_barely_passing_vprk_pnone_glrk3}
			\includegraphics[height=2.5cm]{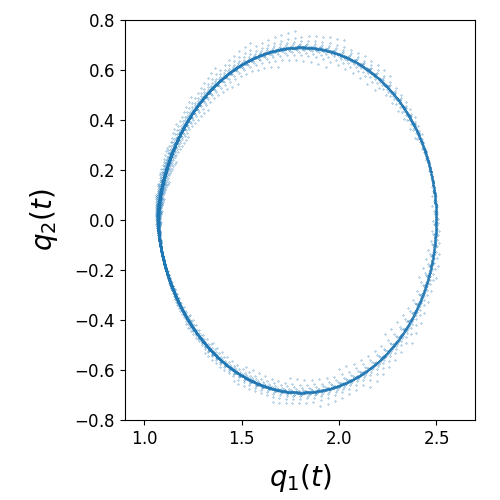}
			\includegraphics[height=2.5cm]{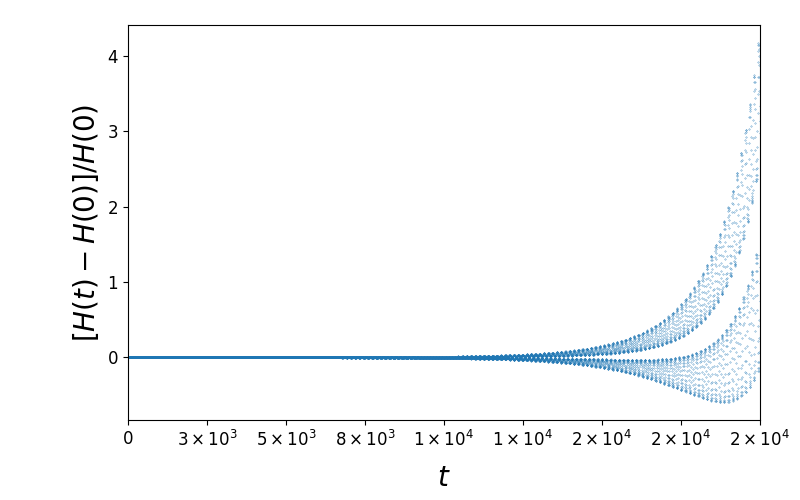}
			\includegraphics[height=2.5cm]{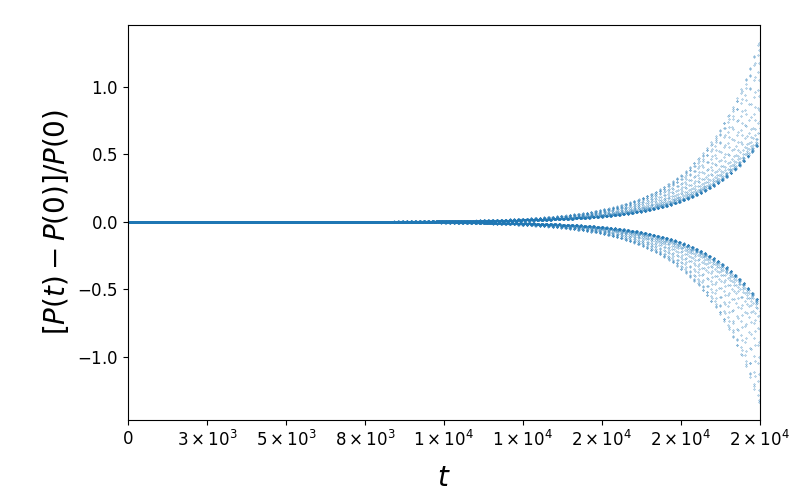}
		}
		
		\subfloat[GLRK4]{\label{fig:guiding_centre_4d_barely_passing_vprk_pnone_glrk4}
			\includegraphics[height=2.5cm]{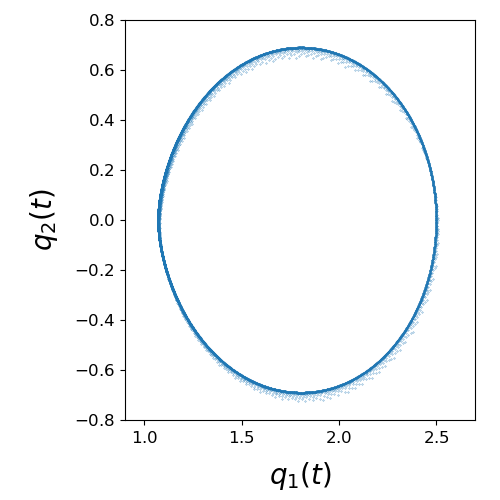}
			\includegraphics[height=2.5cm]{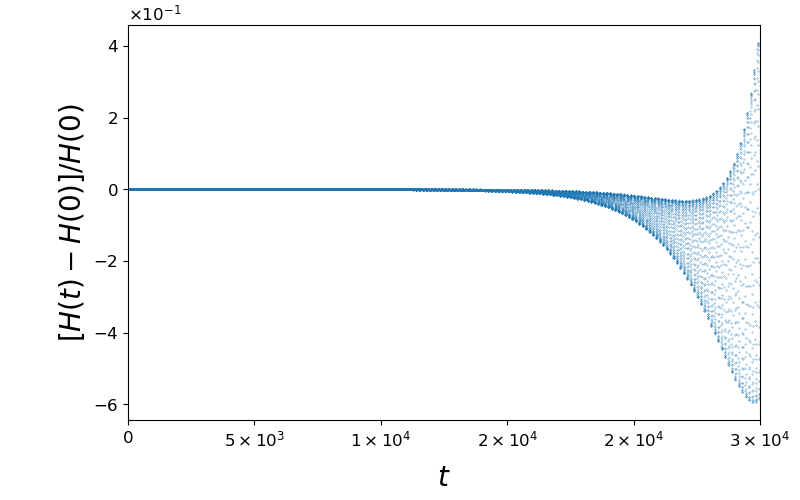}
			\includegraphics[height=2.5cm]{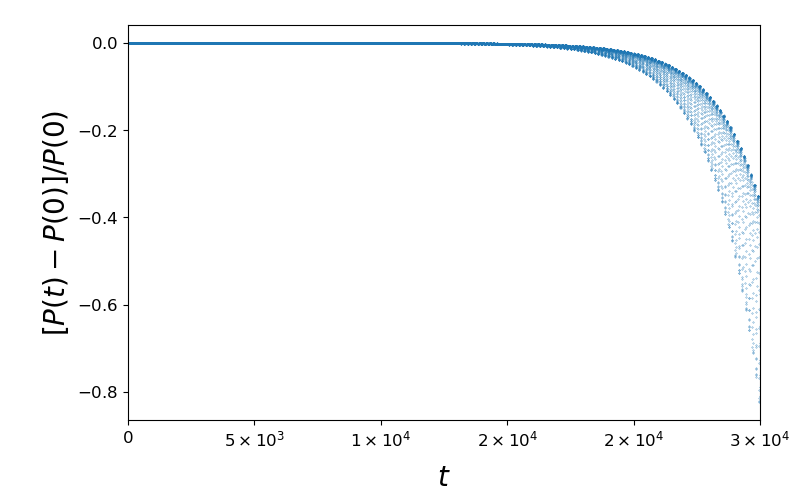}
		}
		
		\subfloat[SRK3]{\label{fig:guiding_centre_4d_barely_passing_vprk_pnone_srk3}
			\includegraphics[height=2.5cm]{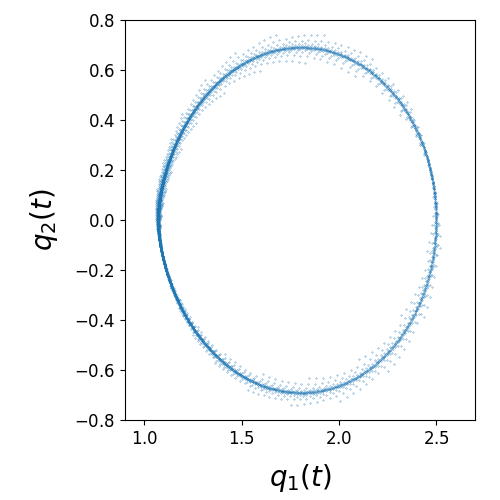}
			\includegraphics[height=2.5cm]{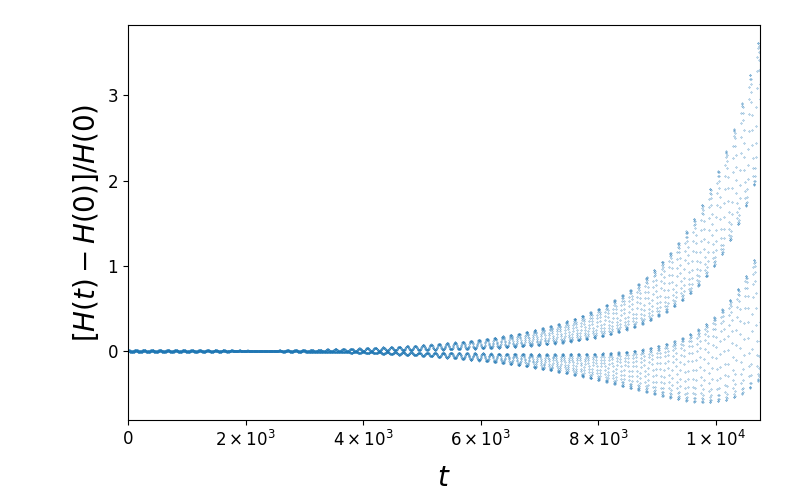}
			\includegraphics[height=2.5cm]{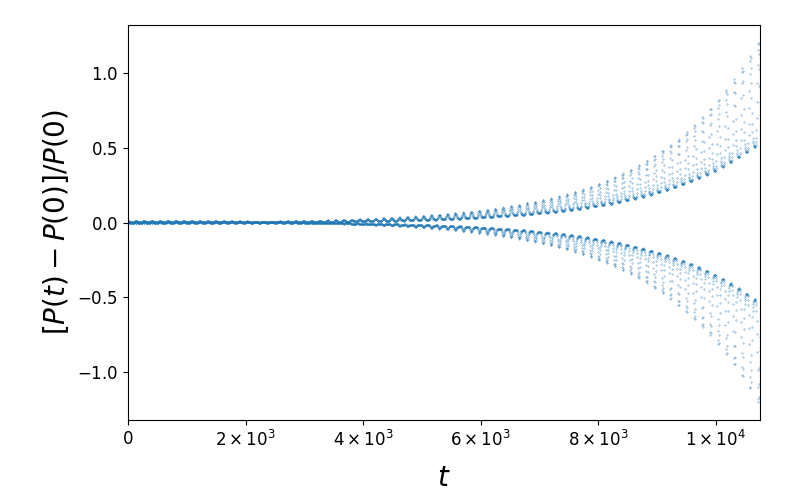}
		}
		
		\subfloat[RadIIA2]{\label{fig:guiding_center_4d_barely_passing_vprk_radIIA2}
			\includegraphics[height=2.5cm]{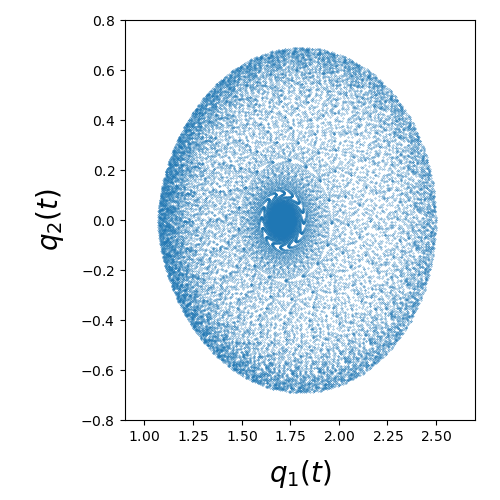}
			\includegraphics[height=2.5cm]{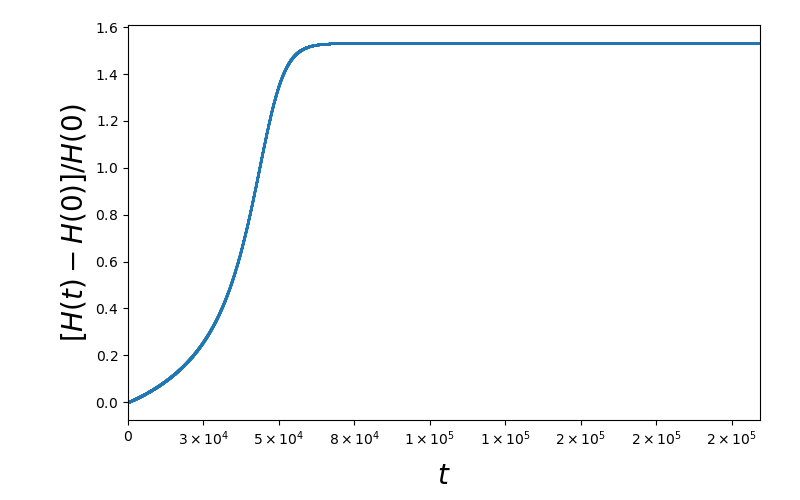}
			\includegraphics[height=2.5cm]{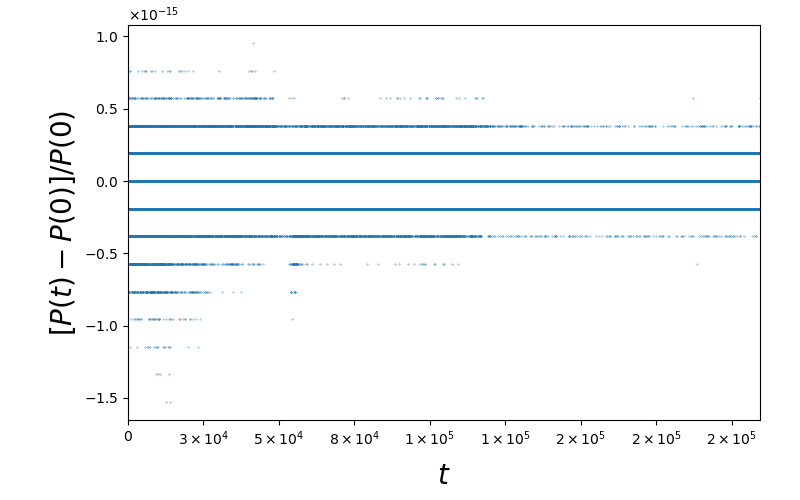}
		}
		
		\subfloat[RadIIA3]{\label{fig:guiding_center_4d_barely_passing_vprk_radIIA3}
			\includegraphics[height=2.5cm]{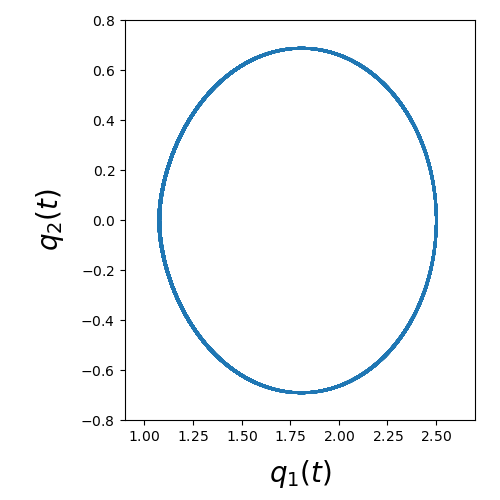}
			\includegraphics[height=2.5cm]{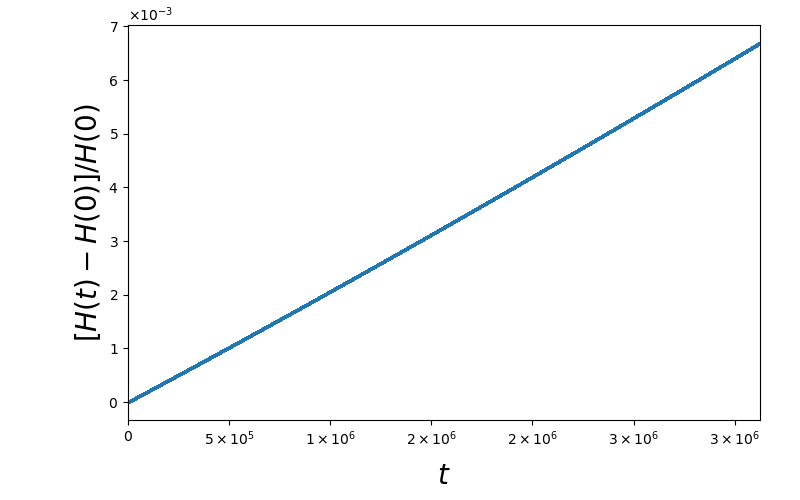}
			\includegraphics[height=2.5cm]{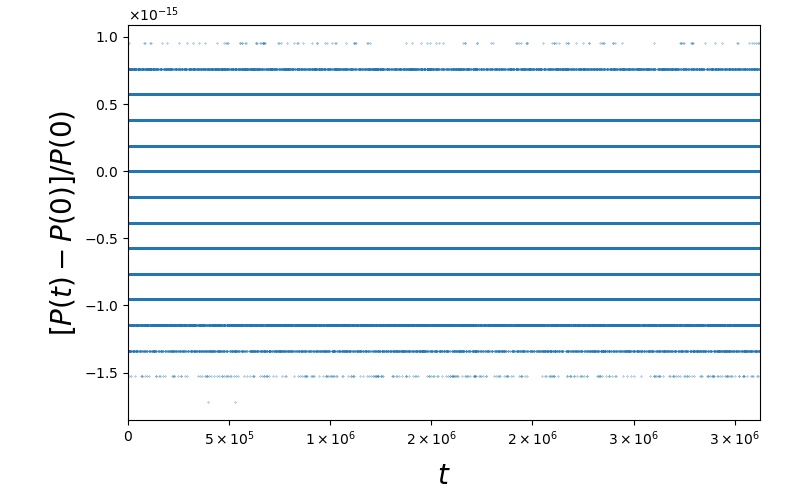}
		}
		
		\caption{Barely passing guiding centre particle with variational Runge--Kutta methods.}
		\label{fig:guiding_centre_4d_barely_passing_vprk_pnone}
	\end{center}
\end{figure}

\begin{figure}[p]
	\begin{center}
		\subfloat[GLRK1]{\label{fig:guiding_centre_4d_barely_trapped_vprk_pnone_glrk1}
			\includegraphics[height=2.5cm]{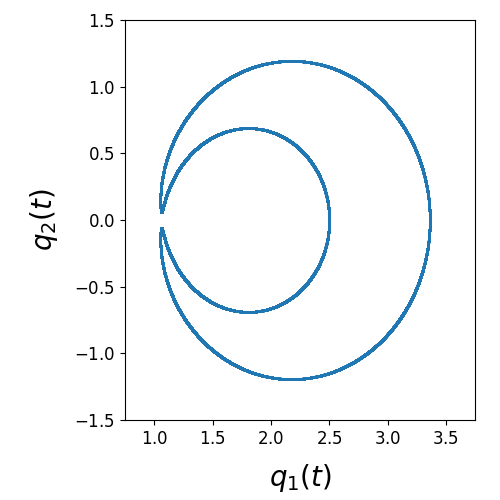}
			\includegraphics[height=2.5cm]{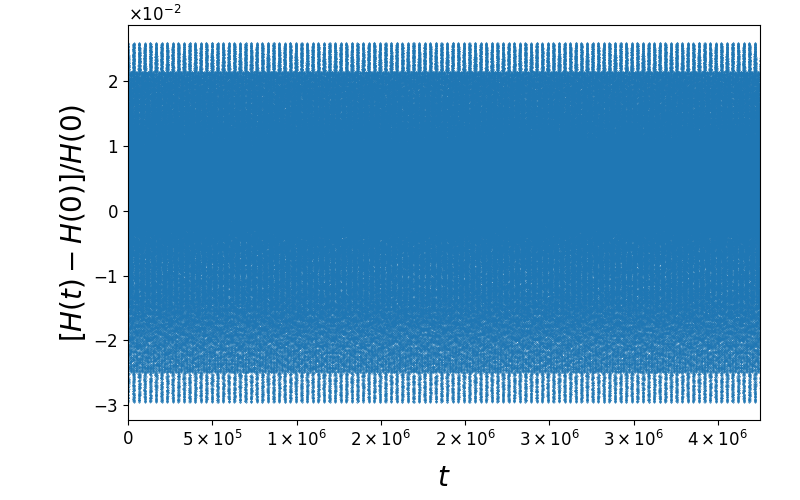}
			\includegraphics[height=2.5cm]{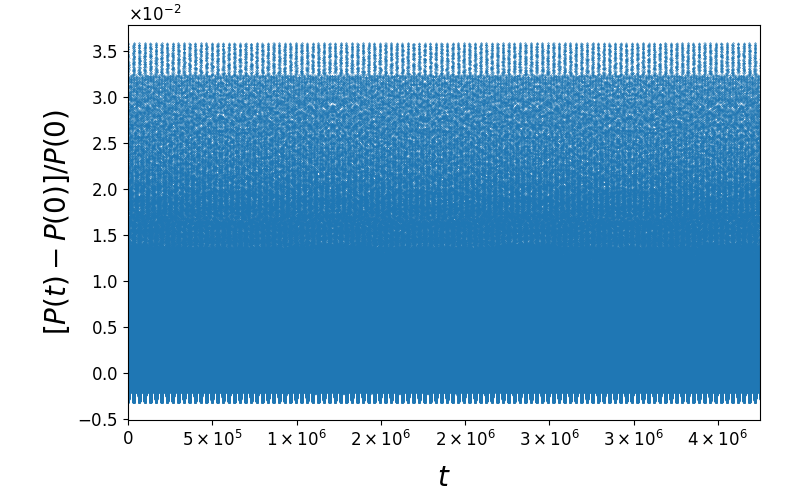}
		}
		
		\subfloat[GLRK2]{\label{fig:guiding_centre_4d_barely_trapped_vprk_pnone_glrk2}
			\includegraphics[height=2.5cm]{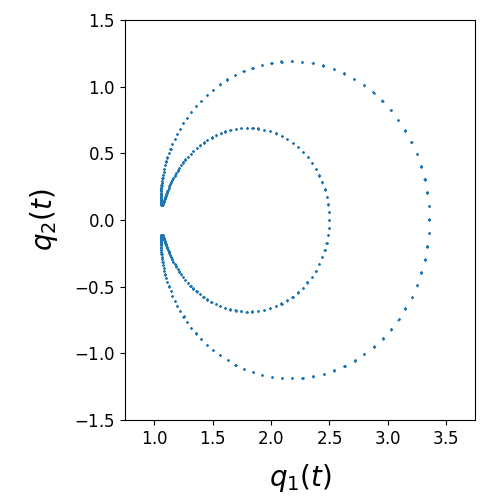}
			\includegraphics[height=2.5cm]{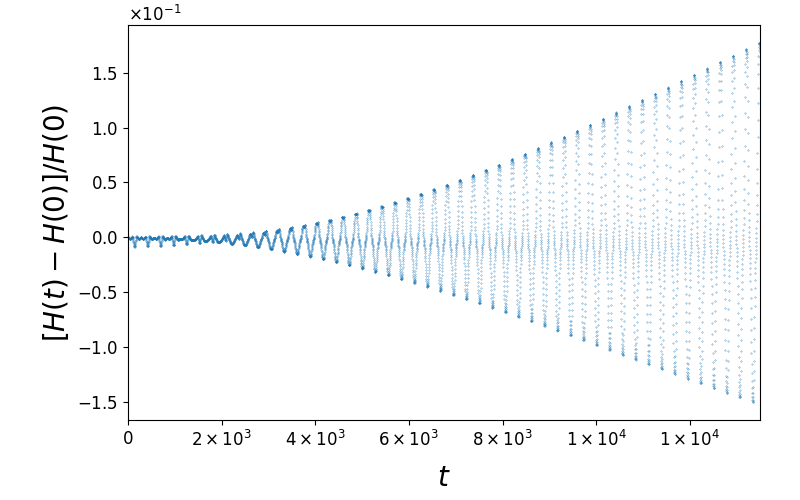}
			\includegraphics[height=2.5cm]{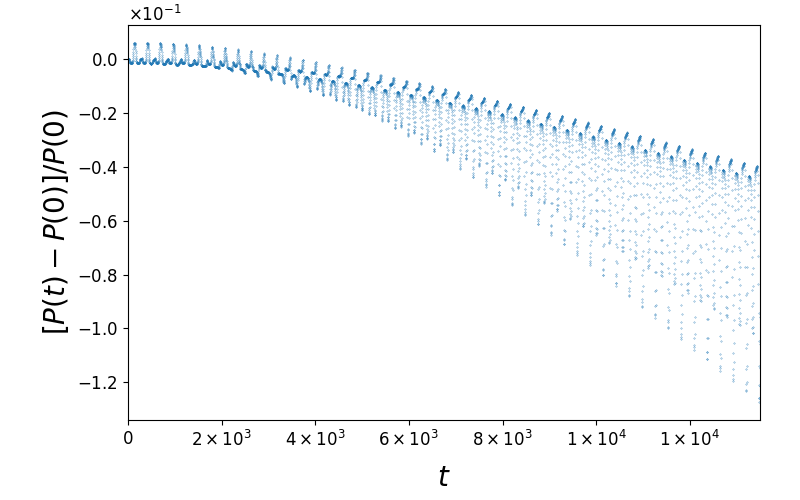}
		}
		
		\subfloat[GLRK3]{\label{fig:guiding_centre_4d_barely_trapped_vprk_pnone_glrk3}
			\includegraphics[height=2.5cm]{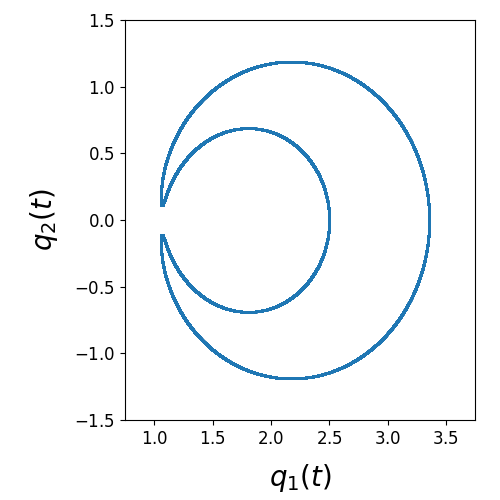}
			\includegraphics[height=2.5cm]{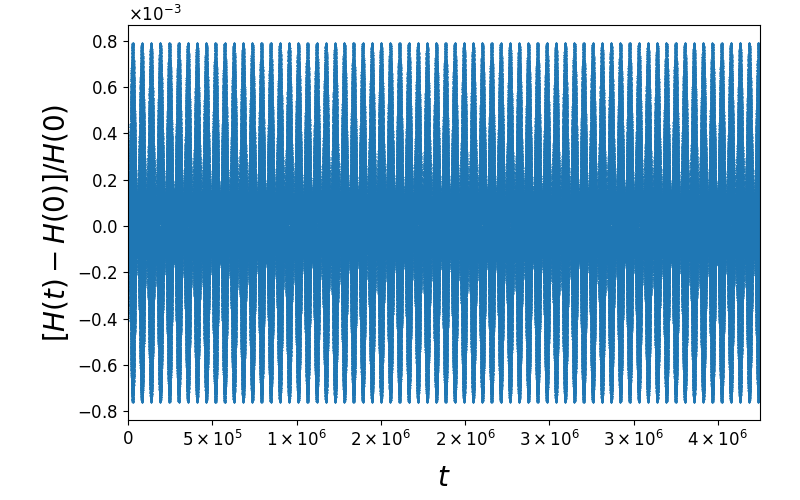}
			\includegraphics[height=2.5cm]{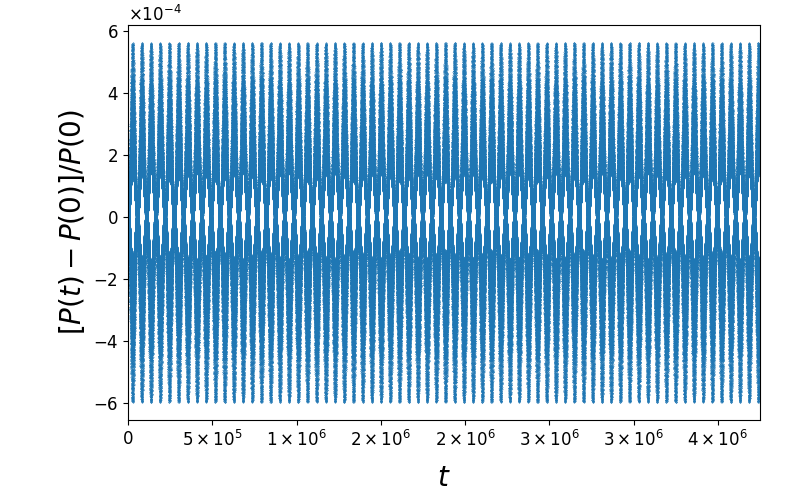}
		}
		
		\subfloat[GLRK4]{\label{fig:guiding_centre_4d_barely_trapped_vprk_pnone_glrk4}
			\includegraphics[height=2.5cm]{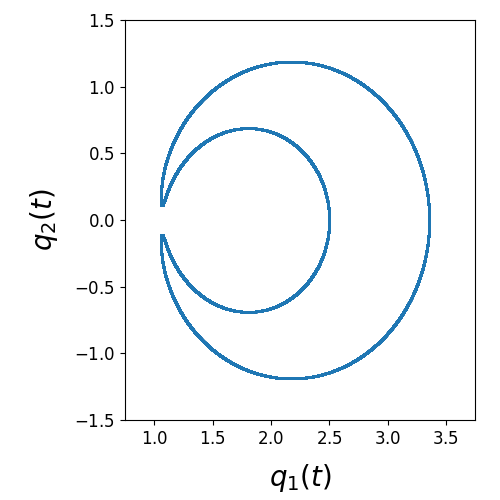}
			\includegraphics[height=2.5cm]{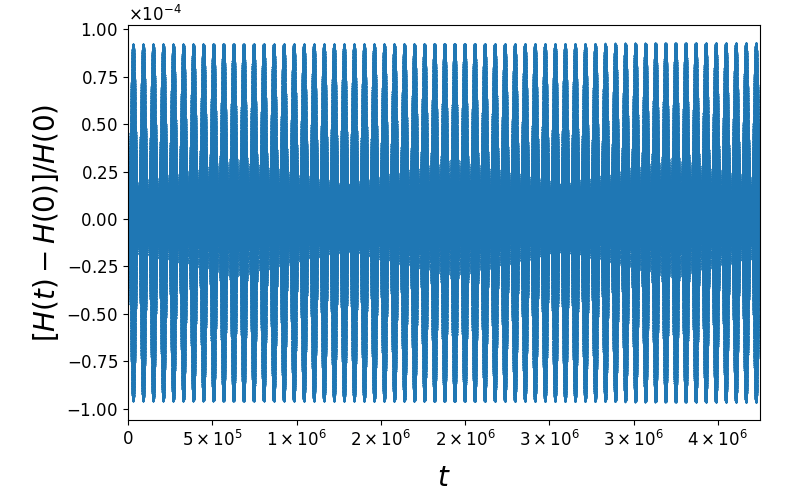}
			\includegraphics[height=2.5cm]{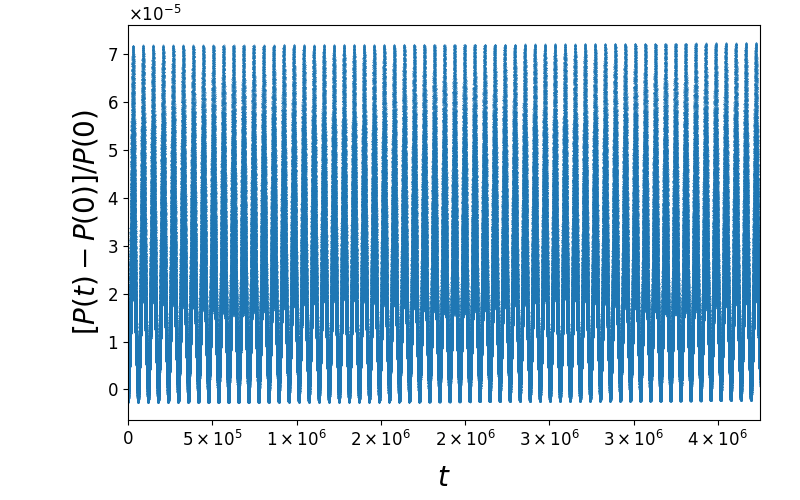}
		}
		
		\subfloat[SRK3]{\label{fig:guiding_centre_4d_barely_trapped_vprk_pnone_srk3}
			\includegraphics[height=2.5cm]{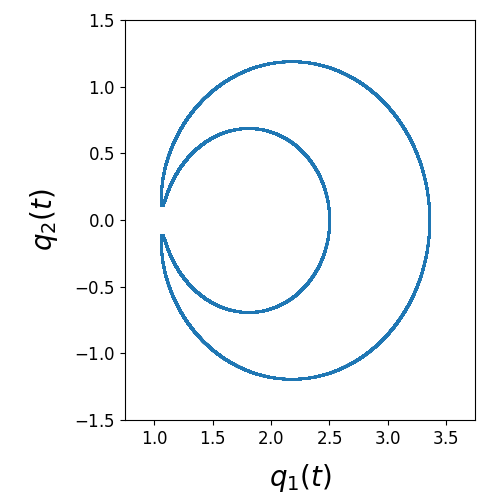}
			\includegraphics[height=2.5cm]{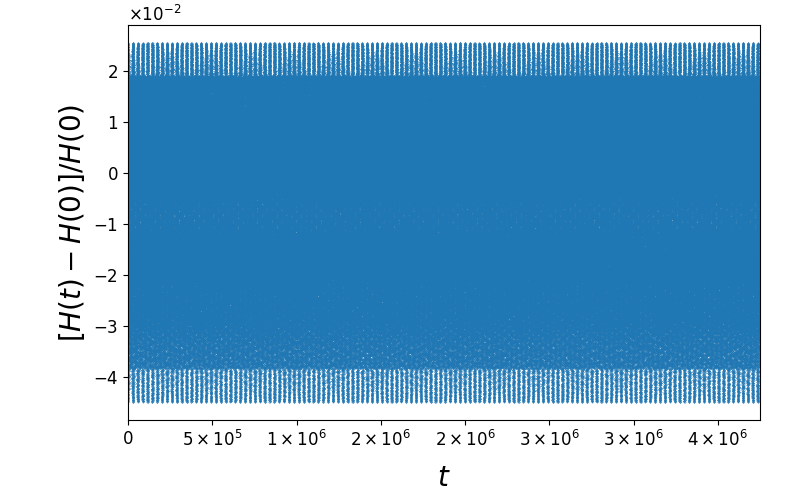}
			\includegraphics[height=2.5cm]{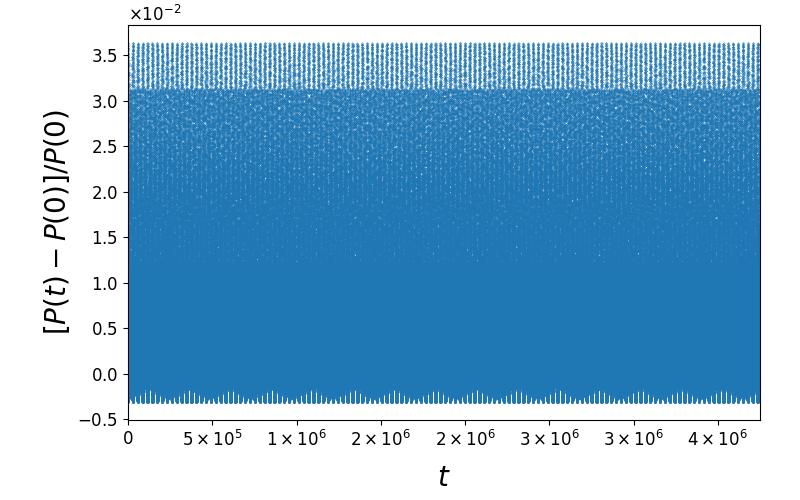}
		}
		
		\subfloat[RadIIA2]{\label{fig:guiding_center_4d_barely_trapped_vprk_radIIA2}
			\includegraphics[height=2.5cm]{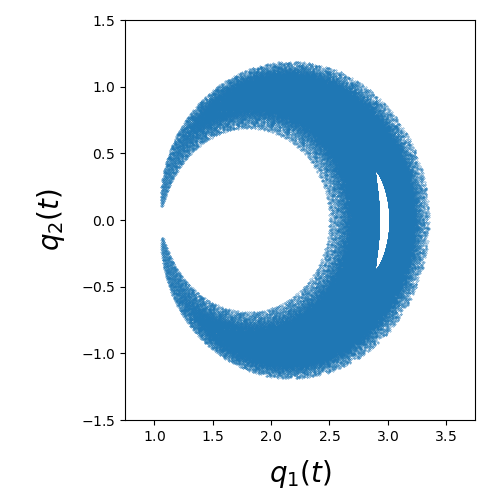}
			\includegraphics[height=2.5cm]{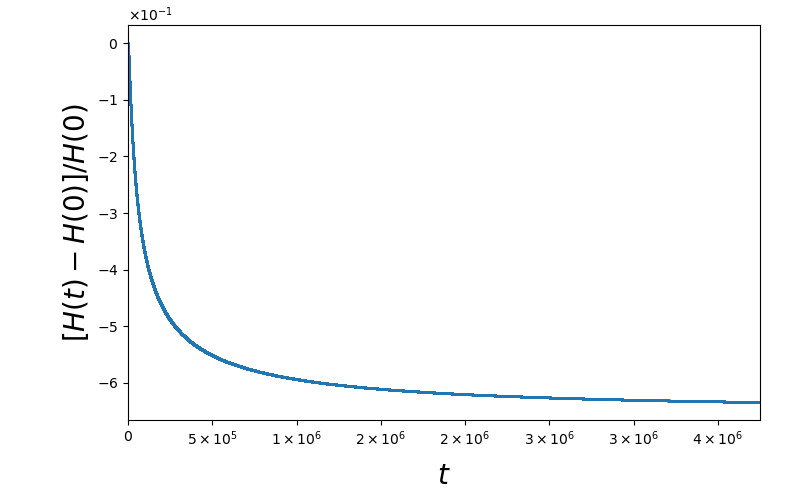}
			\includegraphics[height=2.5cm]{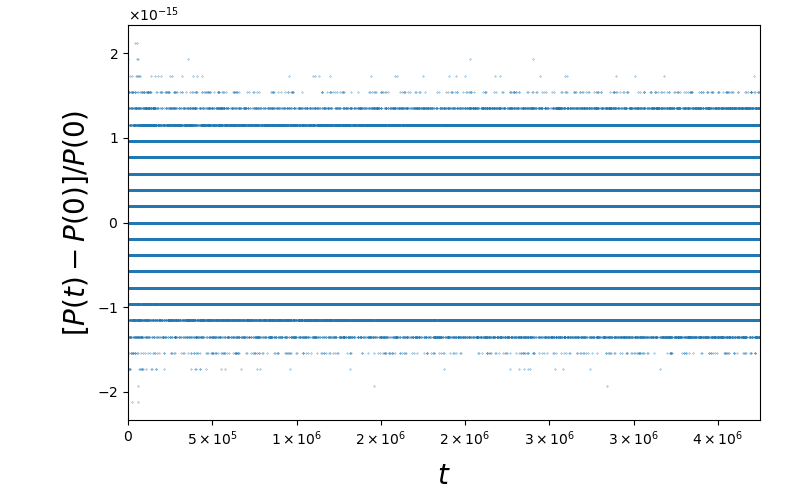}
		}
		
		\subfloat[RadIIA3]{\label{fig:guiding_center_4d_barely_trapped_vprk_radIIA3}
			\includegraphics[height=2.5cm]{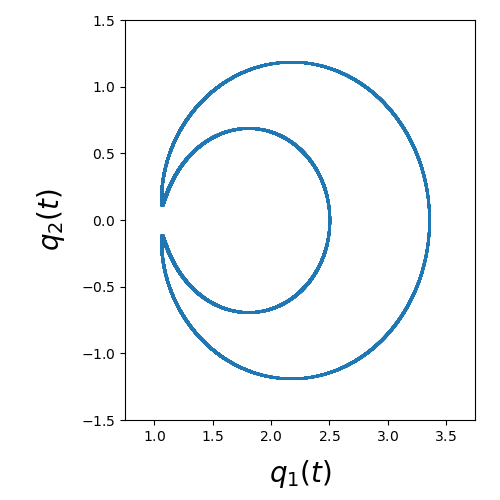}
			\includegraphics[height=2.5cm]{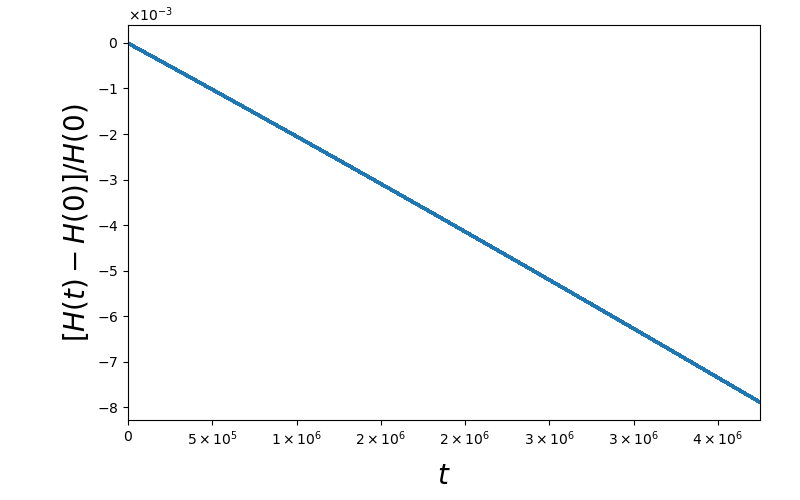}
			\includegraphics[height=2.5cm]{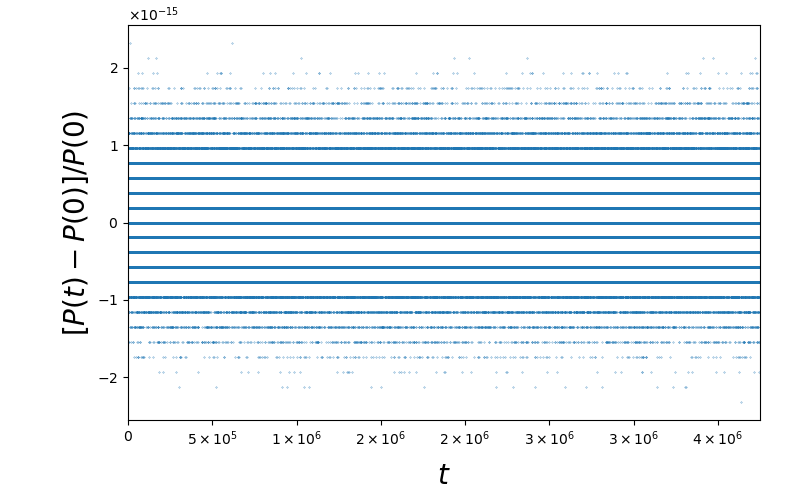}
		}
		
		\caption{Barely trapped guiding centre particle with variational Runge--Kutta methods.}
		\label{fig:guiding_centre_4d_barely_trapped_vprk_pnone}
	\end{center}
\end{figure}

\begin{figure}[p]
	\begin{center}
		\subfloat[GLRK3pStandard]{\label{fig:guiding_centre_4d_barely_passing_vprk_pstandard_glrk3}
			\includegraphics[height=2.5cm]{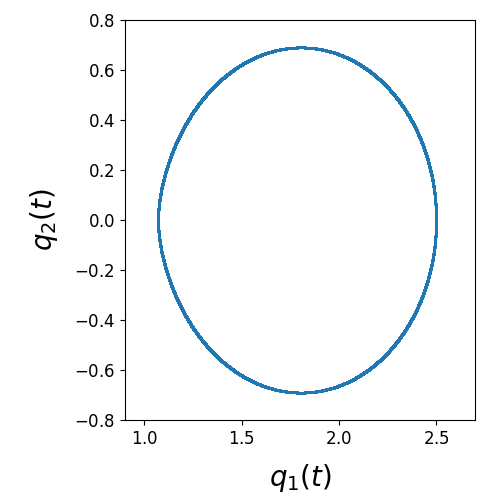}
			\includegraphics[height=2.5cm]{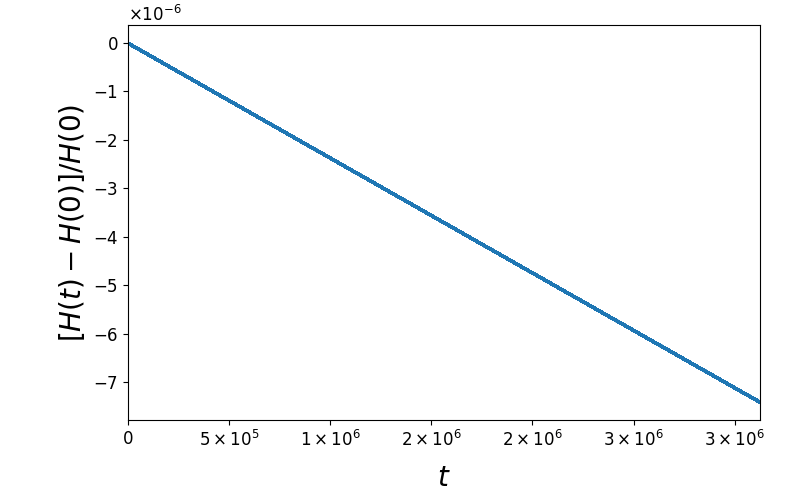}
		}
		\subfloat[GLRK3pSymmetric]{\label{fig:guiding_centre_4d_barely_passing_vprk_psymmetric_glrk3}
			\includegraphics[height=2.5cm]{results/guiding_centre_4d/barely_passing/guiding_center_4d_vprk_glrk3_psymmetric}
			\includegraphics[height=2.5cm]{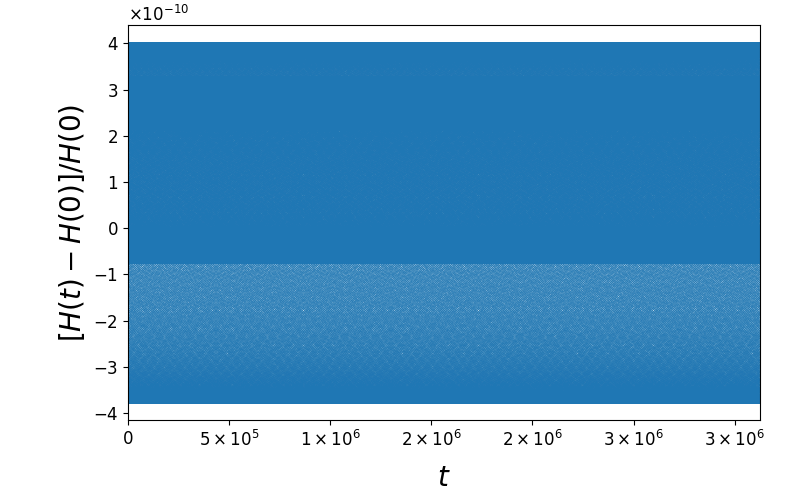}
		}
		
		\subfloat[GLRK4pStandard]{\label{fig:guiding_centre_4d_barely_passing_vprk_pstandard_glrk4}
			\includegraphics[height=2.5cm]{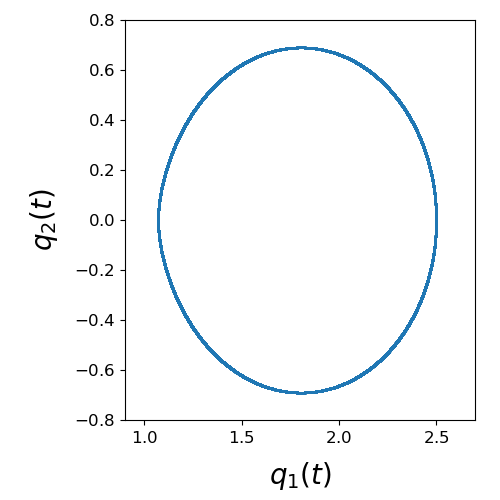}
			\includegraphics[height=2.5cm]{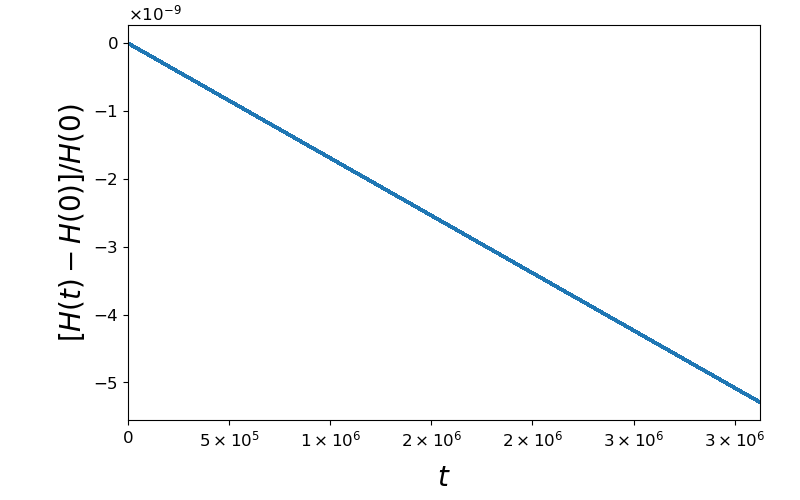}
		}
		\subfloat[GLRK4pSymmetric]{\label{fig:guiding_centre_4d_barely_passing_vprk_psymmetric_glrk4}
			\includegraphics[height=2.5cm]{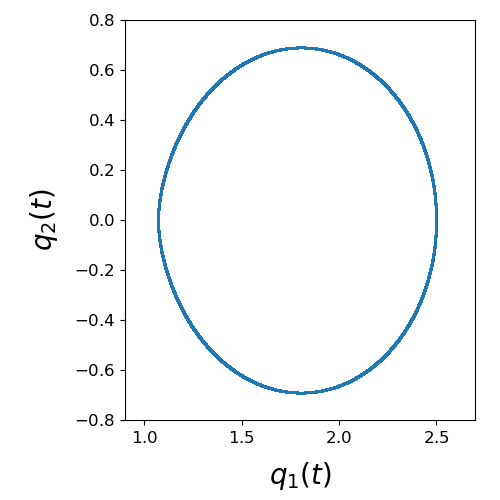}
			\includegraphics[height=2.5cm]{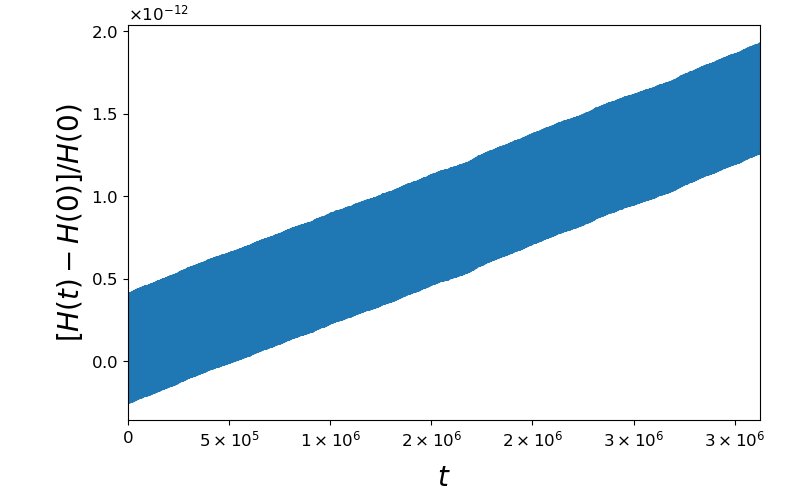}
		}
		
		\subfloat[SRK3pStandard]{\label{fig:guiding_centre_4d_barely_passing_vprk_pstandard_srk3}
			\includegraphics[height=2.5cm]{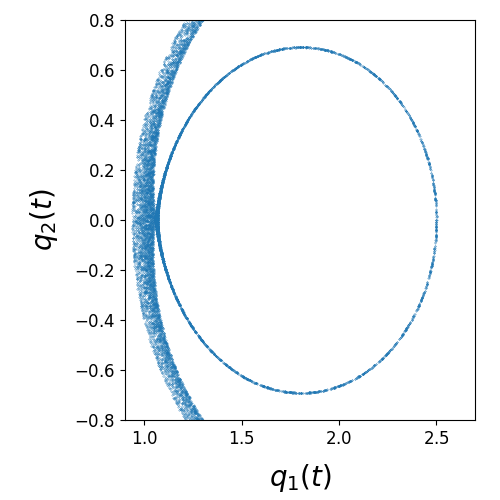}
			\includegraphics[height=2.5cm]{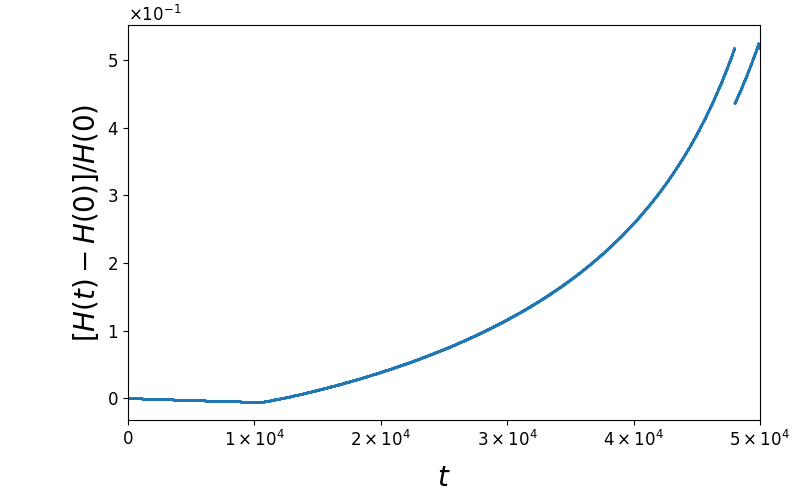}
		}
		\subfloat[SRK3pSymmetric]{\label{fig:guiding_centre_4d_barely_passing_vprk_psymmetric_srk3}
			\includegraphics[height=2.5cm]{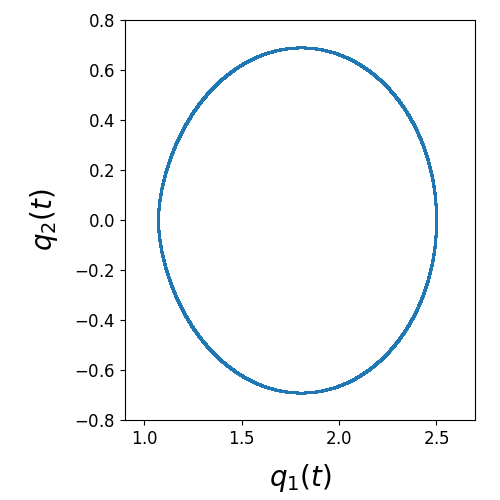}
			\includegraphics[height=2.5cm]{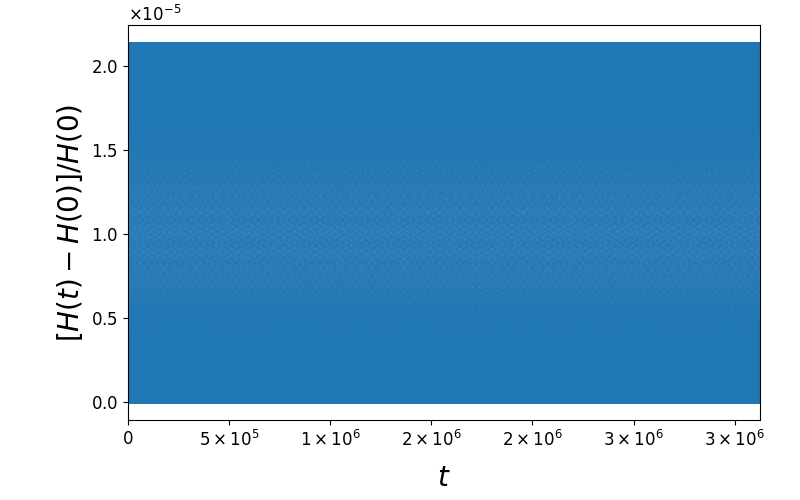}
		}
				
		\caption{Barely passing guiding centre particle with variational Runge--Kutta methods with standard and symmetric projection.}
		\label{fig:guiding_centre_4d_barely_passing_vprk_pstandard}
	\end{center}
\end{figure}

\begin{figure}[p]
	\begin{center}
		\subfloat[GLRK3pStandard]{\label{fig:guiding_centre_4d_barely_trapped_vprk_pstandard_glrk3}
			\includegraphics[height=2.5cm]{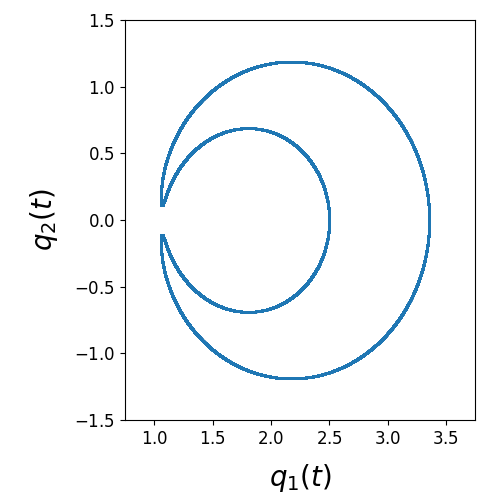}
			\includegraphics[height=2.5cm]{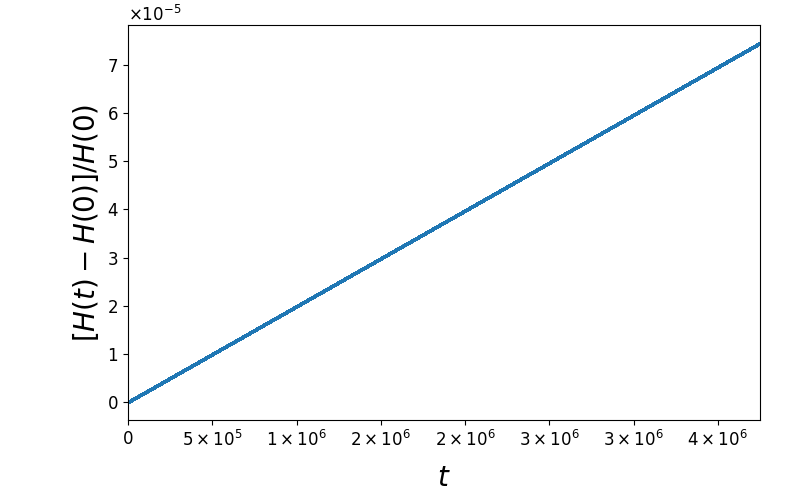}
		}
		\subfloat[GLRK3pSymmetric]{\label{fig:guiding_centre_4d_barely_trapped_vprk_psymmetric_glrk3}
			\includegraphics[height=2.5cm]{results/guiding_centre_4d/barely_trapped/guiding_center_4d_vprk_glrk3_psymmetric}
			\includegraphics[height=2.5cm]{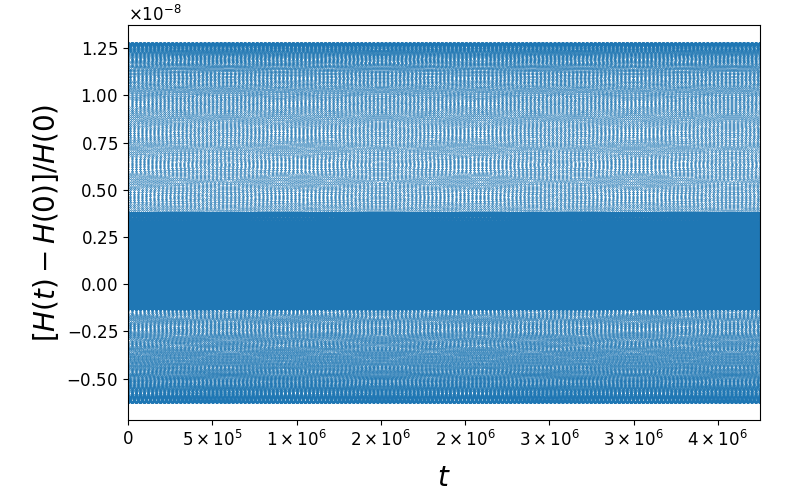}
		}
		
		\subfloat[GLRK4pStandard]{\label{fig:guiding_centre_4d_barely_trapped_vprk_pstandard_glrk4}
			\includegraphics[height=2.5cm]{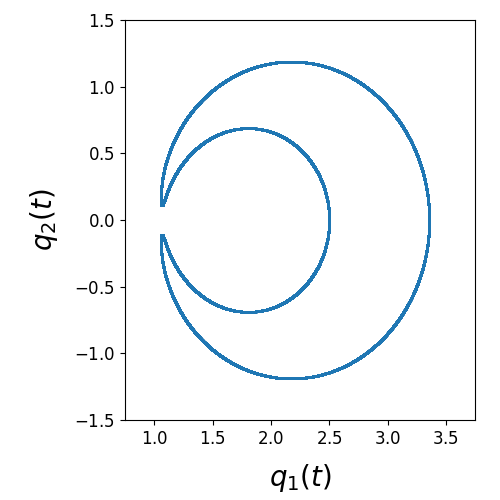}
			\includegraphics[height=2.5cm]{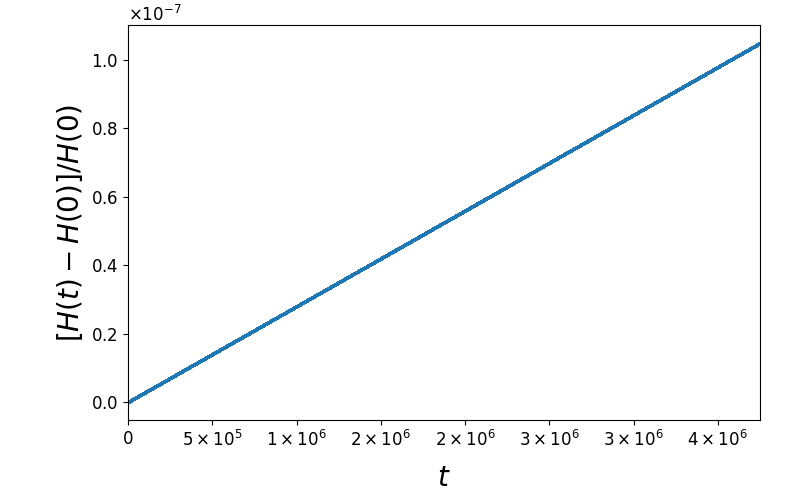}
		}
		\subfloat[GLRK4pSymmetric]{\label{fig:guiding_centre_4d_barely_trapped_vprk_psymmetric_glrk4}
			\includegraphics[height=2.5cm]{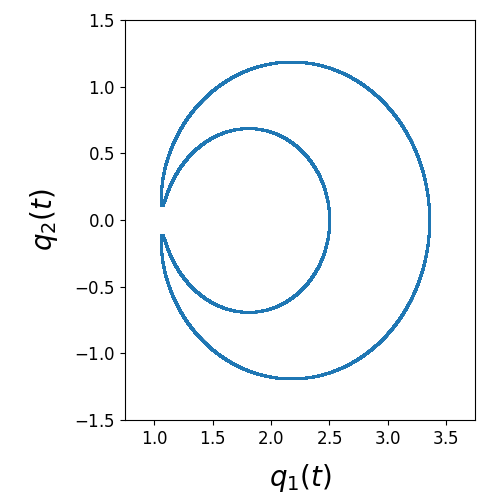}
			\includegraphics[height=2.5cm]{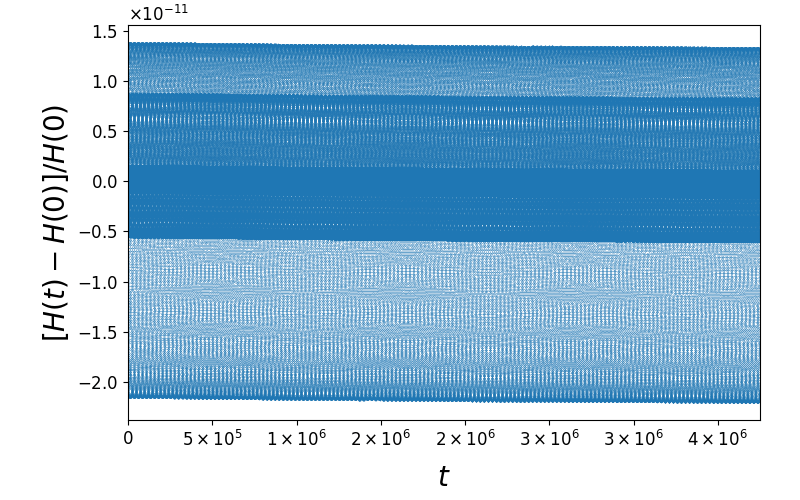}
		}
		
		\subfloat[SRK3pStandard]{\label{fig:guiding_centre_4d_barely_trapped_vprk_pstandard_srk3}
			\includegraphics[height=2.5cm]{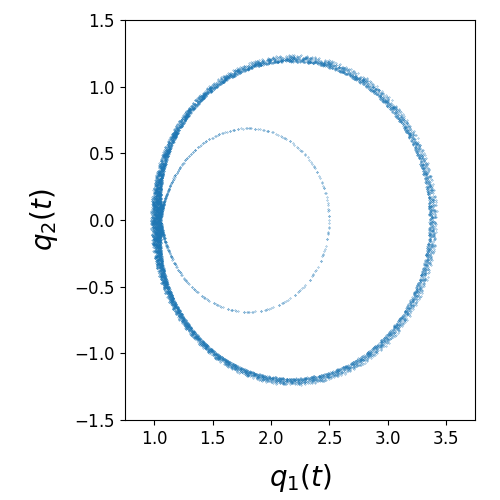}
			\includegraphics[height=2.5cm]{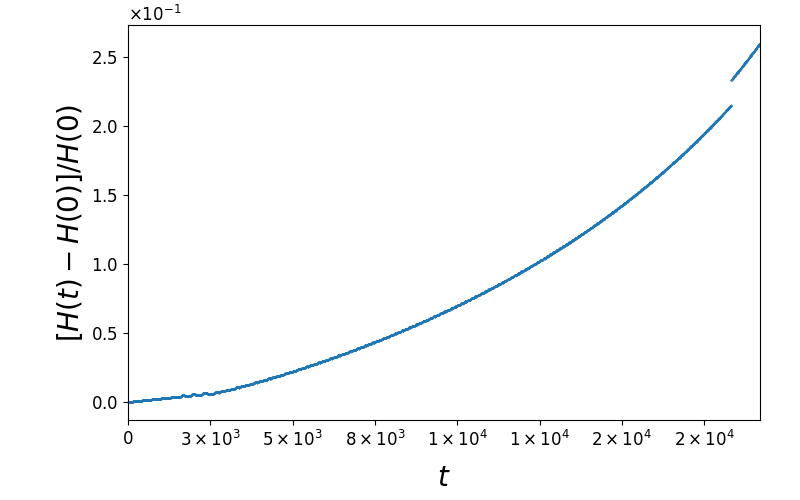}
		}
		\subfloat[SRK3pSymmetric]{\label{fig:guiding_centre_4d_barely_trapped_vprk_psymmetric_srk3}
			\includegraphics[height=2.5cm]{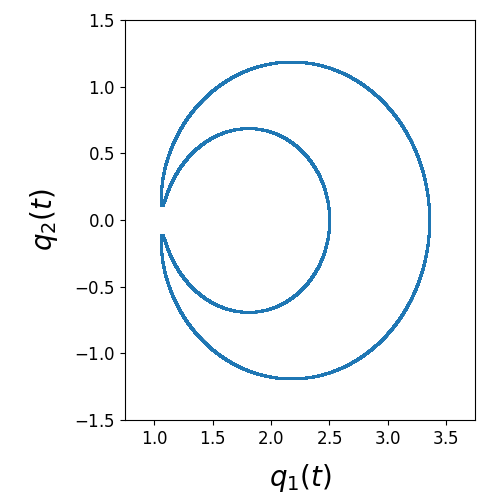}
			\includegraphics[height=2.5cm]{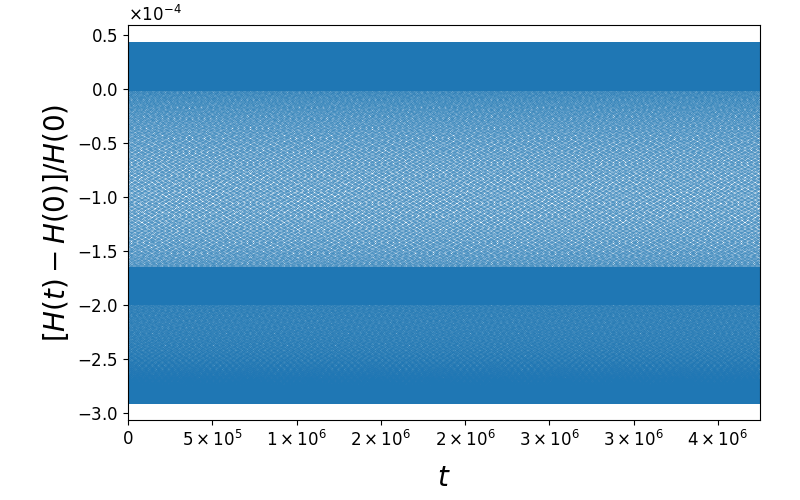}
		}
				
		\caption{Barely trapped guiding centre particle with variational Runge--Kutta methods with standard and symmetric projection.}
		\label{fig:guiding_centre_4d_barely_trapped_vprk_pstandard}
	\end{center}
\end{figure}

\begin{figure}[p]
	\begin{center}
		\subfloat[GLRK3pSymplectic]{\label{fig:guiding_centre_4d_barely_passing_vprk_psymplectic_glrk3}
			\includegraphics[height=2.5cm]{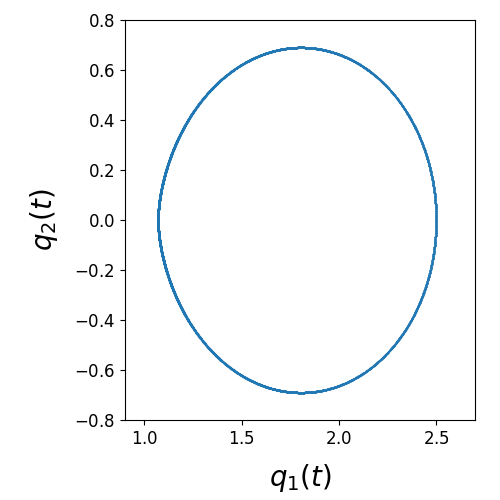}
			\includegraphics[height=2.5cm]{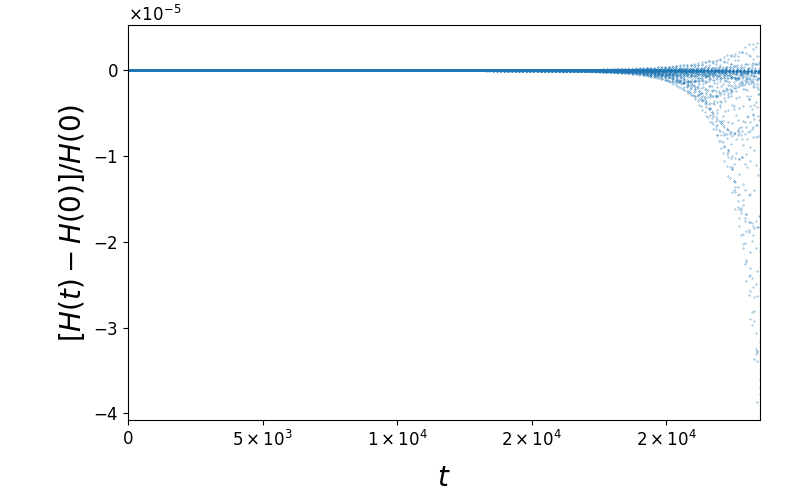}
		}
		\subfloat[GLRK3pMidpoint]{\label{fig:guiding_centre_4d_barely_passing_vprk_pmidpoint_glrk3}
			\includegraphics[height=2.5cm]{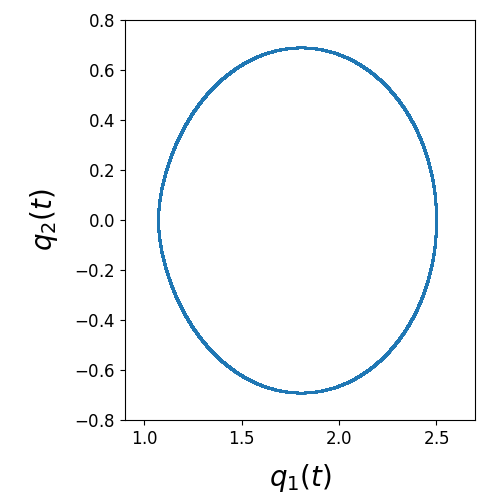}
			\includegraphics[height=2.5cm]{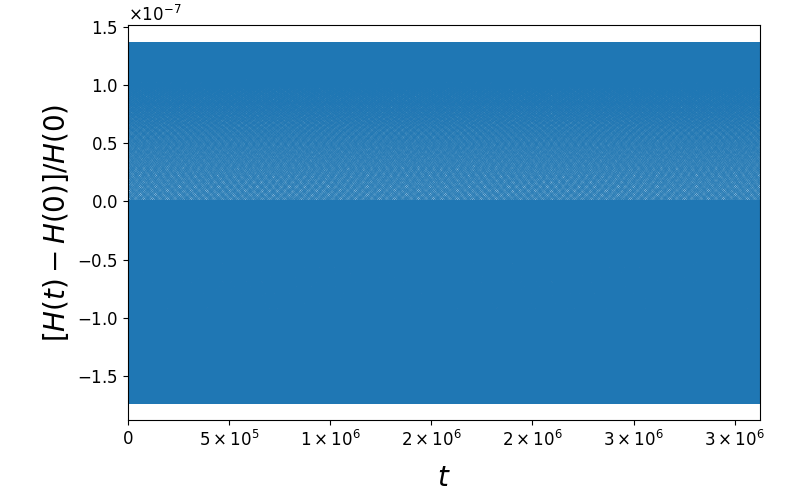}
		}
		
		\subfloat[GLRK4pSymplectic]{\label{fig:guiding_centre_4d_barely_passing_vprk_psymplectic_glrk4}
			\includegraphics[height=2.5cm]{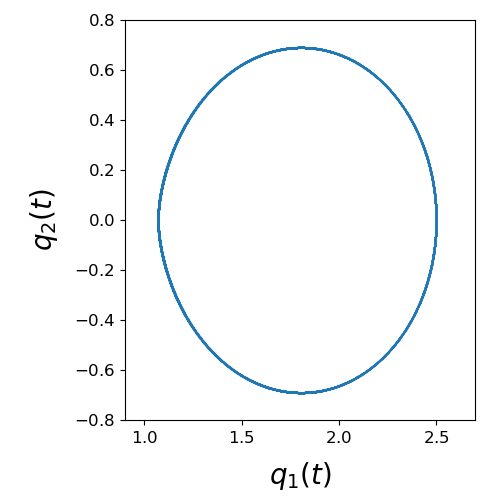}
			\includegraphics[height=2.5cm]{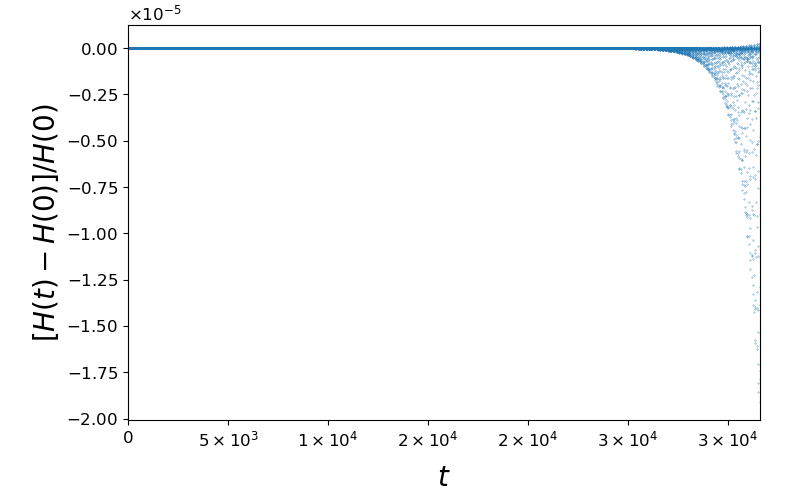}
		}
		\subfloat[GLRK4pMidpoint]{\label{fig:guiding_centre_4d_barely_passing_vprk_pmidpoint_glrk4}
			\includegraphics[height=2.5cm]{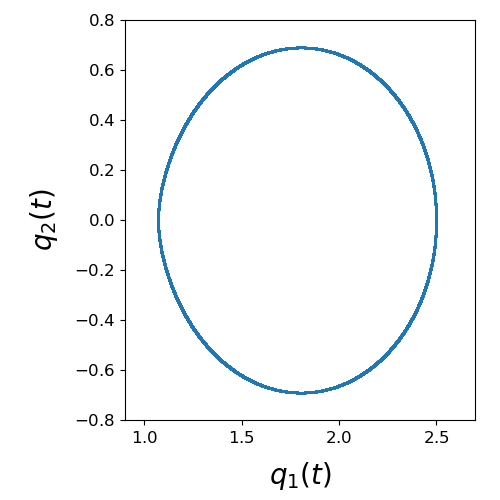}
			\includegraphics[height=2.5cm]{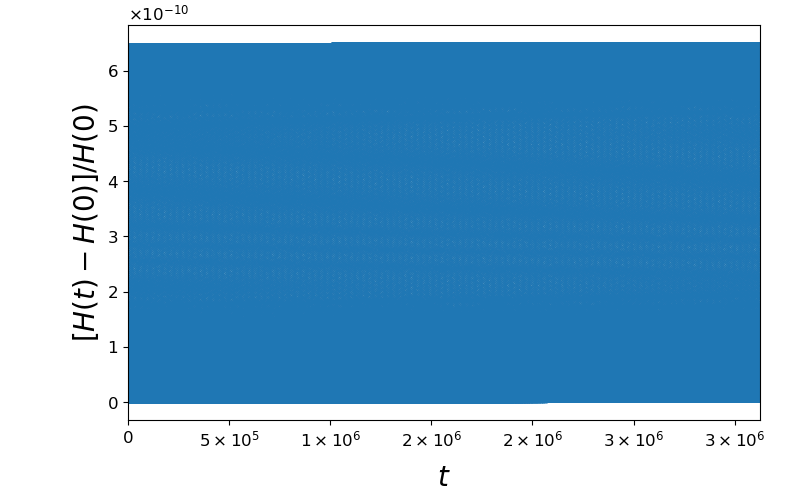}
		}
		
		\subfloat[SRK3pSymplectic]{\label{fig:guiding_centre_4d_barely_passing_vprk_psymplectic_srk3}
			\includegraphics[height=2.5cm]{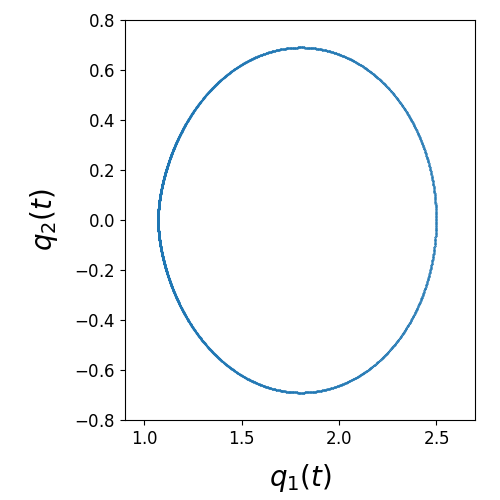}
			\includegraphics[height=2.5cm]{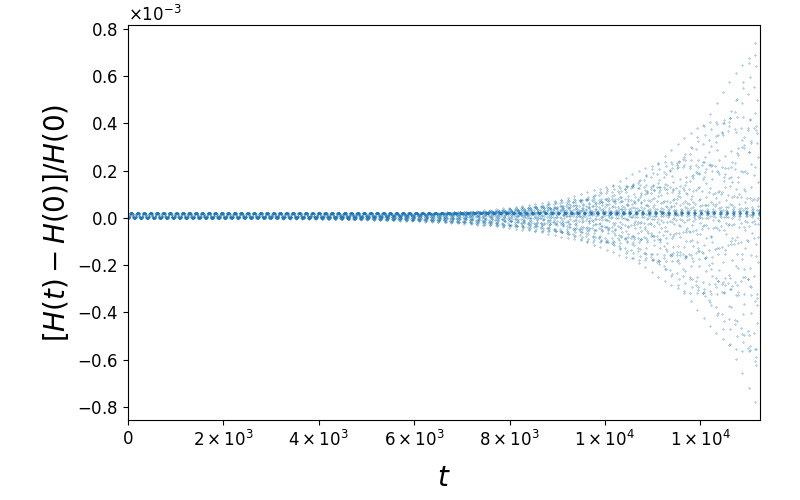}
		}
		\subfloat[SRK3pMidpoint]{\label{fig:guiding_centre_4d_barely_passing_vprk_pmidpoint_srk3}
			\includegraphics[height=2.5cm]{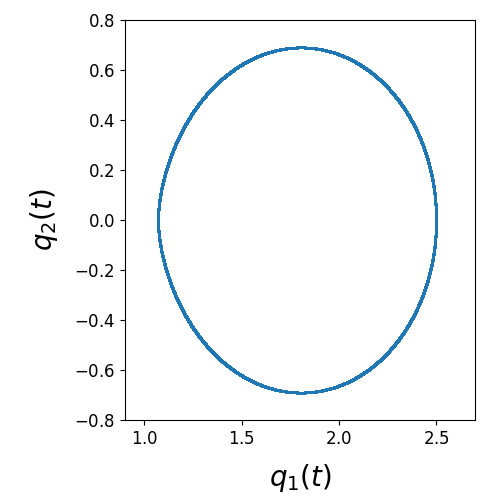}
			\includegraphics[height=2.5cm]{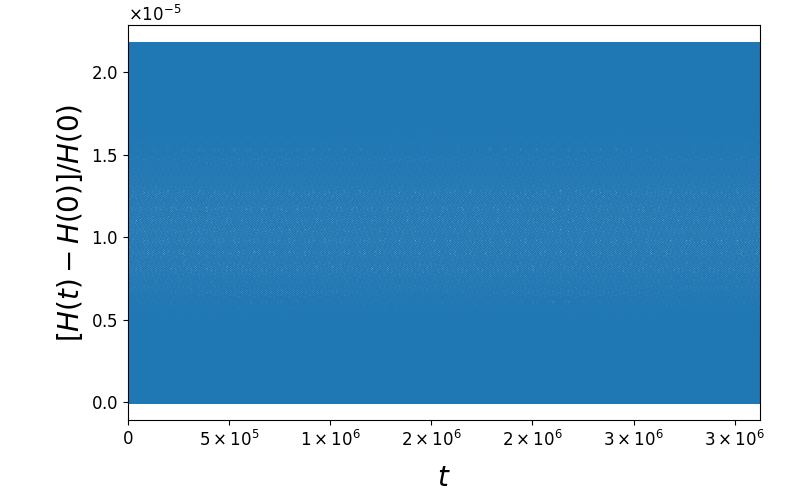}
		}
		
		\caption{Barely passing guiding centre particle with variational Runge--Kutta methods and symplectic and midpoint projection.}
		\label{fig:guiding_centre_4d_barely_passing_vprk_psymplectic}
	\end{center}
\end{figure}

\begin{figure}[p]
	\begin{center}
		\subfloat[GLRK3pSymplectic]{\label{fig:guiding_centre_4d_barely_trapped_vprk_psymplectic_glrk3}
			\includegraphics[height=2.5cm]{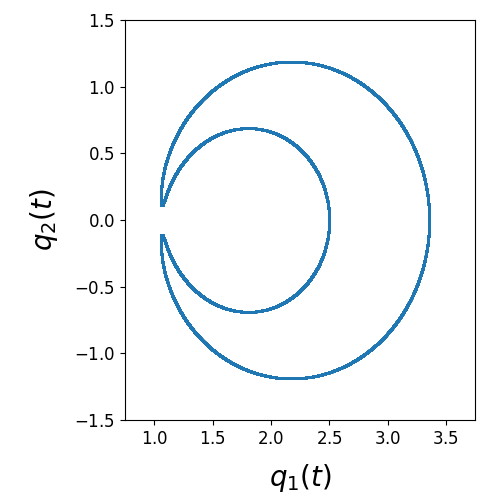}
			\includegraphics[height=2.5cm]{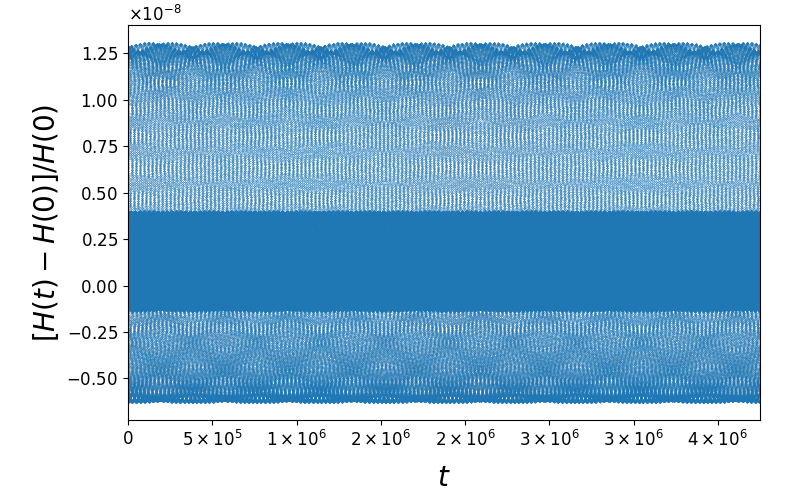}
		}
		\subfloat[GLRK3pMidpoint]{\label{fig:guiding_centre_4d_barely_trapped_vprk_pmidpoint_glrk3}
			\includegraphics[height=2.5cm]{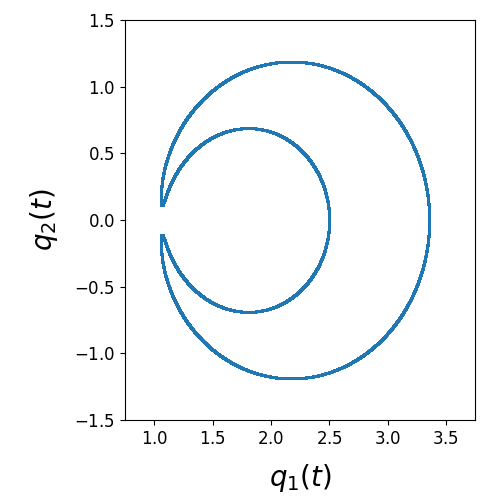}
			\includegraphics[height=2.5cm]{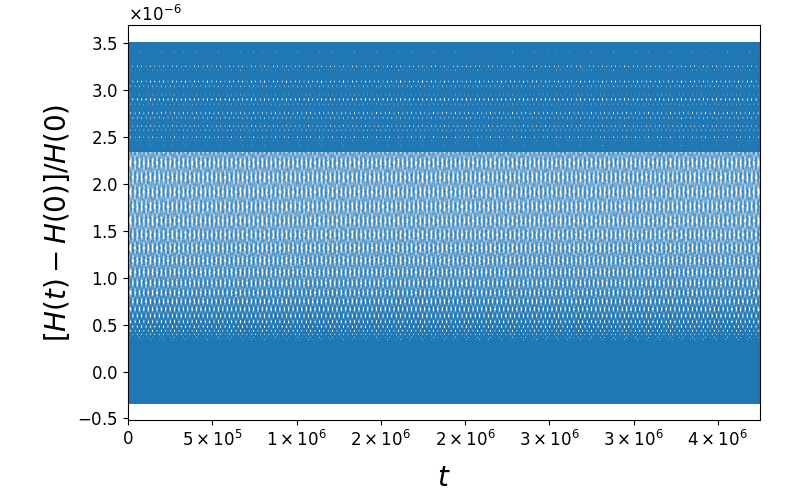}
		}
		
		\subfloat[GLRK4pSymplectic]{\label{fig:guiding_centre_4d_barely_trapped_vprk_psymplectic_glrk4}
			\includegraphics[height=2.5cm]{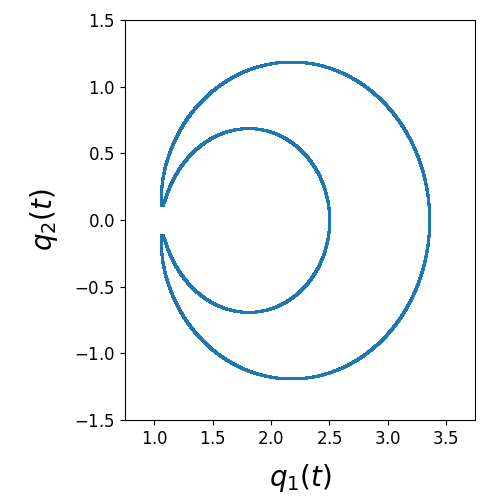}
			\includegraphics[height=2.5cm]{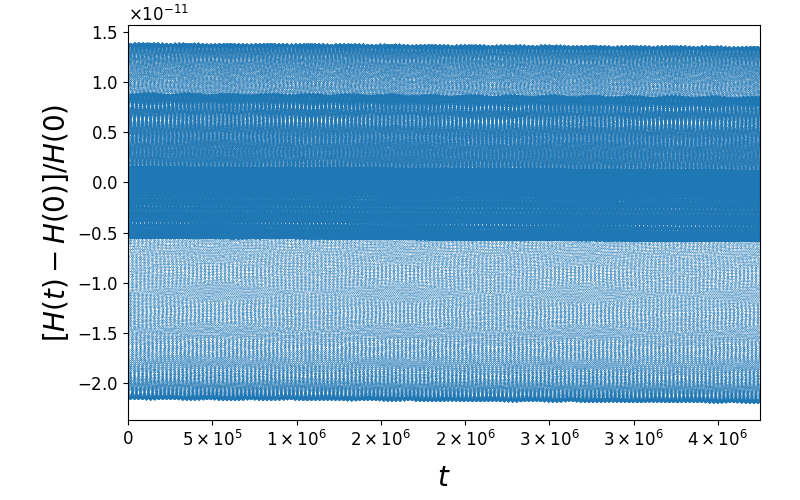}
		}
		\subfloat[GLRK4pMidpoint]{\label{fig:guiding_centre_4d_barely_trapped_vprk_pmidpoint_glrk4}
			\includegraphics[height=2.5cm]{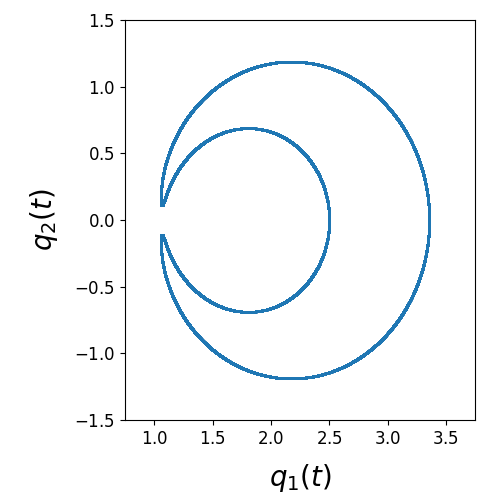}
			\includegraphics[height=2.5cm]{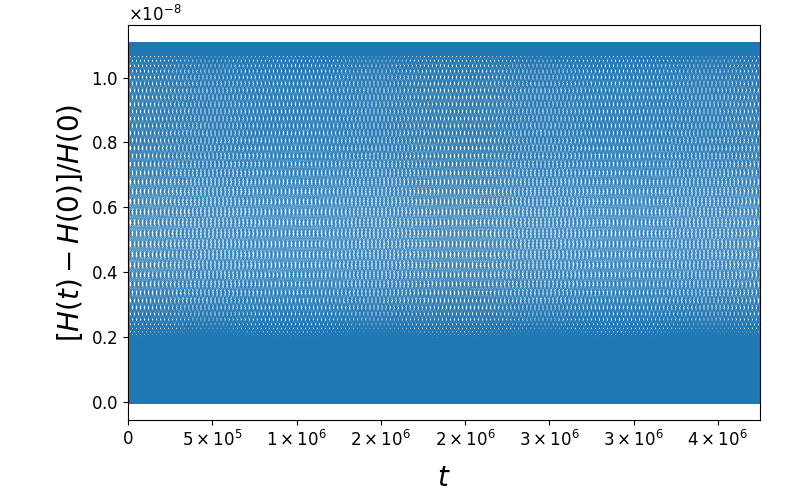}
		}
		
		\subfloat[SRK3pSymplectic]{\label{fig:guiding_centre_4d_barely_trapped_vprk_psymplectic_srk3}
			\includegraphics[height=2.5cm]{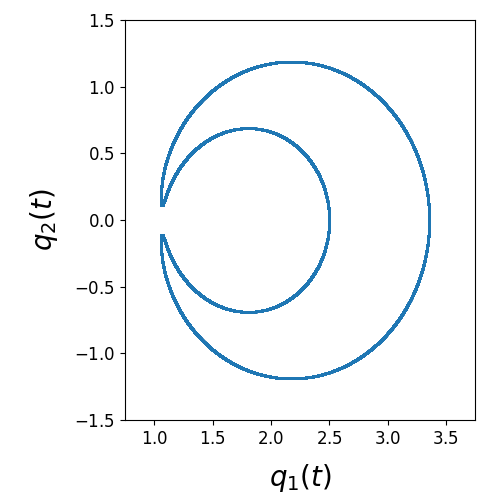}
			\includegraphics[height=2.5cm]{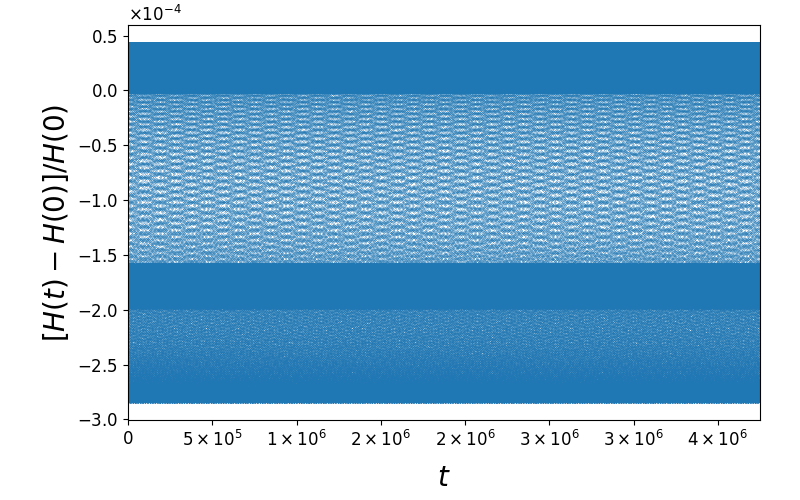}
		}
		\subfloat[SRK3pMidpoint]{\label{fig:guiding_centre_4d_barely_trapped_vprk_pmidpoint_srk3}
			\includegraphics[height=2.5cm]{results/guiding_centre_4d/barely_trapped/guiding_center_4d_vprk_srk3_psymplectic}
			\includegraphics[height=2.5cm]{results/guiding_centre_4d/barely_trapped/guiding_center_4d_vprk_srk3_psymplectic_energy_error}
		}
		
		\caption{Barely trapped guiding centre particle with variational Runge--Kutta methods and symplectic and midpoint projection.}
		\label{fig:guiding_centre_4d_barely_trapped_vprk_psymplectic}
	\end{center}
\end{figure}

We make the following observations:
\begin{itemize}
\item The passing particles (Figure~\ref{fig:guiding_centre_4d_barely_passing_vprk_pnone}) seem to be the more challenging ones of the considered examples. Without projection almost no method is long-time stable and we observe drifts in the energy as well as in the toroidal momentum. The Gauss--Lobatto--IIIA, IIIB and IIIC methods are particularly unstable, crashing after only a few or at best a few hundred time steps. This can be remedied by a smaller time step. Then, however, the computational effort is much larger than with the other methods.
For the Gauss--Legendre Runge--Kutta methods we observe that for the integrators with an odd number of stages, the solution has two branches, while for the integrators with an even number of stages, the solution has only one branch.

\item For the trapped particles (Figure~\ref{fig:guiding_centre_4d_barely_trapped_vprk_pnone}), we obtain good results without projection for the Gauss--Legendre methods and the SRK3 method except for the six-stage Gauss--Legendre method. For the barely trapped particle, this is also true for the two-stage Gauss--Legendre method. The Gauss--Lobatto--IIIA, IIIB and IIIC methods crash quickly, that is after only about one or a few thousand time steps. The Gauss--Lobatto--IIID and IIIE methods work mostly well for the deeply trapped particle but not for the barely trapped particle.

\item The standard projection seems to improve the situation in case of the passing particles (left column of Figure~\ref{fig:guiding_centre_4d_barely_passing_vprk_pstandard}) and worsen the situation in case of the trapped particles (left column of Figure~\ref{fig:guiding_centre_4d_barely_trapped_vprk_pstandard}). In all cases, though, the results are not satisfactory as the solution always exhibits a drift in the energy error.

\item The symplectic projection (left column of Figures~\ref{fig:guiding_centre_4d_barely_passing_vprk_psymplectic} and~\ref{fig:guiding_centre_4d_barely_trapped_vprk_psymplectic}) leads to good results for the trapped particles but not for the passing particles. As before, we observe that in those cases where the unprojected solution is stable, the symplectic projection is also stable, but in those cases where the unprojected solution is unstable, the symplectic projection is also unstable.

\item The symmetric (right column of Figures~\ref{fig:guiding_centre_4d_barely_passing_vprk_pstandard} and~\ref{fig:guiding_centre_4d_barely_trapped_vprk_pstandard}) and midpoint (right column of Figures~\ref{fig:guiding_centre_4d_barely_passing_vprk_psymplectic} and~\ref{fig:guiding_centre_4d_barely_trapped_vprk_psymplectic}) projections lead to good results with almost all methods except for the Gauss--Lobatto--IIIC methods. For the barely trapped particle both projections are unstable for almost all Gauss--Lobatto methods.
For the symmetric projection and the higher-order Gauss--Legendre method we observe a drift in the energy error, which again is only of the order of $10^{-12}$ for 1~250~000 time steps. This drift is not due to the non-symplecticity of the symmetric projection method, as one could expect, but it is caused by round-off errors and disappears in simulations in quadruple precision (c.f., Section~\ref{sec:app_energy_drift}). 

\item With the Radau methods we observe exact conservation of the toroidal momentum, which constitutes one of the components of the Dirac constraint and is thus expected to be preserved, but we also observe dissipation of energy. For the two-stage Radau method these are related to large errors in the solution (Figures~\ref{fig:guiding_center_4d_barely_passing_vprk_radIIA2} and~\ref{fig:guiding_center_4d_barely_trapped_vprk_radIIA2}). For the three-stage Radau method (Figures~\ref{fig:guiding_center_4d_barely_passing_vprk_radIIA3} and~\ref{fig:guiding_center_4d_barely_trapped_vprk_radIIA3}) the errors in the solution are less pronounced.
\end{itemize}

The numerical experiments suggest that the symmetric and midpoint projection lead to good results with all integrators, while the symplectic projection is only stable for those integrators and those examples which are stable even without projection.
Similar to the point vortex example, the solution with the standard projection exhibits a drift in the energy and correspondingly a degradation of the numerical solution, as is expected.

\subsection{Energy Drift}
\label{sec:app_energy_drift}

\begin{figure}[hbt]
	\begin{center}
		\subfloat[GLRK3pSymmetric]{
			\includegraphics[width=.28\textwidth]{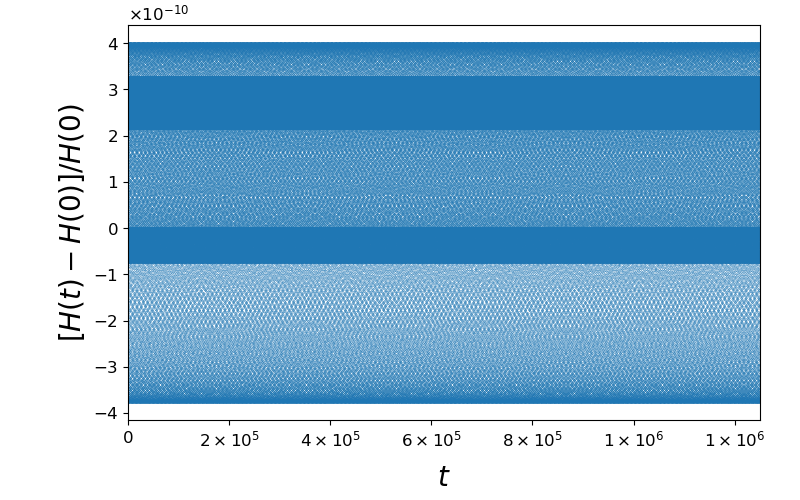}
			\includegraphics[width=.28\textwidth]{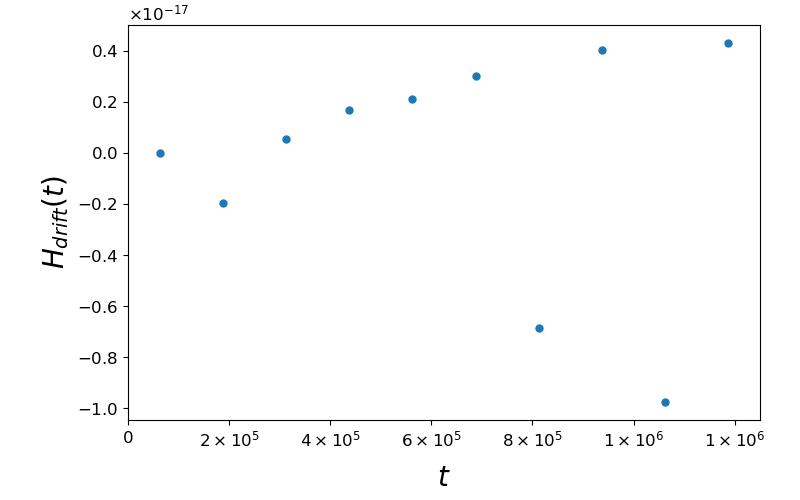}
			\includegraphics[width=.28\textwidth]{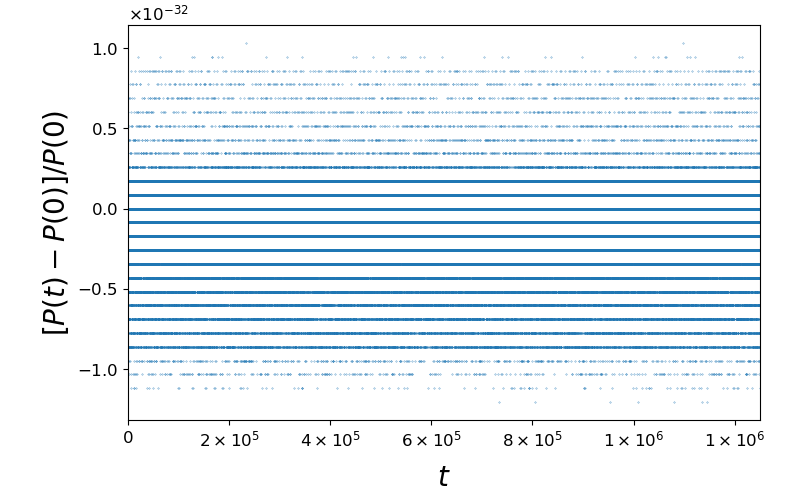}
		}
		
		\subfloat[GLRK4pSymmetric]{
			\includegraphics[width=.28\textwidth]{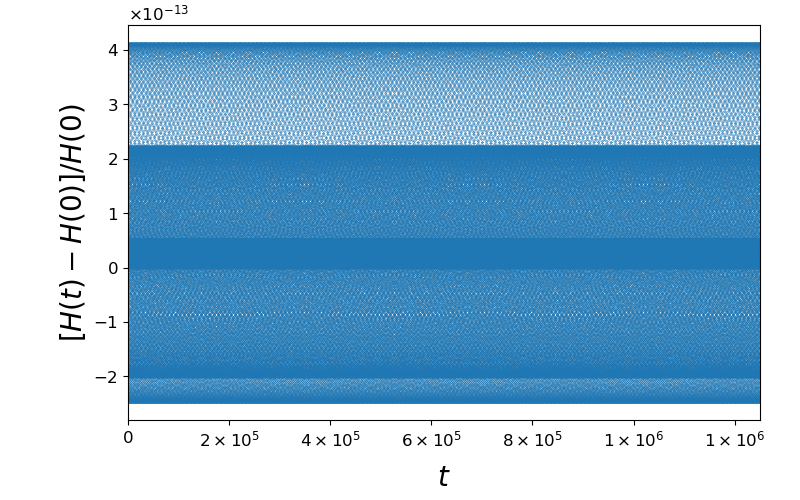}
			\includegraphics[width=.28\textwidth]{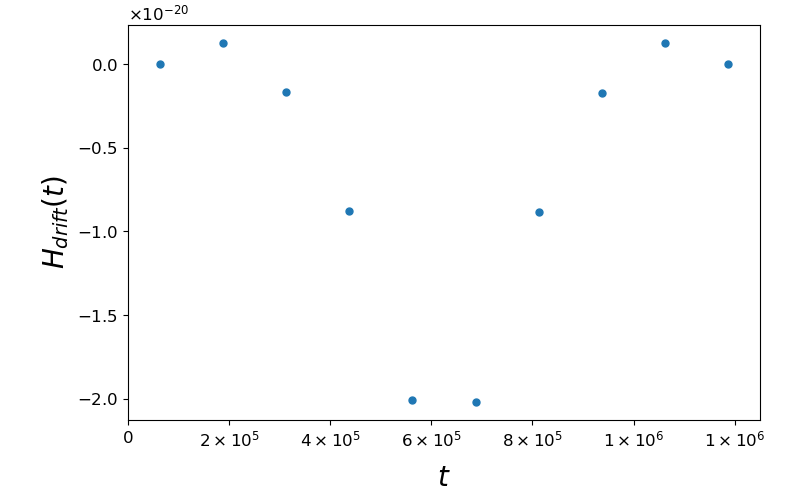}
			\includegraphics[width=.28\textwidth]{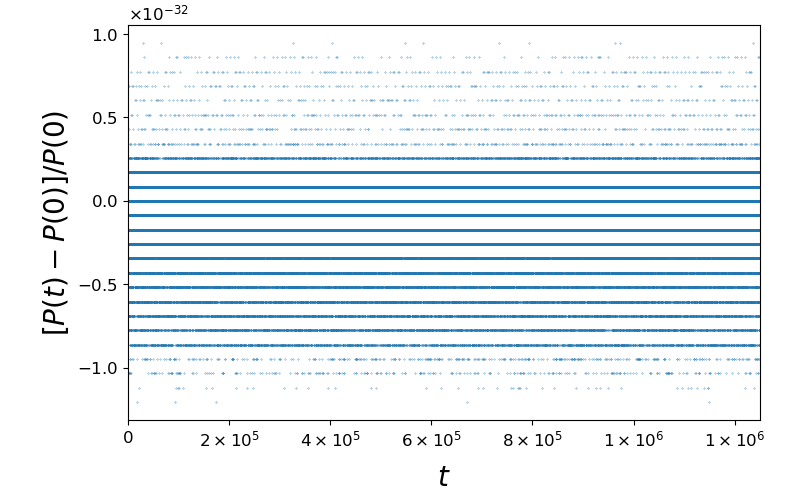}
		}

		\subfloat[GLRK5pSymmetric]{
			\includegraphics[width=.28\textwidth]{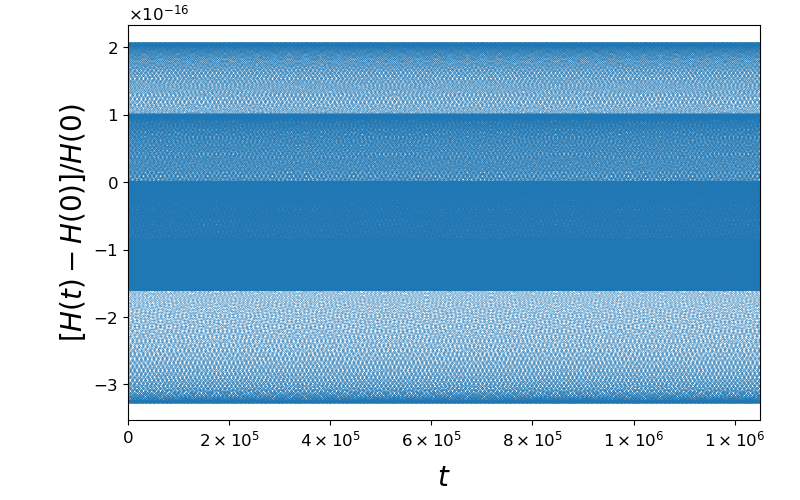}
			\includegraphics[width=.28\textwidth]{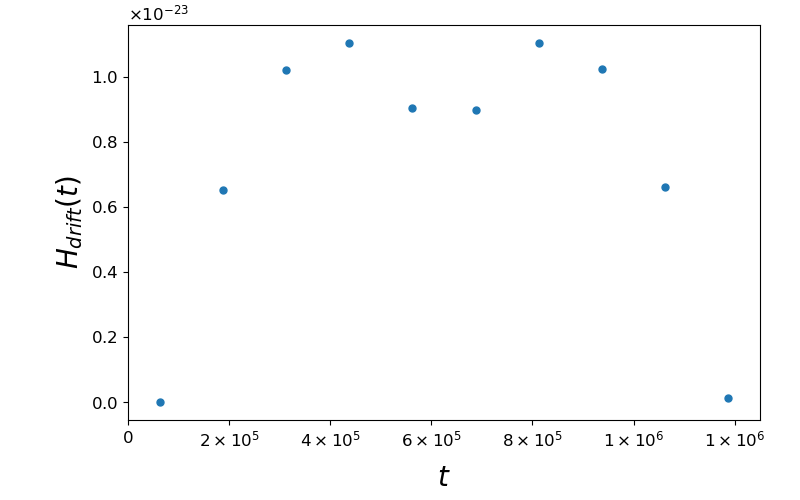}
			\includegraphics[width=.28\textwidth]{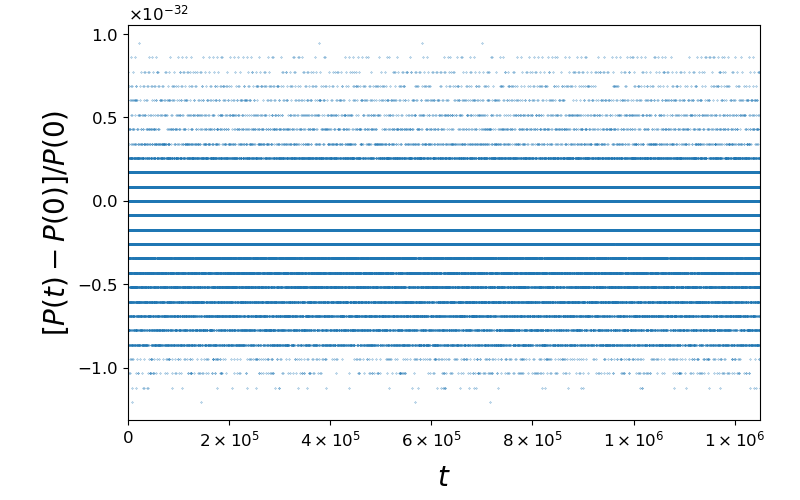}
		}

		\subfloat[GLRK6pSymmetric]{
			\includegraphics[width=.28\textwidth]{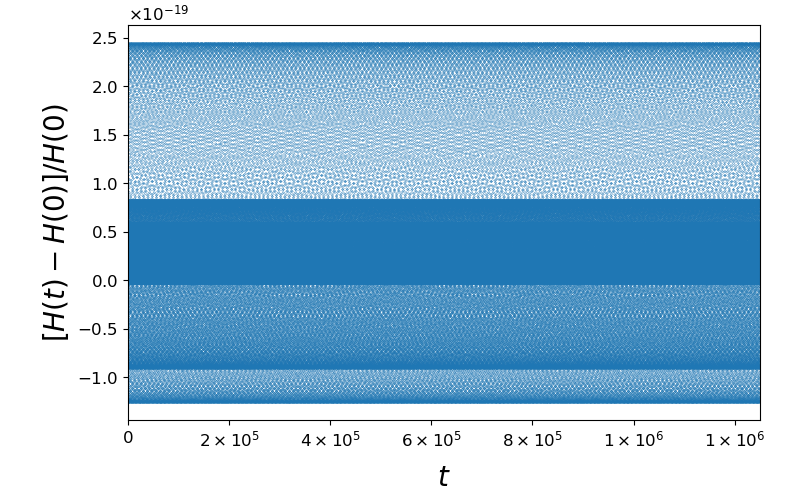}
			\includegraphics[width=.28\textwidth]{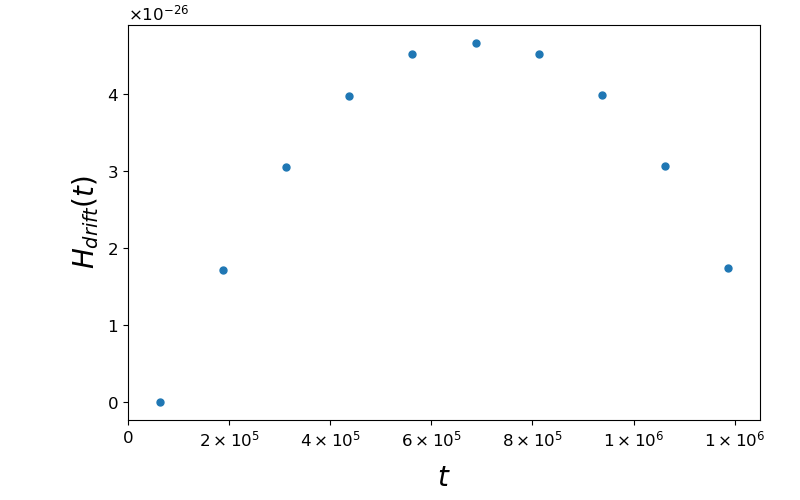}
			\includegraphics[width=.28\textwidth]{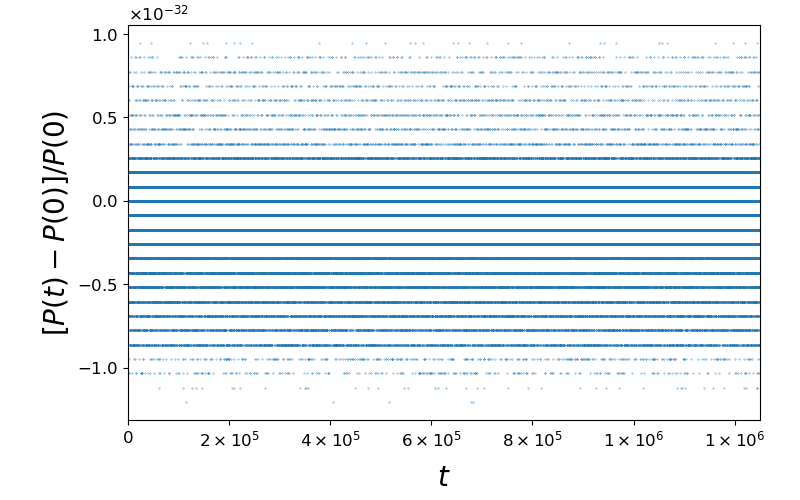}
		}

		\caption{Energy error, drift of the energy error, and toroidal momentum error for the barely passing guiding centre particle with higher-order Gauss--Legendre Runge--Kutta integrators using quadruple precision.}
		\label{fig:guiding_centre_quadruple}
	\end{center}
\end{figure}

For several examples we observed a drift of the energy error with the higher-order variational integrators, e.g., for Gauss--Legendre Runge--Kutta methods with four or more stages.
This behaviour was particularly prominent with the symmetric and midpoint projection when the energy error approached the machine accuracy.
A natural question to ask is if the origin for this phenomenon lies in the non-symplecticity of the two projection methods or if it is just due to round-off errors.
To this end, we repeated some simulations, namely those for a barely passing guiding centre particle, in quadruple precision. The resulting energy error is plotted in Figure~\ref{fig:guiding_centre_quadruple}. We also show the drift of the energy error, which is obtained by splitting the time interval of the simulation in 10 sub-intervals and computing the maximum of the absolute value of the energy error in each interval. As can be seen, there is no trend like linear growth, etc., visible, indicating that all errors due to either the non-symplecticity of the projection methods or due to round-off are much smaller than the energy error in these simulations.

\citet{Hairer:2008} attribute this drift behaviour to inaccuracies in the Runge-Kutta coefficients and weights, leading to an only approximate satisfaction of the symplecticity conditions~\eqref{eq:vprk_coefficients}. This leads to a linear growth of the energy error with time, even though one would expect a growth with the square root of the time as follows from random walk arguments (\emph{Brouwer's law}~\cite{Brouwer:1937}).
In order to reduce the influence of round-off errors for Runge--Kutta methods of high precision, \citet{Hairer:2008} suggest to apply a method inspired by compensated summation, where the Runge-Kutta coefficients and weights are split into two parts,
\begin{align}
b_{i} &= b_{i}^{*} + \otilde{b}_{i} , &
a_{ij} &= a_{ij}^{*} + \otilde{a}_{ij} .
\end{align}
Here, $b_{i}^{*}$ and $a_{ij}^{*}$ are approximations of $b_{i}$ and $a_{ij}$, which are exact to machine precision, and $\otilde{b}_{i}$ and $\otilde{a}_{ij}$ are corrections e.g. to the values of $b_{i}$ and $a_{ij}$ in quadruple precision.
With these coefficients and weights, the definition of the internal stages and the update rule e.g. in~\eqref{eq:vprk} are modified to
\begin{subequations}
\begin{align}
Q_{n,i} &= q_{n} + h \sum \limits_{j=1}^{s} a_{ij}^{*} \, \dot{Q}_{n,j} + h \sum \limits_{j=1}^{s} \otilde{a}_{ij} \, \dot{Q}_{n,j} , \\
q_{n+1} &= q_{n} + h \sum \limits_{i=1}^{s} b_{i}^{*} \dot{Q}_{n,i} + h \sum \limits_{i=1}^{s} \otilde{b}_{i}  \, \dot{Q}_{n,i} , 
\end{align}
\end{subequations}
which comes at practically no cost but allows to recover the missing accuracy in the last few digits of the solution.

Another technique, recently proposed by \citet{Antonana:2017}, consists in choosing the coefficients $\bar{a}_{ij}$ such that the symplecticity condition~\eqref{eq:vprk_coefficients} holds exactly in floating point precision. This means that even for methods like Gauss--Legendre, where $a$ and $\bar{a}$ are usually the same, slightly different coefficients are used.
If this technique is applied, e.g., with Gauss--Legendre Runge--Kutta methods with five and six stages and symmetric projection, the acquired solutions exhibit no drift in the energy error (not shown).

\subsection{Convergence}
\label{sec:app_convergence}

\begin{figure}[p]
	\begin{center}
		\subfloat[No Projection]{
			\includegraphics[width=.3\textwidth]{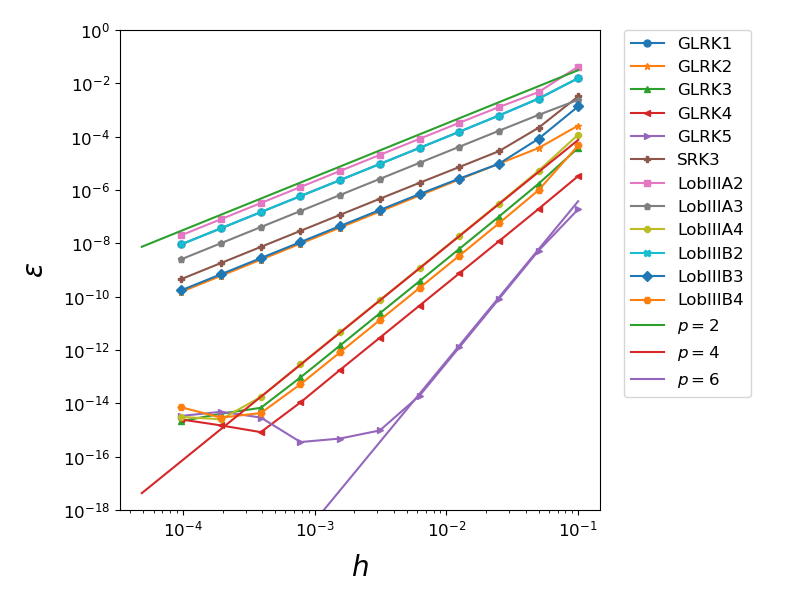}
			\includegraphics[width=.3\textwidth]{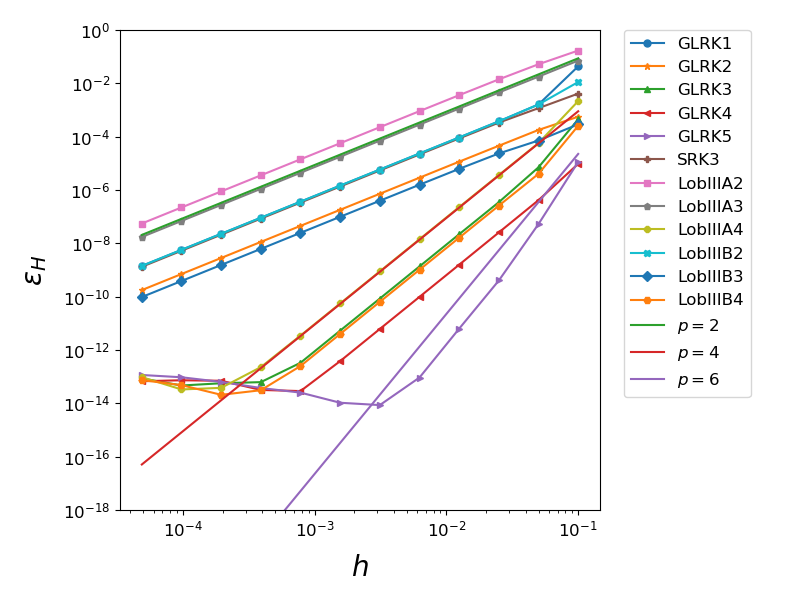}
			\includegraphics[width=.3\textwidth]{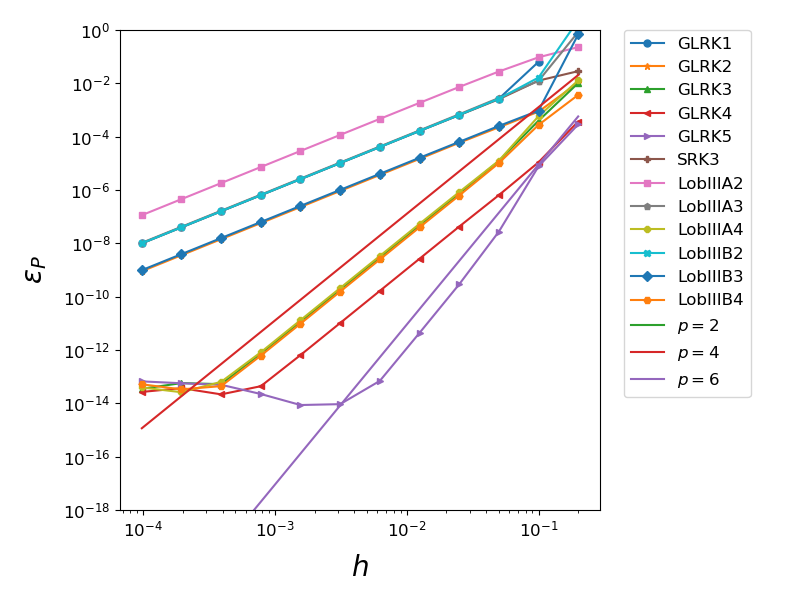}
		}

		\subfloat[Standard Projection]{
			\includegraphics[width=.3\textwidth]{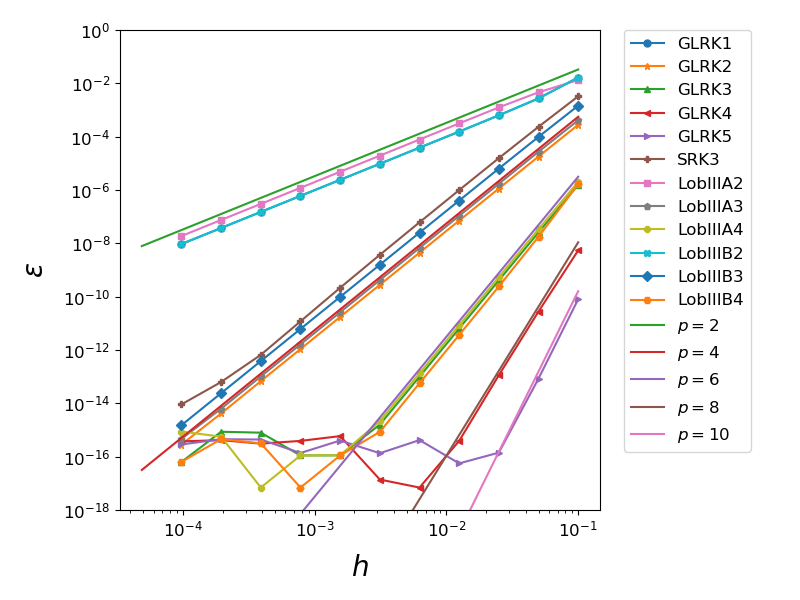}
			\includegraphics[width=.3\textwidth]{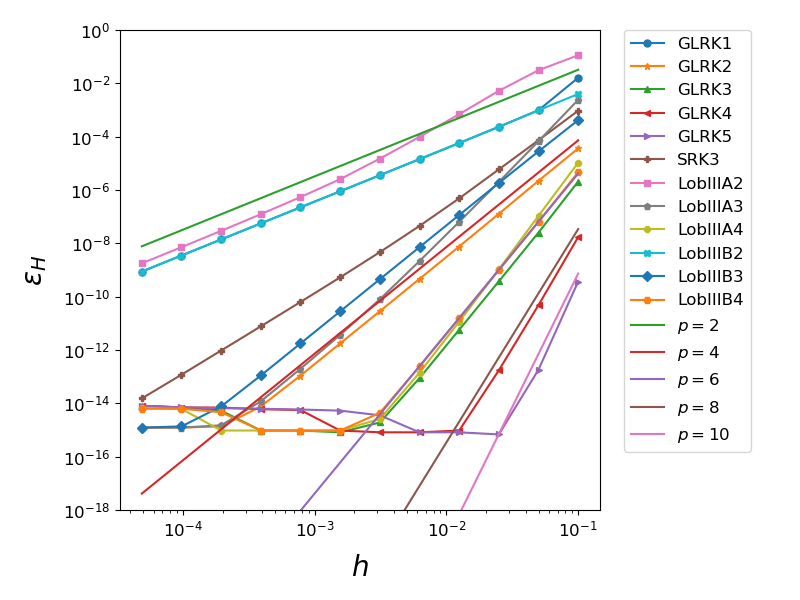}
			\includegraphics[width=.3\textwidth]{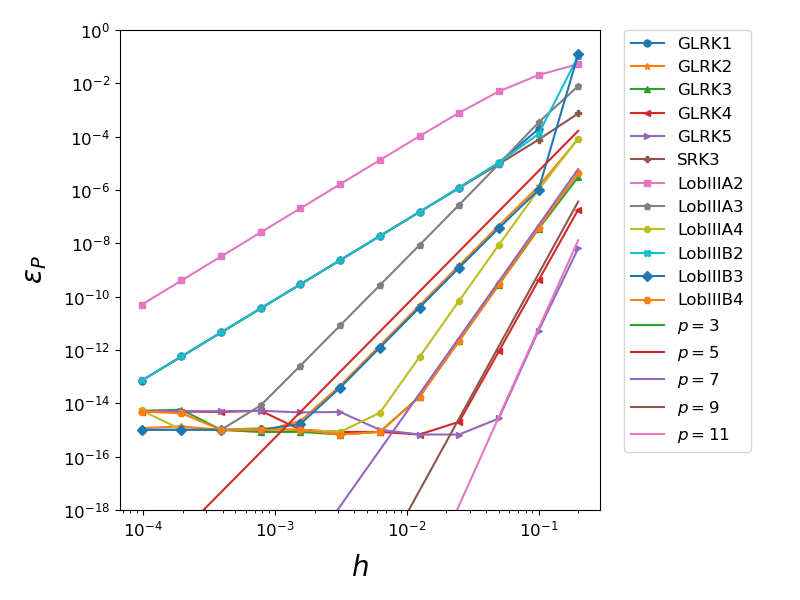}
		}
		
		\subfloat[Symmetric Projection]{
			\includegraphics[width=.3\textwidth]{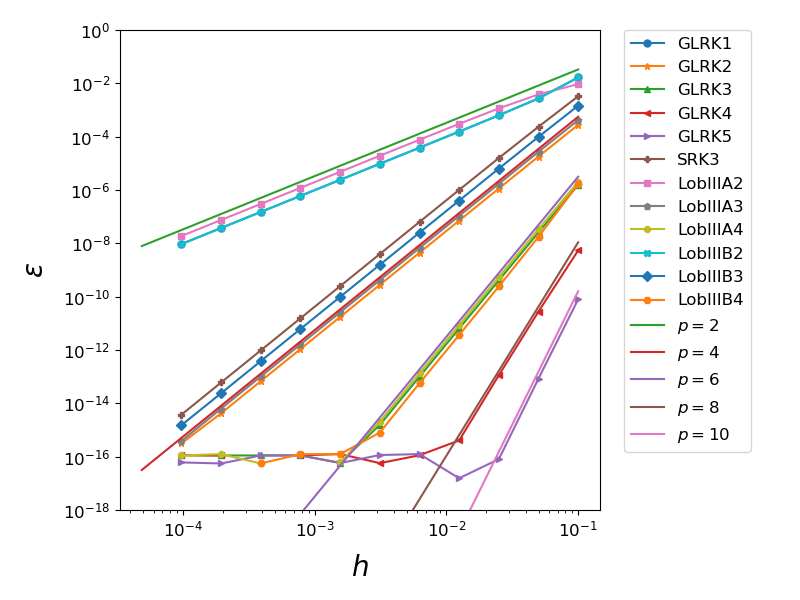}
			\includegraphics[width=.3\textwidth]{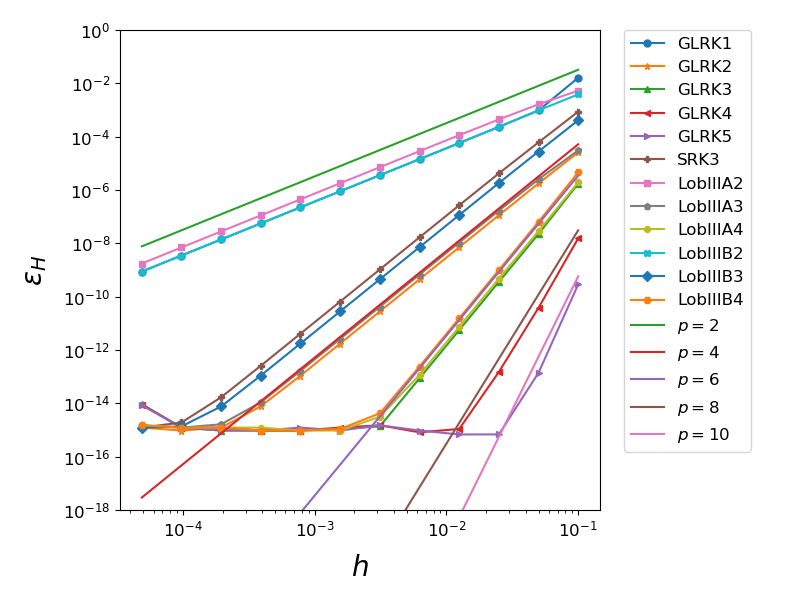}
			\includegraphics[width=.3\textwidth]{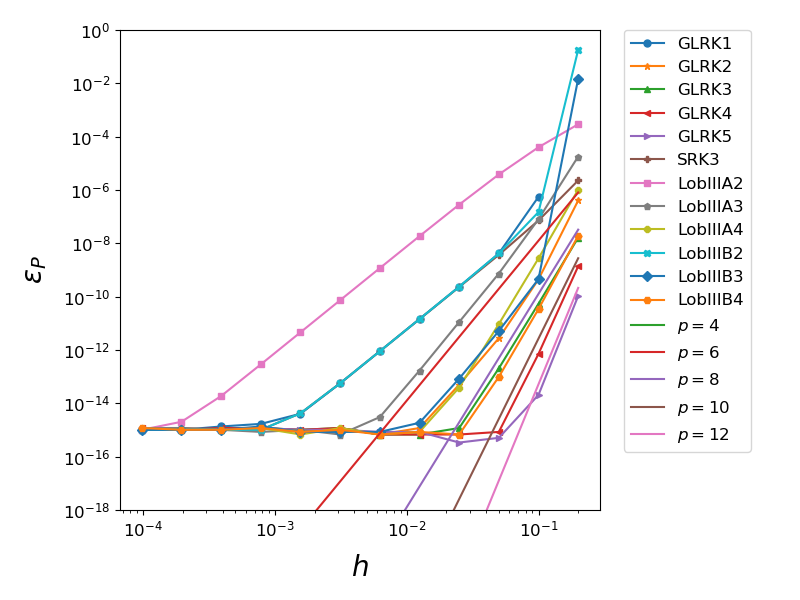}
		}
		
		\subfloat[Symplectic Projection]{
			\includegraphics[width=.3\textwidth]{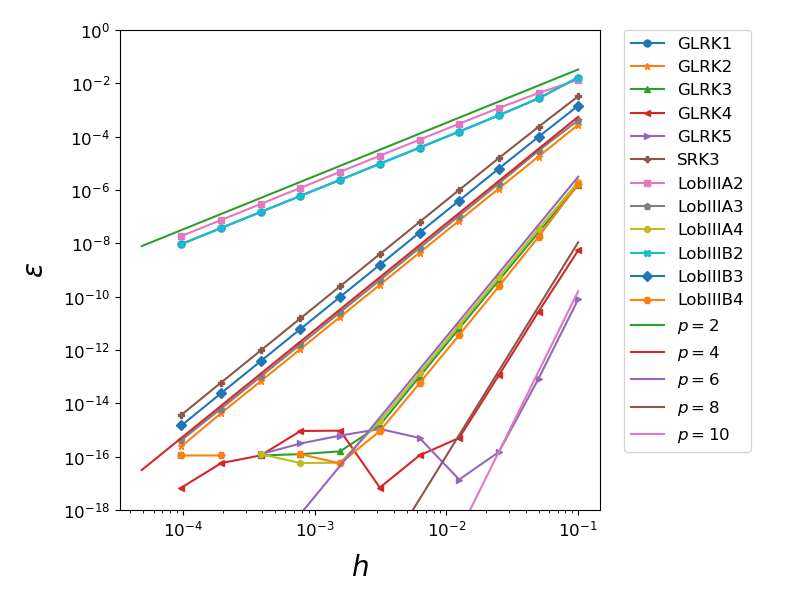}
			\includegraphics[width=.3\textwidth]{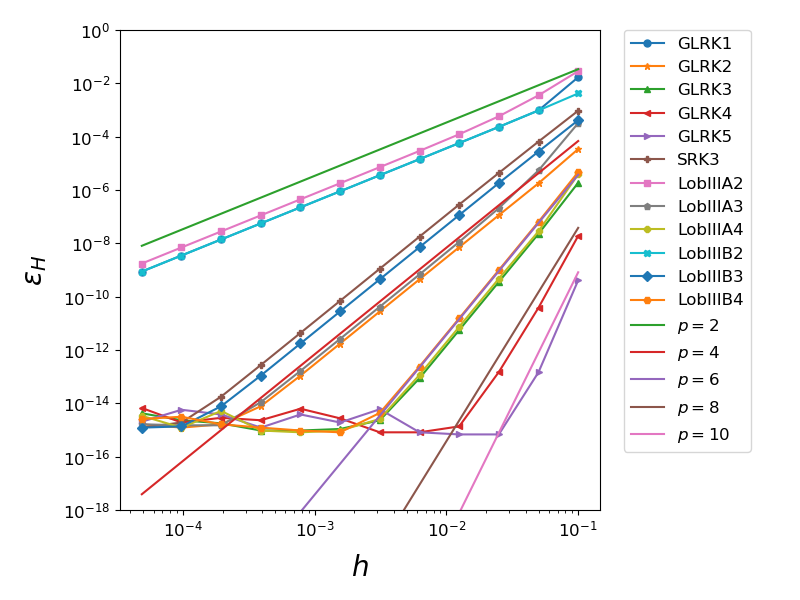}
			\includegraphics[width=.3\textwidth]{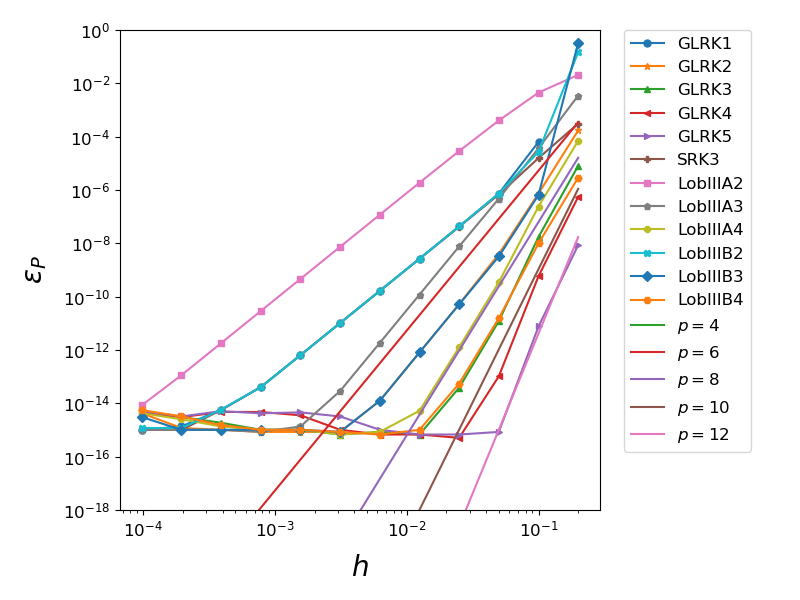}
		}
		
		\subfloat[Midpoint Projection]{
			\includegraphics[width=.3\textwidth]{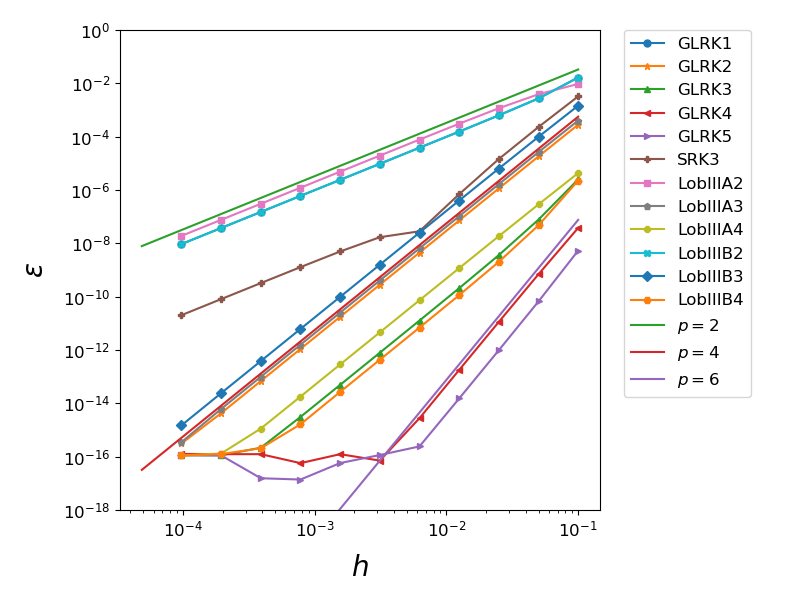}
			\includegraphics[width=.3\textwidth]{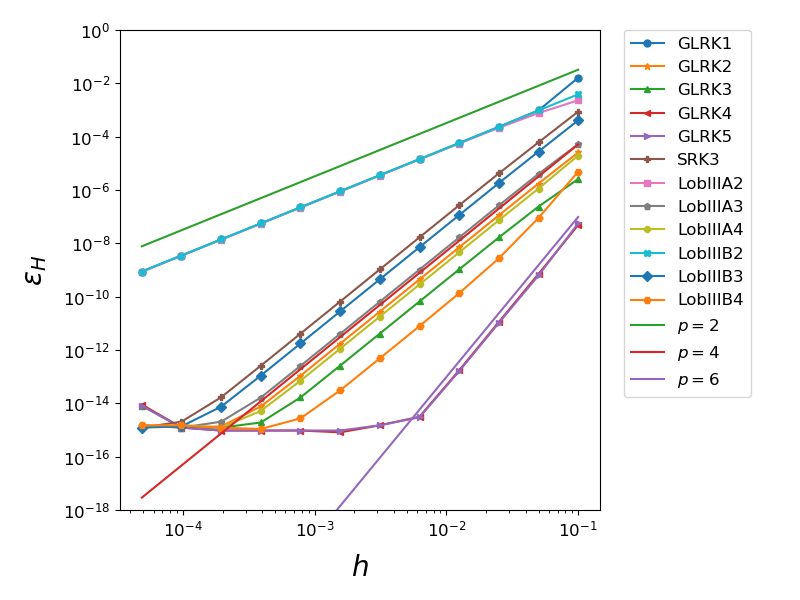}
			\includegraphics[width=.3\textwidth]{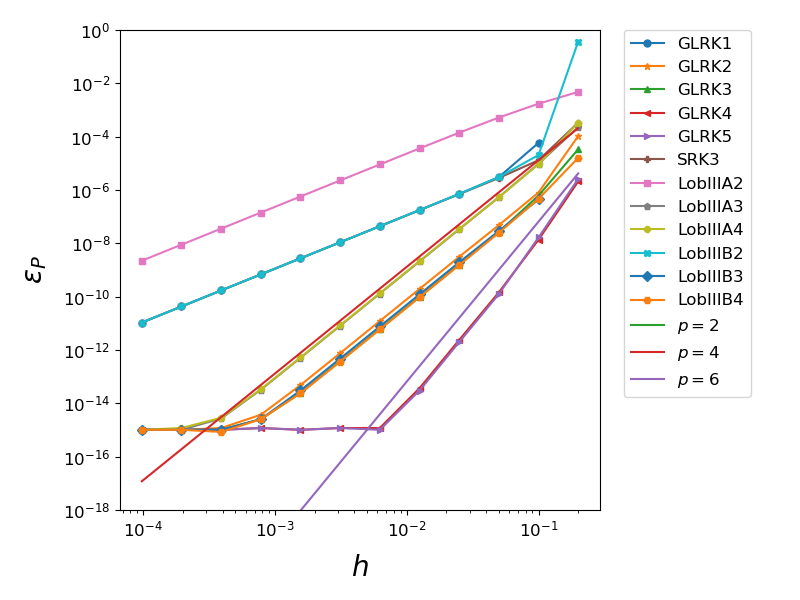}
		}
		
		\caption{Convergence of the solution error $\veps$, the energy error $\veps_{H}$ and the angular momentum error $\veps_{P}$ for the various projection methods and the point vortex example.}
		\label{fig:point_vortices_convergence}
	\end{center}
\end{figure}

We analyse the convergence behaviour of the various integrators together with the various projection methods for the planar point vortex example (c.f., Figure~\ref{fig:point_vortices_convergence}).
The orders for the solution error and the energy error are obtained as listed in Table~\ref{tab:convergence_energy}, while the orders for the angular momentum error are obtained as listed in Table~\ref{tab:convergence_momentum}.

We see that without projection the orders of the integrators are reduced to half the order for ODE problems or half the order plus one in the case of half the order being odd (as all considered methods are symmetric and therefore have even order).
The standard, symmetric and symplectic projection restore the usual order of $2s$ for the Gauss--Legendre methods, $2s-2$ for the Gauss--Lobatto methods and order $4$ for the SRK3 method.
These are mostly well known results (see e.g.~\cite{HairerWanner:1996, Chan:2002, Chan:2004}).
The midpoint projection restores the usual order for the Gauss--Legendre method with one stage and the SRK3 method, that is for the integrators for which the midpoint projection is symplectic, as well as the Gauss--Lobatto methods. For the other Gauss--Legendre methods the order is either $s+1$ or $s+2$ depending on the method being of odd or even $s$, respectively.

Interestingly, for the standard, symmetric and symplectic projection methods, the convergence order of the angular momentum error is increased compared to the convergence order for the solution and the energy error, by one for the standard projection and by two for the symmetric and symplectic projection.
Without projection and with the midpoint projection the convergence of the angular momentum error is the same as the convergence order of the solution and the energy error.
Only the SRK3 method forms an exception, as here the convergence order of the angular momentum error is decreased by one for the standard projection and by two for the midpoint projection with respect to the order of the other errors.

\clearpage

\clearpage

\begin{table}[htb]
	\begin{center}
		\begin{tabular}{|l|c|c|c|c|c|}
		\hline
		 & GLRK   & GLRK   & Lob--IIIA, IIIB & Lob--IIIA, IIIB & SRK3 \\
		 & (odd)  & (even) & (odd)       & (even)      &      \\
		\hline
		No Projection         & $s+1$ & $s$   & $s-1$  & $s$    & $2$ \\
		Standard Projection   & $2s$  & $2s$  & $2s-2$ & $2s-2$ & $4$ \\
		Symmetric Projection  & $2s$  & $2s$  & $2s-2$ & $2s-2$ & $4$ \\
		Symplectic Projection & $2s$  & $2s$  & $2s-2$ & $2s-2$ & $4$ \\
		Midpoint Projection   & $s+1$ & $s+2$ & $2s-2$ & $2s-2$ & $4$ \\
		\hline
		\end{tabular}
	\end{center}
	\caption{Convergence rates for the energy error.}
	\label{tab:convergence_energy}
\end{table}

\begin{table}[htb]
	\begin{center}
		\begin{tabular}{|l|c|c|c|c|c|}
		\hline
		 & GLRK   & GLRK   & Lob--IIIA, IIIB & Lob--IIIA, IIIB & SRK3 \\
		 & (odd)  & (even) & (odd)       & (even)      &      \\
		\hline
		No Projection         & $s+1$  & $s$    & $s-1$  & $s$    & $2$ \\
		Standard Projection   & $2s+1$ & $2s+1$ & $2s-1$ & $2s-1$ & $3$ \\
		Symmetric Projection  & $2s+2$ & $2s+2$ & $2s$   & $2s$   & $4$ \\
		Symplectic Projection & $2s+2$ & $2s+2$ & $2s$   & $2s$   & $4$ \\
		Midpoint Projection   & $s+1$  & $s+2$  & $2s-2$ & $2s-2$ & $2$ \\
		\hline
		\end{tabular}
	\end{center}
	\caption{Convergence rates for the momentum error.}
	\label{tab:convergence_momentum}
\end{table}

\subsection{Poincar\'{e} Integral Invariants}
\label{sec:app_integral_invariants}

It has recently been established that Poincar\'{e} integral invariants provide useful diagnostics for analysing the long-time accuracy of numerical integrators for Hamiltonian dynamics and for distinguishing between symplectic and non-symplectic integrators~\cite{KrausBurby:2017}.
In the following we will apply this diagnostics to the guiding centre example with a simple, symmetric magnetic field.

\begin{figure}[htb]
	\begin{center}
		\subfloat[Trajectories]{
			\includegraphics[height=7cm]{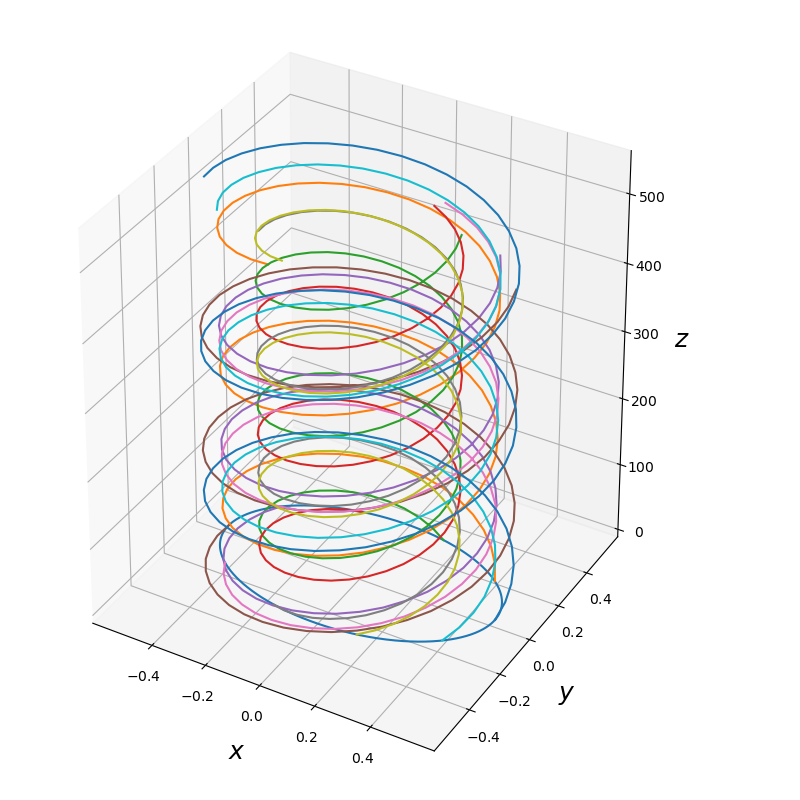}
		}
		\subfloat[Loops]{
			\includegraphics[height=7cm]{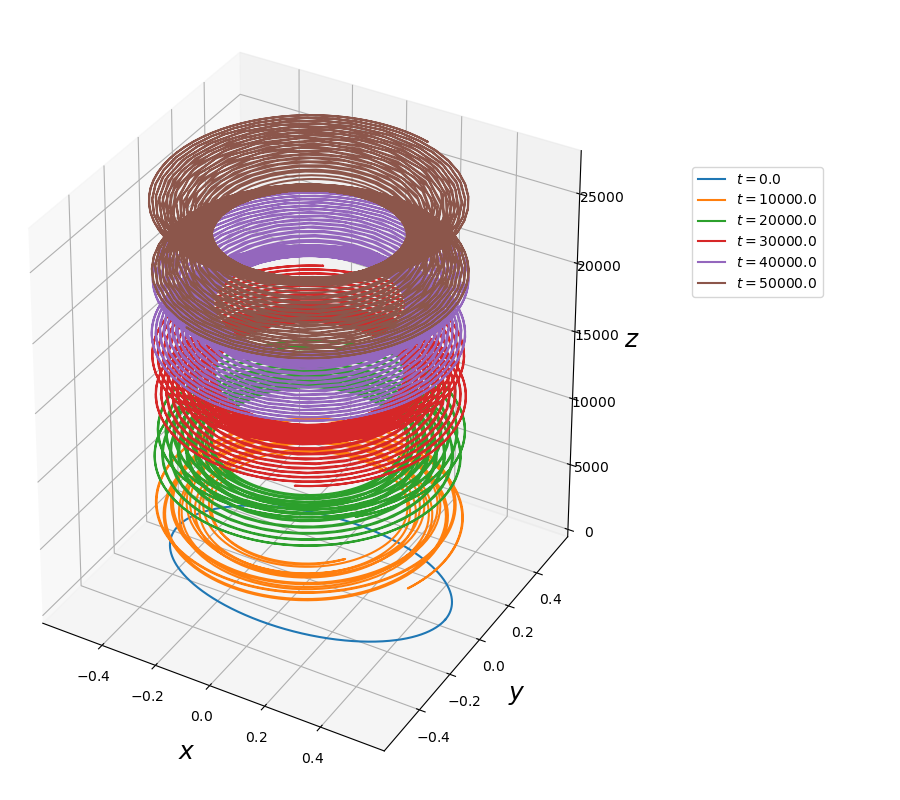}
		}
		
		\caption{Temporal evolution of sample trajectories (left) and a phasespace loop (right) used to compute the first Poincar\'{e} integral invariant for the dynamics of the guiding centre system.}
		\label{fig:guiding_center_4d_poincare1}
	\end{center}
\end{figure}

\begin{figure}[htb]
	\begin{center}
		\subfloat[$t = 0 \, .. \, 1000$]{
			\includegraphics[width=.48\textwidth]{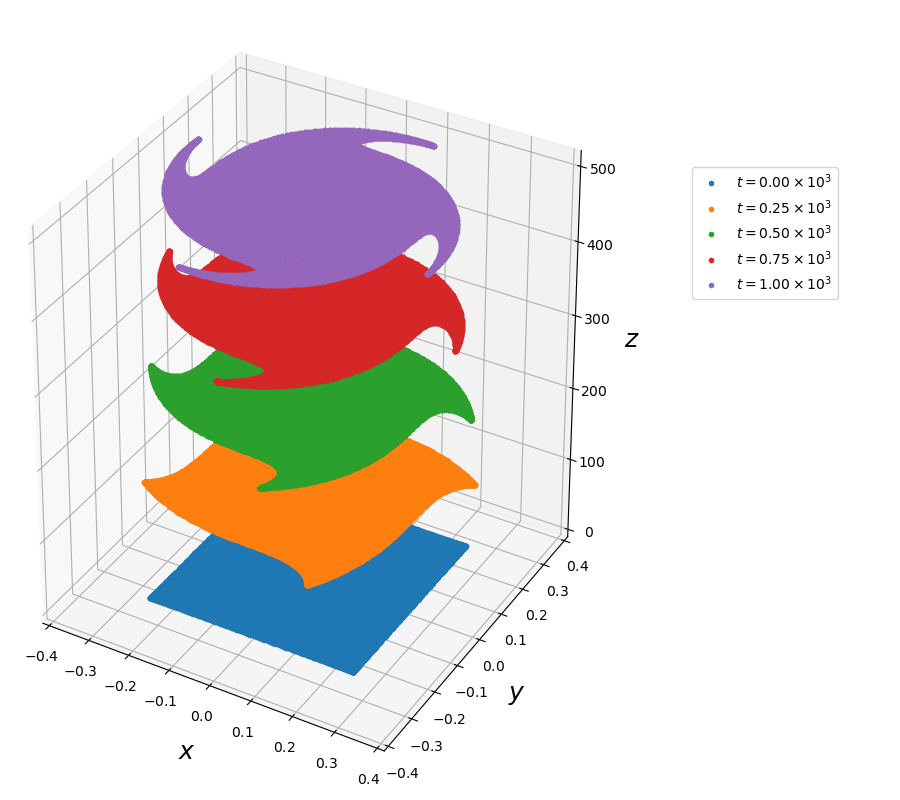}
		}
		\subfloat[$t = 0 \, .. \, 10000$]{
			\includegraphics[width=.48\textwidth]{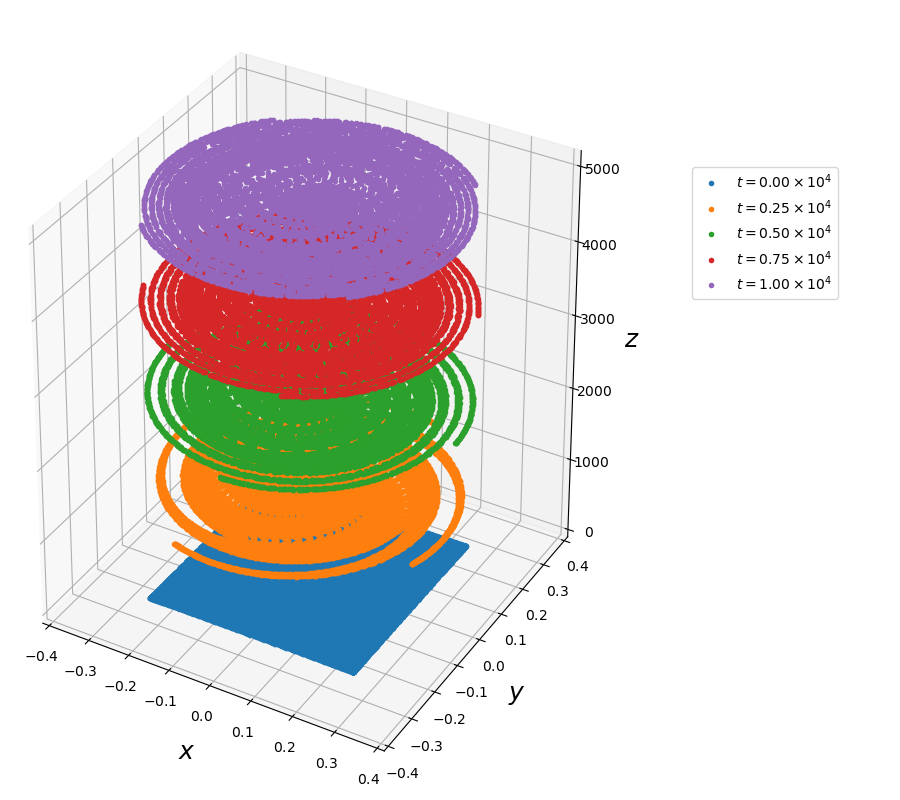}
		}
		
		\caption{Temporal evolution of the phasespace area used to compute the second Poincar\'{e} integral invariant for the dynamics of the guiding centre system.}
		\label{fig:guiding_center_4d_poincare2}
	\end{center}
\end{figure}

In particular, we consider the first and second Poincar\'{e} integral invariant $\vartheta$ and $\obar{\omega}$, which are given by the Lagrangian one- and two-form, respectively.
The one-form $\vartheta$ is a relative integral invariant, which means that the integral
\begin{align}
I_{1} &= \int_{\gamma} \vartheta_{i} (q) \, dq^{i} 
\end{align}
stays constant in time when $\gamma$ is a closed loop in the configuration space $\mf{M}$ (a compact one-dimensional parametrized submanifold of $\mf{M}$ without boundary), that is advected along the solution of the dynamics.
Figure~\ref{fig:guiding_center_4d_poincare1} shows examples of single trajectories of some samples of such a loop, as well as the temporal evolution of the whole loop following the dynamics of the guiding centre system.
The two-form $\obar{\omega}$ is an absolute integral invariant, which means that the integral
\begin{align}
I_{2} &= \int_{S} \obar{\omega}_{ij} (q) \, dq^{i} \, dq^{j} 
\end{align}
stays constant in time when $S$ is any compact two-dimensional parametrized submanifold of $\mf{M}$, advected along the solution of the dynamics.
Figure~\ref{fig:guiding_center_4d_poincare2} shows how an initially rectangular area in phasespace is advected by the dynamics of the guiding centre system.

The loop $\gamma$ is parametrized by $\tau \in [0,1)$, so that
\begin{align}
I_{1} (t) &= \int_{0}^{1} \vartheta_{i} (q_{(\tau)}(t)) \, \dfrac{d q_{(\tau)}^{i}}{d\tau} \, d\tau .
\end{align}
In order to compute this integral, we use $N$ equidistant points in $[0,1)$, so that the derivatives $d q_{(\tau)} / d\tau$ can be efficiently computed via discrete Fourier transforms. The integral is approximated with the trapezoidal quadrature rule, which has spectral convergence on periodic domains~\cite{Trefethen:2014}.
The area $S$ is parametrized by $(\sigma, \tau) \in [0,1]^{2}$, so that
\begin{align}
I_{2} (t) &= \int_{0}^{1} \int_{0}^{1} \obar{\omega}_{ij} (q_{(\sigma, \tau)} (t) \, \dfrac{dq_{(\sigma, \tau)}^{i}}{d\sigma} \, \dfrac{dq_{(\sigma, \tau)}^{j}}{d\tau} \, d\sigma \, d\tau .
\end{align}
Here, we represent the surface in terms of Chebyshev polynomials and thus use Chebyshev points for the discretization of the domain $[0,1]^{2}$. The Chebyshev polynomials allow for an extremely accurate approximation of the surface, even if the latter becomes severely deformed. Moreover, they allow for the use of the \texttt{ApproxFun.jl} package \cite{ApproxFun, Olver:2014, Townsend:2015} for the easy and accurate computation of the derivatives and the integral.
As initial conditions we use
\begin{align}
q_{(\tau)} (0) = (r_{x} \cos (2\pi \tau), \, r_{y} \sin (2\pi \tau), \, r_{z} \sin (2\pi s), \, u_{0} + u_{1} \cos (2\pi s)) 
\end{align}
for the loop $\gamma$ with $r_{x} = 0.5$, $r_{y} = 0.3$, $r_{z} = 0.1$, $u_{0} = 0.5$, $u_{1} = 0.05$ and
\begin{align}
q_{(\sigma, \tau)} (0) = (r_{0} (\sigma-0.5), \, r_{0} (\tau-0.5), \, r_{z} \cos (2\pi \sigma) \cos (2\pi \tau), \, u_{0} + u_{1} \sin (2\pi \sigma) \sin (2\pi \tau) ) 
\end{align}
for the surface $S$ with $r_{0} = 0.1$, $r_{z} = 0.1$, $u_{0} = 0.5$ and $u_{1} = 0.01$. In both cases, the magnetic field is given by $B = (0, \, 0, \, B_{0} \, (1 + x^{2} + y^{2}))$ with $B_{0} = 1$ and we set $\mu = 0.01$. We use a simpler configuration as in Section~\ref{app:guiding_centre_dynamics} as for the magnetic field used there not all integrators are stable.
We use $2000$ points to discretize the loop $\gamma$ and $100 \times 100$ points to discretize the surface $S$. The time step is $h = 10$ in both cases. \\

\begin{figure}[p]
	\begin{center}
		\subfloat[GLRK1pNone]{\label{fig:guiding_centre_4d_poincare1_glrk1_pnone}
			\includegraphics[width=.25\textwidth]{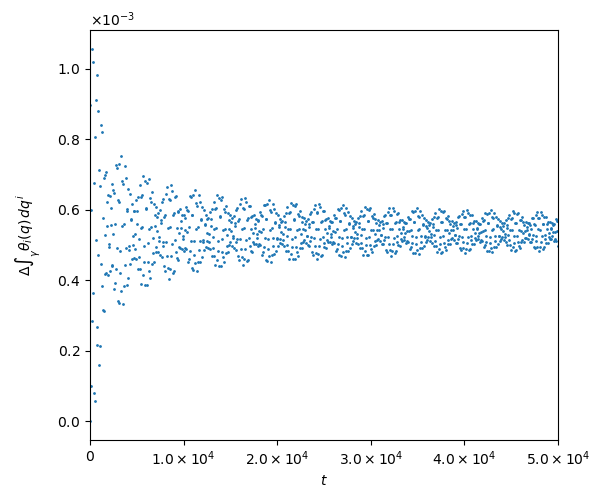}
		}
		\subfloat[GLRK2pNone]{\label{fig:guiding_centre_4d_poincare1_glrk2_pnone}
			\includegraphics[width=.25\textwidth]{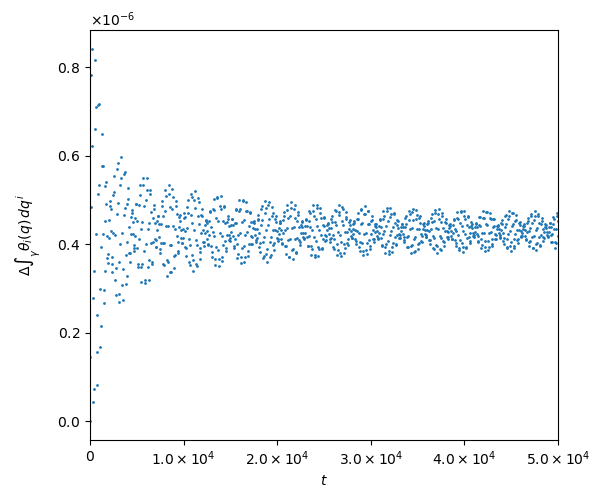}
		}
		\subfloat[SRK3pNone]{\label{fig:guiding_centre_4d_poincare1_srk3_pnone}
			\includegraphics[width=.25\textwidth]{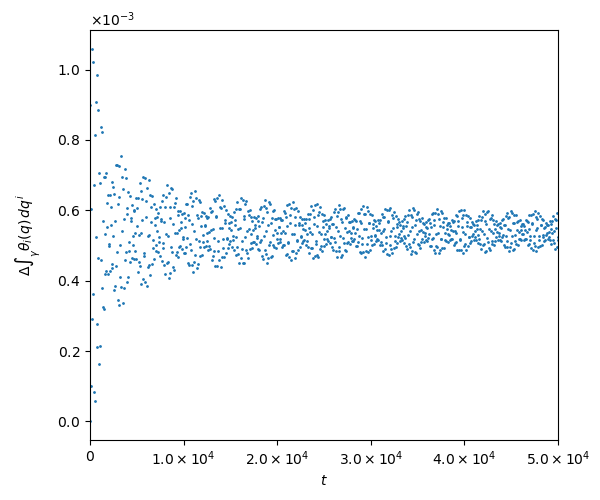}
		}

		\subfloat[GLRK1pStandard]{\label{fig:guiding_centre_4d_poincare1_glrk1_pstandard}
			\includegraphics[width=.25\textwidth]{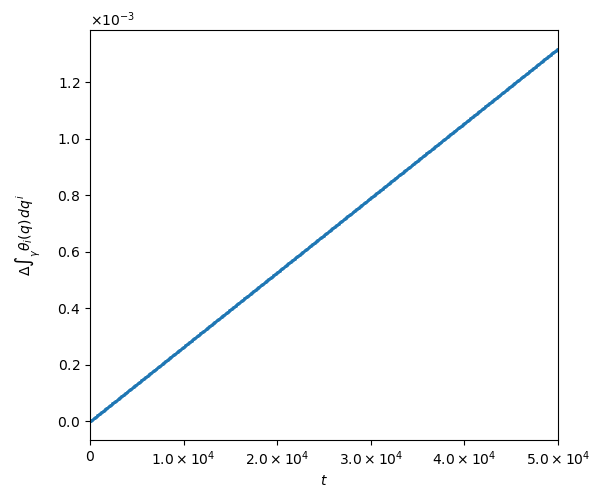}
		}
		\subfloat[GLRK2pStandard]{\label{fig:guiding_centre_4d_poincare1_glrk2_pstandard}
			\includegraphics[width=.25\textwidth]{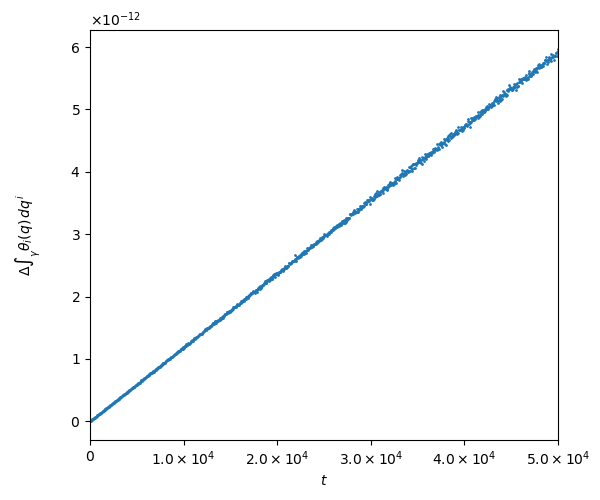}
		}
		\subfloat[SRK3pStandard]{\label{fig:guiding_centre_4d_poincare1_srk3_pstandard}
			\includegraphics[width=.25\textwidth]{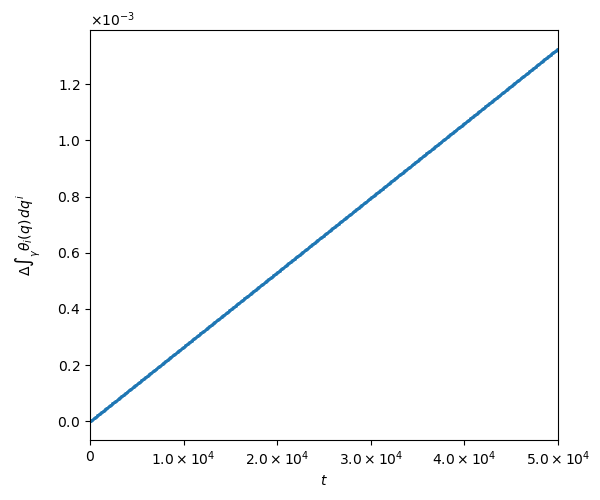}
		}

		\subfloat[GLRK1pSymmetric]{\label{fig:guiding_centre_4d_poincare1_glrk1_psymmetric}
			\includegraphics[width=.25\textwidth]{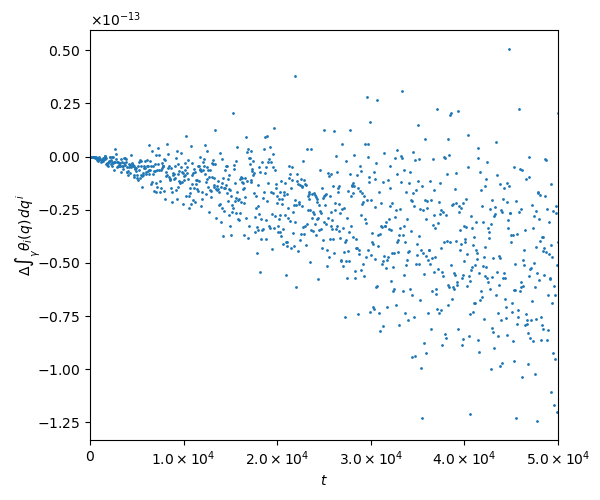}
		}
		\subfloat[GLRK2pSymmetric]{\label{fig:guiding_centre_4d_poincare1_glrk2_psymmetric}
			\includegraphics[width=.25\textwidth]{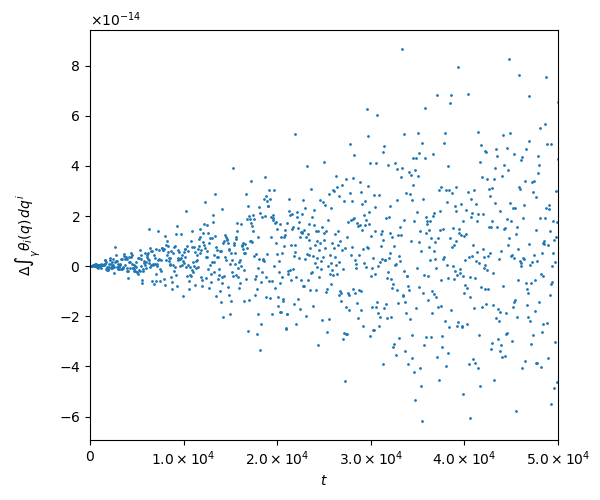}
		}
		\subfloat[SRK3pSymmetric]{\label{fig:guiding_centre_4d_poincare1_srk3_psymmetric}
			\includegraphics[width=.25\textwidth]{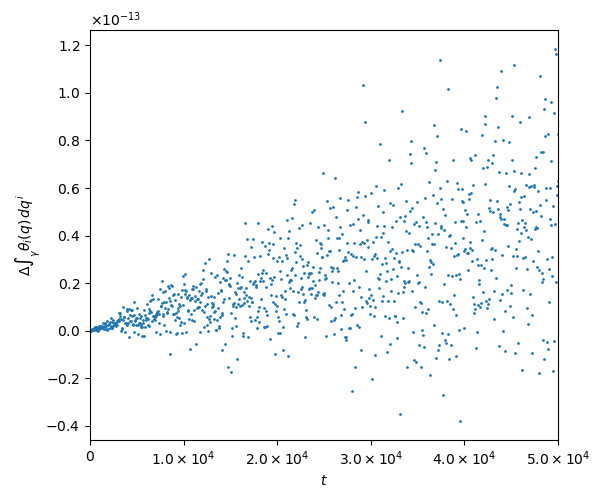}
		}

		\subfloat[GLRK1pSymplectic]{\label{fig:guiding_centre_4d_poincare1_glrk1_psymplectic}
			\includegraphics[width=.25\textwidth]{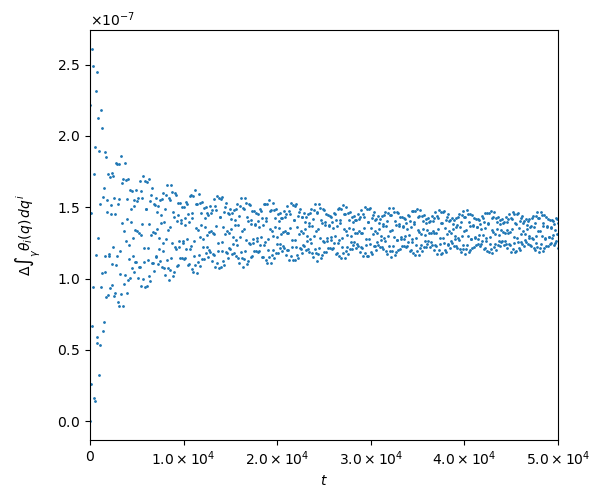}
		}
		\subfloat[GLRK2pSymplectic]{\label{fig:guiding_centre_4d_poincare1_glrk2_psymplectic}
			\includegraphics[width=.25\textwidth]{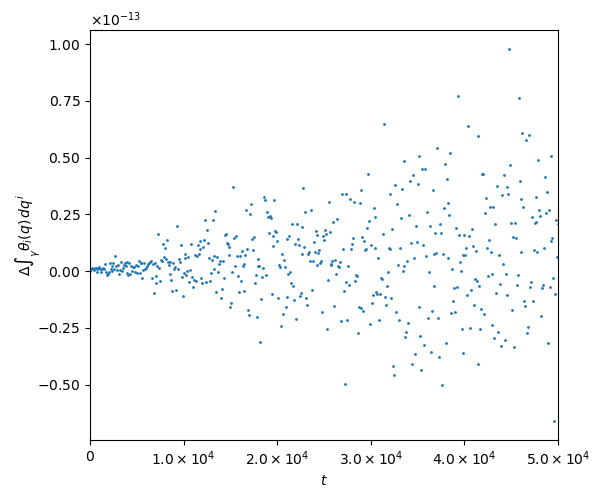}
		}
		\subfloat[SRK3pSymplectic]{\label{fig:guiding_centre_4d_poincare1_srk3_psymplectic}
			\includegraphics[width=.25\textwidth]{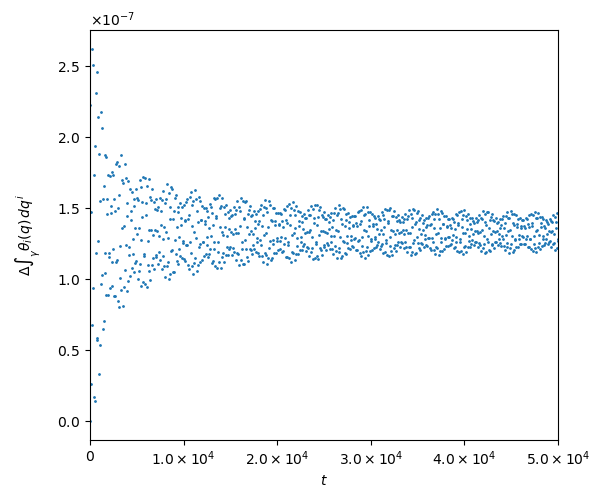}
		}

		\subfloat[GLRK1pMidpoint]{\label{fig:guiding_centre_4d_poincare1_glrk1_pmidpoint}
			\includegraphics[width=.25\textwidth]{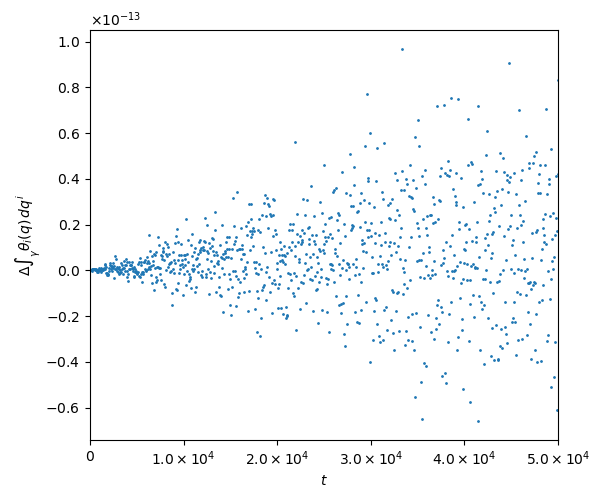}
		}
		\subfloat[GLRK2pMidpoint]{\label{fig:guiding_centre_4d_poincare1_glrk2_pmidpoint}
			\includegraphics[width=.25\textwidth]{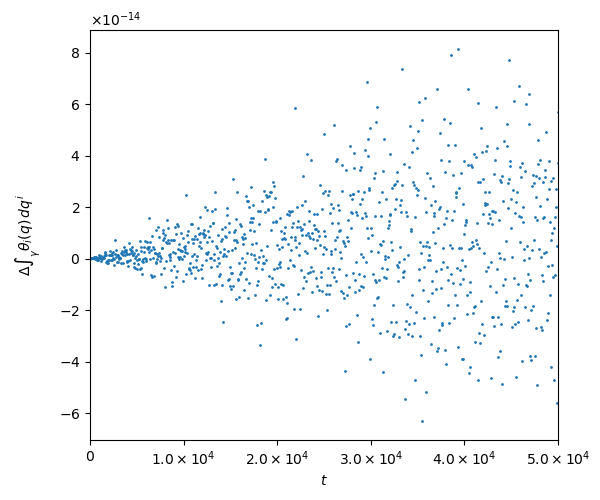}
		}
		\subfloat[SRK3pMidpoint]{\label{fig:guiding_centre_4d_poincare1_srk3_pmidpoint}
			\includegraphics[width=.25\textwidth]{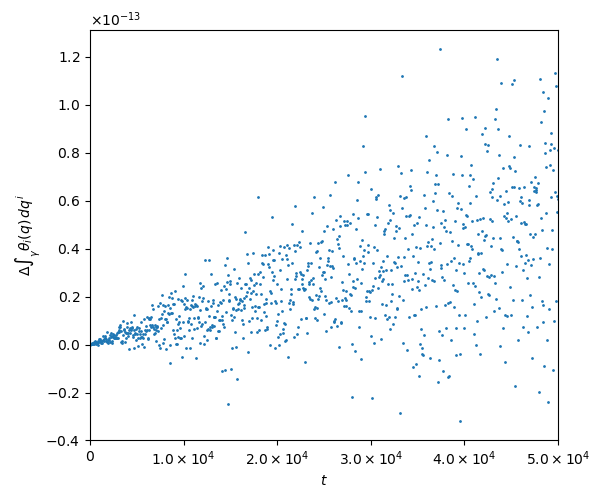}
		}

		\caption{First Poincar\'{e} integral invariant with respect to the noncanonical one-form $\vartheta_i (q) \, \ext q^{i}$ for a guiding centre particle in a symmetric magnetic field with one- and two-stage variational Gauss--Legendre Runge--Kutta and SRK3 methods.}
		\label{fig:guiding_centre_4d_poincare1}
	\end{center}
\end{figure}

\begin{figure}[p]
	\begin{center}
		\subfloat[GLRK1pNone]{\label{fig:guiding_centre_4d_poincare2_glrk1_pnone}
			\includegraphics[width=.24\textwidth]{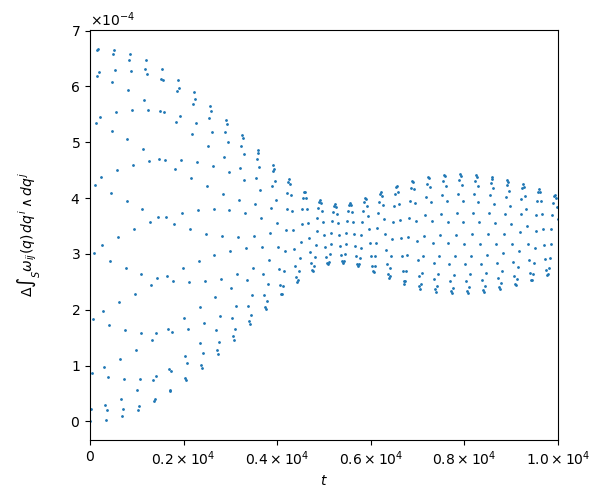}
		}
		\subfloat[GLRK2pNone]{\label{fig:guiding_centre_4d_poincare2_glrk2_pnone}
			\includegraphics[width=.24\textwidth]{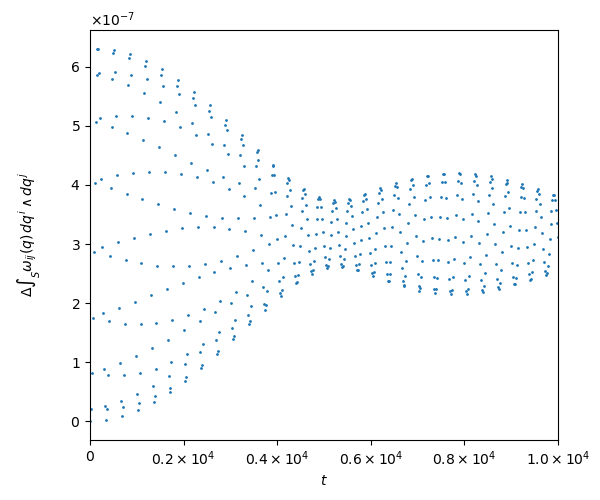}
		}
		\subfloat[SRK3pNone]{\label{fig:guiding_centre_4d_poincare2_srk3_pnone}
			\includegraphics[width=.24\textwidth]{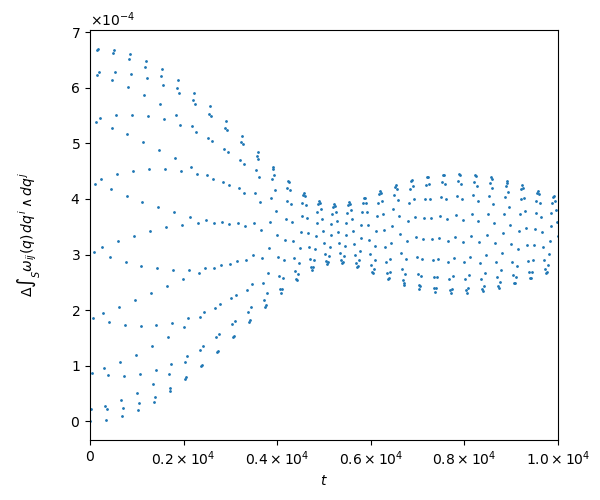}
		}

		\subfloat[GLRK1pStandard]{\label{fig:guiding_centre_4d_poincare2_glrk1_pstandard}
			\includegraphics[width=.24\textwidth]{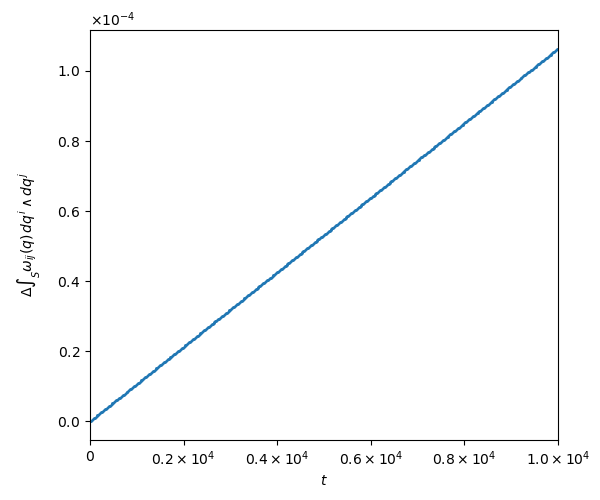}
		}
		\subfloat[GLRK2pStandard]{\label{fig:guiding_centre_4d_poincare2_glrk2_pstandard}
			\includegraphics[width=.24\textwidth]{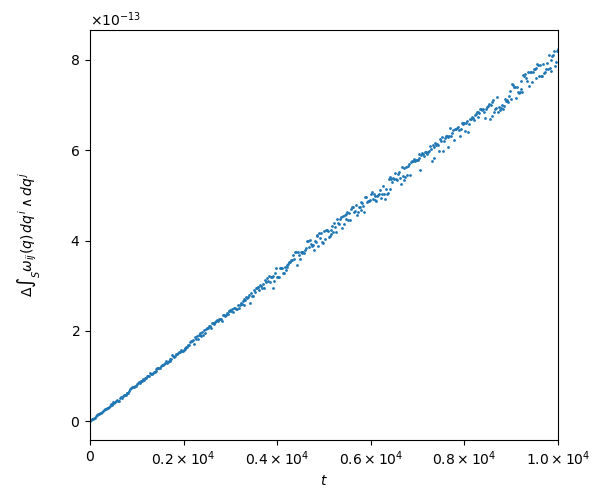}
		}
		\subfloat[SRK3pStandard]{\label{fig:guiding_centre_4d_poincare2_srk3_pstandard}
			\includegraphics[width=.24\textwidth]{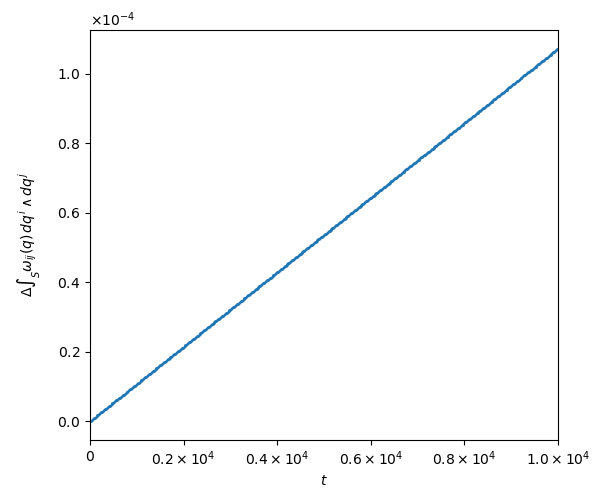}
		}

		\subfloat[GLRK1pSymmetric]{\label{fig:guiding_centre_4d_poincare2_glrk1_psymmetric}
			\includegraphics[width=.24\textwidth]{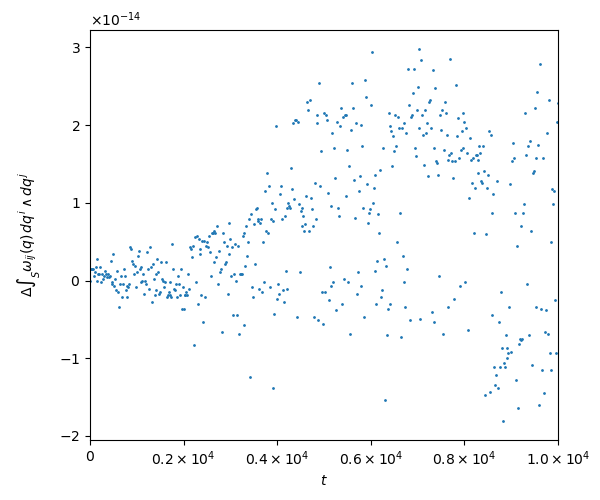}
		}
		\subfloat[GLRK2pSymmetric]{\label{fig:guiding_centre_4d_poincare2_glrk2_psymmetric}
			\includegraphics[width=.24\textwidth]{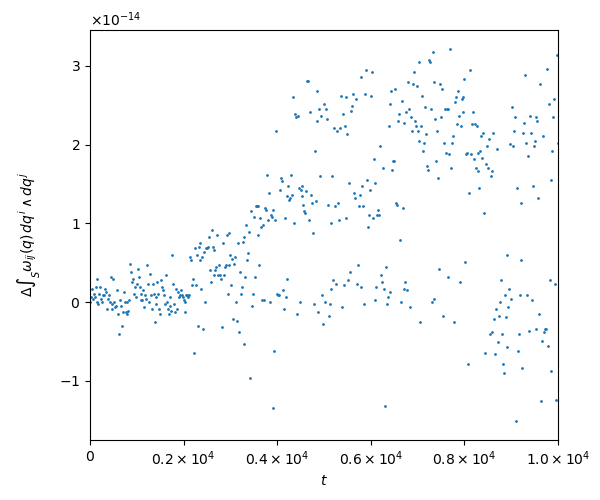}
		}
		\subfloat[SRK3pSymmetric]{\label{fig:guiding_centre_4d_poincare2_srk3_psymmetric}
			\includegraphics[width=.24\textwidth]{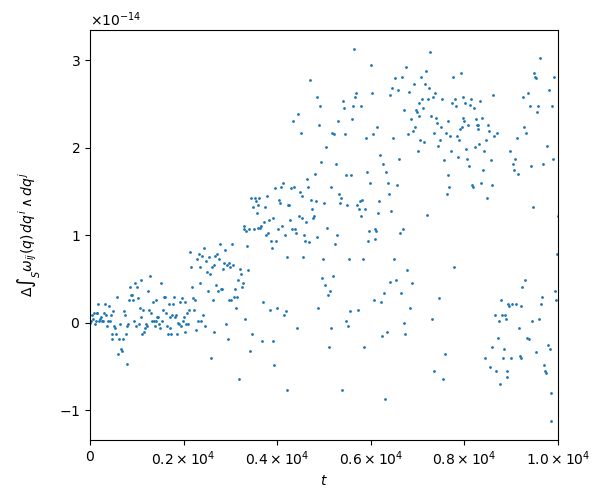}
		}

		\subfloat[GLRK1pSymplectic]{\label{fig:guiding_centre_4d_poincare2_glrk1_psymplectic}
			\includegraphics[width=.24\textwidth]{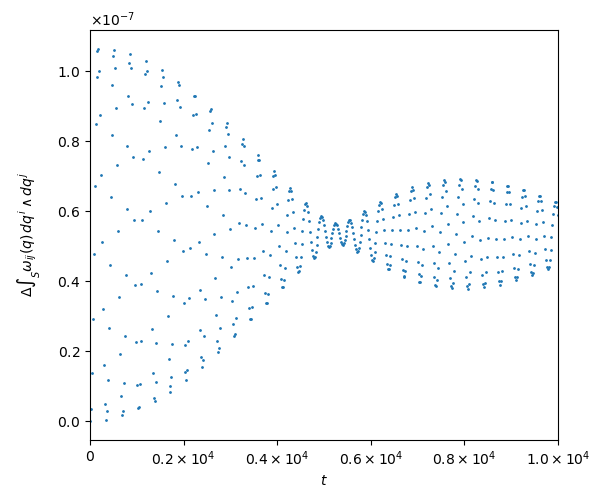}
		}
		\subfloat[GLRK2pSymplectic]{\label{fig:guiding_centre_4d_poincare2_glrk2_psymplectic}
			\includegraphics[width=.24\textwidth]{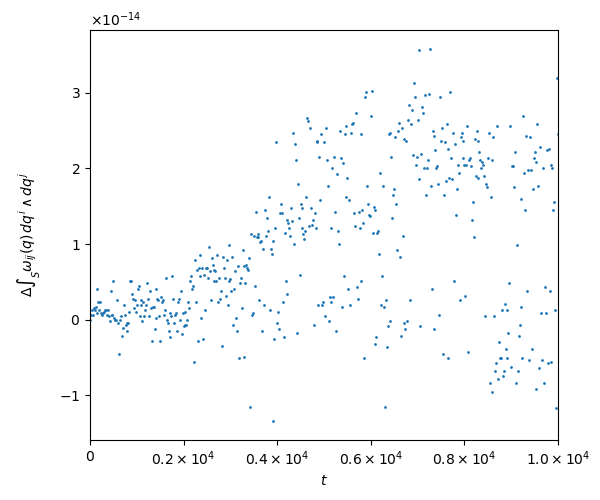}
		}
		\subfloat[SRK3pSymplectic]{\label{fig:guiding_centre_4d_poincare2_srk3_psymplectic}
			\includegraphics[width=.24\textwidth]{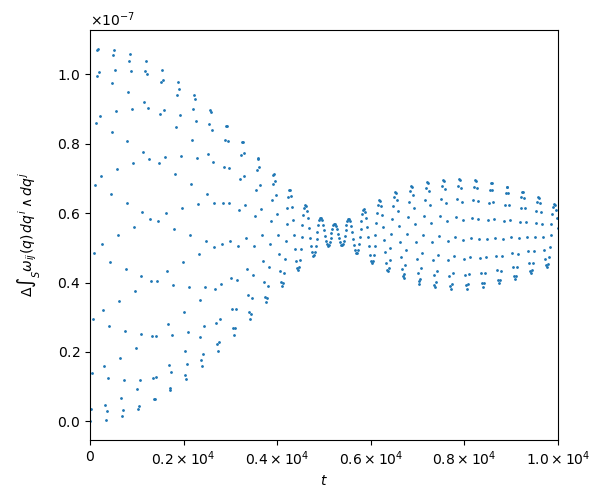}
		}

		\subfloat[GLRK1pMidpoint]{\label{fig:guiding_centre_4d_poincare2_glrk1_pmidpoint}
			\includegraphics[width=.24\textwidth]{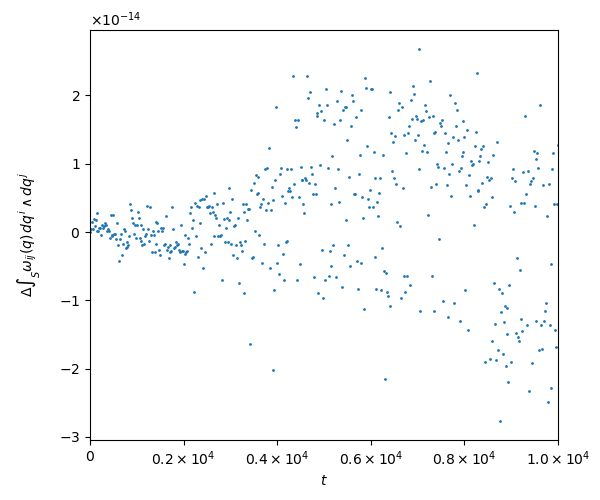}
		}
		\subfloat[GLRK2pMidpoint]{\label{fig:guiding_centre_4d_poincare2_glrk2_pmidpoint}
			\includegraphics[width=.24\textwidth]{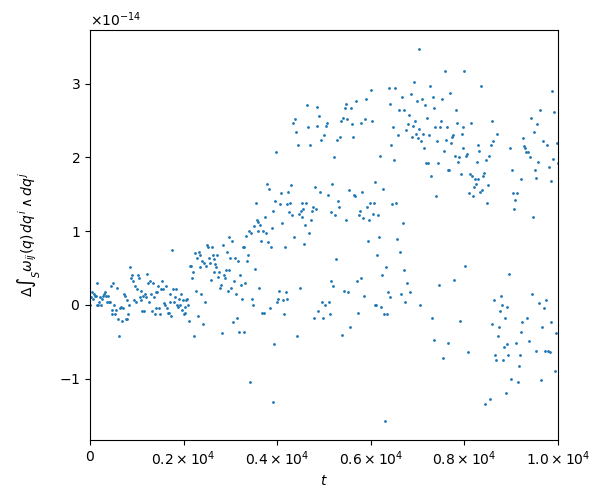}
		}
		\subfloat[SRK3pMidpoint]{\label{fig:guiding_centre_4d_poincare2_srk3_pmidpoint}
			\includegraphics[width=.24\textwidth]{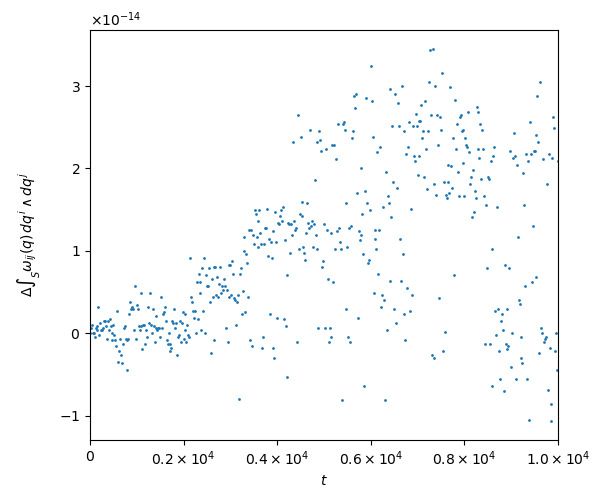}
		}

		\caption{Second Poincar\'{e} integral invariant with respect to the noncanonical two-form $\obar{\omega}_{ij} (q) \, \ext q^{i} \wedge \ext q^{j}$ for a guiding centre particle in a symmetric magnetic field with one- and two-stage variational Gauss--Legendre Runge--Kutta and SRK3 methods.}
		\label{fig:guiding_centre_4d_poincare2}
	\end{center}
\end{figure}

\begin{figure}[p]
	\begin{center}
		\subfloat[GLRK1pNone]{
			\includegraphics[width=.25\textwidth]{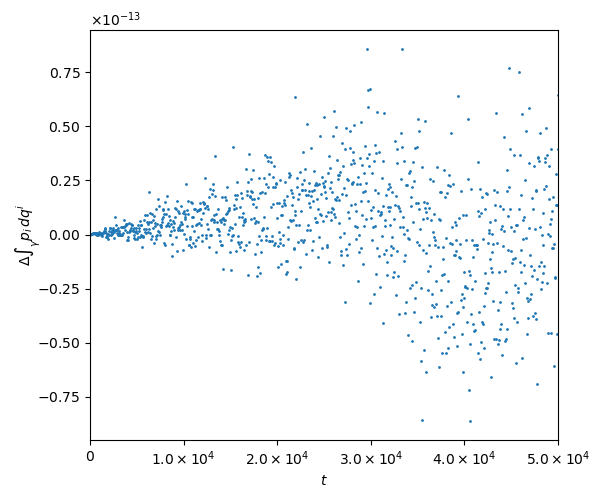}
		}
		\subfloat[GLRK2pNone]{
			\includegraphics[width=.25\textwidth]{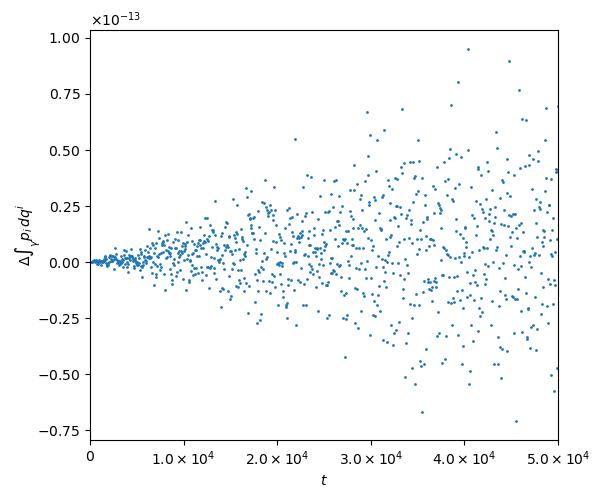}
		}
		\subfloat[SRK3pNone]{
			\includegraphics[width=.25\textwidth]{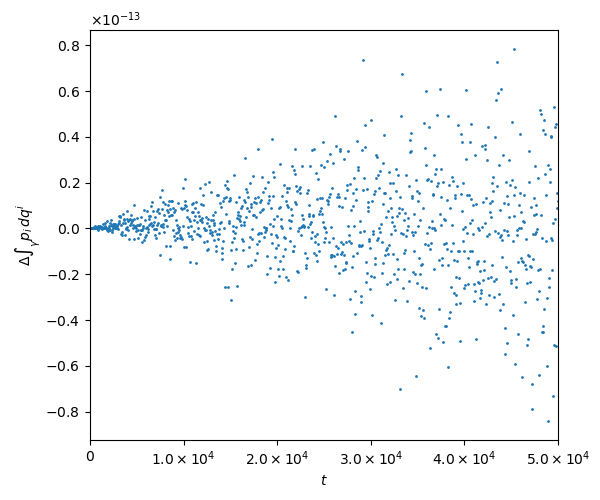}
		}

		\caption{First Poincar\'{e} integral invariant for the canonical one-form $p_{i} \, \ext q^{i}$.}
		\label{fig:guiding_centre_4d_poincare1_canonical}
	\end{center}
\end{figure}

\begin{figure}[p]
	\begin{center}
		\subfloat[GLRK1pNone]{
			\includegraphics[width=.24\textwidth]{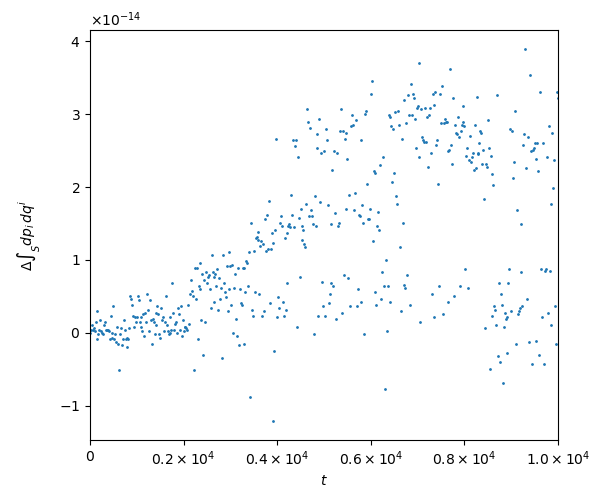}
		}
		\subfloat[GLRK2pNone]{
			\includegraphics[width=.24\textwidth]{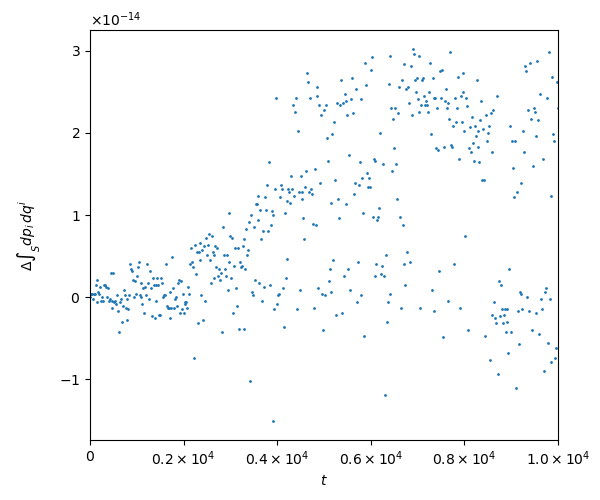}
		}
		\subfloat[SRK3pNone]{
			\includegraphics[width=.24\textwidth]{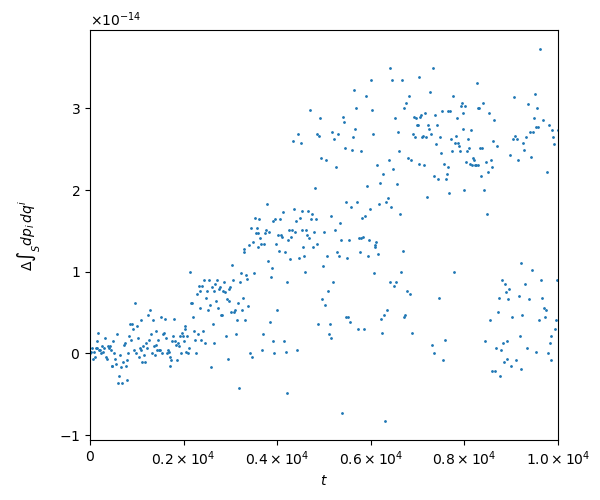}
		}

		\caption{Second Poincar\'{e} integral invariant for the canonical two-form $\ext p_{i} \wedge \ext q^{i}$.}
		\label{fig:guiding_centre_4d_poincare2_canonical}
	\end{center}
\end{figure}

\begin{figure}[p]
	\begin{center}
		\subfloat[GLRK1pSymplectic]{
			\includegraphics[width=.25\textwidth]{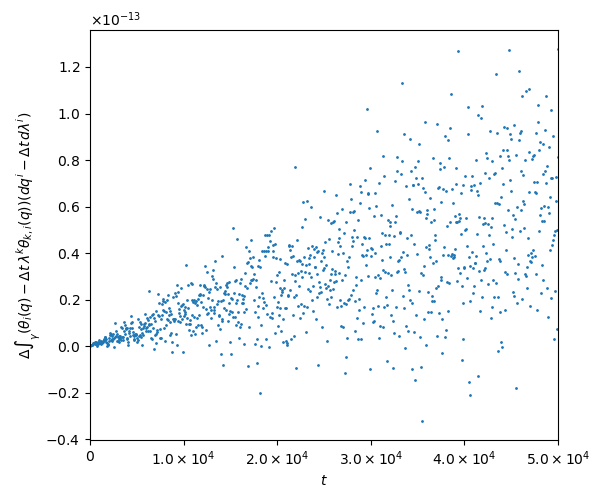}
		}
		\subfloat[GLRK2pSymplectic]{
			\includegraphics[width=.25\textwidth]{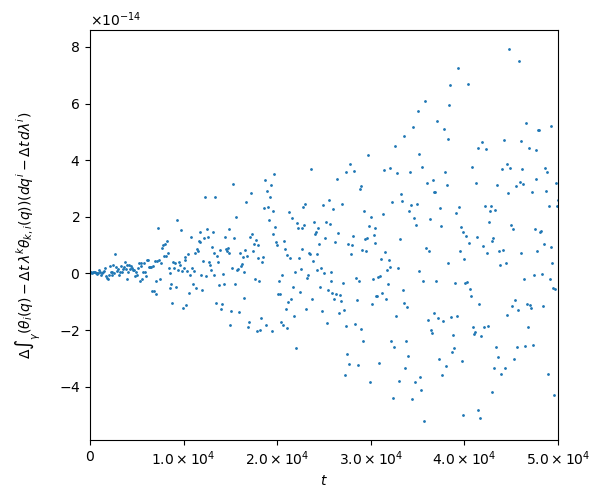}
		}
		\subfloat[SRK3pSymplectic]{
			\includegraphics[width=.25\textwidth]{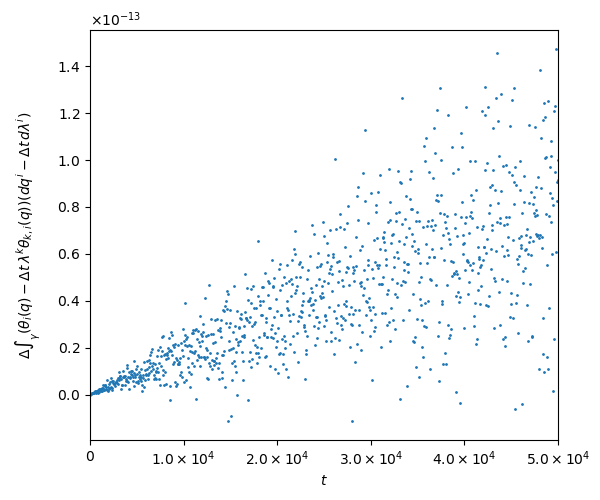}
		}

		\caption{First Poincar\'{e} integral invariant for the noncanonical one-form with correction term, $( \vartheta_{i} (\cq) - h \, \lambda^{k} \vartheta_{k,i} (\cq) ) \, ( \ext \cq^{i} - h \, \ext \cl^{i} )$.}
		\label{fig:guiding_centre_4d_poincare1_corrected}
	\end{center}
\end{figure}

\begin{figure}[p]
	\begin{center}
		\subfloat[GLRK1pSymplectic]{
			\includegraphics[width=.25\textwidth]{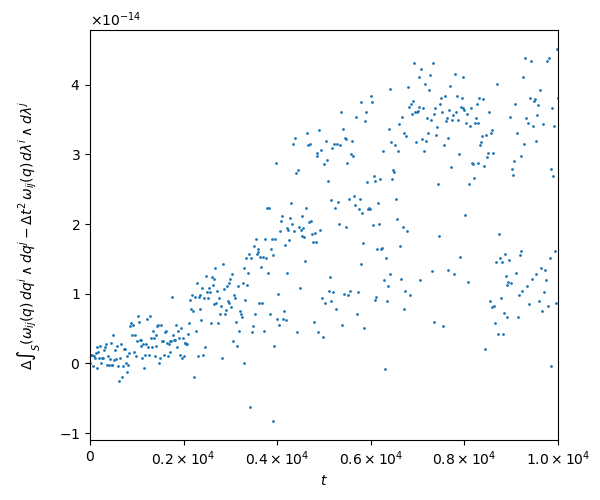}
		}
		\subfloat[GLRK2pSymplectic]{
			\includegraphics[width=.25\textwidth]{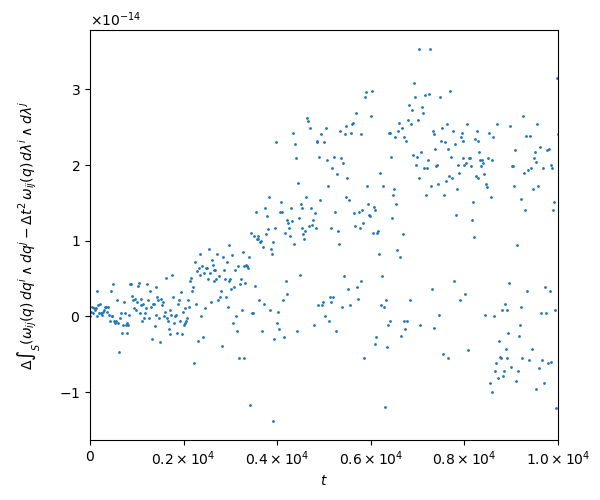}
		}
		\subfloat[SRK3pSymplectic]{
			\includegraphics[width=.25\textwidth]{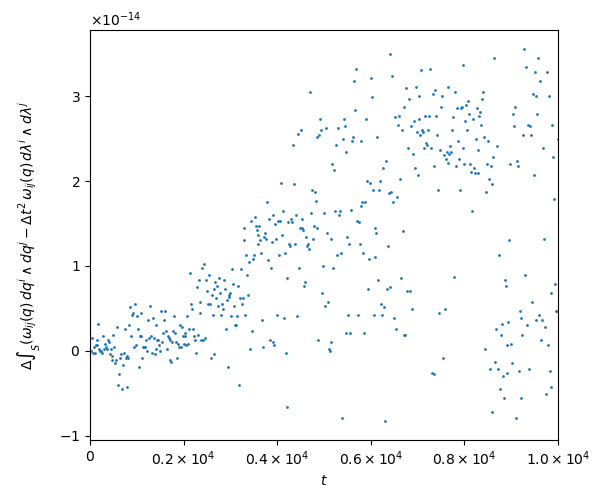}
		}

		\caption{Second Poincar\'{e} integral invariant for the noncanonical two-form with correction term, $\obar{\omega}_{ij} (\cq) \, \ext \cq^{i} \wedge \ext \cq^{j} - h^2 \, \obar{\omega}_{ij} (\cq) \, \ext \cl^{i} \wedge \ext \cl^{j} - h^{2} \cl^{k} \vartheta_{k,ij} (\cq) \, \ext \cq^{i} \wedge \ext \cl^{j}$.}
		\label{fig:guiding_centre_4d_poincare2_corrected}
	\end{center}
\end{figure}

We make the following observations:
\begin{itemize}
\item The unprojected variational integrators preserve the Poincar\'{e} integral invariants with respect to the canonical one-form $p_{i} \, \ext q^{i}$ and two-form $\ext p_{i} \wedge \ext q^{i}$ (Figures~\ref{fig:guiding_centre_4d_poincare1_canonical} and~\ref{fig:guiding_centre_4d_poincare2_canonical}) but not the invariants with respect to the noncanonical one-form $\vartheta_{i} \, \ext q^{i}$ (Figures~\ref{fig:guiding_centre_4d_poincare1_glrk1_pnone}-\ref{fig:guiding_centre_4d_poincare1_srk3_pnone}) and two-form $\obar{\omega}_{ij} (q) \, \ext q^{i} \wedge \ext q^{j}$ (Figures~\ref{fig:guiding_centre_4d_poincare2_glrk1_pnone}-\ref{fig:guiding_centre_4d_poincare2_srk3_pnone}).

This behaviour is expected. The variational integrators preserve the canonical forms by construction but not the noncanonical forms as the solution does not satisfy the Dirac constraint under which the two are equivalent.
However, in cases like this one, where the solution stays close to the constraint submanifold, conservation of the canonical forms implies also approximate conservation of the noncanonical forms. This manifests itself in the error of the integral invariant appearing to be bounded.

In addition, conservation of the canonical forms implies conservation of discrete noncanonical forms, which are obtained by pulling back the canonical forms with the discrete fibre derivatives~\eqref{eq:discrete_fibre_derivative_2}.

\item For the standard projection, we observe a clear drift in both the first and second Poincar\'{e} integral invariant for the one- and two-stage Gauss--Legendre Runge--Kutta integrators (Figures~\ref{fig:guiding_centre_4d_poincare1_glrk1_pstandard}-\ref{fig:guiding_centre_4d_poincare1_glrk2_pstandard} and~\ref{fig:guiding_centre_4d_poincare2_glrk1_pstandard}-\ref{fig:guiding_centre_4d_poincare2_glrk2_pstandard}) as well as the SRK3 method (Figures~\ref{fig:guiding_centre_4d_poincare1_srk3_pstandard} and~\ref{fig:guiding_centre_4d_poincare2_srk3_pstandard}). 
For the higher-order Gauss--Legendre Runge--Kutta integrators (not shown) we see no drift over the runtime of the simulation as the error in the Poincar\'{e} integral invariants introduced by the standard projection is smaller than machine accuracy and therefore not measurable.

\item For the symmetric projection, both the first and second integral invariant is exactly preserved by all three of the considered methods (Figures~\ref{fig:guiding_centre_4d_poincare1_glrk1_psymmetric}-\ref{fig:guiding_centre_4d_poincare1_srk3_psymmetric} and~\ref{fig:guiding_centre_4d_poincare2_glrk1_psymmetric}-\ref{fig:guiding_centre_4d_poincare2_srk3_psymmetric}). 
We only observe some ``diffusion'' of the error. This is expected due to the deformation of the loop and the surface, over which the integrals are evaluated, leading to a degradation of accuracy in the integrals.

\item The symplectic projection preserves neither the first (Figures~\ref{fig:guiding_centre_4d_poincare1_glrk1_psymplectic}-\ref{fig:guiding_centre_4d_poincare1_srk3_psymplectic}) nor the second (Figures~\ref{fig:guiding_centre_4d_poincare2_glrk1_psymplectic}-\ref{fig:guiding_centre_4d_poincare2_srk3_psymplectic}) integral invariant exactly. Instead it shows a behaviour of the error of the invariant similar to the unprojected methods albeit at greatly reduced amplitude (for the two-stage Gauss--Legendre Runge--Kutta integrator the error is already reduced below the machine accuracy, so the invariants appear to be preserved).
This behaviour is expected, though, as the symplectic projection does not preserve the standard noncanonical one- and two-form but the modified forms~\eqref{eq:symplectic_projection_one_form} and~\eqref{eq:symplectic_projection_two_form}. 
Taking the corresponding ``correction terms'' to the standard one- and two-form into account, we find exact conservation of both the first and the second invariant (Figures~\ref{fig:guiding_centre_4d_poincare1_corrected} and~\ref{fig:guiding_centre_4d_poincare2_corrected}).

\item For the midpoint projection, we find exact conservation of both the first and second invariant for all three integrators (Figures~\ref{fig:guiding_centre_4d_poincare1_glrk1_pmidpoint}-\ref{fig:guiding_centre_4d_poincare1_srk3_pmidpoint} and~\ref{fig:guiding_centre_4d_poincare2_glrk1_pmidpoint}-\ref{fig:guiding_centre_4d_poincare2_srk3_pmidpoint}), even for the two-stage Gauss--Legendre Runge--Kutta method, for which the midpoint projection is not symplectic.
\end{itemize}

We have seen that the Poincar\'{e} invariant diagnostics provides a useful tool in experimentally judging the symplecticity of a numerical algorithm.
We see that the unprojected variational integrators only preserve the invariants with respect to the canonical one- and two-form, but not with respect to the noncanonical forms.
With the standard projection all integrators exhibit a clear drift in both invariants.
With the symmetric and the midpoint projection, this drift is essentially gone, even though the former is not symplectic.
Similarly, with the symplectic projection and all the integrators, we find exact conservation of the corrected invariant.
Yet, these results should be taken with a grain of salt. The test case considered here is relatively simple. It is conceivable that in more complicated situations larger errors in the Poincar\'{e} invariants will be observed with the various projection methods. It is still expected, though, that the errors are bounded with all three of the symmetric, midpoint and symplectic projection methods.
If it is conceived necessary, it should always be possible to reduce the errors in the invariants to the level of the machine accuracy by reducing the time step, where in most cases we expect only a mild reduction to be necessary.
When the errors in the invariants are smaller than the machine accuracy, the projected methods essentially behave like any other nonlinearly implicit symplectic integrator, which needs to be solved in finite precision and is thus never exactly symplectic (see \citet{Tan:2005} for details).

\clearpage

\begin{figure}[htb]
	\begin{center}
		\subfloat[$R(\infty)=+1$]{
			\includegraphics[width=.4\textwidth]{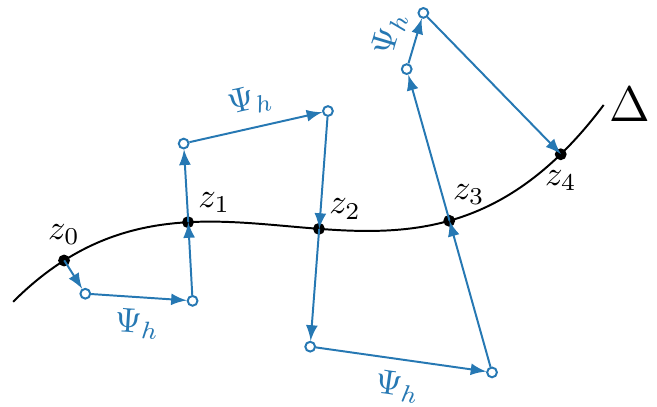}
		}
		\subfloat[$R(\infty)=-1$]{
			\includegraphics[width=.4\textwidth]{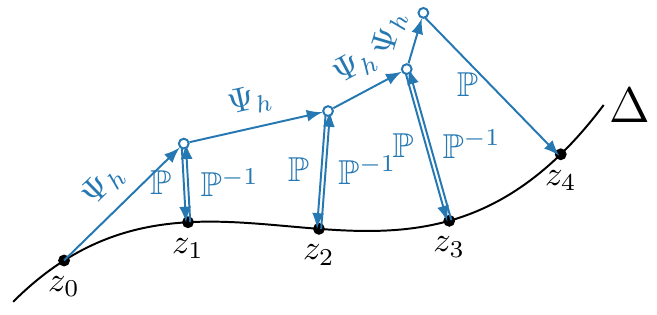}
		}
	
		\caption{Illustration of the symplectic projection method in the unstable case with growing Lagrange multiplier $\lambda$.}
		\label{fig:symplectic-projection}
	\end{center}
\end{figure}

\subsection{Symplectic Projection}
\label{sec:app_symplectic_projection}

\begin{figure}[htb]
	\begin{center}
		\subfloat[$\lambda_1$]{
			\includegraphics[width=.24\textwidth]{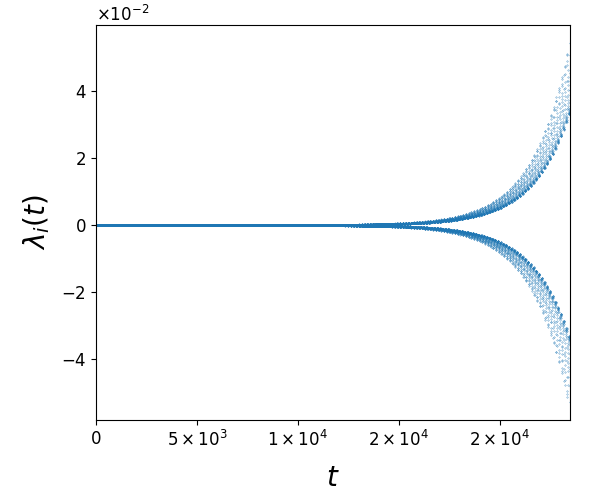}
		}
		\subfloat[$\lambda_2$]{
			\includegraphics[width=.24\textwidth]{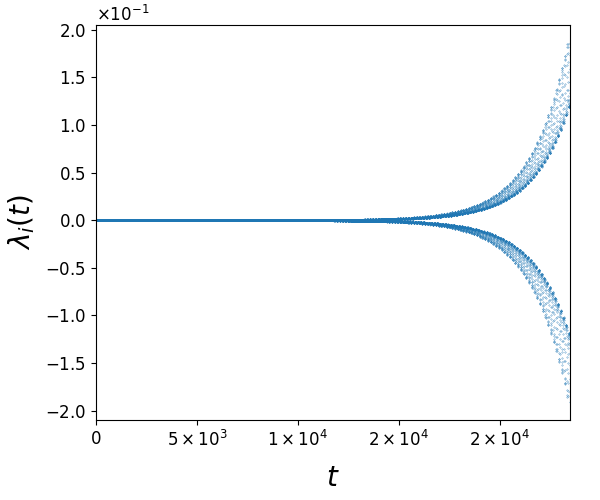}
		}
		\subfloat[$\lambda_3$]{
			\includegraphics[width=.24\textwidth]{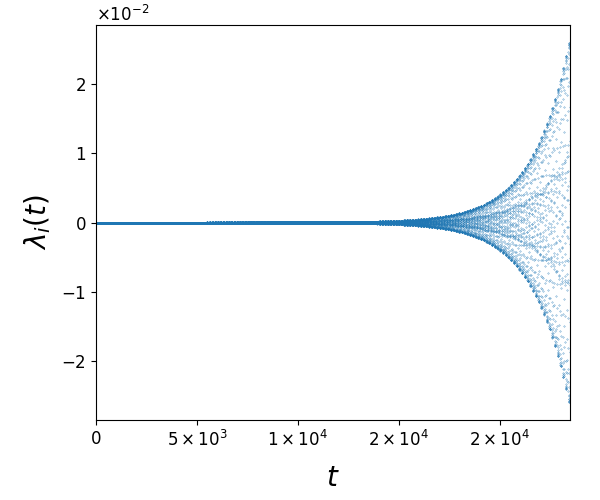}
		}
		\subfloat[$\lambda_4$]{
			\includegraphics[width=.24\textwidth]{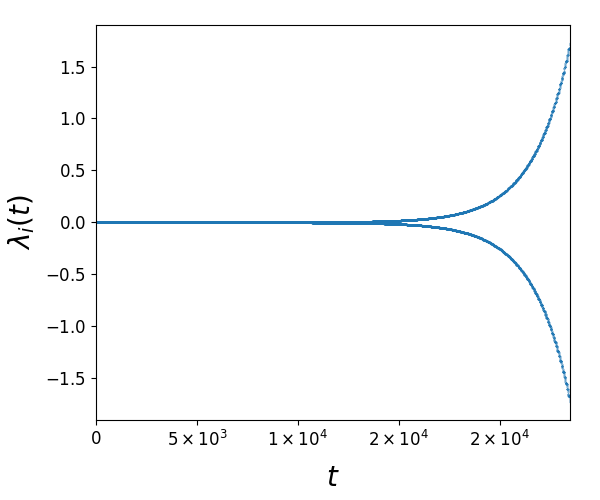}
		}

		\caption{Time evolution of the components of the Lagrange multiplier $\lambda$ for the barely passing guiding centre particle and the three-stage Gauss--Legendre Runge--Kutta method with symplectic projection.}
		\label{fig:guiding_centre_psymplectic_lambda}
	\end{center}
\end{figure}

For many integrators we observed that the symplectic projection is unstable, even when it does not amount to the post projection method.
The reason for this is that in general there are no bounds on the amplitude of Lagrange multiplier $\lambda$ that prevent the unprojected solution from moving further and further away from the constraint submanifold $\Delta$.
This behaviour is illustrated in Figure~\ref{fig:symplectic-projection} and prominently exemplified by the passing guiding centre particles as shown in Figure~\ref{fig:guiding_centre_psymplectic_lambda}. Here, the Lagrange multiplier oscillates between positive and negative values while its amplitude grows without bound.
Eventually, the projection step becomes larger than the integration step and the simulation becomes unstable.
A deeper understanding of this behaviour is very much desirable. One could, for example, apply backward error analysis and try to understand under which conditions the modified equation for the evolution of $\lambda$ is stable. There are also indications of a connection with index reduction of the continuous problem. Such investigations, however, are left for future work.

\section{Summary}

We have devised several projection methods and analyzed their influence on the symplecticity and long-time stability of variational integrators applied to degenerate Lagrangian systems.
The corresponding system of equations constitutes a system of differential-algebraic equations of index two, for which standard symplectic integrators like Gauss--Legendre or Gauss--Lobatto Runge--Kutta methods are well known to deliver poor performance. In particular their order of accuracy is reduced severely compared to their order when applied to ordinary differential equations.

This underperformance can be remedied by the application of appropriate projection methods. In the context of symplectic or variational integrators, this approach raises the question what influence such a projection method has on the symplecticity and long-time stability of the resulting integrator.
While for simple problems all projection methods, and even some of the unprojected integrators, lead to long-time stable simulations, the only universally stable methods have been found to be the symmetric and midpoint projection.
In most examples, the standard projection admits a drift in the solution and thus the energy error, rendering the simulation unstable in finite time. When exactly that happens depends strongly on the order of the underlying numerical integrator and thus the error in the algebraic constraint.
Among the two stable projection methods, the symmetric projection is usually preferable, as the midpoint projection leads to integration methods of reduced order, whereas the symmetric projection restores the order of the underlying numerical integrator.
On the other hand, the midpoint projection is exactly symplectic for the midpoint and SRK3 methods.
Even though the symmetric projection method is not exactly symplectic, we found the error in the symplecticity condition to be so small, that it has no practical influence on the long-time stability of the projected variational integrator.

We also analyzed a modification of the symmetric projection method that is symplectic in a generalized sense, namely in that it is preserving a generalized symplectic structure on a larger space. While for some problems this method leads to similarly good results as the symmetric method, it fails for others. The reason for this failure is that the Lagrange multiplier can grow without bounds. If that happens, the simulation becomes unstable.
However, in those cases where the symplectic projection can be applied, it is preferable over the symmetric projection as it has lower computational cost. For the symplectic projection, the underlying numerical integrator and the projection step can be performed independently, whereas in the symmetric case the whole system has to be solved at once.

\section*{Acknowledgments}

\noindent Helpful discussions with Joshua Burby, Leland Ellison, Melvin Leok and Tomasz Tyranowski are gratefully acknowledged.
Moreover, the author is indebted to Omar Maj and Hiroaki Yoshimura for inspiring conversations as well as reading a draft of the paper. 
The author has received funding from the European Union's Horizon 2020 research and innovation programme under the Marie Sklodowska--Curie grant agreement No 708124. The views and opinions expressed herein do not necessarily reflect those of the European Commission.

\phantomsection
\addcontentsline{toc}{section}{References}
\bibliographystyle{plainnat}
\bibliography{Projected_VIs_for_Degenerate_Lagrangian_Systems}

\end{document}